\documentclass[12pt]{report}

\usepackage{amssymb,amsmath,latexsym}
\voffset
-1.5cm

\usepackage{a4,amscd,graphicx,epsfig}
\usepackage[all]{xy}
\everymath{\displaystyle}

\newtheorem{teo}{Theorem}[section]
\newtheorem{lem}[teo]{Lemma}
\newtheorem{cor}[teo]{Corollary}
\newtheorem{exa}[teo]{Example}
\newtheorem{prop}[teo]{Proposition}
\newtheorem{defi}[teo]{Definition}

\newtheorem{remark}[teo]{Remark}
\newtheorem{remarks}[teo]{Remarks}

\newcommand{\Pm}{\mathbb{P}}
\newcommand{\mr}{\mathbb{R}}
\newcommand{\mc}{\mathbb{C}}
\newcommand{\mz}{\mathbb{Z}}
\newcommand{\mh}{\mathbb{H}}

\newcommand{\mn}{\mathbb{N}}
\newcommand{\mq}{\mathbb{Q}}
\newcommand{\mm}{\mathbb{M}}

\newcommand{\mx}{\mathbb{X}}

\newcommand{\md}{\mathbb{D}}
\newcommand{\ms}{\mathbb{S}}
\newcommand{\mP}{\mathbb{P}}

\newcommand{\Yy}{{\mathcal Y}}
\newcommand{\Aa}{{\mathcal A}}
\newcommand{\Bb}{{\mathcal B}}
\newcommand{\Cc}{{\mathcal C}}
\newcommand{\Dd}{{\mathcal D}}
\newcommand{\Ee}{{\mathcal E}}
\newcommand{\Ff}{{\mathcal F}}
\newcommand{\Gg}{{\mathcal G}}

\newcommand{\Ii}{{\mathcal I}}

\newcommand{\Kk}{{\mathcal K}}
\newcommand{\Ll}{{\mathcal L}}
\newcommand{\Mm}{{\mathcal M}}

\newcommand{\Pp}{{\mathcal P}}
\newcommand{\Qq}{{\mathcal Q}}
\newcommand{\Rr}{{\mathcal R}}
\newcommand{\Ss}{{\mathcal S}}
\newcommand{\Tt}{{\mathcal T}}
\newcommand{\Uu}{{\mathcal U}}
\newcommand{\Ww}{{\mathcal W}}
\newcommand{\Vv}{{\mathcal V}}

\newcommand{\SG}{{\mathfrak S}}

\newcommand{\gG}{{\mathfrak g}}
\newcommand{\sG}{\mathfrak s}
\newcommand{\lG}{\mathfrak l}

\newcommand{\hL}{{\widehat{L}}}

\newcommand{\C}{{\mathbb C}}

\newcommand{\Z}{{\mathbb Z}}
\newcommand{\R}{{\mathbb R}}
\newcommand{\X}{{\mathbb X}}

\newcommand{\D}{{\mathbb D}}
\newcommand{\N}{{\mathbb N}}

\def\ort#1{#1^\perp}

\def\arctgh{\mathrm{arctgh\, }}

\def\Dim{\emph{Proof : }}
\def\cvd{\nopagebreak\par\rightline{$_\blacksquare$}}

\def\E#1#2{\left\langle #1,#2 \right\rangle}  
\def\fut{\mathrm{I}^+}                         
\def\pass{\mathrm{I}^-}                         

\def\SL#1#2{\mathrm{SL}(#1,\,\mathbb{#2})}

\def\OO{\mathrm{O}}
\def\SOO{\mathrm{SO}}

\def\ISO{\mathrm{Isom}}

\def\coom{\mathrm{H}^}
\def\eps{\varepsilon}                        

\def\Ad{\mathrm{Ad}}
\def\tr{\mathrm{tr}}                          

\def\d{\mathrm{d}}
\def\ch{\mathrm{ch\,}}
\def\sh{\mathrm{sh\,}}
\def\tgh{\mathrm{tgh\,}}

\def\arctgh{\mathrm{arctgh\, }}
\def\grad{\mathrm{grad\, }}

\title{Canonical Wick rotations \\ in 3-dimensional gravity}

\author {Riccardo Benedetti and Francesco Bonsante}

\date {}

\begin{document}

\maketitle
\pagestyle{empty}
\noindent Riccardo Benedetti\\
Dipartimento di Matematica, Universit\`a di Pisa,\\
Largo Bruno Pontecorvo, 5 , I-56127 Pisa;\\
benedett@dm.unipi.it\\

\medskip

\noindent Francesco Bonsante\\
Scuola Normale Superiore di Pisa;\\
f.bonsante@sns.it\\ 
\bigskip

\noindent 
The authors have been supported by the INTAS project ``CalcoMet-GT''
03-51-6336 
\bigskip


\begin{abstract}\label{ABS}
We develop a {\it canonical Wick rotation-rescaling theory
in $3$-dimensional gravity}.  This includes
\smallskip

(a) {\it A simultaneous classification}: this shows how maximal
globally hyperbolic spacetimes of arbitrary {\it constant curvature},
which admit a {\it complete Cauchy surface} and {\it canonical
cosmological time}, as well as {\it complex projective structures} on
arbitrary surfaces, are all different materializations of ``more
fundamental'' encoding structures.
\smallskip

(b) {\it Canonical geometric correlations}: this shows how spacetimes
of different curvature, that share a same encoding structure, are
related to each other by {\it canonical rescalings}, and how they can be
transformed by {\it canonical Wick rotations} in {\it hyperbolic
$3$-manifolds}, that carry the appropriate {\it asymptotic} projective
structure. Both Wick rotations and rescalings act along the
canonical cosmological time and have {\it universal rescaling
functions}. These correlations are {\it functorial} with respect to
isomorphisms of the respective geometric categories.
\smallskip

This theory applies in particular to spacetimes with {\it compact}
Cauchy surfaces. By the Mess/Scannell classification, for every fixed
genus $g\geq 2$ of a Cauchy surface $S$, and for any fixed value of
the curvature, these spacetimes are parametrized by pairs $(F,\lambda)
\in\Tt_g\times \Mm\Ll_g$, where $\Tt_g$ is the Teichm\"uller space of
hyperbolic structures on $S$, $\lambda$ is a {\it measured geodesic
lamination} on $F$.  On the other hand, $\Tt_g\times \Mm\Ll_g$ is also
Thurston's parameter space of {\it complex projective structures} on
$S$. The Wick rotation-rescaling theory provides, in particular, a
transparent geometric explanation of this remarkable
coincidence of parameter spaces, and contains a wide generalization of
Mess/Scannell classification to the case of non-compact Cauchy
surfaces.  These general spacetimes of constant curvature are
eventually encoded by a kind of {\it measured geodesic laminations
$\lambda$ defined on some straight convex sets $H$ in $\mh^2$},
possibly in invariant way for the proper action of some discrete
subgroup of $PSL(2,\R)$. We specifically analyze the remarkable
subsectors of the theory made by $\Mm\Ll(\mh^2)$-{\it spacetimes}
($H=\mh^2$), and by $\Qq\Dd$-{\it spacetimes} (associated to $H$
consisting of one geodesic line) that are generated by quadratic
differentials on Riemann surfaces. In particular, these incorporate
the spacetimes with compact Cauchy surface of genus $g\geq 2$, and of
genus $g=1$ respectively.  We analyze {\it broken $T$-symmetry} of AdS
$\Mm\Ll(\mh^2)$-spacetimes and its relationship with {\it earthquake}
theory, beyond the case of compact Cauchy surface.

Wick rotation-rescaling does apply on the {\it ends} of geometrically finite
hyperbolic $3$-manifolds, that hence realize concrete interactions of
their globally hyperbolic {\it ending spacetimes} of constant
curvature. This also provides further ``classical amplitudes"
of these interactions, beyond the volume of the hyperbolic convex
cores.

\end{abstract}


\pagestyle{headings}
\tableofcontents 
\newpage{\pagestyle{empty}\cleardoublepage}
\chapter{General view on themes and contents }\label{INTRO}
A basic fact of $3$-dimensional geometry is that the Ricci tensor
determines the Riemann tensor. This implies that the solutions of pure
3D gravity are the Lorentzian or Riemannian $3$-manifolds of {\it
constant curvature}. Usually the curvature is normalized to be
$\kappa=0,\ \pm 1$. The sign of the curvature coincides with the sign
of the cosmological constant. We stipulate that all manifolds are {\it
oriented} and that the Lorentzian spacetimes are also {\it
time-oriented}. We could also include in the picture the presence of
{\it world lines} of particles. A typical example is given by the {\it
cone manifolds} of constant curvature with cone locus at some embedded
link. The cone angles reflect the ``mass'' of the particles. In the
Lorentzian case we also require that the world lines are of causal
type (see e.g. \cite {tH, BG}(2) and also Chapter \ref{QD} and Section
\ref{particles}).  However in the present paper we shall be mostly concerned
with the matter-free case.

Sometimes gravity is studied by considering separately its different
``sectors'', according to the metric signature (Lorentzian or
Euclidean), and the sign of the cosmological constant. By using the
comprehensive term ``3D gravity'', we propose considering it as
a unitary body, where different sectors actually interact.

The main goal of this paper is to make this suggestion concrete by
fully developing a {\it 3D canonical Wick rotation-rescaling theory}
on $(2+1)$ maximal globally hyperbolic spacetimes of constant
curvature $\kappa = 0,\ \pm 1 \ $, which admit a {\it complete}
spacelike Cauchy surface.  Roughly speaking, we shall see how
spacetimes of arbitrary curvature as well as {\it complex projective
structures} on arbitrary surfaces are all encoded by a common kind of
``more fundamental structures''. Moreover, the theory will establish
explicit {\it canonical geometric correlations}: spacetimes of
different curvature, that share the same instance of encoding
structure, are related to each other via canonical rescalings, and via
canonical Wick rotations can be transformed into hyperbolic
$3$-manifolds that asymptotically carry on the corresponding
projective structure.  In fact such correlations are {\it functorial}
with respect to isometries of spacetimes and isomorphisms of
projective structures on surfaces.

Globally hyperbolic spacetimes with {\it compact} Cauchy surfaces form
a special class to which the theory does apply. This case has been
intensively investigated also in the physics literature (see for
instance \cite{W, Ca2, tH, Mon, A-M-T} and also Section \ref{END}
below). In fact, $(2+1)$ maximal globally hyperbolic spacetimes of
arbitrarily fixed constant curvature $\kappa = 0,\ \pm 1$, having a
compact Cauchy surface $S$, have been classified (up to
Teichm\"uller-like equivalence) by Mess in his pioneering paper
\cite{M}. For $\kappa = 1$ the classification has been completed by
Scannell in \cite{Sc}.  

It turns out that for every fixed genus $g\geq 2$ of $S$ (for $g=1$
see Section \ref{QDST:I}), and for any fixed value of $\kappa$, such
spacetimes are parametrized by the pairs $(F,\lambda)\in \Tt_g\times
\Mm\Ll_g$, where $\Tt_g$ is the classical Teichm\"uller space of
hyperbolic structures on $S$, and $\lambda$ is a {\it measured geodesic
lamination} on $F$ (in fact, in a suitable sense, $\Mm\Ll_g$ only
depends on the topology of $S$).  We could say that all these
spacetimes are different ``materializations'' in 3D gravity of the
same more fundamental structure $\Tt_g\times \Mm\Ll_g$.  On the other
hand, we know that $\Tt_g\times \Mm\Ll_g$ is also Thurston's parameter
space of {\it complex projective structures} on $S$. In particular,
the Wick rotation-rescaling theory will provide a transparent
geometric explanation of this remarkable  coincidence of
parameter spaces. Moreover, the general theory will contain a wide
generalization of Mess-Thurston classification to the case of
non-compact surfaces. This includes the adequate generalization of
$\Mm\Ll_g$. We shall see that suitably defined measured laminations on
{\it hyperbolic surfaces with geodesic boundary}, already introduced
in \cite{Ku} by Kulkarni-Pinkall in order to study complex projective
structures on arbitrary surfaces, furnish the required generalization.
\medskip

{\bf Aknowledgement.} We would like to thank the referee for his
help and suggestions that allows us to substantially improve the
presentation of our work.\\
\medskip

We are going to outline in a bit more detail the main themes and
the contents of the paper.

\section{3-dimensional constant curvature geometry}\label{MOD:I}
Riemannian or Lorentzian $3$-manifolds of constant curvature $\kappa =
0,\pm 1$ have {\it isotropic (local) models}, say $\mx$.  Every
isometry between two open sets of $\mx$ extends to an isometry of the
whole $\mx$. Thus we can adopt the convenient technology of
$(\mx,\Gg)$-{\it manifolds}, where $\Gg = {\rm Isom}(\mx)$, including
{\it developing maps} and ``compatible'' {\it holonomy
representations}.  In the Riemannian case, we will deal mostly with
{\it hyperbolic} manifolds ($\kappa = -1$), so that $\mx = \mh^3$, the
hyperbolic space.  In the Lorentzian case, we will denote the models
by $\mx_\kappa$, called the $3$-dimensional {\it Minkowski} ($\kappa
=0$), {\it de Sitter} ($\kappa =1$) and {\it anti de Sitter} ($\kappa
=-1$) spacetime, respectively. More details on this matter are
collected in Chapter \ref{MOD}, where we also recall some basic facts
about {\it complex projective structures} on surfaces (Section
\ref{SPS}), that are, by definition, $(S^2, PSL(2,\C))$-manifold
structures, where $S^2 = \C\mP^1$ is the Riemann sphere, and
$PSL(2,\C)$ is naturally identified with the group of complex
automorphisms of $S^2$. This will include the notions of {\it
$H$-hull}, {\it canonical stratification} and {\it Thurston metric}.

\section{Wick rotation and rescaling}\label{WICK_ROTATION:I} 
Wick rotation is a very basic procedure for inter-playing Lorentzian
and Riemannian geometry. The simplest example applies to $\R^{n+1}$
endowed with both the standard Minkowski metric $-dx_0^2+\dots
+dx_{n-1}^2 + dx_{n}^2$, and the Euclidean metric $dx_0^2+\dots +
dx_{n}^2$.  By definition (see below), these are related via a Wick
rotation directed by the vector field $\partial/\partial
x_{0}$. Sometimes one refers to it as ``passing to the imaginary
time''. More generally we have:
\begin{defi}\label{WR}{\rm  
Given a manifold $M$ equipped with a Riemannian metric
$g$ and a Lorentzian metric $h$, then we say
that {\it $g,\ h$ are related via a {\rm rough} Wick rotation directed by
$X$}, if:
\smallskip

(1) $X$ is a nowhere vanishing $h$-timelike and future directed
vector field on $M$;
\smallskip

(2) For every $y\in M$, the $g$- and $h$-orthogonal spaces to $X(y)$
coincide and we denote them by $\ort{<X(y)>}$.
\smallskip

\noindent The positive function $\beta$ defined on $M$ by 
$$||X(y)||_g = - \beta(y)||X(y)||_h$$ is called the {\it vertical
rescaling function} of the Wick rotation.

A Wick rotation is said {\it conformal} if there is also a
positive {\it horizontal rescaling function $\alpha$} such that, for
every 
$y\in M$, 
$$g|_{\ort{<X(y)>}} \ =\ \alpha(y)h|_{\ort{<X(y)>}} \ . $$
}
\end{defi}
In fact, all metrics $g,\ h$ as above are canonically related by a
rough Wick rotation: we use $g$ to identify $h$ to a field of linear
automorphisms $h_y \in {\rm Aut}(TM_y)$, and we take as $X(y)$ the
field of $g$-unitary and $h$-future directed eigenvectors of $h_y$,
with negative eigenvalues.  Call $X_{(g,h)}$ this canonical vector
field associated to the pair of metrics $(g,h)$. Any other field $X$
as in Definition \ref{WR} is of the form $X=\lambda X_{(g,h)}$, for some
positive function $\lambda$. 

If we fix a nowhere vanishing vector
field $X$, and two positive functions $\alpha, \ \beta$ 
on $M$ , then the conformal Wick rotation directed by $X$ and with
rescaling functions $(\alpha,\beta)$ establishes a bijection, say
$W_{(X;\beta,\alpha)}$, between the set of Riemannian metrics on
$M$, and the set of Lorentzian metrics which have $X$ as a
 future directed timelike field. In particular, the couple $(g,X)$
encodes part of the global causality of
$h=W_{(X;\beta, \alpha)}(g)$. Clearly $W_{(X;\beta, \alpha)}^{-1}=
W_{(X;\frac{1}{\beta}, \frac{1}{\alpha})}$.
\smallskip

From now on {\it we shall consider only conformal Wick rotation}, so
we do not specify it anymore. The couples $(g,h)$ related via a Wick
rotation, and such that {\it both $g$ and $h$ are solutions of pure
gravity} are of particular interest, especially when the support
manifold $M$ has a non-trivial topology.

\paragraph {Rescaling directed by a vector field.} 
This is a simple operation (later simply called ``rescaling'') on
Lorentzian metrics formally similar to a Wick rotation. Let $h$ and
$h'$ be Lorentzian metrics on $M$.  Let $v$ be a nowhere vanishing
vector field on $M$ as above.  Then $h'$ is obtained from $h$ via a
{\it rescaling directed by $v$, with rescaling functions
$(\alpha,\beta)$}, if
\smallskip

(1) For every $y\in M$, the $h$- and $h'$-orthogonal spaces to $v(y)$
coincide and we denote them by $\ort{<v(y)>}$.
\smallskip

(2) $h'$ coincides with $\beta h$ on the line bundle $<v>$ spanned by $v$.
\smallskip

(3) $h'$ coincides with $\alpha h$ on $\ort{<v>}$.

Again, rescalings which relate different solutions of pure gravity,
possibly with different cosmological constants, are of particular
interest.

\section {Canonical cosmological time}\label{CT:I}
This is a basic notion (see \cite{A}) that will play a crucial r\^ole in
our Wick rotation-rescaling theory. In Section \ref{CT} we will recall
the precise definition and the main properties.  Roughly speaking, for
any spacetime $M$, its {\it cosmological function} gives {\it the
(possibly infinite) proper time that every event $q \in M$ has been in
existence in $M$}. If the cosmological time function is ``regular''
(see \cite{A} or Section ~\ref{CT}) - this means in particular that it
is finite valued for every $q\in M$ - then actually it is a continuous
{\it global time} on $M$.  This canonical {\it cosmological time} (if
it exists) is not related to any specific choice of coordinates in
$M$, is invariant under the automorphisms of $M$, and represents an
intrinsic fundamental feature of the spacetime. In a sense it gives
the Lorentz distance of every event from the ``initial singularity'' of
$M$.  In fact, we are going to deal with spacetimes having rather tame
cosmological time; in these cases the geometry of the initial
singularity will quite naturally arise.

\section {Classification of flat globally hyperbolic \\ 
spacetimes}
\label{FLAT:I}
For the basic notions of global Lorentzian geometry and causality we
refer for instance to \cite{BEE, H-E}. Some specific facts about {\it
(maximal) globally hyperbolic spacetimes} are recalled in Section
\ref{glob:hyp:ST}.
\smallskip

One of our goals is to classify the maximal globally hyperbolic
3-dimensional spacetimes of constant curvature that contain a {\it
complete} Cauchy surface. These properties lift to the respective
locally isometric universal covering spacetimes, so we will classify
the simply connected ones, keeping track of the isometric action of
the fundamental group $\pi_1(S)$, $S$ being any Cauchy surface.
\smallskip

A standard analytic approach to the classification of constant
curvature globally hyperbolic spacetimes $M$ is in terms of solutions
of the Gauss-Codazzi equation at Cauchy surfaces $S$, possibly
imposing some supplementary conditions to such solutions, that
translates some geometric property of the embedding of $S$ into $M$.
For instance a widely studied possibility is to require that the
surface $S$ has constant mean curvature in $M$ (e.g. we refer to
\cite{Mon, A-M-T, B-Z, Kra-Sch}).

Here we follow a rather different approach, initiated by Mess \cite{M}
in the case of compact Cauchy surfaces. By restricting to the
``generic'' case of spacetimes that have cosmological time, we realize
that in a sense each one is determined by the ``asymptotic states'' of
the level surfaces of its canonical cosmological time, rather than the
embedding data of some Cauchy surface $S$. In the case of compact
Cauchy surfaces, the central objects in \cite{M} rather were the
holonomy groups; the r\^ole of the cosmological time was recognized in
\cite{Sc} ($\kappa =1$) and fully stressed (with its asymptotic
states) in \cite{BG}(3) ($\kappa =0$). In general, these asymptotic
states will appear as additional geometric structures on $S$ such as a
hyperbolic structure (possibly with geodesic boundary), and a {\it
measured geodesic lamination} suitably defined on it. The intrinsic
geometry of the level surfaces is determined in terms of them by means
of a grafting-like construction.

The Wick rotation-rescaling mechanism will be based on the 
fact that eventually the intrinsic geometry of these level surfaces
does not depend on the curvature, up to some scale factor.

\medskip 

We will consider at first (Chapter \ref{FGHST}) the {\it flat}
spacetimes (\emph{i.e.} of constant curvature $\kappa = 0$).
\smallskip

In~\cite{Ba}(1), Barbot showed that, except some sporadic cases and
possibly reversing the time orientation, the simply connected maximal
globally hyperbolic flat spacetimes that contain a complete Cauchy
surface coincide with the so called \emph{regular domains} (see below)
$\ \Uu$ of the Minkowski space $\mx_0$. Moreover, when $$\Uu \to
\Uu/\tilde{\Gamma}$$ is a universal covering, $\tilde{\Gamma}\cong
\pi_1(S)$ being a subgroup of $\ISO^+_0(\mx_0)$ (the group of
isometries of $\mx_0$ that preserves the orientations), then we also
have from ~\cite{Ba}(1) informations about $\tilde{\Gamma}$ and its
linear part $\Gamma \subset SO^+(2,1)$ (see below, and
\ref{flat:descr:teo} for the precise statement of these results). As a
corollary we know for example (see corollary \ref{flat:descr:cor})
that:

{\it If $\pi_1(S)$ as above is not Abelian, then the corresponding
universal covering is a regular domain different from the future of a
spacelike geodesic line, and the linear part of the holonomy is a
faithful and discrete representation of $\pi_1(S)$ in $SO(2,1)$.}
\smallskip

A regular domain $\Uu \subset \mx_0$ is a convex domain that coincides
with the intersection of the future of its null support planes. We
also require that there are at least two null support planes. Note
that a regular domain is future complete.
\smallskip

One realizes that all the sporadic exceptions have {\it not} canonical
cosmological time. On the other hand we have (see Proposition
\ref{COSMT_REGD}, Subsection ~\ref{ctrd:sec}, Proposition
\ref{gen:tree:prop}, and also \cite{Ba}(1)).
\begin{prop}\label{flatCT:I}
Every flat regular domain $\Uu$ has canonical cosmological time
$T$. In fact $T$ is a {\rm C}$^{1,1}$-submersion onto
$(0,+\infty)$. Every $T$-level surface $\Uu(a)$, $a\in
(0,+\infty)$, is a complete Cauchy surface of $\Uu$. For every $x\in
\Uu$, there is a unique past-directed geodesic timelike segment
$\gamma_x$ that starts at $x$, is contained in $\Uu$, has finite
Lorentzian length equal to $T(x)$. The other end-point of $\gamma_x$
belongs to the frontier of $\Uu$ in $\mx_0$. The union of these
boundary end-points makes the {\rm initial singularity} $\Sigma_\Uu
\subset \partial \Uu$ of $\Uu$. $\Sigma_\Uu$ is a
spacelike-path-connected subset of $\mx_0$, and this gives it a
natural {\rm $\R$-tree} structure.
\end{prop}
Hence the study of regular domains (and their quotient spacetimes) is
equivalent to the study of maximal globally hyperbolic flat spacetimes
having a complete spacelike Cauchy surface $S$ {\it and} canonical
cosmological time. 
\smallskip

Let us consider the $T$-level surface $\Uu(1)$ of a regular domain
$\Uu$. We have a natural continuous {\it retraction}
$$r: \Uu(1)\to \Sigma_\Uu \ . $$ Moreover, the gradient of $T$ is
a unitary vector field, hence it induces the {\it Gauss map} (here we
are using the standard embedding of the hyperbolic plane $\mh^2$ into
$\mx_0$)
$$N: \Uu(1)\to \mh^2 \ .$$ These two maps are of central importance
for all constructions. We realize that the closure $H_\Uu$ of the
image ${\rm Im}(N)$ of the Gauss map in $\mh^2$ is a {\it straight
convex set} in $\mh^2$, \emph{i.e.} an closed set that is the convex
hull of a {\it ideal} set  contained in the natural
boundary $S^1_\infty$ of $\mh^2$ (Section \ref{RD:sec}). If $\Uu \to
\Uu/\tilde{\Gamma}$ is a universal covering as above, then the
cosmological time is $\tilde{\Gamma}$-invariant, and the action of
$\tilde{\Gamma}$ both extends to an isometric action on the initial
singularity, and to an isometric action on $H_\Uu$ (via the linear
part $\Gamma$ indeed).

Thus we can distinguished two sub-cases:
\smallskip

(a) {\it non-degenerate}: when dim $H_\Uu=2$.  In this case we know
that $\Gamma \cong \tilde{\Gamma}$ is a discrete torsion-free subgroup
of $SO^+(2,1)$ so that $F=\mh^2/\Gamma$ is a complete hyperbolic
surface homeomorphic to $S$; 
\smallskip

(b) {\it degenerate}: when dim $H_\Uu=1$. In such a case $H_\Uu$
consists of one geodesic line of $\mh^2$. Equivalently, the initial
singularity reduces to one complete real line (we say that it is {\it
elementary}).  The group $\tilde{\Gamma}$ is isomorphic to either
$\{0\}$, or $\mz$, or $\mz\oplus\mz$.
\medskip

The degenerate case will be treated in Chapter \ref{QD} (see also
Section \ref{QDST:I} below) in the more general framework of so called
$\Qq\Dd$-{\it spacetimes}. Let us consider here the non-degenerate one.
\smallskip

We take the partition of $\Uu(1)$ given by the closed sets
$r^{-1}(y)$, $y\in \Sigma_\Uu$. Via the retraction we can pullback to
this partition the metric structure of $\Sigma_\Uu$, and (in a
suitable sense) we can project everything onto $H_\Uu$, by means of
the Gauss map.  We eventually obtain a triple 
$$\lambda_\Uu = (H_\Uu, \Ll_\Uu,\mu_\Uu)$$ where $(\Ll_\Uu,\mu_\Uu)$
is a kind of {\it measured geodesic lamination} on $H_\Uu$. The
geometry of the initial singularity is, in a sense, ``dual'' to the
geometry of the lamination.  More precisely, if $r^{-1}(y)$ is
1-dimensional, then it is a geodesic line, so that the union of such
lines makes a ``lamination'' in $\Uu(1)$. We can define on it a 
``transverse measure'' such that the mass of any transverse path
is given by the integral of the Lorentzian norm of the derivative
of $r$. The lamination $\lambda$ is obtained via the
push-forward by $N$ of this lamination on $\Uu(1)$. 

These measured geodesic laminations on straight
convex sets are formalized and studied in Section \ref{laminations}.
The detailed construction of $\lambda_{\Uu}$ is done in
Section \ref{REGD_ML}.  On the other hand, in Section \ref{ML_REGD} we
get the inverse construction. Let us denote by 
$$ \Mm\Ll = \{\lambda = (H,\Ll,\mu)\}$$ the set of such triples, 
that is the set of {\it measured geodesic laminations defined on
some 2-dimensional straight convex set in $\mh^2$}. Let us denote
by $\Rr$ the set of non-degenerate regular domains in $\mx_0$.  Note
that there is a natural left action of $SO(2,1)$ on $\Mm\Ll$, and
of $\ISO^+(\mx_0)$ (hence of the translations subgroup $\mr^3$) on
$\Rr$.
 
Then we construct a map
$$ \Uu^0: \Mm\Ll \to \Rr,\ \ \lambda \to \Uu^0_\lambda $$ such
that, by setting $\Uu=\Uu^0_\lambda$, then $\lambda_{\Uu} =
\lambda$.  Both constructions are rather delicate. Summarizing
they give us the following {\it
classification} theorem for non-degenerate flat regular domains.
\begin{teo}\label{GEN_FLAT_CLASS}
The map $$ \Uu^0: \Mm\Ll \to \Rr,\ \ \lambda \to \Uu^0_\lambda $$
induces a bijection between $\Mm\Ll$ and $\Rr/\mr^3$, and a 
bijection between $\Mm\Ll/SO(2,1)$ and $\Rr/\ISO^+(\mx_0)$. 
\end{teo}
While proving this theorem, in Subsection ~\ref{continuity} we will
also study {\it continuity properties} of the map $\Uu^0$.  The
following corollary is immediate.
\begin{cor}\label{MLH2} 
The regular domains in $\mx_0$ that have {\rm surjective Gauss map}
are parametrized by the measured geodesic laminations (of the most
general type indeed) defined on the whole of $\mh^2$.
\end{cor} 

When $\Uu \to \Uu/\tilde{\Gamma}$ is a universal covering,
$\lambda_{\Uu}$ is $\Gamma$-invariant; $$Z=H_\Uu/\Gamma$$ 
is a {\it straight convex subset of $F$} that is 
convex hyperbolic surface with geodesic boundary embedded in the
complete hyperbolic surface $$F= \mh^2/\Gamma \ .$$ $F$ is
homeomorphic to the interior of $Z$. Finally $\lambda_\Uu$ is
the pull-back of a measured geodesic lamination suitably defined on
$Z$ (see Subsection \ref{Gamma:inv:lam}). Set
 $$ \Mm\Ll^\Ee = \{(\lambda,\Gamma)\}$$ where $\lambda \in \Mm\Ll$,
and $\Gamma$ is a discrete torsion-free subgroup of $SO(2,1)$ acting
on $\lambda$ as above.  Set $$\Rr^\Ee = \{(\Uu,\tilde{\Gamma})\}$$
where $\Uu\in \Rr$ and $\tilde{\Gamma} \subset \ISO^+(\mx_0)$ acts
properly on $\Uu$. The action of $SO(2,1)$ ( $\ISO^+(\mx_0)$) extends
to $ \Mm\Ll^\Ee$ ($\Rr^\Ee$) by conjugation of the subgroup
$\Gamma$ ($\tilde{\Gamma}$). There is a natural {\it equivariant
version} of all constructions, and all structures descend to the
quotient spacetimes. This leads to an extended map
$$ \Uu^0: \Mm\Ll^\Ee \to \Rr^\Ee$$
$$(\lambda,\Gamma) \to (\Uu^0_\lambda,\Gamma^0_\lambda)$$ where
$\Gamma$ is the linear part of $\Gamma^0_\lambda$. Finally (Subsection
\ref{flat:equiv}) we have

\begin{teo}\label{FULLFLAT} The extended map $\Uu^0$ induces
a bijection between $\Mm\Ll^\Ee/SO(2,1)$ and 
$\Rr^\Ee/\ISO^+(\mx_0)$.
\end{teo}

Let us consider {\it marked} maximal globally hyperbolic flat
spacetimes containing a complete Cauchy surface of fixed topological
type $S$, up to {\it Teichm\"uller-like equivalence} of spacetimes.
Assume furthermore that they have cosmological time and non-degenerate
universal covering.  Denote by $MGH_0(S)$ the corresponding
Teichm\"uller-like space.  As a corollary of the previous theorem we
have:
\begin{cor}\label{Teich-like} For every fixed topological type
$S$, $MGH_0(S)$ can be identified with $\Mm\Ll^\Ee_S/SO(2,1)$,
where $\Mm\Ll^\Ee_S$ is the set of $(\lambda,\Gamma)\in \Mm\Ll^\Ee$
such that the complete hyperbolic surface $F=\mh^2/\Gamma$ is
homeomorphic to $S$.
\end{cor}
In other words, such a Teichm\"uller-like space consists of the
couples $(F,\lambda)$ where $F$ is any complete hyperbolic structure
parametrized by the base surface $S$, and $\lambda$ denotes any
measured geodesic lamination on some straight convex set $Z$ in $F$.

Recall that if $\pi_1(S)$ is not Abelian, then the additional
conditions stated before the corollary are always satisfied. So we
have
\begin{cor}\label{non-abel} 
If $\pi_1(S)$ is not Abelian, then the restriction of $\Uu^0$ to
$\Mm\Ll^\Ee_S$ induces a parametrization of maximal globally
hyperbolic flat spacetimes containing a complete Cauchy surface
homeomorphic to $S$.
\end{cor} 

\section{$\Mm\Ll$-spacetimes}\label{ML:I}
Although we adopt a slightly different definition, it turns out that
the above laminations are equivalent to the ones already introduced in
\cite{Ku}.  Since \cite{Ku} it is known that $\Mm\Ll^\Ee/SO(2,1)$ also
parameterizes the 2-dimensional complex projective structures of
``hyperbolic type'' and with ``non-degenerate'' canonical
stratification. This is unfolded in terms of a 3-dimensional
hyperbolic construction (see Section ~\ref{SPS} and Chapter
\ref{HYPE}). Given $(\lambda,\Gamma)\in \Mm\Ll^\Ee$, denote by
$\mathring H$ the interior of the corresponding straight
convex set $H$, $\mathring Z=\mathring H/\Gamma$.  Then we construct
$(D^{\rm hype}_\lambda,h^{\rm hype}_\lambda)$ where
$$D^{\rm hype}_\lambda: \mathring H\times (0,+\infty)\to \mh^3$$ 
is a developing map of a hyperbolic structure $M_\lambda$ on
$\mathring Z \times (0,+\infty)$, 
$$h^{\rm hype}_\lambda: \Gamma \to \ISO^+(\mh^3)$$ is a compatible
holonomy representation.  The map $D^{\rm hype}_\lambda$ extends (in a
$\Gamma$ equivariant way, via $h^{\rm hype}_\lambda$) to $ H\times
\{0\}\cup \mathring H\times \{+\infty\}$. This is in fact the extension to the
{\it completion} of the hyperbolic metric.  The restriction $D^{\rm
hype}_\lambda$ to $H\times \{0\}$ realizes a locally isometric pleated
immersion of $H$ into $\mh^3$, having $\lambda$ as {\it bending}
measured lamination. This gives us the so called {\it hyperbolic
boundary} of $M_\lambda$ (see Section ~\ref{hyp:bend:cocy}, Theorem
~\ref{hyperbolic:compl:teo}).  The restriction of $D^{\rm
hype}_\lambda$ to $\mathring H\times \{+\infty\}$ has values on the boundary
$S^2_\infty$ of $\mh^3$, and is in fact the developing map $D^{{\rm
proj}}$ of a complex projective structure $S_\lambda$ on $\mathring Z$, having
$h^{\rm proj}_\lambda=h^{\rm hype}_\lambda$ as compatible holonomy
representation. This gives the so called {\it asymptotic projective
boundary} of $M_\lambda$ (see Subsection~\ref{proj-bound}).  The
hyperbolic manifold $M_\lambda$ is called the $H$-{\it hull} of
$S_\lambda$.  $(\lambda,\Gamma)$, $(D^{\rm hype}_\lambda,h^{\rm
hype}_\lambda)$ (that is $M_\lambda$), and $(D^{{\rm proj}}_\lambda,
h^{{\rm proj}}_\lambda)$ (that is $S_\lambda$) are determined by each
other, up to natural actions of either $SO(2,1)$ or $\ISO^+(\mh^3)$,
and this provides us with the parametrization mentioned above.
\medskip

Given $(\lambda,\Gamma)\in \Mm\Ll^\Ee$ we construct also suitable
couples $(\Uu^\kappa_\lambda, \Gamma_\lambda^\kappa)$ where
$\Uu^\kappa_\lambda$ is a simply connected maximal globally hyperbolic
spacetimes of constant curvature $\kappa = \pm 1$,
$$\Uu^\kappa_\lambda\to \Uu^\kappa_\lambda/\Gamma_\lambda^\kappa$$ is
a locally isometric universal covering, and
$\Uu^\kappa_\lambda/\Gamma_\lambda^\kappa$ is homeomorphic to $\mathring Z
\times \mr$.  In fact (see Chapter ~\ref{dS}) $\Uu^{1}_\lambda$ is
given in terms of a developing map
$$D^{\rm dS}_\lambda: \mathring H \times \mr  \to \mx_1$$
and a compatible holonomy representation 
$$h^{\rm dS}_\lambda: \Gamma \to \ISO^+(\mx_1)$$ whose construction runs
parallel to the one of $M_\lambda$ (by using the fact, which is
evident in the projective models, that $\mh^3$ and $\mx_1$ share the
same sphere at infinity). It turns out that
$$h^{\rm dS}_\lambda = h^{\rm hype}_\lambda \ .$$ In general these
hyperbolic or de Sitter developing maps are {\it not} injective.
\smallskip 

The construction of $\Uu^{-1}_\lambda$ is based on an AdS version of
the {\it bending procedure} that is carefully analyzed in Section
\ref{bend-ads}. Remarkably, every
$\Uu^{-1}_\lambda$ is a convex domain in $\mx_{-1}$ ({\it i.e.} the
developing map $D^{\rm AdS}$ is an embedding), and
$\Gamma_\lambda^\kappa$ is a subgroup of $\ISO^+(\mx_{-1})$ that acts
properly on it.
\medskip

Hence, for every $\kappa = 0,\pm 1$, we construct a family
$$\Mm\Ll\Ss_\kappa = \{(\Uu^\kappa_\lambda,
\Gamma_\lambda^\kappa)\}/\ISO^+(\mx_\kappa) $$ of maximal globally
hyperbolic spacetimes of constant curvature $\kappa$, sharing the same
``universal'' parameter space $\Mm\Ll^\Ee/SO(2,1)$.  These are
generically called $\Mm\Ll$-{\it spacetimes}.

We roughly collect here some basic properties of the
$\Mm\Ll$-spacetimes. For $\kappa = 1$ see Proposition
\ref{desitter:ct:prop}, for $\kappa = -1$ see Proposition
\ref{memgen:adsresc2:prop}, Corollary ~\ref{memgen:adscompl:cor}, and
Proposition ~\ref{Plambda-past}.

\begin{prop}\label{ML:PROPERTY}
(1) Each $\Uu_\lambda^\kappa$ has canonical cosmological time, say
$T^\kappa_\lambda$, with non-elementary initial singularity.  For
$\kappa =1$ the image of the cosmological time is $(0,+\infty)$, 
whereas for $\kappa = -1$ it is an interval $(0,a_0)$, for some
$\pi/2< a_0 <\pi$.
\smallskip

{\rm {\small [We adopt the following notations. For
every subset $X$ of $(0,+\infty)$, $\Uu_\lambda^\kappa(X) =
(T^\kappa_\lambda)^{-1}(X)$; for every $a\in
T^\kappa_\lambda(\Uu_\lambda^\kappa)$, $\Uu_\lambda^\kappa(a) =
\Uu_\lambda^\kappa(\{a\})$ denotes the corresponding level surface of
the cosmological time. Sometimes we shall also use the notation
$\Uu_\lambda^\kappa(\geq a)$ instead of
$\Uu_\lambda^\kappa([a,+\infty))$, and so on.]}}
\smallskip

(2) For $\kappa =1$, $T^\kappa_\lambda$ is C$^{1,1}$.
 For $\kappa =-1$, it is C$^{1,1}$ on
$$\Pp_\lambda = \Uu_\lambda^{-1}((0,\pi/2))\ .$$ The level surface
$\Uu_\lambda^{-1}(\pi/2)$ is an isometric pleated copy of $\mathring H$
embedded in $\mx_{-1}$, that has $\lambda$ as AdS bending lamination.
\smallskip

(3) Both level surfaces $\Uu_\lambda^{1}(a)$
and $\Uu_\lambda^{-1}(a)$, $a<\pi/2$ are {\it complete}
Cauchy surfaces.
\end{prop}
In Proposition \ref{Plambda-past} we will recognize $\Pp_\lambda$
to be the {\it past part} of $\Uu_\lambda^{-1}$, that is the past
of the ``future boundary of its convex core''.
\smallskip

{\it Some special subfamilies.}

(1) When $H=\mh^2$ and the support of the lamination $\lambda$ is
empty, we say that the corresponding spacetime is {\it static}. This
special case is also characterized by the fact that the initial
singularity of each spacetime $\Uu_\lambda^\kappa$ consists of {\it
  one} point, or that the cosmological times $T_\kappa$ is {\it real
  analytic}. Moreover the cosmological time is also a {\it constant
  mean curvature} (CMC) time.

For a general $H$ we say that it is $H$-static if the support of the
lamination coincides with the boundary geodesics of $H$ (hence they
are all $+\infty$-weighted).  The initial singularity consists now of
one vertex $v_0$ from which a complete half line emanates for each
boundary component of $H$. In the flat case, the portion $r^{-1}(v_0)$
is a homeomorphic deformation retract of the whole spacetime;
it is contained in the $\mh^2$-static spacetime obtained by just
forgetting the lamination $\lambda$, and the respective cosmological
times do agree on such a portion.
\smallskip

(2) As in the above corollary \ref{MLH2}, we point out the
distinguished sub-class of spacetimes, that we call {\it
  $\Mm\Ll(\mh^2)$-spacetimes}, obtained by imposing that $H$ consists
of the whole hyperbolic plane $\mh^2$. It is convenient and
instructive to analyze specific aspects of our Wick rotation-rescaling
theory in such a case. The remarkable spacetimes with compact Cauchy
surfaces of genus $g\geq 2$ belong to this sector of the theory. In
fact the only straight convex set in a compact closed hyperbolic
surface $F$ is the whole of $F$.  In a sense, this cocompact
$\Gamma$-invariant case has {\it tamest features}. For
instance, it implies strong constraints on the measured geodesic
laminations on $\mh^2/\Gamma$. Throughout the paper, we shall focus on
these special features, against the different phenomena that arise for
general $\Mm\Ll(\mh^2)$-spacetimes, even in the {\it finite coarea},
but not cocompact case (see Section \ref{3cusp}).  A key point here is
that we can work with geodesic laminations on $F=\mh^2/\Gamma$ that do
not necessarily have compact support. For a first account of the
cocompact case one can see also Section \ref{END} below.

We will also realize that the $\Mm\Ll(\mh^2)$-spacetimes of curvature
$\kappa = - 1$ have interesting characterizations among general
$\Mm\Ll$ ones. This is related to the interesting behaviour of the
spacetimes with respect to the $T$-{\it symmetry} obtained by
reversing the time orientation (see Section \ref{BROKEN:I}). Moreover,
at some points of the AdS treatment it is technically convenient to
deal first with the $\Mm\Ll(\mh^2)$ subcase (see
Section~\ref{resc-class}).
\medskip

Our linked goals  consist in: 
\smallskip

(a) Pointing out natural and explicit geometric correlations between
the hyperbolic manifolds $M_\lambda$ and the $\Mm\Ll$-spacetimes
$(\Uu^\kappa_\lambda, \Gamma_\lambda^\kappa)$ that share a same
encoding $\lambda$.
\smallskip

(b) Eventually obtain also for $\kappa =\pm 1$, an {\it intrinsic}
characterization of $\Mm\Ll$-spacetimes $\Mm\Ll\Ss_\kappa$, similarly
to the classification already outlined for $\kappa =0$.

\section{Canonical Wick rotations and rescalings}\label{CANWR:I}
The correlations evoked in (a), will be either canonical Wick
rotations or rescalings {\it directed by the gradient of the
cosmological time} and with {\it universal rescaling functions}. This
means that these are constant on each level surface (on which they are
defined) and their value only depends on the corresponding value of
the cosmological time. We stress that they do not depend on
$\lambda$. A delicate point is that all this happens in fact {\it up
to determined ${\rm C}^1$ diffeomorphism}. Note that any ${\rm C}^1$
isometry between Riemannian metrics induces at least an isometry of
the underlying {\it length spaces}. In the Lorentzian case, it
preserves the global causal structure. Wick rotations and rescalings
have higher regularity (they are real analytic indeed) only in the
case of $\mh^2$-static spacetimes.
\smallskip

We are going to summarize (in somewhat rough way) our main results.

\begin{teo}\label{DSRESC:I}
A canonical rescaling directed by the gradient of (the restriction of)
the cosmological time $T^0_\lambda$, with universal rescaling
functions, converts (in equivariant way)
$(\Uu_\lambda^0(<1),\Gamma^0_\lambda)$ into
$(\Uu_\lambda^1,\Gamma^1_\lambda)$. The inverse rescaling satisfies
the same properties with respect to the cosmological time
$T^{1}_\lambda$. The rescaling extends to an (equivariant) isometry
between the respective initial singularities.
\end{teo}
All this is proved in Chapter ~\ref{dS}.

\begin{teo}\label{ADSRESC:I}
A canonical rescaling directed by the gradient of the cosmological
time $T^0_\lambda$, with universal rescaling functions, converts (in
equivariant way) $(\Uu_\lambda^0,\Gamma^0_\lambda)$ into the past part
$\Pp_\lambda$ of $(\Uu_\lambda^{-1},\Gamma^{-1}_\lambda)$.  The
inverse rescaling satisfies the same properties with respect to the
cosmological time $T^{-1}_\lambda$ restricted to $\Pp_\lambda$.  The
rescaling extends to an (equivariant) isometry between the respective
initial singularities.
\end{teo}
This is proved in Chapter \ref{AdS}.

\begin{teo}\label{WR:I} 
(1) A canonical Wick rotation directed by the gradient of (the
restriction of) the cosmological time $T^0_\lambda$, with universal
rescaling functions, converts (in equivariant way)
$(\Uu_\lambda^0(>1), \Gamma^0_\lambda)$ into the hyperbolic
$3$-manifold $M_\lambda$. The inverse Wick rotation satisfies the same
properties with respect to the distance function from the hyperbolic
boundary of $M_\lambda$. The intrinsic spacelike metric on the level
surface $\Uu_\lambda^0(1)$ coincides with the {\rm Thurston metric}
associated to the asymptotic complex projective structure whose
$M_\lambda$ is the $H$-hull.  The {\rm canonical stratification}
associated to such projective structure coincides with the
decomposition of $\Uu_\lambda^0(1)$ given by the fibers of the
retraction on the initial singularity.
\smallskip

(2) This Wick rotation can be transported onto the AdS slab
$\Uu_\lambda^{-1}((\pi/4,\pi/2))$ by means of the rescaling of Theorem
\ref{ADSRESC:I}. This extends continuously to the boundary of this
slab. The boundary component $\Uu^{-1}_\lambda (\pi/4)$ maps
homeomorphically onto the asymptotic projective boundary of
$M_\lambda$, and its intrinsic spacelike metric is again the
associated Thurston metric.  The restriction to the boundary component
$\Uu^{-1}_\lambda (\pi/2)$ is an isometry between such an AdS pleated
surface and the hyperbolic one that makes the hyperbolic boundary of
$M_\lambda$.
\smallskip

(3) The de Sitter rescaling (Theorem \ref{DSRESC:I}) on
$(\Uu_\lambda^0(<1),\Gamma^0_\lambda)$, and the Wick rotation on
$\Uu_\lambda^0(>1)$ ``fit well'' at the level surface
$\Uu_\lambda^0(1)$, and give rise to an immersion of the whole of
$\Uu_\lambda^0$ in $\mP^3$ (by using the projective Klein models).
\end{teo}
The statements in (1) are proved in Chapter \ref{HYPE}.  Note that the
canonical Wick rotations {\it cut the initial singularity off}.  For
(2) see Chapter~\ref{AdS}; (3) is proved in Chapters ~\ref{HYPE},
~\ref{dS}.

\section{Full classification}\label{FULL_CLASS:I}
Finally we get the following {\it full classification result} of
maximal globally hyperbolic spacetimes of constant curvature
containing a complete Cauchy surface.

\begin{teo}\label{FULL_CLASS}   For every $\kappa = 0, \pm 1$, the family
of $\Mm\Ll$-spacetimes $\Mm\Ll\Ss_\kappa$ coincides with the family
$MGH_\kappa$ of maximal globally hyperbolic spacetimes of constant
curvature $\kappa$, with cosmological time, non-elementary initial
singularity, and containing a complete Cauchy surface, considered up
to Teichm\"uller-like equivalence (by varying the topological type).
Hence, these spacetimes are parametrized by $\Mm\Ll^\Ee/SO(2,1)$, for
every value of $\kappa$.  Wick rotations and rescalings establish
canonical bijections with the set of non-degenerate surface complex
projective structures of hyperbolic type.
\end{teo}
The cases $\kappa = 0,\ 1$ are treated in Chapters
~\ref{FGHST},~\ref{HYPE} and ~\ref{dS}; $\kappa = -1$ in
~\ref{AdS}. The proofs are rather demanding, especially in the AdS
case. The reader can find in the introduction of Chapter \ref{AdS} a
more detailed outline of this matter. As a by-product of the
classification we shall see that every maximal globally hyperbolic AdS
spacetime that has a complete Cauchy surface, does actually admit canonical
cosmological time. A delicate point in order to show that every such
spacetime actually belongs to $\Mm\Ll\Ss_{-1}$, consists in the proof
that the level surfaces contained in its past part are in fact
complete Cauchy surfaces (see Proposition \ref{memgen:adscompl:prop},
and also \cite{Ba}(2),(3) for a different proof).

\section {The other side of $\Uu^{-1}_\lambda$ - (Broken) $T$-symmetry}
\label{BROKEN:I} 
We will show (see Section \ref{T-symm}):
\begin{prop}\label{TSYM:PROP} Reversing the time orientation produces
an involution on $\Mm\Ll\Ss_{-1}$ (hence on $\Mm\Ll^\Ee/SO(2,1)$):
$$(\Uu^{-1}_\lambda,\Gamma^{-1}_\lambda) \to
(\Uu^{-1}_{\lambda^*},\Gamma^{-1}_{\lambda^*})$$ (for some $\lambda^*
\in \Mm\Ll$).
\end{prop}

This is called $T$-{\it symmetry} and is studied in Section
\ref{T-symm}.  Here is a few of its features. The isometry group of
the spacetime $\mx_{-1}$ is isomorphic to $PSL(2,\R)\times PSL(2,\R)$
(see Chapter \ref{MOD}). So $\Gamma^{-1}_\lambda$ is given by an
ordered pair $(\Gamma_L,\Gamma_R)$ of representations of $\Gamma$ with
values in $PSL(2,\R)$.  Then $\Gamma^{-1}_{\lambda^*}$ simply
corresponds to $(\Gamma_R,\Gamma_L)$.  On the other hand, we will see
while proving the classification Theorem \ref{FULL_CLASS}, that
$\Uu^{-1}_\lambda$ is determined by its {\it curve at infinity}
$C_\lambda$, contained in the ``boundary'' of $\mx_{-1}$.  This
boundary is canonically diffeomorphic to $S^1_\infty\times S^1_\infty$
and has a natural {\it causal structure}, actually depending on the
fixed orientations; $C_\lambda$ is a {\it nowhere timelike} embedded
curve homeomorphic to $S^1$. In fact, $\Uu^{-1}_\lambda$ is the
(interior of the) {\it Cauchy development} of $C_\lambda$ in
$\mx_{-1}$ (see Section \ref{MGHADS}).  The spacetime
$\Uu^{-1}_{\lambda^*}$ is obtained similarly by taking the Cauchy
development of the curve $C_\lambda^*$, that is the image of
$C_\lambda$ via the homeomorphism $(x,y)\to (y,x)$ of
$\partial\mx_{-1}$.
\smallskip

It is interesting to investigate whether distinguished sub-families of
spacetimes are closed or not under the $T$-symmetry.
\smallskip

{\it T-symmetry in the cocompact case.}  

\noindent Assume that $H=\mh^2$, and $F=\mh^2/\Gamma$ is a compact
hyperbolic surface.  It is known since \cite{M},
that in such a case $(\Gamma_L,\Gamma_R)$ is a couple of faithful
cocompact representations of the same genus of $F$, and that any such a
couple uniquely determines an AdS spacetime of this kind.  It follows
that cocompact $\Mm\Ll(\mh^2)$-spacetimes of curvature $\kappa = -1$
are closed under the $T$-symmetry.  We can also say that the initial
and {\it final} singularities of $\ \Uu^{-1}_\lambda$ have in this
case the same kind of structure. The same facts hold if $F$ is not
necessarily compact, but we confine ourselves to consider only
laminations with compact support.
\medskip

{\it Broken $T$-symmetry for general $\Mm\Ll(\mh^2)$-spacetimes.}  

\noindent We have (see Proposition \ref{mh2-case}):

\begin{prop}\label{-1mlh2}
$(\Uu^{-1}_\lambda,\Gamma^{-1}_\lambda)\in \Mm\Ll\Ss_{-1}$ is a
$\Mm\Ll(\mh^2)$-spacetime if and only if the special level surface
$\Uu_\lambda^{-1}(\pi/2)$ of the cosmological time ({\it i.e.} the future
boundary of the convex core) is a {\rm complete} Cauchy surface.
\end{prop}

We will show that $\Mm\Ll(\mh^2)$-spacetimes (even of finite co-area
but non cocompact) contained in $\Mm\Ll\Ss_{-1}$ are {\it not} closed
under the $T$-symmetry. In fact we will show that the characterizing
property of Proposition \ref{-1mlh2} is not preserved by the symmetry:
in general the level surface $\Uu^{-1}_{\lambda^*}(\pi/2)$ is not
complete (and it is only \emph{future} Cauchy). We could also say that
the initial and final singularities of $\Uu^{-1}_\lambda$ are not
necessarily of the same kind (for instance ``horizons censoring black
holes'' can arise). See the examples in Section ~\ref{3cusp} and
Remark~\ref{multiBH}.
\smallskip

It is an intriguing problem to characterize AdS
$\Mm\Ll(\mh^2)$-spacetimes (and broken $T$-symmetry) purely in terms
of the curve at infinity $C_\lambda$.  This also depends on a subtle
relationship between these spacetimes and Thurston's Earthquake Theory
\cite{Thu2}, {\it beyond the cocompact case} already depicted in
\cite{M}. In Section \ref{ge:quake} we get some partial results in
that direction. We recall that \cite{M} actually contains a new
``AdS'' proof of the ``classical'' Earthquake Theorem in the cocompact
case (see Proposition~\ref{earth:class:prop}).  On the other hand, in
~\cite{Thu2} there is a formulation of the Earthquake Theorem that
strictly generalizes the cocompact case. In Section~\ref{ge:quake}, we
study the relations between generalized earthquakes defined on
straight convex sets of $\mh^2$ and general Anti de Sitter
spacetimes. As a corollary, we will point out an ``AdS'' proof of such
a general formulation, and we also show that the holonomy of any
maximal globally hyperbolic Anti de Sitter spacetime containing a
complete Cauchy surface and with non-Abelian fundamental group is
given by a pair of discrete representations (Proposition
\ref{discrete-rep}). We study in particular the case of those achronal
curves $C$ that are {\it graphs of some homeomorphism of $S^1$} (as it
happens for instance in the cocompact case).  We show that the
boundary curve $C_\lambda$ of an ADS $\Mm\Ll$-spacetime is the graph
of a homeomorphism of $S^1$ iff the lamination $\lambda$ generates
surjective earthquakes onto $\mh^2$ (Proposition
\ref{earth:homeo:prop}). Moreover, we give examples (see \ref{earthex}
and Section \ref{3cusp}) showing that there is no logical implication
between being a $\Mm\Ll(\mh^2)$-spacetime and having $C$ graph of some
homeomorphism.
\section {$\Qq\Dd$-spacetimes}\label{QDST:I}
In Chapter \ref{QD} we develop the sector of our theory based on the
simplest flat regular domains \emph{i.e.} the future $\fut(r)$ of a
spacelike geodesic line $r$ of $\mx_0$. This is the degenerate case
when the image $H_\Uu$ of the Gauss map just consists of one geodesic
line of $\mh^2$. It is remarkable that the theory on $\fut(r)$ can be
developed in a completely explicit and self-contained way, eventually
obtaining results in complete agreement (for instance, for what
concerns the universal rescaling functions) with what we have done for
the $\Mm\Ll$-spacetimes. Combining them and $\Mm\Ll$ results we get in
particular the ultimate classification theorem, just by removing
``non-degenerate'' or ``non-elementary'' in the above statements.
\smallskip

Moreover, the theory over $\fut(r)$ extends to so called $\Qq\Dd$-{\it
spacetimes}. In contrast with the $\Mm\Ll$-spacetimes, these present
in general world-lines of {\it conical singularities}, and the
corresponding developing maps (even for the flat ones) are {\it not}
injective.  So they are also a first step towards a generalization of
the theory in the presence of ``particles'' (see also Section
\ref{particles}). The globally hyperbolic $\Qq\Dd$-spacetimes are
``generated'' by {\it meromorphic quadratic differentials} on $\Omega
=\ms^2,\ \C,\ \mh^2$ (possibly invariant with respect to the proper
action of some group of conformal automorphisms of $\Omega$).

Flat $\Qq\Dd$-spacetimes are locally modeled on $\fut(r)$, and had
been already considered in \cite{BG}(3). In particular, the {\it
quotient spacetimes} of $\fut(r)$ with compact Cauchy surface, realize
all {\it non-static} maximal globally hyperbolic flat spacetimes with
{\it toric} Cauchy surfaces.  Via canonical Wick rotation, we get the
kind of non-complete hyperbolic structures on $(S^1\times S^1)\times
\R$ that occur in Thurston's {\it Hyperbolic Dehn Filling} set up (see
\cite{Thu, BP}).  By canonical AdS rescaling of the quotient
spacetimes of $I^+(r)$ homeomorphic to $(S^1\times\mr) \times \R$, we
recover the so called {\it BTZ black holes} (see \cite{BTZ, Ca}).  In
contrast with the $\Mm\Ll$-spacetimes, these non-static quotient
$\Qq\Dd$-spacetimes have {\it real analytic} cosmological time, and
this is even a CMC time.

When $F=\mh^2/\Gamma$ is compact, there is a natural bijection between
the space $\Mm\Ll(F)$ of measured geodesic laminations and the space
$\Qq\Dd(F)$ of {\it holomorphic} quadratic differentials on $F$ (see
Remark \ref{QDcompact} of Chapter \ref{QD}, and also Section \ref{END}
below). It is remarkable that the ``same'' parameter space
$\Mm\Ll(F)\cong \Qq\Dd(F)$ gives rise to {\it different} families of
globally hyperbolic spacetimes with compact Cauchy surfaces, belonging
to the $\Mm\Ll(\mh^2)$ and $\Qq\Dd$ sectors of the theory
respectively. In fact $\Mm\Ll(\mh^2)$ or $\Qq\Dd$ spacetimes that
share the same encoding parameter have the same initial singularity.

We will also show that general $\Qq\Dd$-spacetimes can realize
arbitrarily complicated topologies and causal structures.

So, although the basic domain $\fut(r)$ is extremely simple, the
resulting sector of the theory is far from being trivial.

\section {Along rays of spacetimes}\label{MLRAY:I}
Given $(\lambda, \Gamma)$, $\lambda =(H,\Ll,\mu)$, $\mu$ being the
transverse measure, we can consider the {\it ray} $(t\lambda,
\Gamma)$, $t\lambda=(H,\Ll,t\mu)$, $t\in [0,+\infty)$, where it is
natural to stipulate that for $t=0$ we have the lamination just
supported by the geodesic boundary of $H$. So we have the
corresponding $1$-parameter families of spacetimes
$(\Uu^\kappa_{t\lambda },\Gamma^\kappa_{t\lambda })$ and hyperbolic
manifolds $M_{t\lambda}$, emanating from the static case at $t=0$.
The study of these families (made in Section \ref{der}, together with
further complements) gives us interesting information about the Wick
rotation-rescaling mechanism.

We study the ``derivatives'' at $t=0$ of the spacetimes
$\Uu^\kappa_{t\lambda }$, and of holonomies and ``spectra'' of the
quotient spacetimes.  In particular, let us denote by
$\frac{1}{t}\Uu^\kappa_{t\lambda }$ the spacetime obtained by
rescaling the Lorentzian metric of $\Uu^\kappa_{t\lambda }$ by the
constant factor $1/t^2$. So $\frac{1}{t}\Uu^\kappa_{t\lambda }$ has
constant curvature $\kappa_t=t^2\kappa$.  Then, we shall prove
(using a suitable notion of convergence) that (see \ref{derivateST})
$$ \lim_{t\to 0}\ \ \frac{1}{t}\Uu^\kappa_{t\lambda } = \Uu^0_\lambda
 \ .$$ We stress that this convergence is at the level of
 Teichm\"uller-like classes; here working up to reparametrization
 becomes important.

For every $\kappa$, the space $\Mm\Ll\Ss_\kappa(S)$ of maximal
 globally hyperbolic spacetimes of constant curvature $\kappa$ and
 fixed topological support $S\times\mr$ has a natural ``Fuchsian''
 locus, corresponding to $H$-static spacetimes. In a sense, thanks to
 the above limit, the space $\Mm\Ll\Ss_0(S)$ could be considered as
 the normal bundle of this Fuchsian locus. In the case of compact $S$,
 an algebraic counterpart of this fact is recalled in Section
 \ref{END}.  Then the canonical rescaling produces a map
 $\Mm\Ll\Ss_0(S)\rightarrow\Mm\Ll\Ss_\kappa(S)$ that could be regarded
 as a sort of exponential map. As a by-product, we find some formulae
 relating interesting classes of spectra associated to each
 $\Uu^\kappa_\lambda$. In fact these formulae are proved in a
 different context also in~\cite{Go-Ma}.

In Section~\ref{more_cocompact} some specific applications in
the case of compact Cauchy surfaces are given. We consider the family
$\Uu^{-1}_{(t\lambda)^*}$ as in Section \ref{BROKEN:I}.  In such a case
(see for instance Section \ref{END} below), the set of
$\Gamma$-invariant measured laminations has an $\R$-linear structure,
and it makes sense to consider $-\lambda$. Then we shall show (see Section 
\ref{derivate*})
$$ \lim_{t\to 0}\ \ \frac{1}{t} \Uu^{-1}_{(t\lambda_t)^*} =
\Uu^0_{-\lambda} \ .$$ 

Let $Q_t$ be a smooth family of homeomorphic quasi-Fuchsian manifolds
such that $Q_0=\mh^3/\Gamma$ is Fuchsian and $\lambda_t,\lambda^*_t$
be the bending loci of the boundary of the convex core of
$Q_t$. Bonahon \cite{Bon}(2) proved that the family of measured
geodesic laminations $\lambda_t/t$ and $\lambda^*_t/t$ converge to
geodesic laminations $\lambda_0, \lambda^*_0$ such that
$\lambda^*_0=-\lambda_0$ with respect to the linear structure of
$\Mm\Ll(\mh^2/\Gamma)$.  Notice that this result is strongly similar
to that one we get in Anti de Sitter setting. Roughly speaking, we can
conclude that bending in Anti de Sitter space is the same as bending
in hyperbolic space at infinitesimal level (that is, the boundary
components of the convex core are the same).  On the other hand on the
large scale the bending behaviour in the two contexts is very
different (see Section ~\ref{big-t}).

Finally, giving a partial answer to a question
of \cite{M}, we establish formulae relating the volume of the past of
a given compact level surface of the canonical time, its area, and the
{\it length} of the associated measured geodesic lamination on $F$. We
recover in the present set up a simple proof of a continuity property
of such a length function (see Section \ref{volume}).

Similarly, for $\Qq\Dd$-spacetimes we will consider the Wick
rotation-rescaling behaviour along lines of quadratic differentials
$t^2\omega$.

\section{QFT and ending spacetimes}\label{END}
This Section contains a rather long expository digression. Shortly,
this could be summarized by the following sentences:

{\it The Wick rotation-rescaling theory applies on the ends of
geometrically finite hyperbolic $3$-manifolds. Hence these manifolds
realize concrete interactions of their globally hyperbolic ending
spacetimes of constant curvature. This provides natural geometric
instances of morphisms of a $(2+1)$ bordism category suited to support
a quantum field theory pertinent to 3D gravity. Moreover, the finite
volume of the slabs of the AdS spacetimes that support the Wick
rotations ({\rm see Section~\ref{more_cocompact}}) together with the
volume of the hyperbolic convex cores, furnish classical
``amplitudes'' of these interactions.  }
\smallskip

\noindent If satisfied with this, the reader can skip to the next
Chapter. However, we believe that the digression shall display
important underlying ideas and a background that have motivated this
work.
\smallskip

We want to informally depict a few features of a {\it quantum field
theory} QFT that would be pertinent to 3D gravity, by taking
inspiration from \cite{W}, and having as model the current
formalizations of {\it topological quantum field theories}.

\paragraph{Bits of axiomatic QFT.}
Following \cite{At, Tur}, by $(2+1)$ QFT we mean any functor,
satisfying a demanding pattern of axioms, from a $(2+1)$-bordism
category to the tensorial category of finite dimensional complex
linear spaces.
\smallskip

The {\it objects} of the bordism category are (possibly empty) finite
union of suitably {\it marked} connected compact closed oriented 
surfaces $\Sigma$.  Every marking includes (at least) an oriented
parametrization $\phi: S \to \Sigma$, by some {\it base} surface $S$.

Every {\it morphism} is a compact oriented $3$-manifold $Y$ with marked
and bipartited boundary components; hence $Y$ realizes a ``tunnel'',
a transition from its {\it input} boundary object $\partial_-$ towards
the {\it output} one $\partial_+$.

A QFT functor associates to every object $\alpha$ a complex linear
space $V(\alpha)$ and to every marked bordism $Y$ a linear map $Z_Y:
V(\partial_-) \to V(\partial_+)$, the {\it tunneling amplitude}.  This
is functorial with respect to composition of bordisms on one side, and
usual composition of linear maps on the other. Amplitudes can be
considered as a generalization of time evolution operators (where $Y$
is a cylinder).

If an object $\alpha$ is union of connected components $\Sigma_j$'s,
then $V(\alpha)$ is the tensor product of the
$V(\Sigma_j)$'s. $V(-\Sigma) = V^*(\Sigma)$, where $-\Sigma$ denotes
the surface with the opposite orientation, $V^*$ is the dual of $V$.
$V(\emptyset)=\C$; hence if the boundary of $Y$ is empty, then $Z_Y\in
\C$ is a scalar. This numerical invariant of the $3$-manifold $Y$ is
usually called its {\it partition functions}. If $\partial_-
=\emptyset$, then $Z_Y$ is a vector in $V(\partial_+)$ (the {\it
vacuum state} of $Y$); if $\partial_+ =\emptyset$, then $Z_Y$ is a
functional on $V(\partial_-)$.

The amplitudes are sensitive to the action of the {\it mapping
class groups} on the markings.

A crucial feature of any QFT is that we can express any amplitude by
using (infinitely many) different decompositions of the given bordism,
associated for instance to different Morse functions for the triple
$(Y,\partial_-,\partial_+)$.

Possibly a QFT is {\it not purely topological}, \emph{i.e.} the marked
$3$-manifolds can carry more structure and we deal with a bordism category
of such equipped manifolds. Both $V$ and $Z$ possibly depend also on
the additional structure.

Pertinence to 3D gravity starts arising by specializing the additional
structure on each marked surface $\Sigma$, and showing later that
classical gravity naturally furnishes a wide set of bordisms in the
appropriate category.

\paragraph{Matter-free Witten phase space.}
At first, it seems reasonable to select a sector of 3D gravity by
fixing the signature (Lorentzian or Euclidean) of the $3$-dimensional
metrics, and the sign $\sigma(\Lambda)$ of the cosmological
constant. The basic idea is to give $\Sigma$ the structure of a {\it
spacelike surface} embedded in some {\it universe} $U$ of constant
curvature $\kappa = \sigma(\Lambda)$ (for the Euclidean signature we
just consider embedded surfaces).

Let $\mx$ be any model of constant curvature geometry, and $\Gg$
denote its group of isometries. In \cite{W} a reformulation of the
corresponding sector of pure classical 3D gravity is elaborated as a
theory with {\it Chern-Simons
action}, for which the relevant fields are the connections on principal
$\Gg$-bundles on $3$-manifolds, up to gauge transformations, rather
than the metrics.  The ``classical phase space'' finally associated to every
connected base surface $S$ is the space of {\it flat connections} on
principal $\Gg$-bundles over $S$, up to {\it gauge}
transformations. Equivalently, this consists of the {\it space of
representations} $$ {\rm Hom}(\pi_1(S),\Gg)$$ on which the group $\Gg$
acts by conjugation. To each spacelike surface $\Sigma$ in some
universe $U$ as above, one associates a holonomy representation of its
{\it domain of dependence} $D(\Sigma)$ in $Y$ (see
\cite{H-E}). $D(\Sigma)$ is globally hyperbolic, and $\Sigma$ is a
Cauchy surface, hence $D(\Sigma)$ is homeomorphic to $\Sigma \times
\R$, and $\pi_1(D(\Sigma))=\pi_1(S)$. For the Euclidean signature we
just take an open tubular neighbourhood of $\Sigma$ in $U$. More
generally, as additional structure on a marked surface $\Sigma$, we
just take any element of such a ``Witten phase space''.

Natural instances of morphisms in the corresponding bordism category
arise in classical gravity. They should be compact submanifolds with
boundary $Y$ of some universe $U$ of constant curvature $\kappa$, the
boundary being made by (bipartited) spacelike surfaces.  The holonomy
of the whole spacetime ${\rm Int}(Y)$, induces the one of the
``ending'' globally hyperbolic spacetimes.  Such a transition would
realize, in particular, a change of topology from the set of input
spacetimes towards the set of output ones. The prize for it consists
in a severe weakening of the global causal structure of $Y$: normally
${\rm Int}(Y)$ contains closed timelike curves.

\paragraph{Mess-Thurston  phase space.}
Mess/Scannell and Thurston parameterizations in the case of compact
surfaces of genus $g\geq 2$ (if $g=1$ see Section \ref{QDST:I} above)
leads to a somewhat different and unified way to look at the
classical phase space. That is we would look for a QFT built on the
``universal'' parameter space $\Tt_g\times \Mm\Ll_g$, and that should
deal simultaneously with spacetimes of arbitrary constant curvature
and complex projective surfaces.
 
Some comments are in order to establish some point of contact with the
former paragraph.  We use the notations of the previous Sections.

(1) $\Tt_g\times \Mm\Ll_g$ should be considered as a {\it trivialized}
fiber bundle $\Bb \to \Tt_g$ over the classical Teichm\"uller space
$\Tt_g$ of hyperbolic structures on $S$. The fiber over any $F\in
\Tt_g$ is $\Mm\Ll(F)$. Given $F$ and $F'$ in $\Tt_g$, we fix the
canonical topological bijection between $\Mm\Ll(F)$ and $\Mm\Ll(F')$
that identifies two laminations if and only if they share the same
``marked spectrum'' of measures. This trivialization respects the {\it
  ray structure} considered in Subsection \ref{MLRAY:I}. This induces
a trivialized fiber bundle $\Tt_g \times \mP^+(\Mm\Ll_g)\cong
\Tt_g\times S^{6g-7}$, where each fiber is a copy of the {\it
  Thurston's boundary of $\Tt_g$}.

(2) {\it Fuchsian slices.}  The relevant Lorentzian isometry groups
are: $\Gg_0 = ISO^+(2,1)$ \emph{i.e.} the Poincar\'e group of affine
isometries of the Minkowski space $\mx_0$ that preserve the
orientations; $\Gg_{-1} = PSL(2,\R)\times PSL(2,\R)$; $\Gg_{1}=
PSL(2,\C)$ (see Chapter \ref{MOD}). For every $\kappa$, there is a
canonical embedding of $PSL(2,\R)$ into $\Gg_\kappa$. The subgroup
$SO^+(2,1)\subset \ISO^+(2,1)$ of linear isometries is canonically
isomorphic to $PSL(2,\R)\cong {\ISO}^+(\mh^2)$ (by using both the
hyperboloid and half-plane models of the hyperbolic plane).  For
$\kappa=-1$ we have the {\it diagonal} embedding, for $\kappa=1$ we
take the {\it real part} of $PSL(2,\C)$.  Denote by $\Ff\Rr_g$ the
subset of ${\rm Hom}(\pi_1(S),PSL(2,\R))/PSL(2,\R)$ of Fuchsian
representations.  By using the above embeddings, this determines the
{\it Fuchsian slice} of each ${\rm Hom}(\pi_1(S),\Gg_k)/\Gg_\kappa$.
$\Tt_g$ can be identified with (a connected component of) $\Ff\Rr_g$,
so that each Fuchsian slice corresponds to the ``$0$-section'' of the
bundle $\Tt_g \times \Mm\Ll_g$ which parameterizes the static
spacetimes. So, for every $\kappa$, any ray in $\Mm\Ll(F)$ over a
given $F\in \Tt_g$ can be consider as a $1$-parameter family of {\it
  deformations} of the static spacetime associated to $F$.

(3) {\it Coincidence of infinitesimal deformations.} Let us denote by
$\gG_\kappa$ the Lie algebra of $\Gg_\kappa$. By using the adjoint
representation restricted to the embedded copy of $PSL(2,\R)$, each
$\gG_\kappa$ can be considered as a $PSL(2,\R)$-module. It is a fact
that these modules are canonically isomorphic to each other. This is
immediate for $\kappa=\pm 1$. For $\kappa=0$, this follows from the
canonical linear isometry between $\sG\lG(2,\mr)$, endowed with its
killing form, and the Minkowski space $(\mr^3,\E{\cdot}{\cdot})$.

For every Fuchsian group $\Gamma$ as above, the ``infinitesimal
deformations'' of $\Gamma$ in $\Gg_\kappa$, are parametrized by
$H^1(\Gamma,\gG_\kappa)$, through such a $PSL(2,\R)$-module
structure. Thus we can say that the infinitesimal deformations of
$\Gamma$ considered as a static spacetime of {\it arbitrary} constant
curvature, as well as a Fuchsian projective structure on $S$ actually
coincide. Canonical
Wick rotation theory realizes, in a
sense, a full ``integration'' of such an infinitesimal coincidence.

(4) {\it Holonomy pregnancy.} Mess's work for $\kappa=0$ includes an
identification of $\Tt_g \times \Mm\Ll_g$ with the subset of ${\rm
  Hom}(\pi_1(S),\Gg_0)/\Gg_0$ made by the representations with
Fuchsian linear part; for $\kappa = -1$, with the subset $\Ff\Rr_g
\times \Ff\Rr_g$ of ${\rm Hom}(\pi_1(S),\Gg_{-1})/\Gg_{-1}$ (via a
subtle Lorentzian revisitation of Thurston {\it earthquake}
theory). This means, in particular, that for $\kappa =0, -1$ the
maximal hyperbolic spacetimes with compact Cauchy surface are
determined by their holonomy, so that the Witten phase spaces properly
contain parameter spaces of these spacetimes. On the other hand de
Sitter spacetimes and compact complex projective surfaces are no
longer determined by the respective holonomies. It happens in fact
that on a same ray, spacetimes (projective structures) with injective
({\it quasi-Fuchsian}) developing maps (corresponding to small values
of the ray parameter), and others having non-injective and surjective
developing maps (for big values of the parameter) share the same
holonomy - see \cite{Go}(1)).

(5) {\it Linear structures on $\Mm\Ll(F)$.}  Recall that $\Tt_g$ is
homeomorphic to the open ball $B^{6g-6}$ and the fiber $\Mm\Ll_g$ to
$\R^{6g-6}$. In fact, the bundle $\Bb \to \Tt_g$ can be endowed with
natural vector bundle structures that identify it as the {\it
  cotangent bundle} $T^*(\Tt_g)$ (hence we have a honest classical
phase space, with $\Tt_g$ as {\it configuration space}).  For example,
given $F=\mh^2/\Gamma$ as above, one establishes a natural bijection
between $\Mm\Ll(F)$ and $H^1(\Gamma,\R^3)$ (where $\R^3$ is here the
subgroup of translations of $ISO^+(2,1)$), and this gives each fiber
of $\Bb$ the required linear structure. This induces the usual ray
structure. Another linear structure is through complex analysis. The
fiber of the complex cotangent bundle $T^*(\Tt_g)$ over $F$ is
identified with the space $\Qq\Dd(F)$ of holomorphic quadratic
differentials on $F$, considered now as a Riemann surface. On the
other hand, there is a natural bijection between $\Qq\Dd(F)$ and
$\Mm\Ll(F)$ (through horizontal measured foliations of quadratic
differentials - see \cite{Ke} and also Chapter \ref{QD}).

However, the canonical (topological) trivialization mentioned in (1)
above is {\it not} compatible with any such natural vector bundle
structures.

In \cite{Mon}, the flat case $\kappa = 0$, for every genus $g \geq 1$,
is treated as a Hamiltonian system over $\Tt_g$ (considered here as a
space of Riemann surface structures on $S$), in such a way that
$T^*(\Tt_g)$ (with its complex vector bundle structure) is the phase
space of this system. In this analytic approach the relevant global
time fills each spacetime by the evolution of CMC Cauchy surfaces.
Recall that $T^*(\Tt_g)$ is also the phase space for a family of
globally hyperbolic $\Qq\Dd$-spacetimes. Moncrief's and $\Qq\Dd$
theories coincide exactly when $g=1$.
\medskip
 
Keeping both the basic idea of Witten's approach of dealing with flat
connections as fundamental fields, and the unifying viewpoint
underlying the Wick rotation-rescaling theory, we would adopt the 
space of representations
$$ {\rm Hom}(\pi_1(S),PSL(2,\C))$$ (on which the group $PSL(2,\C)$
acts by conjugation) as a reasonable ``universal'' phase space,
although we are aware - point (4) - that it is strictly weaker than
the universal parameter space $\Tt_g\times \Mm\Ll_g$. A nice fact is
that hyperbolic $3$-manifolds furnish a wide natural class of
morphisms in the corresponding bordism category.

\paragraph{Hyperbolic  bordisms.}
For the notions of hyperbolic geometry used in this paragraph we shall
refer, for instance, to~\cite{Thu,Ka, Ot}.

Recall that a (complete) hyperbolic $3$-manifold $Y=\mh^3/G$ (where
$G$ is a Kleinian group isomorphic to $\pi_1(Y)$) is {\it
topologically tame} if it is homeomorphic to the interior ${\rm int}\
M$ of a compact manifold $M$. Roughly speaking, each boundary
component of $M$ corresponds to an {\it end} of $Y$. By using the
holonomy of $Y$, we can associate to each boundary component a point
in our phase space.  Any bipartition of the boundary components gives
rise to a bordism in the appropriate bordism category. This can be
considered as a transition from the set of input towards the set of
output ends.

In general the asymptotic geometry on the ends is a rather subtle
stuff.  This is much simpler for the important subclass of
geometrically finite manifolds. Recall that $Y$ is {\it geometrically
finite} if any $\epsilon$-neighbourhood of its {\it convex core}
$C(Y)=C(G)$ is of {\it finite volume}. For simplicity, we also assume
that there is no {\it accidental parabolic} in $G$.

Let us assume furthermore that $C(Y)$ is {\it compact}, and that $Y$
is not compact. Then:

- The  group $G$ does not contain parabolic elements.

- $C(Y)$ is a compact manifold with non-empty boundary; $Y$ is
homeomorphic to the interior of $C(Y)$; for every end $E(S)$
corresponding to a given boundary component $S$ of $C(Y)$, there is a
unique connected component of $Y\setminus C(Y)$ (also denoted $E(S)$),
characterized by the fact that $\partial \overline{E}(S)=S$. $E(S)$ is
homeomorphic to $S\times (0,+\infty)$. The asymptotic behavior of (the
universal covering of) $E(S)$ at the boundary $S^2_\infty$ of $\mh^3$,
determines a {\it quasi-Fuchsian} complex projective structure on
$S$. In fact $E(S)$ is its $H$-hull, and $S$ inherit also the
intrinsic hyperbolic structure for being the hyperbolic boundary of
$E(S)$.

- We can apply on each $E(S)$ the Wick rotation-rescaling theory
developed for $\Mm\Ll(\mh^2)$-spacetimes. Thus we can convert $E(S)$
into the future domain of dependence of the level $1$ surface of the
cosmological time for a determined maximal globally hyperbolic {\it
  flat} spacetime. By canonical rescaling we can convert it in a
suitable slab of such a spacetime of constant curvature $\kappa=-1$.
\begin{remark}\label{TUNNEL} 
{\rm Notice that the canonical Wick rotation on each end $E(S)$ {\it
cuts off} the initial singularity of the associated ending spacetime.
In a $(2+1)$ set up, this is clearly reminiscent of the basic
geometric idea underlying the so called Hartle-Hawking {\it
no-singularity proposal} and the related notion of {\it real tunneling
geometries} (see \cite{Gi, Gi2}). However, note that, in our
situation, the surface $S$ is not in general totally geodesic in $Y$.
This happens exactly when the associated ending spacetime is static,
or, equivalently, if we require furthermore that the canonical
cosmological time is twice differentiable (in fact real analytic)
everywhere, or that the initial singularity consists of one point. So
these (equivalent) assumptions are not mild at all as they ``select''
very special configurations.}
\end{remark}

\paragraph{Including particles.} In contrast with the matter-free case,
the parameter spaces of globally hyperbolic spacetimes of constant
curvature, with compact Cauchy surface, and {\it coupled with
particles} are, as far as we know, not yet well understood (see for
instance \cite{BG}(2), and also Section \ref{particles}). Moreover, even the
reformulation of gravity as a Chern-Simons theory should be not so
faithful in this case (see \cite{Mat}). Nevertheless, by slightly
generalizing the above discussion, we can propose a meaningful ``phase
space'' also in this case, and a relative bordism category.

We consider now compact surfaces $S$ as above with a fixed finite
non-empty set $V$ of {\it marked points} on it. Set $S'=S\setminus V$.
Then we consider the space of representations
$${\rm Hom}(\pi_1(S'),PSL(2,\C)) \ .$$ In building the bordism category
we have to enhance the set up by including $1$-dimensional tangles
$L$, properly embedded in $M$, such that $\partial L$ consists of the
union of distinguished points on the components of $\partial M$.  Set
$M'= M\setminus L$, then an additional structure on $M$ is just an
element of $ {\rm Hom}(\pi_1(M'),PSL(2,\C)) $. The tangle components
mimic the particle world lines. 

It is useful to ``stratify'' ${\rm Hom}(\pi_1(S'),PSL(2,\C))$, by
specializing the {\it peripheral representations} on the loops
surrounding all distinguished points.  For instance we get the
``totally parabolic peripheral stratum'' by imposing that they are all
of {\it parabolic type}.  We can also specialize elliptic peripheral
representations, by fixing the cone angles of the local quotient
surfaces, and so on. Note that the stratum where these representations
are all {\it trivial} is not really equivalent to the matter-free
case, because we keep track of the marked points; in particular we act
with a different mapping class group.

Let $Y=\mh^3/G$ be again a geometrically finite hyperbolic
$3$-manifold.  Assume now that $G$ contains some parabolic
element. Assume also, for simplicity, that $G$ does not contain
subgroups isomorphic to $\Z \oplus \Z$, \emph{i.e.} that $Y$ has no
{\it toric cusps}.  Then:

- Any boundary component $\Sigma$ of the convex core $C(Y)$ either is
compact or is homeomorphic to some $S' = S\setminus V$ as above.  In
any case, $\Sigma$ inherits an intrinsic complete hyperbolic structure
of {\it finite area}. The holonomy of $Y$ at every loop surrounding
the punctures of $S'$ is of parabolic type, and every parabolic
element of $G$ arises in this way (up to conjugation).

- There are a compact $3$-manifold $M$, with non-empty boundary
$\partial M$, and a properly embedded tangle $L$ made by arcs
($\partial L = L\cap \partial M$) such that $Y$ is homeomorphic to
$M\setminus (\partial M \cup L)$. Moreover, there is a natural
identification between the boundary components of $M\setminus L$ and
the ones of $C(Y)$. Note that $M$ is obtained by adding $2$-handles to
the boundary of a suitable compact manifold $W$ (the co-core of
the handles just being the arcs of $L$). Then $Y$ is homeomorphic to
$\mathring W$. Usually, the so-called {\it geometric ends} of $Y$
correspond to the connected components of $\partial W$.  On the other
hand, we can consider also the ``ends'' of $Y \cong M \setminus
(\partial M \cup L)$, each one homeomorphic to $\Sigma \times \R$, for
some boundary component of $C(Y)$. We stipulate here to take these
last as ends of $Y$. Then, the asymptotic behavior of (the universal
covering of) an end of $Y$ at the boundary of $\mh^3$, determines a
{\it quasi-Fuchsian} projective structure on the corresponding surface
$\Sigma$.

- We can apply on each end of $Y$ the Wick rotation-rescaling theory
developed for $\Mm\Ll(\mh^2)$-spacetimes. Note that in this case the
involved geodesic laminations have always compact support.
\smallskip

These considerations can be enhanced by allowing also toric cusps, and
replacing $Y$ by $Y'$, obtained by removing from $Y$ the cores of the
Margulis tubes of a given {\it thick-thin} decomposition (for
instance, the one canonically associated to the Margulis constant). In
general $Y'$ is no longer complete, and possibly presents further
ends, homeomorphic to $(S^1\times S^1)\times \R$, either corresponding
to toric cusps or tubes of the thin part.  On these ends the
$\Qq\Dd$-Wick rotation-rescaling sector does apply.
\smallskip

\noindent Summarizing, geometrically finite hyperbolic $3$-manifolds
furnish natural examples of bordisms in the pertinent category, and by
Wick rotation-rescaling we show that they eventually realize geometric
transitions of {\it the ending spacetimes of constant curvature}.
\smallskip

We have considered geometrically finite manifolds for the sake of
simplicity. In particular, the $\Mm\Ll(\mh^2)$ sector of the theory
suffices in this case. The whole theory should allow to treat in the
same spirit more general tame manifolds.
\smallskip

\paragraph{Some realized ``exact" QFT.}
We say that a QFT is {\it exact} if the amplitudes are expressed by
exact formulae, based on some effective encoding of the (equipped) 
marked bordisms.

It has been argued (see for instance \cite{F-K}) that Turaev-Viro
state sum invariants $TV_q(.)$ of compact closed $3$-manifolds
\cite{T-V} are partition functions of a QFT pertinent to 3D gravity
with {\it Euclidean signature and positive cosmological constant}. In
fact they can be considered as a countable family of {\it
regularizations}, obtained by using the quantum groups
$U_q(sl(2,\C))$, of the Ponzano-Regge calculus based an the classical
unitary group $SU_2$.  In fact $TV_q(.)=|W_q(.)|^2$, where $W_q(.)$
denotes the Witten-Reshetikhin-Turaev invariant and a complete
topological QFT has been developed that embodies these last partition
functions (see \cite{Tur}).
\smallskip

In the papers \cite{BB} a so called {\it quantum hyperbolic field
theory} QHFT was developed by using the bordism category (including
particles) that we have depicted above. In fact this is a countable
family of exact QFT's, indexed by odd integers $N\geq 1$. The building
blocks are the so called {\it matrix dilogarithms}, which are
automorphisms of $\C^N\otimes \C^N$, associated to oriented hyperbolic
ideal tetrahedra encoded by their cross-ratio moduli and equipped with
an additional decoration. They satisfy fundamental {\it five-term
identities} that correspond to all decorated versions of the basic
$2\to 3$ move on $3$-dimensional triangulations. For $N=1$, it is
derived from the classical ``commutative'' Rogers dilogarithm; for
$N>1$ they are derived from the cyclic representations theory of a
Borel quantum subalgebra of $U_\zeta(sl(2,\C))$, where $\zeta
=\exp(2i\pi/N)$. As for $TV_q$, the exact amplitude formulae are state
sums, supported by manifold decorated triangulations, satisfying non
trivial {\it global constraints}.  A very particular case of partition
functions equal Kashaev's invariants of links in the $3$-sphere
\cite{K1}, later identified by Murakami-Murakami \cite{MM} as special
instances of colored Jones invariants.

\paragraph{Classical invariants and ``Volume Conjectures''.}
A fundamental test of pertinence to 3D gravity of the above exact QFT
should consist in recovering, by suitable ``asymptotic expansions'',
some fundamental classical geometric invariants. One can look in this
spirit at the activity about asymptotic expansions of quantum
invariants of knots and $3$-manifolds (see for instance Chapter 7 of
\cite{Oht} for an account of this matter).

One usually refers to such challenging general problem about the
asymptotic behaviour of QHFT partition functions as ``Volume
Conjectures'', adopting in general the name of Kashaev's germinal one
for the special case of hyperbolic knots $K$ in $S^3$ \cite{K2, MM}.
\smallskip

{\it QHFT classical state sums.}  A nice feature of QHFT is that as
partition functions of the classical member of the family ($N=1$) one
computes 
$${\rm Vol}(.) + i{\rm CS}(.)$$ Vol$(.)$ and CS$(.)$ being
respectively the volume and the Chern-Simons invariant of both
hyperbolic $3$-manifolds of finite volume (compact and cusped ones)
and of principal flat $PSL(2,\C)$-bundles on compact closed manifolds
(see \cite{N}, \cite{}).
\smallskip

{\it Volume rigidity.} If $W$ is a compact closed hyperbolic
$3$-manifold, its volume coincides with the one of its holonomy, and
this last is the {\it unique maximum} point of the volume function
defined on ${\rm Hom}(\pi_1(W),PSL(2,\C))/PSL(2,\C)$ (see
\cite{Dun}). So this ``volume rigidity'' looks like a geometric
realization of the minimal action principle in 3D Euclidean gravity
with negative cosmological constant; Vol$(.) + i$CS$(.)$ can be
regarded as a natural complexification of that action, very close
to the spirit of \cite{W}.

{\it QHFT ``Volume Conjectures''.} We can find in \cite{BB} some
instances of ``Volume Conjecture'' (involving cusp manifolds,
hyperbolic Dehn filling and the convergence to cusp manifolds ....)
which, roughly speaking, predict that: 
\smallskip

{\it ``For $N\to \infty$, suitable ``quantum'' partition functions
asymptotically recover classical ones''.}
\smallskip

\noindent These conjectures would be ``geometrically well motivated''
because both classical and quantum state sums are basically computed
on the very same geometric support.  As far as we know, all known
instances of Volume Conjecture are open.
\smallskip

{\it Wick rotation-recaling and classical tunneling invariants.}  Let
$Y$ be a geometrically finite hyperbolic $3$-manifold as above.  The
finite volume of the convex core $C(Y)$ is a basic invariant of
$Y$. Let $Y$ be non-compact, and assume that it is the support of a
tunnel from a given input set of ends $E_-$ toward an output set
$E_+$. The Wick rotation-rescaling theory suggests other natural
invariants that are sensitive to the tunneling and not only to the
geometry of the support $3$-manifold: if $E(S)\in E_{-}$ we can
associate to it the volume of the slab of the AdS $\Mm\Ll(\mh^2)$
ending spacetime $\Uu^{-1}$ that supports the Wick rotation; if
$E(S)\in E_{+}$ we do the same by replacing $\Uu^{-1}$ with
$(\Uu^{-1})^*$, via the $T$-symmetry.  Recall that the $T$-symmetry
holds because we are dealing with laminations with compact support. In
fact these slabs are of finite volume, and volume formulas are derived
at the end of Section \ref{more_cocompact}.


\newpage{\pagestyle{empty}\cleardoublepage} 
\chapter{Geometry models}\label{MOD}

\section{Generalities on $(X,G)$-structures}\label{XG}
A nice feature of 3D gravity is that we have very explicit (local)
{\it models} for the manifolds of constant curvature, which we usually
normalize to be $\kappa=0$ or $\kappa=\pm 1$. 

In the Riemannian case, these are the models $\mr^3$, $\ms^3$ and
$\mh^3$ of the fundamental $3$-dimensional isotropic geometries: {\it
flat}, {\it spherical} and {\it hyperbolic}, respectively.  These are
the central objects of Thurston's geometrization program, which has
dominated the $3$-dimensional geometry and topology on the last
decades.  We will deal mostly with hyperbolic geometry. We refer, for
instance, to \cite{Thu, BP, Rat}, and we assume that the reader is
familiar with the usual concrete models of the hyperbolic plane
$\mh^2$ and space $\mh^3$. We just note that the {\it hyperboloid}
model of $\mh^n$, embedded as a spacelike hypersurface in the
Minkowski space $\mm^{n+1}$, establishes an immediate relationship
with Lorentzian geometry. The restriction to $\mh^n$ of the natural
projection of $\mm^{n+1}\setminus \{0\} $ onto the projective space,
gives the {\it Klein projective} model of $\mh^n$.  On the other hand, the
{\it Poincar\'e disk (or half-space)} model of $\mh^3$ concretely
shows its natural boundary at infinity $S^2_\infty = \mc\Pm^1$, and the
identification of the isometry group ${\rm Isom}^+(\mh^3)$ with the
group $PSL(2,\C)$ of projective transformations of the Riemann sphere
($\mh^3$ is oriented in such a way that the boundary orientation
coincides with the complex one).

Also in the Lorentzian case, for every $\kappa = 0,\pm1$, we will
present an {\it isotropic model}, that is a Lorentzian manifold, of
constant curvature equal to $\kappa$, such that the isometry group
acts transitively on it and the stabilizer of a point is the group
$O(2,1)$. This is denoted by $\mx_\kappa$, and called the $3$-dimensional
{\it Minkowski}, {\it de Sitter} and {\it Anti de Sitter} spacetime,
respectively. 

An interesting property of an isotropic manifold $\mx$ is that every
isometry between two open sets of $\mx$ extends to an isometry of the
whole $\mx$. Moreover, each of our concrete models is real analytic and
the isometry groups is made by analytic automorphisms.  Thus, we can
adopt the very convenient technology of $(\X,\Gg)$-{\it manifolds},
\emph{i.e.}  manifolds equipped with (maximal) {\it special atlas}
(see e.g. \cite{Thu} or {Chapter B of \cite{BP} for more details). We
recall that $\X$ denotes the model (real analytic) manifold, $\Gg$ is
a group of analytic automorphisms of $\X$ (which possibly preserve the
orientation). A special atlas has charts with values onto open sets of
$\X$, and any change of charts is given by the restriction to each
connected component of its domain of definition of some element $g\in
\Gg$. For every $(\X,\Gg)$-manifold $M$, a very general analytic
continuation-like construction, gives pairs $(D,h)$, where $ D:
\tilde{M} \to \X$ is a {\it developing map} defined on the universal
covering of $M$, $ h: \pi_1(M)\to \Gg$ is a {\it holonomy
representation} of the fundamental group of $M$.  Moreover, we can
assume that $D$ and $h$ are {\it compatible}, that is, for every
$\gamma \in \pi_1(M)$ we have $ D(\gamma (x))= h(\gamma)(D(x))$, where
we consider the natural action of the fundamental group on
$\tilde{M}$, and the action of $\Gg$ on $\X$, respectively.  The map
$D$ is a local isomorphism (a local isometry in our concrete cases),
and it is unique up to post-composition by elements $g\in \Gg$.  The
holonomy representation $h$ is unique up to conjugation by $g\in \Gg$. An
$(\X,\Gg)$-structure on $M$ lifts to a locally isomorphic structure on
$\tilde{M}$, and these share the same developing maps. In many
situations it is convenient to consider this lifted structure, keeping
track of the action of $\pi_1(M)$ on $\tilde{M}$.  

With this terminology, oriented (non necessarily complete) hyperbolic
$3$- manifolds coincide with $(\mh^3,{\rm Isom}^+(\mh^3))$-manifolds,
as well as any $2+1$ spacetime of constant curvature $\kappa$ is just
a $(\mx_\kappa, {\rm Isom}^+(\mx_\kappa))$-manifold.

We will apply this technology also to deal with (complex) {\it
projective structures} on surfaces, that are by definition
$(S^2_\infty, PSL(2,\C))$-manifold structures. In every case, the data
$(D, h)$ determine the isomorphism class of such manifolds.

Although geometry models are a very classical matter, for the
convenience of the reader we are going to treat somewhat diffusely the
Lorentzian models, giving more details for the perhaps less commonly
familiar de Sitter and anti de Sitter ones.  At the end of the Chapter
we will also discuss complex projective structures on surfaces,
including the notions of {\it $H$-hull}, {\it canonical
stratification} and {\it Thurston metric}.

\section{Minkowski space}\label{MinMod}
We denote by $\mx_0$ the $3$-dimensional Minkowski space $\mm^3$, that
is $\mr^3$ endowed with the flat Lorentzian metric
\[
h_0=-\mathrm dx_0^2+\mathrm dx_1^2+\mathrm dx_2^2 \ .
\]
Geodesics are straight lines and totally geodesic planes are affine
planes.

The orthonormal frame
\[
   \frac{\partial\,}{\partial x_0},\,\frac{\partial\,}{\partial
     x_1},\,\frac{\partial\,}{\partial x_2}
\]
gives rise to an identification of every tangent space $T_x\mx_0$ with
$\mr^3$ provided with the Minkowskian form
\[
   \E{v}{w}=-v_0w_0+v_1w_1+v_2w_2 \ .
\]
This (ordered) framing also determines an orientation of $\mx_0$ and a
time-orientation, by postulating that $\frac{\partial\,}{\partial
x_0}$ is a future timelike vector. The isometries of $\mx_0$ coincide
with the affine transformations of $\mr^3$ with linear part preserving
the Minkowskian form. The group $ISO^+(\mx_0)$ of the isometries that
preserve both the orientation of $\mx_0$ and the time orientation,
coincides with $\mr^3\rtimes SO^+(2,1)$, where $SO^+(2,1)$ denotes the
group of corresponding linear parts. 

There is a standard isometric embedding of $\mh^2$ into $\mx_0$ which
identifies the hyperbolic plane with the set of future directed
unitary timelike vectors, that is
\[
   \mh^2=\{v\in\mr^3| \E{v}{v}=-1\textrm{ and } v_0>0\}.
\]
Clearly $SO^+(2,1)$ acts by isometries on $\mh^2$, this action is
faithful and induces an isomorphism between $SO^+(2,1)$ and the whole
group of orientation preserving isometries of $\mh^2$. Many of the
above facts extend to Minkowski spaces $\mm^n$ of arbitrary dimension.

\section{De Sitter space}\label{dSMod}
Let us consider the $4$-dimensional Minkowski spacetime
$(\mm^4,\E{\cdot}{\cdot})$ and set
\[
  \hat\mx_1=\{v\in\mr^4|\E{v}{v}=1\} \ .
\]
It is not hard to show that $\hat\mx_1$ is a Lorentzian sub-manifold
of constant curvature $1$. Moreover the group $\OO(3,1)$ acts on it by
isometries.  This action is transitive and the stabilizer of a point
is $\OO(2,1)$. It follows that $\hat\mx_1$ is an isotropic Lorentzian
spacetime and $\OO(3,1)$ coincides with the full isometry group of
$\hat\mx_1$.  Notice that $\hat\mx_1$ is orientable and
time-orientable. In particular $\SOO^+(3,1)$ is the group of
time-orientation and orientation preserving isometries whereas
$\SOO(3,1)$ (resp. $\OO^+(3,1)$) is the group of orientation
(resp. time-orientation) preserving isometries.

The projection of
$\hat\mx_1$ into the projective space $\mathbb P^3$ is a local
embedding onto an open set that is the exterior of the Klein model of
$\mh^3$ in $\mathbb P^3$ (that is a regular neighbourhood of $\mathbb
P^2$ in $\mathbb P^3$).  We denote this set by $\mx_1$.  Now the
projection $\pi:\hat\mx_1\rightarrow\mx_1$ is a $2$-fold covering and
the automorphism group is $\{\pm Id\}$. Thus the metric on $\hat\mx_1$
can be pushed forward to $\mx_1$. In what follows we consider always
$\mx_1$ endowed with such a metric and we call it the {\it Klein model
of de Sitter space}. Notice that it is an oriented spacetime
(indeed it carries the orientation induced by $\mathbb P^3$) but it is
not time-oriented (automorphisms of the covering
$\hat\mx_1\rightarrow\mx_1$ are not time-orientation preserving).\par
Since the automorphism group $\{\pm Id\}$ is the center of the
isometry group of $\OO(3,1)$, $\mx_1$ is an
isotropic Lorentz spacetime. The isometry group of $\mx_1$ is
$\OO(3,1)/\pm Id$. Thus the projection
$\OO^+(3,1)\rightarrow\ISO(\mx_1)$ is an isomorphism.

\begin{figure}
\begin{center}
\input GR2005_fig_dual_ds.pstex_t
\caption{{\small The plane in $\mh^3$ dual to a point
$p\in\mx_1$.}}\label{dual:ds:fig}
\end{center}
\end{figure}

We give another description of $\mx_1$. Given a point $v\in \hat \mx_1$,
the plane $v^\perp$ cuts $\mh^3$ along a totally geodesic plane
$P_+(v)$.  In fact we can consider on $P_+(v)$ the orientation induced
by the half-space $U(v)=\{x\in\mh^3|\E{x}{v}\leq 0\}$.  In this way
$\hat\mx_1$ parameterizes the oriented totally geodesic planes of
$\mh^3$.  If we consider the involution given by changing the
orientation on the set of the oriented totally geodesic planes of
$\mh^3$, then the corresponding involution on $\hat\mx_1$ is simply
$v\mapsto-v$. In particular, $\mx_1$ parameterizes the set of
(un-oriented) hyperbolic planes of $\mh^3$.  For every $v\in\mx_1$,
$P(v)$ denotes the corresponding plane.  For $\gamma\in\OO^+(3,1)$ we
have
\[
    P(\gamma x)=\gamma (P(x))
\]
(notice that $\OO^+(3,1)$ is the isometry group of both $\mx_1$ and
$\mh^3$).\\

Just as for $\mh^3$, the geodesics in $\hat\mx_1$ are obtained by
intersecting $\hat\mx_1$ with linear $2$-spaces. Thus geodesics in
$\mx_1$ are projective segments. It follows that, given two points
$p,q\in\mx_1$, there exists a unique complete geodesic passing through
them.\par A geodesic line in $\hat\mx_1$ is spacelike (resp. null,
timelike) if and only if it is the intersection of $\hat\mx_1$ with a
spacelike (resp. null, timelike) plane.  For $x\in\hat\mx_1$ and a
vector $v$ tangent to $\hat\mx_1$ at $x$ we have
\[
\begin{array}{ll}
\textrm{if }\E{v}{v}=1  &  \exp_x(tv)=\cos t\, x + \sin t\, v\\
\textrm{if }\E{v}{v}=0  &  \exp_x(tv)=x+tv\\
\textrm{if }\E{v}{v}=-1 &  \exp_x(tv)=\ch t\, x + \sh t\, v.
\end{array}
\]
This implies that a complete geodesic line in $\mx_1$ is
spacelike (resp. null, timelike) if it is a complete projective line
contained in $\mx_1$ (resp. if it is a projective line tangent to $\mh^3$,
if it is a projective segment with both the end-points in
$\partial\mh^3$).  Spacelike geodesics have finite length equal to
$\pi$. Timelike geodesics have infinite Lorentzian length.\\

Take a point $x\in\mh^3$ and a unit vector $v\in T_x\mh^3=x^\perp$.
Clearly we have $v\in\hat\mx_1$ and $x\in T_v\hat\mx_1$. Notice that the
projective line joining $[x]$ and $[v]$ in $\mathbb P^3$ intersects
both $\mh^3$ and $\mx_1$ in complete geodesic lines $c$ and $c^*$.
They are parametrized in the following way
\[
\begin{array}{l}
c(t)=[\ch t\, x +\sh t\, v] \ , \\
c^*(t)=[\ch t\, v + \sh t\, x] \ . 
\end{array}
\]
We say that $c^*$ is the continuation to $c$. They have the same
end-points on $S^2_\infty$ that are $[x+v]$ and $[x-v]$.  Moreover
if $c'$ is the geodesic \emph{ray} starting from $x$ with speed $v$
the continuation ray $(c')^*$ is the geodesic \emph{ray} on $c^*$
starting at $v$ with the same limit point on $S^2_\infty$ as
$c'$.\\

\section{Anti de Sitter space}\label{AdSMod}
As for the other cases, we will recall some general features of the
AdS local model that we will use later.  In particular, both spacetime
and time orientation will play a subtle r\^ole, so it is important to
specify them carefully.

Let $\mathrm M_2(\mr)$ be the space of $2\times2$ matrices with real
coefficients endowed with the scalar product $\eta$ induced by the
quadratic form \[ q(A)=-\det A\ .\] The signature of $\eta$ is
$(2,2)$. The group\[ \SL{2}{R}=\{A|q(A)=-1\}\] is a Lorentzian
sub-manifold of $\mathrm M_2(\mr)$, that is the restriction of $\eta$
on it has signature $(2,1)$.  Given $A,\ B\in\SL{2}{R}$, we have that
\[
   q(AXB)=q(X) \qquad\textrm{ for }X\in\mathrm M_2(\mr)
\]
Thus,  the left action of $\SL{2}{R}\times\SL{2}{R}$ on $\mathrm M_2(\R)$ 
given by
\begin{equation}\label{ads:act:eq}
   (A,B)\cdot X=AXB^{-1}
\end{equation}
preserves $\eta$. In particular, the restriction of $\eta$ on
$\SL{2}{R}$ is a bi-invariant Lorentzian metric, that actually
coincides with its Killing form (up to some multiplicative factor). 
For $X,Y\in\sG\lG(2,\mr)$ we have the usual formula
\[
    \tr XY=2\eta(X,Y) \ .
\]
 We denote by $\hat\mx_{-1}$ the pair $(\SL{2}{R}, \eta)$.  Clearly
$\hat\mx_{-1}$ is an orientable and time-orientable spacetime.  Hence,
the above rule~(\ref{ads:act:eq}) specifies a transitive isometric
action of $\SL{2}{R}\times\SL{2}{R}$ on $\hat\mx_{-1}$.

The stabilizer of $Id\in \hat\mx_{-1}$ is the diagonal group
$\Delta\cong\SL{2}{R}$.  The differential of isometries corresponding
to elements in $\Delta$ produces a surjective representation
\[      
\Delta\rightarrow\SOO^+(\sG\lG(2,\mr),\eta_{Id}) \, .
\]
It follows that $\hat\mx_{-1}$ is an isotropic Lorentzian spacetime
and the isometric action on $\hat\mx_{-1}$ induces a surjective
representation
\[   
 \hat\Phi:\SL{2}{R}\times\SL{2}{R}\rightarrow\ISO_0(\hat\mx_{-1})\,.
\]
Since $\ker\hat\Phi=(-Id, -Id)$, we obtain
\[   
  \ISO_0(\hat\mx_{-1})\cong\SL{2}{R}\times\SL{2}{R}/(-Id,-Id)\,.
\]
The center of $\ISO_0(\hat\mx_{-1})$ is generated
by $[Id,-Id]=[-Id,Id]$. Hence, $\eta$ induces
on the quotient\[PSL(2,\R)=\SL{2}{R}/\pm Id\]
an isotropic Lorentzian structure. 
We denote by $\mx_{-1}$ this spacetime and call it the 
{\it Klein model} of Anti de Sitter spacetime. 
Notice that left and right translations are isometries and
the above remark implies that the induced representation
\[   \Phi:PSL(2,\R)\times PSL(2,\R)\rightarrow\ISO_0(\mx_{-1})\]
is an isomorphism.\\
\smallskip

\paragraph {The boundary of $\mx_{-1}$}
Consider the topological closure $\overline{PSL(2,\R)}$
of $PSL(2,\R)$ in $\mathbb P^3=\mathbb P(\mathrm M_2(\mr))$.  Its
boundary is the quotient of the set
\[    \{X\in\mathrm M_2(\mr)\setminus \{0\}|q(X)=0\}\]
that is the set of rank $1$ matrices. In particular,
$\partial PSL(2,\R)$ is the image of the {\it Segre embedding}
\[   \mathbb P^1\times\mathbb P^1\ni([v],[w])
\mapsto[v\otimes w]\in\mathbb P^3.\] Thus $\partial PSL(2,\R)$ is a
torus in $\mathbb P^3$ and divides it in two solid tori. In
particular, $PSL(2,\R)$ is topologically  a solid torus.\par
The action of $PSL(2,\R)\times PSL(2,\R)$ extends to the whole of
$\overline\mx_{-1}$.  Moreover, the action on $\mathbb
P^1\times\mathbb P^1$ induced by the Segre embedding is simply
\[   (A,B)(v,w)=(Av,B^*w)\]
where we have set $B^*= (B^{-1})^T$ and considered the natural action
of $PSL(2,\R)$ on $\mathbb P^1=\partial\mh^2$.  If $E$ denotes the
rotation by $\pi/2$ of $\mr^2$, it is not hard to show that
\[    EAE^{-1}=(A^{-1})^T\]
for $A\in PSL(2,\R)$.

It is convenient to consider the following modification of Segre
embedding
\[   S:\mathbb P^1\times\mathbb P^1\ni([v],[w])
\mapsto [v\otimes(Ew)]\in\mathbb P^3\] With respect to such a new
embedding, the action of $PSL(2,\R)\times PSL(2,\R)$ on
$\partial\mx_{-1}$ is simply
\[   (A,B)(x,y)=(Ax, By).\]
In what follows, we will consider the identification of the boundary
of $\mx_{-1}$ with $\mathbb P^1\times\mathbb P^1$ given by $S$. \\

The product structure on $\partial\mx_{-1}$ given by $S$ is preserved
by the isometries of $\mx_{-1}$. This allows us to define a {\it
conformal Lorentzian structure} (\emph{i.e.} a {\it causal structure}) on
$\partial\mx_{-1}$.  More precisely, we can define two foliations on
$\partial\mx_{-1}$.  The left foliation is simply the image of the
foliation with leaves
\[    l_{[w]}=\{([x],[w])|[x]\in\mathbb P^1\}\]
and a leaf of the right foliation is the image of
\[    r_{[v]}=\{([v],[y])|[y]\in\mathbb P^1\}.\]
Notice that left and right leaves are projective lines in $\mathbb P^3$.
Exactly one left and one right leaves
pass through a given point. On the other hand, 
given right leaf and left leaf meet each other at one
point. Left translations preserve leaves of left foliation, whereas
right translations preserve leaves of right foliation.\par

\begin{figure}
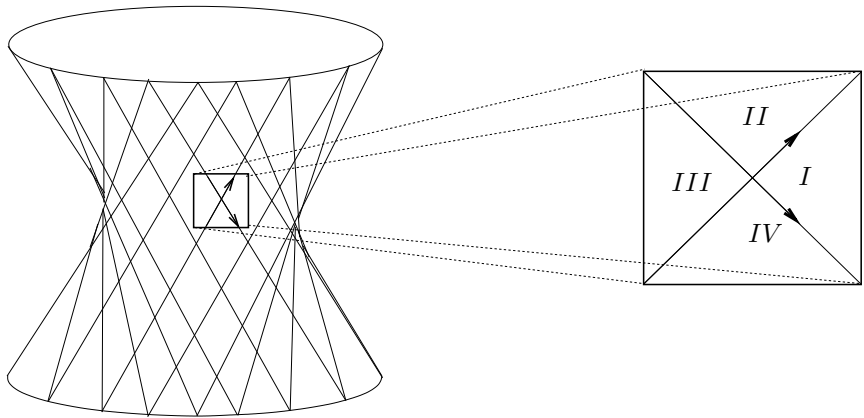

\begin{center}
\input GR2005_fig_ads_bound.pstex_t
\caption{{\small The product structure on
$\partial\mx_{-1}$.}}\label{ads:bound:fig}
\end{center}
\end{figure}

The leaves of the right and left foliations are oriented as the
boundary of $\mh^2$.  Thus if we take a point $p\in\partial\mx_{-1}$,
the tangent space $T_p\partial\mx_{-1}$ is divided by the tangent
vector of the foliations in four quadrants.  By using leaf
orientations we can enumerate such quadrants as shown in figure
~\ref{ads:bound:fig}.  Thus we consider the $1+1$ cone at $p$ given by
choosing the second and fourth quadrants. We make this choice because
in this way the causal structure on $\partial\mx_{-1}$ is the
``limit'' of the causal structure on $\mx_{-1}$ in the following
sense.

Suppose $A_n$ to be a sequence in $\mx_{-1}$ converging to
$A\in\partial\mx_{-1}$, and suppose $X_n\in T_{A_n}\mx_{-1}$ to be a
sequence of timelike vectors converging to $X\in T_A\partial\mx_{-1}$,
then $X$ is non-spacelike with respect to the causal structure of the
boundary.

Notice that oriented left (resp. right) leaves are
homologous non-trivial simple cycles on $\partial\mx_{-1}$, so they
determine non-trivial elements of $\coom1(\partial\mx_{-1})$ that we
denote by $c_L$  and $c_R$.

\paragraph{ Geodesic lines and planes}
Geodesics in $\mx_{-1}$ are obtained by intersecting projective lines
with $\mx_{-1}$.

A geodesic is timelike if it is a projective line entirely contained
in $\mx_{-1}$; its Lorentzian length is $\pi$.  In this case it is a
non-trivial loop in $\mx_{-1}$ (a core).  Take $x\in\mx_{-1}$ and
a unit timelike vector $v\in T_x\mx_{-1}$ .  If $\hat x$ is a pre-image
of $x$ in $\hat\mx_{-1}$ and $\hat v\in T_{\hat x}\hat\mx_{-1}$ is a
pre-image of $v$, then we have
\[    
 \exp_x tv = [\cos t\, \hat x+\sin t\, \hat v] \, .
\]

A geodesic is null if it is contained in a projective line tangent to
$\partial\mx_{-1}$.  Given $x\in\mx_{-1}$, and a null vector $v\in
T_x\mx_{-1}$, if we take $\hat x$ and $\hat v$ as above we have
\[
\exp_x tv=[\hat x+t\hat v] \, .
\]

Finally, a geodesic is spacelike if it is contained in a projective
line meeting $\partial\mx_{-1}$ at two points; its length is infinite.
Given $x\in\mx_{-1}$ and a unit spacelike vector $v$ at $x$, fixed
$\hat x$ and $\hat v$ as above, we have
\[ 
\exp_x tv=[\ch t\, \hat x+\sh t\,\hat v] \, .
\]

Geodesics passing through the identity are $1$-parameter subgroups.
Elliptic subgroups correspond to timelike geodesics, parabolic
subgroups correspond to null geodesics and hyperbolic subgroups are
spacelike geodesics.\\

Totally geodesic planes are obtained by intersecting projective planes
with $\mx_{-1}$.

If $W$ is a subspace of dimension $3$ of $\mathrm M_2(\mr)$ and the
restriction of $\eta$ to it has signature $(m_+,m_-)$, then the
projection $P$ of $W$ in $\mathbb P^3$ intersects $PSL(2,\R)$ if and
only if $m_- >0$. In this case the signature of $P\cap\mx_{-1}$ is
$(m_+,m_--1)$.  Since $\eta$ restricted to $P$ is a flat metric we
obtain that
\begin{enumerate}

\item
If $P\cap\mx_{-1}$ is a Riemannian plane, then it is isometric to
$\mh^2$.
\item If $P\cap\mx_{-1}$ is a Lorentzian plane, then it is a Moebius
band carrying an Anti de Sitter metric.
\item If $P\cap\mx_{-1}$ is a null plane, then $P$ is tangent to
$\partial\mx_{-1}$.
\end{enumerate}

Since every spacelike plane cuts every timelike geodesic at one point,
spacelike planes are compression disks of $\mx_{-1}$.  The boundary of
a spacelike plane is a spacelike curve in $\partial\mx_{-1}$ and it is
homologous to $c_L+c_R$. Every Lorentzian plane is a Moebius band.
Its boundary is homologous to $c_L- c_R$. Every null plane is a
pinched band. Its boundary is the union of a right and a left
leaf.\\

\paragraph{Duality in $\mx_{-1}$}
The form $\eta$ induces a duality in $\mathbb P^3$ between points and
planes, and between projective lines. Since the isometries of
$\mx_{-1}$ are induced by linear maps of $\mathrm M_2(\mr)$ preserving
$\eta$, this duality is preserved by isometries of $\mx_{-1}$.

\begin{figure}
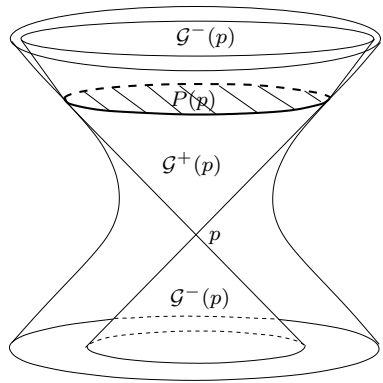

\begin{center}
\input GR2005_fig_dual_ads.pstex_t
\caption{{\small The dual plane of a point $p\in\mx_{-1}$.}}
\label{dual:ads:fig}
\end{center}
\end{figure}

If we take a point in $\mx_{-1}$ its dual projective plane defines a
Riemannian plane in $\mx_{-1}$ and, conversely, Riemannian planes are
contained in projective planes dual to points in $\mx_{-1}$.  Thus a
bijective correspondence between points and Riemannian planes exists.
Given a point $x\in\mx_{-1}$, we denote by $P(x)$ its dual plane and,
conversely, if $P$ is a Riemannian plane, then $x(P)$ denotes its dual
point. If we take a point $x\in\mx_{-1}$ and a timelike geodesic $c$
starting at $x$ and parametrized in Lorentzian arc-length we have that
$c(\pi/2)\in P(x)$.  Moreover, this intersection is orthogonal.
Conversely, given a point $y$ in $P(x)$, there exists a unique
timelike geodesic passing through $x$ and $y$ and such a geodesic is
orthogonal to $P(x)$.  By using this characterization, we can see that
the plane $P(Id)$ consists of those elliptic transformations of
$\mh^2$ that are the rotation by $\pi$ at their fixed points.  In this
case an isometry between $P(Id)$ and $\mh^2$ is simply obtained by
associating to every $x\in P(Id)$ its fixed point in
$\mh^2$. Moreover, such a map
\[    
   I:P(Id)\rightarrow\mh^2
\] 
is natural in the following sense.  The isometry group of $P(Id)$ is
the stabilizer of the identity, that we have seen to be the diagonal
group $\Delta\subset PSL(2,\R)\times PSL(2,\R)$.  Then we have\[
I\circ(\gamma,\gamma)=\gamma\circ I.\] The boundary of $P(Id)$ is the
diagonal subset of $\partial\mx_{-1}=\mathbb P^1\times\mathbb P^1$ that
is
\[   
\partial P(Id)=\{(x,x)\in\partial\mx_{-1}|x\in\mathbb P^1\} \, .
\]
The map $I$ extends to $\overline{P(Id)}$, by sending the
point $(x,x)\in\partial P(Id)$ to $x\in\mathbb P^1=\partial\mh^3$.\\

\begin{remark}\label{mod:ort:rm}
\emph{ The plane $P(id)$ could be regarded as the set of projective
classes of timelike vectors of $(\sG\lG(2,\mr),\eta_{Id})$.}  \emph{
The isometry $I:P(id)\rightarrow\mh^2$ extends to a linear isometry
$\bar I:\sG\lG(2,\mr)\rightarrow\mx_0$ that is $PSL(2,\mr)$
equivariant, where $PSL(2,\mr)$ acts on $\sG\lG(2,\mr)$ via the
adjoint representation and on $\mx_{0}$ via the canonical isomorphism
$PSL(2,\mr)\cong SO(2,1)$.  } \emph{ If $l$ is an oriented geodesic of
$\mh^2=P(id)$ whose end-points are projective classes of null vectors
$x^-,x^+$, then the infinitesimal generator of positive translations
along $l$, say $X(l)\in\sG\lG(2,\mr)$ is sent by $\bar I$ to the unit
spacelike vector $v\in\mx_0$ orthogonal to $l$ such that $x^-,x^+,v$
form a positive basis of $\mx_0$.  }
\end{remark}

The dual point of a null plane $P$ is a point $x(P)$ on the
boundary $\partial\mx_{-1}$.  It is the intersection point of the left
and right leaves contained in the boundary of that plane.  Moreover,
the plane is foliated by null-geodesics tangent to $\partial\mx_{-1}$
at $x(P)$.  Conversely, every point in the boundary is dual to the
null-plane tangent to $\partial\mx_{-1}$ at $x$.

Finally the dual line of a spacelike line $l$ is a spacelike line
$l^*$.  Actually $l^*$ is the intersection of all $P(x)$ for $x\in l$.
$l^*$ can be obtained by taking the intersection of null planes dual
to the end-points of $l$. In particular, if $x_-$ and $x_+$ are the
end-points of $l$, then the end-points of $l^*$ are obtained by
intersecting the left leaf through $x_-$ with the right leaf through
$x_+$, and the right leaf through $x_-$ with the left leaf through
$x_+$, respectively.

There is a simple interpretation of the dual spacelike geodesic
for a hyperbolic $1$-parameter subgroup $l$. In this case $l^*$ is
contained in $P(Id)$ and is the inverse image through $I$ of the axis
fixed by $l$ in $\mh^2$.  Conversely, geodesics in $P(Id)$ correspond
to hyperbolic $1$-parameter subgroups.

\paragraph{Orientation and time-orientation of $\mx_{-1}$}
In order to define a time-orientation it is enough to define a time
orientation at $Id$.  This is equivalent to fixing an orientation on the
elliptic $1$-parameter subgroups.  We know that such a subgroup
$\Gamma$ is the stabilizer of a point $p\in\mh^2$.  Then we stipulate
that an infinitesimal generator $X$ of $\Gamma$ is future directed if
it is a positive infinitesimal rotation around $p$.

A spacelike surface in a oriented and time-oriented spacetime is
oriented by means of the rule: {\it first the normal future-directed
vector field}.  So, we choose the orientation on $\mx_{-1}$ that
induces the orientation on $P(Id)$ that makes $I$ an
orientation-preserving isometry.

Clearly, orientation and time-orientation on $\mx_{-1}$ induce 
orientation and a time-orientation on the boundary. 
Choose a future-directed unit timelike vector $X_0\in T_{Id}\mx_{-1}$ and
consider the $1$-parameter group of isometries of $\mx_{-1}$ given by 
\[
    u_t=(\exp(tX_0),1)
\]
The corresponding Killing vector field is the right-invariant field
$X(A)=X_0A$, that is future-directed. Such a field
extends to $\overline\mx_{-1}$.  Moreover, since for every
$x,y\in\mathbb P^1$ we have
\[
    u_t S(x,y)=S(\exp(tX_0)x,y)
\]
we see that on $\partial\mx_{-1}$ the vector field $X$ is a
positive generator of the left foliation (positive with respect  to the
identification of left leaves with $\mathbb P^1$ given by $S$).
So, a positive generator of the left foliation is \emph{future-oriented}.

An analogous computation shows that a positive generator of the right
foliation is \emph{past oriented}.

Let $e_L$ and $e_R$ be respectively positive generators of the left
and right foliations on $\partial\mx_{-1}$. The positive generator of
$\partial P(Id)$ (positive with respect to the orientation induced by
$P(Id)$) is a positive combination of $e_L$ and $e_R$.  It easily
follows that a positive basis of $T\partial\mx_{-1}$ is given by
$(e_R, e_L)$.

\section{Complex projective structures on surfaces}\label{SPS}
A complex projective structure on a oriented connected surface $S$ is
a $(S^2_\infty, PSL(2,\mc))$-structure (respecting the orientation).

We will often refer to a parameterization of complex
projective structures given in~\cite{Ku}.  That work is a
generalization (even in higher dimension) of a previous classification
due to Thurston when the surface is assumed to be compact. Another
very similar treatment can be found in~\cite{Ap}.\par In this
section we will recall the main constructions of~\cite{Ku}, but we
will omit any proof, referring to that work.

Let us take a projective structure on a surface $S$ and consider a
developing map
\[
    D:\tilde S\rightarrow S^2_\infty\,.
\]
Pulling back the standard metric of $S^2_\infty$ on $\tilde S$ is not a
well-defined operation, as it depends on the choice of the
developing map.  Nevertheless, by the compactness of $S^2_\infty$, the
completion $\hat S$ of $\tilde S$ with respect to such a metric is
well-defined.  By looking at $\hat S$ we can focus on $3$ cases that
lead to very different descriptions:
\medskip\par\noindent 1) $\tilde S$ is complete: in this case $D$ is a
homeomorphism so that $S$ is $S^2_\infty$ (with the standard
structure). We say that $S$ is of \emph{elliptic} type.
\medskip\par\noindent
2) $\hat S\setminus\tilde S$ consists only of one point: in this case $\tilde
S$ is projectively equivalent to $\mr^2$ and the holonomy action preserves the
standard Euclidean metric (in particular $S$ is equipped with a Euclidean
structure). We say that $S$ is of \emph{parabolic} type.
\medskip\par\noindent 3) $\hat S\setminus\tilde S$ contains at least
$2$ points: in this case we say that $S$ is of \emph{hyperbolic} type.
\medskip\par\noindent Clearly, the most interesting case is the third
one (for instance, it includes the case when $\pi_1(S)$ is not
Abelian).  In this case, by following an idea of Thurston, Kulkarni
and Pinkall constructed a canonical stratification of $\tilde S$.  Let
us quickly explain their procedure.\par A \emph{round disk} in $\tilde
S$ is a set $\Delta$ such that $D|_\Delta$ is injective and the image
of $\Delta$ is a round disk in $S^2_\infty$ (this notion is well
defined because maps in $PSL(2,\mc)$ send round disks onto round
disks).  Given a \emph{maximal} disk $\Delta$ (with respect to the
inclusion), we can consider its closure $\overline\Delta$ in $\hat
S$.\par $\overline\Delta$ is sent by $D$ to the closed disk
$\overline{D(\Delta)}$.  In particular if $g_\Delta$ denotes the
pull-back on $\Delta$ of the standard {\it hyperbolic} metric on
$D(\Delta)$, we can regard the boundary of $\Delta$ in $\hat S$ as its
ideal boundary.\par Since $\Delta$ is supposed to be maximal,
$\overline\Delta$ is not contained in $\tilde S$. So, if
$\Lambda_\Delta$ denotes the set of points in
$\overline\Delta\setminus\tilde S$, let $\hat \Delta$ be the convex
hull in $(\Delta, g_\Delta)$ of $\Lambda_\Delta$ (by maximality
$\Lambda_\Delta$ contains at least two points).\par
\begin{prop}\cite{Ku}
For every point $p\in\tilde S$ there exists a unique maximal disk
$\Delta$ containing $p$ such that $p\in\hat\Delta$.
\end{prop}
So, $\{\hat\Delta|\Delta\textrm{ is a maximal disk}\}$ is a partition
of $\tilde S$. We call it the \emph{canonical stratification} of
$\tilde S$. Clearly the stratification is invariant under the action of
$\pi_1(S)$.

Let $g$ be the Riemannian metric on $\tilde S$ that coincides at $p$
with the metric $g_\Delta$, where $\Delta$ is the maximal disk such
that $p\in\hat\Delta$.  It is a conformal metric, in the sense that it
makes $D$ a conformal map.  It is $\mathrm C^{1,1}$ and is invariant
under the action of $\pi_1(S)$. So, it induces a metric on $\tilde S$. We
call it the \emph{Thurston metric} on $\tilde S$.
\begin{center}
\begin{figure}
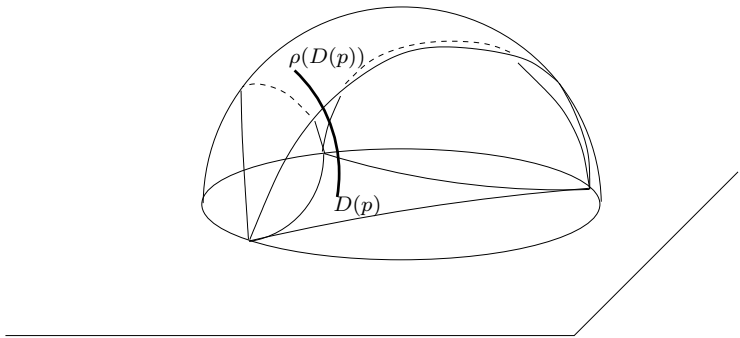

\input GR2005_fig_hhull.pstex_t
\caption{{\small The construction of the $H$-hull.}}\label{hhull:fig}
\end{figure}
\end{center}
By means of the canonical stratification, we are going to construct a
hyperbolic structure on $S\times (0,+\infty)$. In particular, we  
construct an $h$-equivariant local homeomorphism
\[
  dev:\tilde S\times(0,+\infty)\rightarrow\mh^3
\]
(where $h$ is the holonomy of $\tilde S$).  For $p\in\tilde S$ let
$\Delta(p)$ denote the maximal disk such that $p\in\hat\Delta(p)$. The
boundary of $D(\Delta(p))$ can be regarded as the boundary of a plane
$P(p)$ of $\mh^3$. Denote by $\rho:\overline \mh^3\rightarrow P$ the
nearest point retraction. Then $dev(p,\cdot)$ parameterizes in
arc-length the geodesic ray of $\mh^3$ with end-points $\rho(D(p))\in
P$ and $D(p)$ (see Fig.~\ref{hhull:fig}).
\begin{prop}\cite{Ku}
 $dev$ is a $\mathrm C^{1,1}$ developing map for a hyperbolic structure on
 $S\times (0,+\infty)$.
Moreover it extends to a map
\[
    \overline {dev}:\tilde S\times (0,+\infty]\rightarrow\overline\mh^3
\]
such that $dev|_{\tilde S\times\{+\infty\}}$ is a developing map for
the complex projective structure on $S$.
\end{prop}
We call such a hyperbolic structure the \emph{H-hull} of $S$ and
denote it by $H(S)$.  Now $H(\tilde S)$ is never a complete hyperbolic
manifold. If $\overline{H(\tilde S)}$ is the completion of $H(\tilde
S)$ let us set $P(\tilde S)=\overline{H(\tilde S)}\setminus H(\tilde
S)$. Now, $P(\tilde S)$ takes a well-defined distance, induced by the
distance of $\mh^3$ through the developing map. In~\cite{Ku} it is
shown that $P(\tilde S)$ is isometric to a straight convex set of
$\mh^2$. Moreover the developing map (regarding $P$ as boundary of
$H(\tilde S)$)
\[
   dev: P(\tilde S)\rightarrow\mh^3
\]
is the bending of $P(\tilde S)$ along a measured geodesic lamination
$\lambda(\tilde S)$ on it (see Section~\ref{laminations} for the rigorous
definition of measured geodesic lamination on a straight convex
set).\par

\begin{teo}\cite{Ku}
The correspondence
\[
   \tilde S\mapsto(P(\tilde S), \lambda(\tilde S))
\]
induces a bijection among the complex projective structures on a disk
of hyperbolic type (up to projective equivalence) and the set of
measured geodesic laminations on straight convex sets (up to
$PSL(2,\mr)$-action).
\end{teo}    
\begin{center}
\begin{figure}
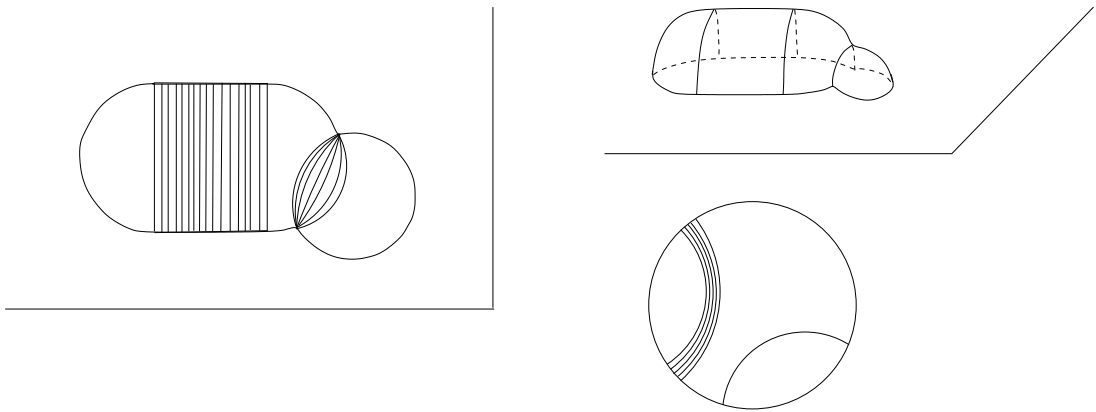

\input GR2005_fig_qf.pstex_t
\caption{{\small The $H$-hull and the canonical stratification associated to a
  simply connected domain of $\mc$.}}\label{qf:fig}
\end{figure}
\end{center}

\begin{remark}
\emph{ When the developing map is an embedding $\tilde S\rightarrow
S^2_\infty$ and the boundary of $\tilde S$ in $S^2_\infty$ is a Jordan
curve, we can give a simpler description of Kulkarni-Pinkall
constructions. In this case, we can consider the convex hull $K$ of
$\partial\tilde S$ in $\mh^3$ (see Fig.~\ref{qf:fig}). Then,}
\medskip\par\noindent \emph{ 1. $H(\tilde S)$ is the component of
$\mh^3\setminus K$ that is close to $\tilde S$.}
\medskip\par\noindent
\emph{
2. $P(\tilde S)$ is the component of $\partial K$ facing towards $\tilde S$.}
\medskip\par\noindent \emph{ 3. The canonical stratification of
$\tilde S$ is obtained by taking the inverse images of the faces of
$P(\tilde S)$ through the nearest point retraction.}
\medskip\par\noindent \emph{ 4. The lamination associated to $\tilde
S$ is obtained by depleating $P(\tilde S)$.}
\end{remark}

If $S$ has non-Abelian fundamental group, then the projective
structures on $S$ are of hyperbolic type. Moreover the pleated set
$P(\tilde S)$ is not a single geodesic (\emph{i.e.}, the interior of
$P(\tilde S)$ is non-empty).  By looking at the construction of the
$H$-hull and of $P(\tilde S)$ we can see that there exists a natural
retraction
\[
   \rho:\overline{H(\tilde S)}\rightarrow P(\tilde S)\,.
\]
In particular $\pi_1(S)$ acts free and properly discontinuously on the
interior of $P(\tilde S)$, and the quotient $P(\tilde S)/\pi_1(S)$ is
homeomorphic to $S$.

\begin{cor}
Let $S$ be a surface with non-Abelian fundamental group. For a
projective structure $F$ on $S$ denote by $h_F$ the hyperbolic
holonomy of $P(F)$ and by $\lambda_F=(P(\tilde F),\Ll_F,\mu_F)$ the
measured geodesic lamination associated to the universal covering
space.  Then the map
\[
    F\mapsto[h_F,\lambda_F]
\]
induces a bijection between the set of complex projective
structures on $S$ and the set $\Mm\Ll^\Ee_S$ (defined in Corollary~\ref{Teich-like}).
\end{cor}


\chapter{Flat globally hyperbolic spacetimes}\label{FGHST}
The main goal in this chapter is to prove the results on the
classification of 3-dimensional maximal globally hyperbolic flat
spacetimes stated in Section~\ref{FLAT:I}. Before doing it we will
recall some general facts about globally hyperbolic spacetimes and the
cosmological time.
 
\section {Globally hyperbolic spacetimes}\label{glob:hyp:ST}

A spacelike hypersurface $S$ in a spacetime $M$ is a \emph{Cauchy
surface} if any inextensible causal curve of $M$ intersects $S$
exactly in one point.

A spacetime $M$ is \emph{globally hyperbolic} if it contains a Cauchy
surface.

\begin{remark}\emph{
The usual definition of global hyperbolcity is rather different. On
the other hand by a theorem of Geroch \cite{H-E} it is equivalent to
that given above, which is more expressive for our purpose.
}\end{remark}

A globally hyperbolic spacetime $M$ satisfies some strong properties:
\medskip\par\noindent
1) It is ``topologically simple'', \emph{i.e.} it is homeomorphic to
   $S\times\mr$, $S$ being a Cauchy surface.
\smallskip\par\noindent 2) There exists a real-valued function $\tau$
on $M$ such that the gradient of $\tau$ is a unit timelike vector and
level surfaces of $\tau$ are Cauchy surfaces.
\smallskip
\par
\noindent  
3) If $S$ is a Cauchy surface for a spacetime
$M$, then its lifting, say $\tilde S$, on the isometric universal
covering $\tilde M$, is a Cauchy surface for $\tilde M$. In
particular, if $M$ is globally hyperbolic, so is $\tilde M$.

If $M$ is a spacetime of constant curvature $\kappa$, and $S$ is a
spacelike slice of $M$, then its first and second fundamental forms
obey to a differential equation, said Gauss-Codazzi equation.

Conversely, by a classical result of Choquet-Bruhat and Geroch
\cite{Cho-Ge}, given a scalar product $g$ and a symmetric bilinear
form $b$ on $S$ satisfying the Gauss-Codazzi equation, there exists a
unique (up to isometries) maximal globally hyperbolic Lorentzian
structure of constant curvature $\kappa$ on $S\times\mr$, such that
\medskip\par\noindent
- $S\times\{0\}$ is a Cauchy surface;
\smallskip

\par\noindent - The first and the second fundamental form of
$S\times\{0\}$ are respectively $g$ and $b$.
 
\begin{remark}\emph{
Let us make precise what}  maximal \emph{means in this context.}

\emph{ A constant curvature globally hyperbolic spacetime $M$ is said}
maximal \emph{ if every isometric embedding of $M$ into a constant
curvature spacetime $M'$} sending a Cauchy surface of $M$ onto a
Cauchy surface of $M'$ \emph{ is an isometry.}

\emph{ For instance the future of $0$ in $\mx_0$ is a maximal globally
hyperbolic spacetime (in the above sense), even if it is embedded in a
bigger globally hyperbolic spacetime (the whole $\mx_0$).  }
\end{remark}

A standard way to treat the classification problem of constant
curvature globally hyperbolic spacetimes is then to solve the
Gauss-Codazzi equation on $S$.

A problem in such an approach is that different solutions can lead to
the same spacetime: in fact solutions of Gauss-Codazzi equations are
in $1$-to-$1$ correspondence with pairs $(M,S)$, $M$ being a maximal
globally hyperbolic spacetime of constant curvature $\kappa$ and $S$
being a Cauchy surface embedded in $M$.

A possible way to overcome this difficulty is to impose some
supplementary condition to the space of solutions that translates some
geometric property on $S$ into $M$.

For instance a widely studied possibility is to require the trace of
$b$ with respect to $g$ to be constant, that is that the surface $S$
has constant mean curvature in $M$ (e.g. we refer to \cite{Mon, A-M-T, B-Z,
Kra-Sch}).

The approach we follow in this paper is rather different.  We will
give a global description of constant curvature spacetimes in terms of
some parameters, that are a hyperbolic structure (possibly with
geodesic boundary) on $S$ and a measured geodesic lamination.  In some
sense, such parameters encode the intrinsic geometry of the asymptotic
states of the so called \emph{cosmological time} rather than the
embedding data of some Cauchy surface.

\section {Cosmological time}\label{CT}
We refer to \cite{A} for a general and careful treatment of this
matter.  Here we limit ourselves to recalling the main features of
this notion.

Let $M$ be any spacetime. The {\it cosmological function} of $M$
$$\tau: M \to (0,+\infty]$$ is defined as follows: let $C^-(q)$ be the
set of past-directed causal curves in $M$ that start at $q\in M$. Then
$$\tau (q) = \sup\{L(c)|\ c\in C^-(q)\}\ ,$$ where $L(c)$ denotes the
Lorentzian length of $c$.  In general, the cosmological function
can be very degenerate; for example, on the Minkowski space $\tau$ is
the $+\infty$-constant function.  We say that $\tau$ is {\it regular}
if:
\smallskip

(1) $\tau(q)$ is finite valued for every $q\in M$;

(2) $\tau \to 0$ along every past-directed inextensible causal curve.
\medskip

It turns out that if $\tau$ is regular, it is a continuous {\it global
time} on $M$ that we call its {\it canonical cosmological time}. 

Having canonical cosmological time has strong consequences for the
structure of $M$, and $\tau$ itself has stronger properties than the
simple continuity; we just recall that:
\smallskip

- $\tau$ is locally Lipschitz and twice differentiable almost
  everywhere;

- each $\tau$-level surface is a {\it future} Cauchy surface
({\it i.e.} each inextensible causal curve that intersects the future
of the surface actually intersects it once); 

- $M$ is globally hyperbolic;

- For every $q\in M$, there exists a future-directed time-like unit
speed geodesic ray $\gamma_q : (0,\tau(q)]\to M$ such that:
$$ \gamma_q(\tau(q))=q\ , \ \ \ \tau(\gamma_q(t))=t \ .$$

The equivalence classes of these rays up to a suitable past-asymptotic
equivalence can be interpreted as forming the {\it initial
singularity} of $M$ (which of course is not contained in the
spacetime). This set could be also endowed with a stronger geometric
structure. In fact, we are going to deal with spacetimes having rather
tame canonical cosmological time; in these cases the geometry of the
initial singularity will quite naturally arises.

\section{Regular domains}\label{RD:sec}

A {\it flat regular domain} $\Uu$ in the Minkowski
space $\mx_0$, is a convex domain that coincides with the intersection
of the future of its null support planes. We also require that there
are at least two null support planes. Note that a regular domain is
future complete. We state here Barbot's results in 3 dimensions.

\begin{teo}\label{flat:descr:teo} \cite{Ba}(1)
Let $M$ be a maximal globally hyperbolic spacetime containing a {\rm
complete} Cauchy surface $S$. Then any developing map is an embedding,
so that a universal covering $\tilde M$ can be identified with a
domain in $\mx_0$. Moreover, up to changing the time orientation, one
of the following sentences holds
\medskip\par
\noindent 1. $\tilde M=\mx_0$ and the holonomy
group acts by spacelike translations on $\mx_0$ (either
$\pi_1(M)$ is isomorphic to $\{0\}$, or to $\mz$, or to $\mz\oplus\mz$);
\medskip\par
\noindent 2. $\tilde M$ is the future of a null plane $P$, and the
holonomy acts by spacelike translations (in particular, $\pi_1(M)$ 
either is $\{0\}$ or $\mz$);
\medskip\par\noindent 3. $\tilde M=\fut(P)\cap\pass(Q)$, where $P$
and $Q$ are parallel null planes: in this case the holonomy group is
isomorphic to either $\{0\}$ or $\mz$;
\medskip\par
\noindent 4. $\tilde M$ is the future of a spacelike
geodesic line, and the holonomy group is isomorphic to either $\{0\}$,
or $\mz$, or $\mz\oplus\mz$;
\medskip\par
\noindent 5. $\tilde M$ is a regular domain different from
the future of a spacelike geodesic line and the linear holonomy
$\pi_1(M)=\pi_1(S)\rightarrow SO^+(2,1)$ is a faithful and discrete
representation.
\end{teo}

\begin{cor}\label{flat:descr:cor}
Let the Cauchy surface $S$ have non-Abelian fundamental group. Then
the universal covering $\tilde M$ of $M$ (as in the previous Theorem)
is a regular domain different from the future of a spacelike geodesic
line and the linear holonomy is a faithful and discrete representation
$h_L:\pi_1(S)\rightarrow SO(2,1)$.
\end{cor}

\begin{figure}
\begin{center}
\input GR2005_fig_rdex.pstex_t
\caption{{\small Examples of regular domains.}}\label{rdex:fig}
\end{center}
\end{figure}

Note that in the first 3 cases of the above theorem, $\tilde M$ does
not have canonical cosmological time. On the other hand, we have:
\begin{prop}\label{COSMT_REGD} 
Every flat regular domain $\Uu$ has canonical cosmological time
$T$. In fact $T$ is a {\rm C}$^{1,1}$-submersion onto
$(0,+\infty)$. Every $T$-level surface $\Uu(a)$, $a\in
(0,+\infty)$, is a complete Cauchy surface of $\Uu$. For every $x\in
\Uu$, there is a unique past-directed geodesic timelike segment
$\gamma_x$ that starts at $x$, is contained in $\Uu$, has finite
Lorentzian length equal to $T(x)$. 
\end{prop}
\Dim We only give a sketch of the proof of this theorem, referring
to~\cite{Ba, Bon} for a complete proof.

A first very simple remark is that the cosmological time
function $T$ is finite-valued. In fact, since $\Uu$ is convex, the
Lorentzian distance between two time-related points is realized by the
geodesic segment between them. Since $\pass(p)\cap\overline\Uu$ is
compact for every $p\in\Uu$, the maximum of the distance from $p$ is
attained by some point on the boundary.

In $\pass(p)$ the level surfaces of the distance from $p$ are strictly
convex in the past, whereas $\Uu$ is convex in the future. It follows
that, if the maximum of the distance from $p$ is attained at $r$, then
the level surface through $r$ and the boundary of $\Uu$ meet only at
$p$.

Eventually, there exists only a point $r=r(p)\in\partial\Uu$ such that
$T(p)$ is the length of the timelike segment $[r(p),p]$.  The map
\[
  p\mapsto r(p)
\]
turns to be Lipschitz, so that $T$ is a continuous function.

The point $r(p)$ can be characterized as the unique point $r$ on
$\partial\Uu$ such that the plane through it, orthogonal to the
segment $[r,p]$ is a spacelike support plane.  In particular, we get
that $r(q)=r(p)$ for every point $q$ on the geodesic line passing
through $p$ and $r(p)$.

$T$ is $\mathrm C^{1,1}$ and the gradient of $T$ can be expressed by
the following formula
\begin{equation}\label{CTderivative:eq}
  \mathrm{grad}\, T(p)=-\frac{1}{T(p)}(p-r(p))\,.
\end{equation}
Let us sketch how this formula can be proved.  Given $p\in\Uu$ let
$\Uu_-$ be the future of the point $r(p)$ and $\Uu_+$ be the future of
the spacelike plane through $r(p)$ orthogonal to $p-r(p)$.  We have
\[
   \Uu_-\subset\Uu\subset\Uu_+\,.
\]
Let $T_\pm$ denote the cosmological time on $\Uu_\pm$. These functions
are smooth and we have
\[
\begin{array}{l}
    T_-\leq T\leq T_+\\
    T_-(p)=T(p)=T_+(p)\\
    \mathrm{grad}\, T_-(p)=\mathrm{grad}\, T_+(p)=-\frac{1}{T(p)}(p-r(p))\,.
\end{array}
\]
Formula~(\ref{CTderivative:eq}) easily follows.\\

Finally, $\mathrm{grad}\, T$ is a past directed unit timelike vector
field such that its integral lines are geodesics. On the other hand,
it is proved in \cite{Ba}(1) that the level sets of $T$ are
\emph{complete Cauchy surfaces}.

\cvd

We have associated to every regular domain
\medskip\par\noindent
1) The \emph{cosmological time} $T$;
\medskip\par\noindent
2) The \emph{retraction} in the past $r:\Uu\rightarrow\partial\Uu$;
\medskip\par\noindent
3) The \emph{Gauss map} $N:\Uu\ni p\mapsto -\mathrm{grad}\, T(p)\in\mh^2$.
\medskip\par\noindent
By the formula of the gradient we obtain the following (very convenient)
decomposition of points in $\Uu$
\begin{equation}\label{flatdec}
    p=r(p)+T(p)N(p)\qquad\textrm{ for every }p\in\Uu\,.
\end{equation}

If $f$ is an isometry of $\Uu$ the following invariance conditions hold
\[
   \begin{array}{l}
    T(f(p))=T(p)\\
    r(f(p))=f(r(p))\\
    N(f(p))=L(f)N(p)
   \end{array}
\]
where $L(f)$ is the linear part of $f$. Let us remark that the name
``Gauss map'' is appropriate, because it coincides with the Gauss map
of the level surface of $T$.

The image of the Gauss map can be interpreted as the set of spacelike
directions orthogonal to some spacelike support plane for $\Uu$. So,
it is a convex set. Since $\Uu$ is a regular domain, it follows that
its closure $H_\Uu$ is a \emph{straight convex set in $\mh^2$}, that
is a closed subset obtained as the the convex hull of a set of points
in $S^1_\infty = \partial\mh^2$ (corresponding in this case to the
null support planes). We say that such a straight convex set (as well
as the associated regular domain) is {\it non-degenerate} if it is
2-dimensional.

Given $p$ in $\Uu$ we have seen that $\Uu$ is contained in the future
of the plane passing through $r(p)$ and orthogonal to $N(p)$. So we
obtain the following inequality, that we shall often use throughout
this work:
\begin{equation}\label{fund-ineq}
    \E{N(p)}{r(q)-r(p)}\leq 0
\end{equation}
for every $p,q\in\Uu$. The equality holds if and only if the segment
$[r(p),r(q)]$ is contained in $\partial\Uu$.

Let us point out another meaningful inequality that immediately
descends from~(\ref{fund-ineq})
\[
  \E{N(p)-N(q)}{r(p)-r(q)}\geq 0\qquad\textrm{ for every }p,q\in\Uu \,.
\]
Again the equality holds if and only if the segment $[r(p),r(q)]$ is
contained in $\partial\Uu$.\\

By means of this inequality we may prove that the restriction of the
Gauss map on the level surface $\Uu(a)$ is $1/a$-Lipschitz. Indeed, if
$c(t)$ is a Lipschitz path on $\Uu(a)$ then we have
\[
   c(t)=r(t)+aN(t)
\]
where $r(t)=r(c(t))$ and $N(t)=N(c(t))$ are Lipschitz.  By
inequality~(\ref{fund-ineq}) we see that $\E{\dot r}{\dot N}\geq 0$
almost everywhere. So we deduce that
\[
   |\dot c|\geq a |\dot N|
\]
almost everywhere. This inequality shows that
$N:\Uu(a)\rightarrow\mh^2$ is $1/a$-Lipschitz.\\

The \emph{initial singularity} $\Sigma$ is the image of $r$. It
coincides with the set of points in the boundary of $\Uu$ admitting a
spacelike support plane. Since $\Uu$ is a regular domain, for every
point $p\in\partial\Uu$ there exists a future complete null ray
starting from $p$ and contained in $\partial\Uu$.

For $r_0\in\Sigma$ the set $\Ff(r_0)=N(r^{-1}(r_0))$ can be regarded
as the set of directions orthogonal to some spacelike support plane
through $r_0$. In fact also $\Ff(r_0)$ turns to be a straight convex
set. Moreover, if $r_1,r_2$ are two points in $\Sigma$, then $r_2-r_1$
is a spacelike vector and the orthogonal plane cuts $\mh^2$ along a
geodesic that divides $\Ff(r_1)$ from $\Ff(r_2)$ (this follows
from~(\ref{fund-ineq})).  In particular, we have that $H_\Uu$ is
decomposed in a union of straight convex sets that do not meet each
other transversally. 
\smallskip

We describe now a natural length-space structure on $\Sigma$. Since
the retraction $r:\Uu(1)\rightarrow\Sigma$ is a Lipschitz map,
$\Sigma$ is connected by spacelike Lipschitz paths. Thus, we can
consider the (a priori) pseudo-distance $\delta$ of $\Sigma$ defined
by
\[
   \delta(x,y)=\inf\{\ell(c)|\ c \textrm{ Lipschitz path in
   }\Sigma\textrm{ joining }x\textrm{ to }y\}
\]
where $\ell(c)$ denotes the length of $c$.\par In fact it is a
\emph{distance}, and $(\Sigma,\delta)$ can be regarded as the limit of
the level surfaces $\Uu(a)$ of $T$ as $a\rightarrow 0$. More precisely
we have

\begin{teo}\cite{Bo}
The function $\delta$ is a distance. Moreover if $\delta_a$ denotes
the distance on $\Uu(a)$ we have that $\delta_a\rightarrow\delta$ as
$a\rightarrow 0$ in the following sense.\par Given $p\in\Uu(1)$ let us
set $r_a(p)=r(p)+aN(p)$ the we have
\[
   \delta_a(r_a(p), r_a(q))\rightarrow \delta(r(p),r(q))
\]
and the convergence is uniform on the compact sets of $\Uu(1)$.
\end{teo}

In Section \ref{more:initial:sing}, we will analyze more carefully
this length-space structure of the initial singularity. 

\section{Measured geodesic laminations on straight convex sets}
\label{laminations}

\subsection{Laminations}
A {\it geodesic lamination} $\Ll$ on a complete Riemannian surface $S$
with geodesic boundary is a {\it closed} subset $L$ (the {\it support}
of the lamination), which is foliated by complete geodesics
(\emph{leaves} of the lamination).  More precisely

1) $L$ is covered by boxes $B$ with a product structure $B \cong
[a,b]\times [c,d]$, such that $L\cap B$ is of the form $X\times
[c,d]$, and for every $x\in X$, $\{x\}\times [c,d]$ is a geodesic arc.

2) The product structures are compatible on the intersection of two
   boxes.

3) The boundary of $S$ is a subset of $L$ and each boundary component is a
   leaf.

Each leaf admits an arc length parametrization defined on the whole
real line $\R$. Either this parametrization is injective and we call
its image a {\it geodesic line} of $S$, or its image is a simple
closed geodesic.  In both cases, we say that they are {\it simple}
(complete) geodesics of $S$.

The leaves together with the connected components of $S\setminus L$
make a {\it stratification} of $S$.
\smallskip

We specialize the discussion to geodesic laminations on 
non-degenerate straight convex 
sets $H$ in $\mh^2$. These sets  can be equivalently
characterized (up to isometry) as the simply-connected complete
hyperbolic surfaces with geodesic boundary.   

We claim that if a closed subset $L$ of $H\subset
\mh^2$ is the disjoint union of complete geodesics, say
$\Ll=\{l_i\}_{i\in I}$, and the boundary components of $H$ are in
$\Ll$, then $\Ll$ is a foliation of $L$ in the above sense.

In fact, fix a point $p_0\in L$ and consider a geodesic arc $c$
transverse to the leaf $l_0$ through $p_0$.  There exists a
neighbourhood $K$ of $p$ such that if a geodesic $l_i$ meets $K$ then
it cuts $c$. Orient $c$ arbitrarily and orient any geodesic $l_i$
cutting $c$ in such a way that respective positive tangent vectors at
the intersection point form a positive base.  Now for $x\in L\cap K$
define $v(x)$ as the unitary positive tangent vector of the leaf
through $x$ at $x$.  The following lemma ensures
that $v$ is a $1$-Lipschitz vector field on $L\cap K$ (see \cite{Ep-M}
for a proof).
\begin{lem}
Let $l,l'$ be disjoint geodesics in $\mh^2$. Take $x\in l$ and $x'\in
l'$ and unitary vectors $v,v'$ respectively tangent to $l$ at $x$ and
to $l'$ at $x'$ pointing in the same direction. Let $\tau(v')$ the
parallel transport of $v'$ along the geodesic segment $[x,x']$ then
\[
     || v-\tau(v') ||< d_{\mh}(x,x')
\]
where $d_{\mh}$ is the hyperbolic distance.
\end{lem}
\cvd 
\smallskip

\noindent
Thus there exists a $1$-Lipschitz vector field $\tilde v$ on
$K$ extending $v$. The flow $\Phi_t$ of this field allows us to
find a box around $p_0$.  Indeed for $\eps$ sufficiently small the
map
\[
  F: c\times (-\eps, \eps)\ni (x,t)\mapsto\Phi_t(x)\in\mh^2
\]
creates a box around $p_0$.
\smallskip

\begin{exa}\label{lamrd}
\emph{ Let $\Uu$ be a regular domain of $\mx_0$, $H=H_\Uu$ be the closure
of the image of the Gauss map $N$ of $\Uu$, and $\Sigma$ be the initial
singularity. Assume that dim $H = 2$. In the previous section
we have seen that for $r_0\in\Sigma$ the set $\Ff(r)=N(r^{-1}(r_0))$
is a straight convex set contained in $H$ and $\Ff(r)$ and $\Ff(s)$ do
not meet transversally. Thus geodesics that are either boundary
components of $H$, or boundary components of some $\Ff(r)$ or that
coincide with some $\Ff(r)$ are pairwise disjoint and provide an
example of a geodesic lamination, say $\Ll_\Uu$ of $H$.}
\end{exa}


\subsection{Transverse measures} \label{laminazioni:meas:sec}
Given a geodesic lamination $\Ll$ on a complete Riemannian surface $S$
with geodesic boundary, a rectifiable arc $k$ in $S$ is {\it
transverse} to the lamination if for every point $p\in k$ there exists
a neighbourhood $k'$ of $p$ in $k$ that intersects each leaf in at
most a point and each $2$-stratum in a connected set.  A {\it
transverse measure $\mu$} on $\Ll$ is the assignment of a positive
measure $\mu_k$ on each rectifiable arc $k$ transverse to $\Ll$ (this
means that $\mu_k$ assigns a non-negative {\it mass} $\mu_k(A)$ to
every Borel subset of the arc, in a countably additive way) in such a
way that:

(1) The support of $\mu_k$ is $k\cap L$;
\smallskip

(2)  If $k' \subset k$, then $\mu_{k'} = \mu_k|_{k'}$;
\smallskip

(3) If $k$ and $k'$ are homotopic through a family of arcs transverse to
$\Ll$, then the homotopy sends the measure $\mu_k$ to $\mu_{k'}$;
\smallskip

(4) $\mu_k(k)=+\infty$ if and only if $k\cap\partial S\neq\varnothing$.
\medskip

A {\it measured geodesic lamination on $S$} is a pair $(\Ll,\mu)$,
where $\Ll$ is a geodesic lamination and $\mu$ is a transverse measure
on $\Ll$. 

Let us specialize it again to  non-degenerate straight convex sets 
$H$ in $\mh^2$.

\begin{remark}
{\rm Notice that if $k$ is an arc transverse to a lamination of $H$
 there exists a transverse piece-wise geodesic arc homotopic to $k$
 through a family of transverse arcs.  Indeed there exists a finite
 subdivision of $k$ in sub-arcs $k_i$ for $i=1,\ldots,n$ such that
 $k_i$ intersects a leaf in a point and a $2$-stratum in a sub-arc. If
 $p_{i-1},p_i$ are the end-points of $k_i$ it is easy to see that each
 $k_i$ is homotopic to the geodesic segment $[p_{i-1},p_i]$ through a
 family of transverse arcs. It follows that a transverse measure on a
 lamination of $H$ is determined by the family of measures on
 transverse geodesic arcs.}\par
\end{remark}

\begin{remark}
{\rm While a geodesic lamination on $H$ can be eventually regarded as
  a particular lamination on $\mh^2$, condition $(4)$ ensures that a
  {\it measured} geodesic lamination on $H$ cannot be extended beyond $H$.}
\end{remark}

\begin{remark}\label{deg-lam}
{\rm We could include in the picture the {\it degenerate} straight
convex sets $H$ formed by a single geodesic; in this case the measured
lamination consists of a single $+\infty$-weighted leaf ($\Ll$
coincides with the whole of $H$). This {\it degenerate lamination} can
be regarded as the ``limit'' of measured geodesic laminations on
non-degenerate convex sets $H_n$, when $H_n$ tends to a geodesic. As
we are going to see this terminology is convenient because many
constructions we will implement on measured laminations on non-degenerate
straight convex sets easily extend to the degenerate lamination. However, the
degenerate case will be fully treated in Chapter \ref{QD}.}
\end{remark}

There is an equivalent way to describe a transverse measure on a
lamination $\Ll$ in terms of boxes of $\Ll$.

For each box $B=[a,b]\times [c,d]$ consider a positive measure $\mu_B$
on $[a,b]$ such that

(1) The support of $\mu_B$ is the set of $t$'s such that
    $\{t\}\times[c,d]$ is a leaf of $\Ll$.

(2) The total mass of $[a,b]$ is $+\infty$ iff either $\{a\}\times[c,d]$ or
    $\{b\}\times [c,d]$ are boundary leaves.

(3) If $B'=[a',b']\times[c',d']$ is a sub-box of $B=[a,b]\times[c,d]$, then
    $\mu_{B'}=\mu_B|_{[a',b']}$.

The proof that this definition is equivalent to the
original one is left to the reader.
\smallskip

The simplest example of a {\it measured lamination} $
(\Ll,\mu)$ on $\mh^2$ is given by any {\it finite} family of disjoint
geodesic lines $l_1,\dots,l_s$, each one endowed with a {\it real
positive weight}, say $a_j$. A relatively compact subset $A$ of an arc
$k$ transverse to $\Ll$ intersects it at a finite number of
points, and we set $\mu_k(A)=\sum_i a_i|A\cap l_i|$. 

The simplest example of a \emph{measured lamination on $H$} is given
by the boundary of $H$ such that each leaf carries the weight $+\infty$
(that is $\mu_k(A)=0$ except if $A\cap\partial H\neq\varnothing$ and
in that case $\mu_k(A)=+\infty$). Notice that by condition $(4)$ the
weight of each boundary curve is necessarily $+\infty$. In fact
boundary curves carry an infinite weight whenever the lamination is
locally finite.

More generally if $(\Ll,\mu)$ is a lamination on $H$, then for
each boundary curve, say $l$, two possibilities can happen: given a
geodesic arc, say $k$, with an end-point in the interior of $H$ and an
end-point $p_0\in l$ then either $\mu_k(k\setminus\{p_0\})=+\infty$ or
$\mu_k(k\setminus\{p_0\})<+\infty$. In the latter case we say that $l$
is a \emph{weighted boundary leaf} (of weight $+\infty$).

Let us point out two interesting subsets of $\Ll$ associated to a
measured geodesic lamination $(\Ll,\mu)$ on $H$. The {\it
simplicial part} $\Ll_S$ of $\Ll$ consists of the union of the
isolated leaves of $\Ll$. $\Ll_S$ does not depend on the measure
$\mu$. In general this is not a sub-lamination, that is its
support $L_S$ is not a closed subset of $H$.  

 A leaf, $l$, is called \emph{weighted} if there exists a transverse
arc $k$ such that $k\cap l$ is an atom of $\mu_k$.  By property (3) of
the definition of transverse measure, if $l$ is weighted then for
every transverse arc $k$ the intersection of $k$ with $l$ consists of
atoms of $\mu_k$ whose masses are equal to a positive number $A$
independent of $k$. We call this number the weight of $l$.  The {\it
weighted part} of $\lambda$ is the union of all the weighted leaves.
It depends on the measure and it is denoted by $\Ll_W=\Ll_W(\mu)$.

Since every compact set $K\subset\mathring{H}$ ($\mathring{H}$ being
the interior part of $H$) intersects finitely many weighted leaves
with weight bigger than $1/n$, it follows that $\Ll_W$ is a countable
set.  On the other hand, it is not in general a sub-lamination of
$\Ll$. For instance, consider the case $H=\mh^2$ and take the set of
geodesics $\Ll$ with a fixed end-point $x_0\in S^1_\infty$. Clearly
$\Ll$ is a geodesic lamination of $\mh^2$ and its support is the whole
of $\mh^2$.  In the half-plane model suppose $x_0=\infty$ so that
geodesics in $\Ll$ are parametrized by $\mr$. Let $l_t$
denote the geodesic in $\Ll$ with end-points $t$ and $\infty$.  If we
choose a dense sequence $(q_n)_{n\in\mn}$ in $\mr$ it is not difficult
to construct a measure on $\Ll$ such that $l_{q_n}$ carries the weight
$2^{-n}$.  For that measure $L_W$ is a dense subset of $\mh^2$.
 
As $L$ is the support of $\mu$, then we have the inclusion
$\Ll_S\subset\Ll_W(\mu)$. The previous example shows that in general 
this inclusion is strict.
\begin{remark}\emph{
The word  \emph{simplicial} 
mostly refers to the ``dual'' geometry of the initial
singularity of the spacetimes that we will associate to the measured
geodesic laminations on $H$.  In fact, it turns out that
when $\lambda$ coincides with its simplicial part, then  
the initial singularity is a simplicial metric tree.
}\end{remark}

\subsection{Standard finite approximation}
Given a measured geodesic lamination $\lambda=(\Ll,\mu)$ on $H$, and
a compact set $K\subset\mathring{H}$, we say that a sequence of measured
geodesic lamination $\lambda_n$ converges to $\lambda$ on $K$ if for
any arc $k\subset K$ transverse to $\lambda$ and with no end-point in
$L_W(\mu)$ we have

(1) $k$ is transverse to $\lambda_n$ for big $n$;

(2) for any continuous function $\varphi$ on $k$ 
\[
    \lim_{n\rightarrow+\infty}\int_k\varphi \mathrm
    d(\mu_n)_k=\int_k\varphi\mathrm d\mu_k\,.
\]

In this work we will need to construct a sequence of finite
laminations converging to a given lamination $\lambda$ on some compact
set $K$.

The following construction ensures that such a sequence exists if $K$ is a box
of $\lambda$.

Let us fix a box $B=[a,b]\times[c,d]$ for $\lambda$ contained
in the interior of $H$.

Fix $n$ and subdivide $[a,b]$ into the union of intervals $c_1,\ldots,
c_r$ such that $c_i$ has length less than $1/n$ and all the end-points
of $c_i$, but $a$ and $b$, are not in $L_W(\mu)$.  For every $c_j$ let
us set $\alpha_j=\mu_{k}(c_j)$. If $\alpha_j>0$, then choose a
leaf $l_j$ of $\Ll$ that cuts $c_j$, with the only restriction that if
$a$ (resp. $b$) is an atom of $\mu$, then $l_1$ (resp. $l_r$) is the
weighted leaf along it. Thus consider the finite lamination
$\Ll_n=\{l_j | a_j>0\}$ on $\mh^2$ and associate to every $l_j$ the
weight $\alpha_j$.  In such a way we define a measure $\mu_n$ transverse to
$\Ll_n$.

\begin{lem}
$\lambda_n$ converges to $\lambda$ on $B$.
\end{lem}
\Dim 
It is clear that any arc $k\subset B$ transverse to $\lambda$
is transverse to $\lambda_n$ for any $n$. So (1) is verified.

On the other hand any arc $k\subset B$ transverse to $\lambda$ is a
finite composition of transverse arcs $k_1,\ldots, k_N$ such that each
$k_i$ has no end-point on $L_W(\mu)$ and is homotopic through a family
of transverse arcs to a horizontal arc. Thus it is sufficient to check
(2) for a horizontal arc with no endpoint on $L_W(\mu)$ and this is
straightforward.  \cvd 

The sequence $\lambda_n$ is called a
{\it standard approximation of $\lambda$}

\begin{remark}
{\rm There is a natural topology on the set of measured geodesic
laminations on $\mh^2$, obtained by considering a measured geodesic
lamination as a \emph{geodesic current} on $\mh^2$ \cite{Bon}(3).
With respect to this topology, the sequence $(\lambda_n)$ of standard
approximation of $\lambda$ on $B$ converges to the measured geodesic
lamination,say $\lambda|_B$, obtained by considering only leaves of
$\lambda$ that intersect $B$.}
\end{remark}

\subsection {$\Gamma$-invariant  measured geodesic laminations}
\label{Gamma:inv:lam}
Let us generalize the notion of straight convex set.  Let
$F=\mh^2/\Gamma$ be a complete hyperbolic surface.
A {\it straight convex set $Z$ in $F$} is a closed convex surface with
geodesic boundary embedded in $F$. Recall that a subset $K$ of $F$ is
said convex if every geodesic segment with end-points in $K$ is
contained in $K$.  It is not hard to see that $Z$ is a straight convex
set of $F$ if and only if the pull-back $H$ of $Z$ in $\mh^2$ is a
$\Gamma$-invariant straight convex set of $\mh^2$. Conversely every
$\Gamma$-invariant straight convex set of $\mh^2$ projects to a
straight convex set of $F$. Moreover, the interior of $Z$ is
homeomorphic to $F$. Measured geodesic laminations on straight convex
sets of $F$ lift to $\Gamma$-invariant measured geodesic laminations
of straight convex sets of $\mh^2$ and this correspondence is actually
bijective.

\begin{remark}
{\rm In general, $F$ contains several straight convex sets $Z$; all of
  them contain the convex core of $F$. On the other hand, if $Z$ has
  finite area and all the boundary components are closed geodesics,
  then it coincides with the convex core of $F$, so in that case $Z$
  is determined by $\Gamma$.}
\end{remark}

\paragraph{ Cocompact $\Gamma$}
The simplest example of a measured geodesic lamination on a compact closed
(\emph{i.e.} without boundary) hyperbolic surface
$F=\mh^2/\Gamma$ is a finite family of disjoint, weighted simple
closed geodesics on $F$.  This lifts to a $\Gamma$-invariant measured
lamination of $\mh^2$ made by a countable family of weighted geodesic
lines, that do not intersect each other on the {\it whole}
$\overline{\mh}^2= \mh^2 \cup S^1_\infty$.  The measure is defined
like in the case of a finite family of weighted geodesics. We call
such special laminations {\it weighted multi-curves}.
\smallskip

When $F=\mh^2/\Gamma$ is compact closed , the $\Gamma$-invariant measured
geodesic laminations $(\Ll,\mu)$ on $\mh^2$ have particularly
good features, that do not hold in general. We limit ourselves to
recall of a few of them:
\begin{enumerate}

\item
The lamination $\Ll$ is determined by its support $L$.
The support $L$ is a no-where dense set of null area.

\item 
The simplicial part $\Ll_S$ and the weighted part $\Ll_W$ 
actually coincide; moreover $\Ll_S$ is the maximal weighted
multi-curve sub-lamination of $\lambda$.

\item
Let $\Mm\Ll(F)$ denote the space of measured geodesic
laminations on $F$. It is homeomorphic to $\mr^{6g-6}$, $g\geq 2$
being the genus of $F$. Any homeomorphism $f: F\to F'$ of hyperbolic
surfaces, induces a natural homeomorphism $f_\Ll: \Mm\Ll(F)\to
\Mm\Ll(F')$; if $f$ and $f'$ are homotopic, then $f_\Ll=f'_\Ll$. This
means that $\Mm\Ll(F)$ is a topological object which only depends on
the genus of $F$, so we will denote it by $\Mm\Ll_g$.
Varying $[F]$ in the Teichm\"uller space $\Tt_g$, the above
considerations allows us to define a trivialized fiber bundle 
$\Tt_g\times \Mm\Ll_g$ over $\Tt_g$ with fiber $\Mm\Ll_g$.
\end{enumerate}

\begin{exa}\label{genlam}
{\rm (1) Mutata mutandis, similar facts hold more generally when
$F=\mh^2/\Gamma$ is homeomorphic to the interior of a compact surface
$S$, possibly with non empty boundary, providing that the lamination
on $F$ has {\it compact support}. However, even when $F$ is of {\it
finite area} but non compact, we can consider laminations that do not
necessarily have compact support (see Section \ref{3cusp}).
\smallskip

(2) Let $\gamma$ be either a geodesic line or a horocycle in $\mh^2$.
  Then, the geodesic lines orthogonal to $\gamma$ make, in the
  respective cases, two different geodesic foliations both having the
  whole $\mh^2$ as support. We can also define a transverse measure
  $\mu$ which induces on $\gamma$ the Lebesgue one.}
\end{exa}

Recall the sets $$\Mm\Ll = \{\lambda = (H,\Ll,\mu)\}$$ and
$$\Mm\Ll^\Ee = \{(\lambda,\Gamma)\}$$ already defined in Section
\ref{FLAT:I} of Chapter ~\ref{INTRO}.  Roughly, the first consists of
the {\it measured geodesic laminations defined on some non-degenerate
straight convex set of $\mh^2$}; the second covers the set of {\it
measured geodesic laminations defined on some straight convex set ($Z=
H/\Gamma$) in some complete hyperbolic surface ($F=\mh^2/\Gamma$)}.
Recall that there are natural actions of $SO(2,1)$ on $\Mm\Ll$ and
$\Mm\Ll^\Ee$.

\section {From measured geodesic laminations towards 
flat regular domains}\label{ML_REGD} Recall that $\Rr$ denotes the set
of non-degenerate regular domains in $\mx_0$.  The group
$\ISO^+(\mx_0)$ (hence the translations subgroup $\mr^3$) naturally
acts on it. In this section we show how a general construction
produces a map
$$ \Uu^0: \Mm\Ll \to \Rr,\ \ \lambda \to \Uu^0_\lambda $$ 
that induces maps (for simplicity we keep the same notation)
$$ \Uu^0: \Mm\Ll \to \Rr/\mr^3$$
and
$$\Uu^0: \Mm\Ll/SO(2,1) \to \Rr/\ISO^+(\mx_0)\ . $$ 

\begin{figure}
\begin{center}
\input{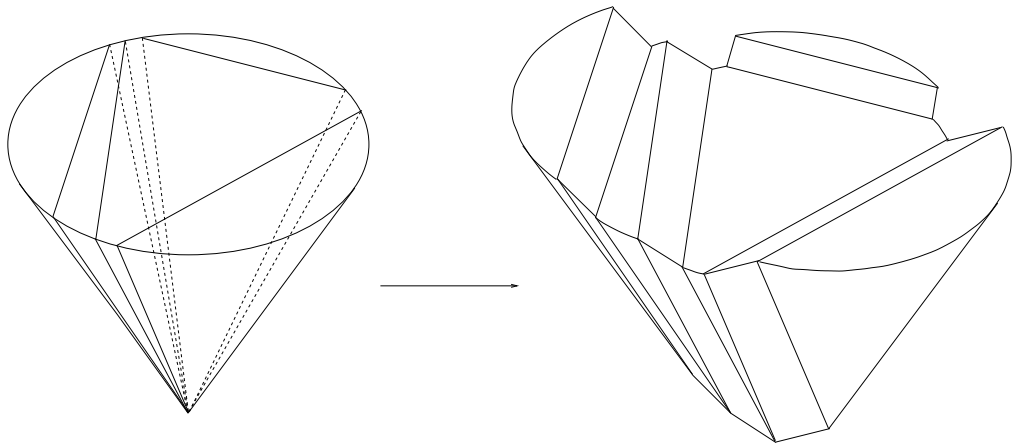}
\caption{{\small An example of regular domain associated to a finite
measured lamination.}}
\end{center}
\end{figure}
For the case $H=\mh^2$ see \cite{M}, \cite{Bo1,Bo}; in fact the
construction extends with minor changes to arbitrary $H$.  Here we
sketch this construction.

\subsection{Construction of regular domains}\label{construction:rd:ssec}
Fix a base-point $x_0\in\mathring{H} \setminus L_W$.  For every
$x\in\mathring H \setminus L_W$ choose an arc $c$ transverse to $\Ll$ with
end-points $x_0$ and $x$. For $t\in c\cap L$, let $v(t)\in \R^3$
denote the unitary spacelike vector tangent to $\mh^2$ at $t$,
orthogonal to the leaf through $t$ and pointing towards $x$. For $t\in
c\setminus L$, let us set $v(t)=0$.  In this way we define a function
\[
    v:c\rightarrow\mr^3
\]
that is continuous on the support of $\mu\,$.
We can define
\[
\rho(x)=\int_{c}v(t)\d\mu(t).
\]
It is not hard to see that $\rho$ does not depend on the path $c$.
Moreover, it is constant on every stratum of $\lambda$ and it is a
continuous function on $\mathrm{H}\setminus L_W$.  

The domain $\Uu^0_\lambda$
can be defined in the following way
\[
  \Uu^0_\lambda=\bigcap_{x\in\mathrm{ H}\setminus
  L_W}\fut(\rho(x)+\ort{x})\ .
\]

Firstly let us prove that $\Uu^0_\lambda$ coincides with the
intersection of its null support planes.

Since $\rho$ is constant on every stratum $F$, for every null vector
$v$ that represents some accumulation point of $F$ in $\overline\mh^2$
the plane $\rho(x)+\ort v$ is a support plane for $\Uu^0_\lambda$
(where $x$ is any point in $F$).  More precisely, if $\partial_\infty
F$ denotes the set of accumulation points of $F$ in $\partial\mh^2$ ,
and $x_F$ is a point in $F$ we have that
\begin{equation}\label{nullplanes}
   \Uu^0_\lambda=\bigcap_{F\textrm{ stratum of
   }\Ll}\ \ \ \bigcap_{[v]\in\partial_\infty F}\fut(\rho(x_F)+\ort v)\,.
\end{equation}
The inclusion $(\subset)$ follows from the above remark. To show that
also the opposite inclusion holds, notice that every $x\in F$ is the
barycentric combination of some null vectors representing points in
$\partial_\infty F$. Thus
\[
    \bigcap_{[v]\in\partial_\infty F}\fut(\rho(x)+\ort v)\subset
    \fut(\rho(x)+\ort x)
\]
The intersection on all $x$'s of both sides shows the inclusion $(\supset)$.

Eventually, in order to show that $\Uu^0_\lambda$ is a regular domain it is
sufficient to prove that it is not empty.

Given $x,y\in\mathring{H}\setminus L_W$ we have
\[
\rho(y)-\rho(x)=\int_{[x,y]}v(t)\mathrm d\mu(t)\,.
\]
Since $\E{v(t)}{y}\geq 0$ for every $t\in [x,y]$,
the following inequality holds
\[
   \E{\rho(x)-\rho(y)}{y}\leq 0
\]
and the equality holds iff $x$ and $y$ lie in a same stratum of $\Ll$
(see~\cite{Bon}).

This inequality implies that $\rho(x)\in\partial\Uu^0_\lambda$ for every
$x\in\mathring{H}\setminus L_W$.

\subsection{Cosmological time on $\Uu^0_\lambda$}\label{ctrd:sec}
First let us determine the image of the Gauss map $N$ of
$\Uu^0_\lambda$.

The definition of $\Uu^0_\lambda$ implies that $\rho(x)+\ort x$ is a
support plane for $\Uu^0_\lambda$, for every $x\in\mathring H\setminus
L_W$. Thus $\mathring H\setminus L_W$ is contained in ${\rm Im} N$. Since
${\rm Im} N$ is a convex set it actually contains the whole $\mathring H$.
More precisely suppose $x$ belongs to a weighted leaf $l$ contained in
the interior of $H$, with weight $A$. Then we can consider a geodesic
arc $c$ starting from $x_0$ (the base point) and passing through
$x$. Then there exist left and right limits
\[
   \begin{array}{l}
    \rho_-(x)=\lim_{t\rightarrow x^-}\rho(t)\\
    \rho_+(x)=\lim_{t\rightarrow x^+}\rho(t)
    \end{array}
\]
and the difference $\rho_+(x)-\rho_-(x)$ is the spacelike vector with
norm equal to $A$ orthogonal to $l$ pointing as $c$. Notice that the vector
$\rho_+(x)-\rho_-(x)$ depends only on the leaf through $x$.
Clearly the plane passing through $\rho_-(x)$ and orthogonal to $x$ is a
support plane for $\Uu^0_\lambda$ (and in fact such a plane contains also
$\rho_+(x)$).

Now take $x\in\partial H$. If the leaf through $x$, say $l$, is a
weighted leaf then for $t\in[x_0,x]$ going to $x$ we have
\[
   \rho(t)\rightarrow\rho(x):= \int_{[x_0,x)}v(s)\mathrm d\mu(s)\,.
\] 
Thus $\rho(x)+\ort x$ is a support plane for $\Uu^0_\lambda$ and
$x\in {\rm Im} N$.

If $l$ is not weighted, we have
$|\rho(t)|\rightarrow+\infty$ as $t\rightarrow x$. Indeed since for
$t,s\in[x_0,x)\cap\Ll$ the geodesics orthogonal to $v(s)$ and $v(t)$
are disjoint,the plane generated by $v(s)$ and $v(t)$ is not spacelike
so the reverse of the Schwarz inequality holds, that is
$\E{v(t)}{v(s)}\geq 1$ (the scalar product of $v(s)$ and $v(t)$ is
positive because they point in the same direction).  As a consequence,
we have
\begin{equation}\label{divrd}
  |\rho(t)|^2=\E{\int_{[x_0,t]} v(s)}{\int_{[x_0,t]}v(s)}=
   \int_{[x_0,t]}\int_{[x_0,t]}\E{v(s)}{v(\sigma)}\geq(\mu([x_0,t]))^2\,.
\end{equation}
Suppose there exists a point $\rho\in\partial\Uu^0_\lambda$ such that
$\rho+\ort x$ is a support plane. By~(\ref{divrd}), the segment
$[\rho, \rho(t)]$ converges to a ray $R$ starting at $\rho$. More
precisely, by~(\ref{fund-ineq}) we have
\[
\E{\rho-\rho(t)}{x}<0 \qquad\qquad  \E{\rho-\rho(t)}{t}>0
\]
so the ray $R$ is $\rho+\mr_{\leq 0} u$ where $u$ is a tangent vector
of $\mh^2$ at $x$, orthogonal to $l$ and pointing outside $H$.  This
gives a contradiction: in fact $\E{u}{y}<0$ for $y\in H$, so
$\E{\rho+tu}{y}\rightarrow+\infty$ as $t\rightarrow-\infty$.

An analogous argument shows that points outside $H$ do not belong to ${\rm Im} N$. Eventually, the image of $N$ is $\mathring{H}\cup(\partial
H\cap L_W)$.

Now given $x\in {\rm Im} N$ we are going to determine $N^{-1}(x)$.  By the
characterization of the retraction given in
Proposition~\ref{COSMT_REGD}, we have that $N^{-1}(x)$ is the union of
timelike rays parallel to $x$ starting at points in
$\partial\Uu^0_\lambda\cap P_x$, $P_x$ being the spacelike support
plane orthogonal to $x$.

The following sentences can be proved as  in Lemma 8.7 of~\cite{Bon} 

1) If $x\in\mathring{H}\setminus L_W$ then
   $P_x\cap\partial\Uu^0_\lambda=\{\rho(x)\}$.

2) If $x\in\mathring{H}\cap L_W$ then
   $P_x\cap\partial\Uu^0_\lambda=[\rho_-(x),\rho_+(x)]$.

3) If $x\in\partial H\cap L_W$ then
   $P_x\cap\partial\Uu^0_\lambda=\rho(x)+\mr_{\geq 0} u(x)$, where $u(x)$ is
   the unit outward normal vector.

\begin{figure}
\begin{center}
\input{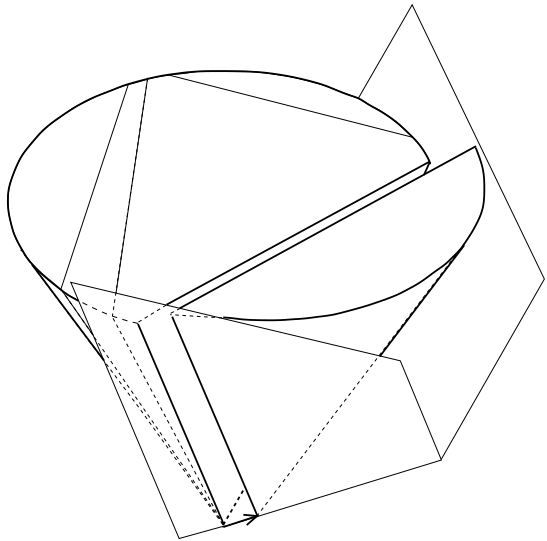}
\caption{{\small The construction of the domain associated to a finite
lamination.}}
\end{center}
\end{figure}

By (\ref{flatdec}) we have
\[
   \begin{array}{l}
   \Uu^0_\lambda=\{ax+\rho(x)|x\in H\setminus(\partial H\cup L_W)\}\cup\\
   \cup\{ax+t\rho_+(x)+(1-t)\rho_-(x)|x\in L_W\setminus\partial H,\ t\in[0,1],\
   a>0\}\cup\\
   \cup\{ax+\rho(x)+tu(x)|x\in L_W\cap\partial H,\ t>0\,\ a>0\}
   \end{array}
\]

Moreover, for $x\in H\setminus(\partial H\cup
L_W)$ we have
\[
\begin{array}{l}
    T(ax+\rho(x))=a\\
    N(ax+\rho(x))=x\\
    r(ax+\rho(x))=\rho(x)\,.
\end{array}
\]
For $x\in L_W\setminus\partial H$ we have
\[
\begin{array}{l}
    T(ax+t\rho_+(x)+(1-t)\rho_-(x))=a\\
    N(ax+t\rho_+(x)+(1-t)\rho_-(x))=x\\
    r(ax+t\rho_+(x)+(1-t)\rho_-(x))=t\rho_+(x)+(1-t)\rho_-(x)\,.
\end{array}
\]
Finally for $x\in\partial H\cap L_W$ we have
\[
\begin{array}{l}
    T(ax+\rho(x)+tu(x))=a\\
    N(ax+\rho(x)+tu(x))=x\\
    r(ax+\rho(x)+tu(x))=\rho(x)+tu(x)\,.
\end{array}
\]

\begin{remark}\emph{
The lamination associated to $\Uu^0_\lambda$ in 
Example~\ref{lamrd} is the support of $\lambda$. In fact we have that
$\Ff(\rho(x))$ is the stratum of $\lambda$ through $x$.}
\end{remark}

\begin{remark}\emph{
The domain associated to the degenerate lamination $\lambda_0$ (Remark
\ref{deg-lam}) is the future of the spacelike line dual to the
geodesic of $\lambda_0$.  With our convention the definition of
$\Uu^0_\lambda$ can be applied also in this case. }\end{remark}

\subsection {Measured geodesic laminations on the $T$-level surfaces}
\label{lam:supliv}
Let us fix a level surface $\Uu^0_\lambda(a)=T^{-1}(a)$ of the
cosmological time of $\Uu^0_\lambda$.  We want to show that
$\Uu^0_\lambda(a)$ also carries a natural measured geodesic lamination
on. Consider the closed set $\hL_a=N^{-1}(L)\cap \Uu^0_\lambda(a)$.
First we show that $\hL_a$ is foliated by geodesics.

If $l$ is a non-weighted leaf of $L$ we have that $\hat
l_a:=N^{-1}(l)\cap\Uu^0_\lambda(a)=a l+\rho(x)$ where $x$ is any point
on $l$. Since the map $N:\Uu^0_\lambda(a)\rightarrow\mh^2$ is
$1/a$-Lipschitz, $\hat l_a$ is a geodesic of $\Uu^0_\lambda(a)$.

If $l$ is a weighted leaf contained in $\mathring{H}$, 
then $N^{-1}(l)\cap\Uu^0_\lambda(a)= al +
[\rho_-(x),\rho_+(x)]$ where $x$ is any point on $l$. Thus we have
that $N^{-1}(l)$ is a Euclidean band foliated by geodesics $\hat
l_a(t):=al +t\rho_-(x)+(1-t)\rho_+(x)$ where $t\in[0,1]$.  

Finally if $l$ is a weighted boundary leaf then
$N^{-1}(l)\cap\Uu^0_\lambda(a)=a l+\rho(x)+\mr_{\geq 0} u$, where $x$
is any point of $l$ and $u=u(x)$ is the vector outward normal vector
at $x$. Thus $N^{-1}(l)$ is a Euclidean band (with infinite width)
foliated by geodesics $\hat l_a(t)=al+\rho(x)+tu$ for $t\geq0$.
 
Thus the set
\[
  \hL_a= \bigcup_{l\subset L\setminus L_W} \hat l_a 
\cup\bigcup_{l\subset L_W}\bigcup_{\ t\in [0,1]} \hat l_a(t)\cup
\bigcup_{l\subset L_W\cap\partial H}\bigcup_{\ t\geq 0}\hat l_a(t)
\]
is a geodesic lamination of $\Uu^0_\lambda(a)$.  

Given a Lipschitz path 
$c(t)$ transverse to this lamination we know that $r(t)=r(c(t))$ is a
Lipschitz map. Thus we have that $r$ is differentiable almost
every-where and
\[
  r(p)-r(q)=\int_c \dot r(t)\mathrm d t
\]
where $p$ and $q$ are the endpoints of $c$.  It is not hard to see
that $\dot r(t)$ is a spacelike vector. Thus we can define a measure
$\hat\mu_c=\E{\dot r(t)}{\dot r(t)}^{1/2}\mathrm dt$.  If $N(p)$ and
$N(q)$ are not in $L_W$ then
\begin{equation}\label{www}
  \hat\mu (c)=\mu(N(c)) \ .
\end{equation}
By this identity we can deduce that $(\hL_a,\hat\mu)$ is a measured
geodesic lamination on $\Uu^0_\lambda(a)$. The measure $\hat\mu_c$
defined on every rectifiable transverse arc is absolutely continuous
with respect to the Lebesgue measure of $c$. Moreover, inequality
~(\ref{fund-ineq}) implies that $\E{\dot r(t)}{\dot r(t)}\leq\E{\dot
c(t)}{\dot c(t)}=1$.

Hence the total mass of $\hat\mu_c$ is bounded by the length
$c$. Moreover, the density of $\hat\mu$ with respect the Lebesgue
measure of $c$ is bounded by $1$.

\begin{lem}\label{geo:est:lem}
Let us fix $p,q\in\Uu^0_\lambda(1)$ such that $N(p)$ and $N(q)$ are
not in $L_W$. If $p$ and $q$ are not in the same stratum, then the
geodesic segment in $\Uu^0_\lambda$ connecting them, say $\hat c$, is
a transverse path and its mass $\hat\mu(\hat c)$ is equal to the mass
(with respect to $\mu$) of the geodesic segment, say $c$ of $\mh^2$
connecting $N(p)$ to $N(q)$.  Moreover the following inequalities hold
\[
\begin{array}{l}
\hat\mu(\hat c)\leq \ell(\hat c)\\
\ell_{\mh^2}(c)\leq \ell(\hat c).
\end{array}
\]
\end{lem}

\Dim Since $\Uu^0_\lambda(1)$ is convex, its curvature is
non-positive. Thus there exists a unique geodesic segment $\hat c$
joining $p$ to $q$. Clearly it cannot be tangent to any leaf of
$\lambda$ (since leaves are geodesics). The same argument shows that
it intersects each leaf at most once. Since each leaf of $\hat\lambda$
disconnects $\Uu^0_\lambda$ in two half-planes, $\hat c$ intersects
only the leaves that disconnects $p$ from $q$. Thus $N(\hat c)$ is
homotopic to $c$ through a family of transverse arcs.  This proves
that $\hat\mu(\hat c)=\mu(c)$.

Since the density of $\hat\mu$ is bounded by $1$, the first inequality
follows. The second inequality is a consequence of the fact that the
map $N:\Uu^0_\lambda(1)\rightarrow\mh^2$ is $1$-Lipschitz.  \cvd

\subsection {Continuous dependence of $\ \Uu^0_\lambda$}\label{continuity} 

We discuss  how the construction of $\Uu^0_\lambda$ 
continuously depends on $\lambda$ (see \cite{Bo1, Bo}).
 
Fix a compact connected domain $K\subset\mathring H$ containing the base point
$x_0$.  Suppose that $\lambda_n$ is a sequence of measured geodesic
laminations such that $\lambda_n\rightarrow\lambda_\infty$ on $K$.\\
We shall denote by $\Uu_n$ (resp. $\Uu_\infty$) the domain associated
to $\lambda_n$ (resp. $\lambda_\infty$) and by $T_n,r_n,N_n$
(resp. $T_\infty,r_\infty, N_\infty$) the corresponding cosmological
time, retraction and Gauss map.

\begin{prop}\label{piatto:conv:prop}
For any pair of positive numbers $a<b$ let $U(K;a,b)$ be the set of
points $x$ in $\Uu_\infty$ such that $a<T_\infty(x)<b$ and
$N_\infty(x)\in K$. We have
\begin{enumerate}
\item
  $U(K;a,b)\subset\Uu_n$ for $n\gg 0$;
\item
  $T_n\rightarrow T_\infty$ in $\mathrm C^1(U(K;a,b))$;
\item
  $N_n\rightarrow N_\infty$ and $r_n\rightarrow r_\infty$ uniformly on
  $U(K;a,b)$.
\end{enumerate}    
\end{prop}
\Dim  
For any $x\in K \setminus
(L_\infty)_W$
\begin{equation}\label{zzz}
   \int_{c_{x_0,x}}v_n(t)\d\mu_n(t)\rightarrow\rho_\infty(x)
\end{equation}
where $c_{x_0,x}$ is any transverse path contained in $K$ joining
$x_0$ to $x$ and $v_n(t)$ is the orthogonal field of $\Ll_n$. Indeed such a
field is $C$-Lipschitz on $L_n\cap c_{x_0,x}$, for some $C$ that
depends only on $K$. Thus we can extend $v_n|_{L_n\cap c_{x_0,x}}$ to a
$C$-Lipschitz field $\tilde v_n$ on $c_{x_0,x}$. Clearly
\[
    \int_{c_{x_0,x}}v_n(t)\d\mu_n(t)=\int_{c_{x_0,x}}\tilde
    v_n(t)\d\mu_n(t) \ .
\]
Possibly up to passing to a subsequence, we have that $\tilde
v_n\rightarrow \tilde v_\infty$ on $\mathrm C^0(c_{x_0,x})\in
\mr^3$. Since $\tilde v_\infty=v_\infty$
on $L_\infty\cap c_{x_0,x}$, (\ref{zzz}) follows. 

By this fact we can deduce the following result.
\begin{lem}\label{piatto:conv1:lem}
Let us take $p\in\Uu_\infty(a)$ such that $N_\infty(p)\in K$.  There
exists a sequence $p_n\in\Uu_n(a)$ such that $p_n\rightarrow
p_\infty$.
\end{lem}
\Dim
First assume that $x=N(p)\notin (L_\infty)_W$. Then we know that 
\[
  p= a x +\rho_\infty(x) \ .
\]
Hence $p_n= a x +\int_{c_{x_0,x}}v_n(t)\d\mu_n(t)$ works.

Now assume that $x=N(p)\in (L_\infty)_W$ so $p$ lies in a band
$al+[\rho_-(x),\rho_+(x)]$.  We can take two points $y,z\notin
(L_\infty)_W$ such that $||y-x||<\eps$ and $||z-x||<\eps$ and
$||\rho(y)-\rho_-(x)||<\eps$ and $||\rho(z)-\rho_+(x)||<\eps$ (where
$||\cdot||$ is the Euclidean norm).  If we put $q^-=a y+\rho(y)$ and
$q^+=a z+\rho(z)$, the Euclidean distance of $p$ from the
segment $[q^-, q^+]$ is less than $2(1+a)\eps$.  Now let us set
$q^-_n=ay+\rho_n(y)$ and $q^+_n=az+\rho_n(z)$ and choose $n$
sufficiently large such that $||q^\pm-q^\pm_n||<\eps$.  We have that
the distance of $p$ from $[q^-_n, q^+_n]$ is less than
$2(2+a)\eps$. On the other hand since the support planes for the
surface $\Uu_{n}(a)$ at $q^-_n$ and at $q^+_n$ are close it is easy to
see that the distance of any point on $[q^-_n, q^+_n]$ from $\Uu_n(a)$
is less then $\eta(\eps)$ and $\eta\rightarrow 0$ for $\eps\rightarrow
0$.  It follows that we can take a point $p_n\in\Uu_n(a)$ arbitrarily
close to $p$ for $n$ sufficiently large.  \cvd

Choose coordinates $(y_0,y_1,y_2)$ such that the coordinates of the
base point $x_0$ are $(1,0,0)$.  We have that the surface $\Uu_n(a)$
(resp. $\Uu_\infty(a)$) is the graph of a positive function
$\varphi_n^a$ (resp. $\varphi_\infty^a$) defined over the horizontal
plane $P=\{y_0=0\}$.  Moreover, $\varphi_n^a$ is a
$1$-Lipschitz convex function and $\varphi_n^a(0)=a$.  Thus
Ascoli-Arzel\`a Theorem implies that $\{\varphi_n^a\}_{n\in\mn}$ is a
pre-compact family in $\mathrm C^0(P)$.  Up to passing to a
subsequence, there exists a function $\varphi$ on $P$ such that
$\varphi_n^a\rightarrow\varphi$ as $n\rightarrow +\infty$.  Consider
the compact domain of $P$
\[
    P(K,a)=\{p\in P | N_\infty(\varphi_\infty^a(p),p)\in K\} \ .
\]
By Lemma~\ref{piatto:conv1:lem} it is easy to check that
$\varphi=\varphi_\infty^a$ on $P(K,a)$.  Thus we can deduce
\begin{equation}\label{piatto:conv2:eq}
  \varphi_n^a|_{P(K,a)}\rightarrow\varphi_\infty^a|_{P(K,a)} \ .
\end{equation}
Fix $b>a>\alpha$. The domain $U(K;a,b)$ is contained in the
future of the portion of surface $N_\infty^{-1}(K)\cap\Uu_\infty(\alpha)$.
By (\ref{piatto:conv2:eq}) we see that $U(K;a,b)$ is contained in the future
of $\Uu_n(\alpha)$ for $n$ sufficiently large. Thus we have
\begin{equation}\label{piatto:conv3:eq}
    U(K;a,b)\subset\Uu_n(\geq\alpha)\qquad\textrm{ for }n\gg 0 \ .
\end{equation}
Since we are interested in the limit behaviour of functions
$T_n,N_n,r_n$ we can suppose that $U(K;a,b)$ is contained in
$\Uu_n(\geq\alpha)$ for any $n$.\\
Hence, $T_n,N_n, r_n$ are defined on $U(K;a,b)$ for any $n$.
Moreover notice that 
\[
\begin{array}{l}
    T_n(\xi,0,0)=\xi \ ,\\
    N_n(\xi,0,0)=(1,0,0) \ ,\\
    r_n(\xi,0,0)=0 \ .
\end{array}
\]
Thus we have that $r_n(p)$ lies in the half-space $P^+=\{y_0>0\}$ for
every $p\in\Uu_n$.  Since $U(K;a,b)$ is compact then there exists a
constant $C$ such that for every $p\in U(K;a,b)$ and for every past
directed vector $v$ such that $p+v$ is in $P^+$ then $||v||<C$.  Since
$r_n(p)=p-T_n(p)N_n(p)$ we have that
$$||T_n(p)N_n(p)||<C$$ for every $n\in\mn$ and for every $p\in
U(K;a,b)$.  Since $T_n(p)\geq\alpha$ the following
property follows.
\begin{lem}\label{piatto:Gauss:lem}
The family $\{N_n\}$ is bounded in $\mathrm C^0(U(K;a,b);\mh^2)$.
\end{lem}
\cvd Since $N_n(p)=-\mathrm{grad} T_n(p)$ Lemma~\ref{piatto:Gauss:lem}
implies that the family $\{T_n\}$ is equi-continuous on $U(K;a,b)$. On
the other hand since $||N_n(p)||\geq 1$ we have that $|T_n(p)|<C$ for
every $p\in U(K;a,b)$.  Thus the family $\{T_n\}$ is pre-compact in
$\mathrm C^0(U(K;a,b))$.  On the other hand by 
Lemma~\ref{piatto:conv1:lem},  $T_n\rightarrow
T_\infty$.\\ Finally the same argument as in Proposition 6.5 of
\cite{Bo} shows that $N_n\rightarrow N_\infty$ in $\mathrm
C^0(U(K;a,b);\mh^2)$.  The proof of Proposition \ref{piatto:conv:prop}
easily follows.

\section {From flat regular domains towards measured 
geodesic laminations}\label{REGD_ML}
In this section we will construct the inverse maps of 
$$ \Uu^0: \Mm\Ll \to \Rr/\mr^3, \ \ \ \Uu^0: \Mm\Ll/SO(2,1) \to
\Rr/\ISO^+(\mx_0) $$ eventually proving the classification of
non-degenerate flat regular domains stated in Theorem
\ref{GEN_FLAT_CLASS} of Chapter \ref{INTRO}.

\begin{figure}
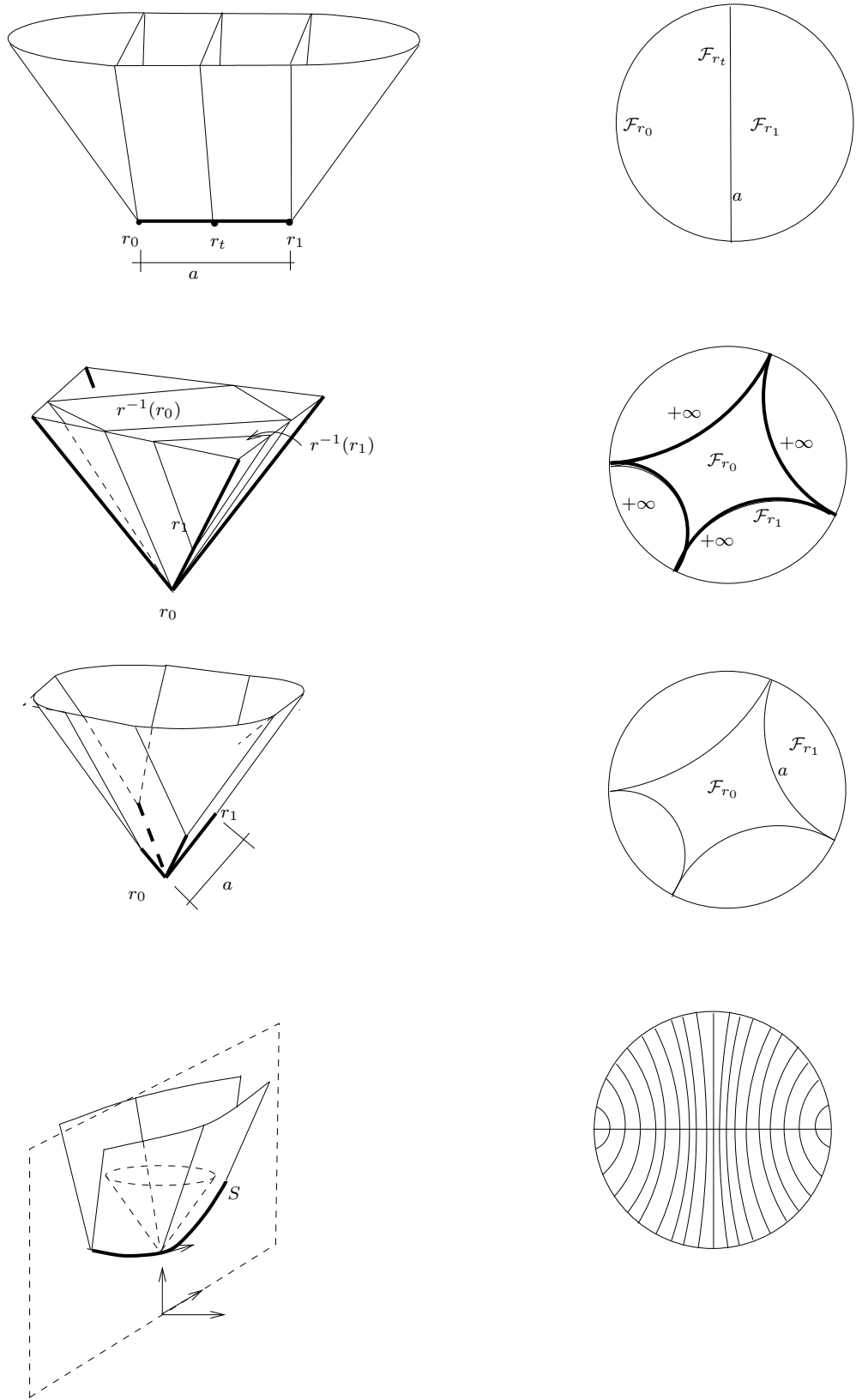

\begin{center}
\input GR2005_fig_class.pstex_t
\caption{{\small The laminations corresponding to some regular
domains.}}\label{class:fig}
\end{center}
\end{figure}

Let us fix a regular domain $\Uu$.  Thus for $r_0\in\Sigma$ let us put
$\Ff(r_0)=N(r^{-1}(r_0))$ (as usual $N$ denotes the Gauss map, $r$ the
retraction on the initial singularity $\Sigma$). We have seen in
Example~\ref{lamrd} that
\[
   \Ll=\bigcup_{\Ff(r):\dim \Ff(r)=1}\Ff(r)\ \cup\
   \bigcup_{\Ff(r):\dim\Ff(r)=2}\partial\Ff(r)\ \cup\ \partial H
\]
is a geodesic lamination on $H=\overline{{\rm Im} N}$.

In the remaining part of this section, we are going to construct a
transverse measure $\mu$ on $\Ll$, such that $\Uu=\Uu^0_\lambda$ where
$\lambda=(\Ll,\mu)$.

In order to explain the idea to construct $\mu$, let us consider the
following situation. Suppose $\Uu=\Uu^0_\lambda$.  Let $c(t)$ be a
geodesic arc, and suppose for simplicity that it does not meets the
weighted part of $\lambda$.  Let $\rho(t)$ be the unique point of
$\partial\Uu$ admitting a support plane orthogonal to $c(t)$.  By
construction we have
\[
\rho(t)-\rho(t_0)=\int_{[c(t_0),c(t)]}v(t)\d\mu_c(t)\,.
\]
Notice that $\rho$ turns to be a rectifiable and its length is exactly
the total mass of $c$. Since the length of $\rho$ is determined by the
geometry of $\Uu$, the total mass of $c$ is determined by $\Uu$.

To make this argument to work in the general case, we need the following
Proposition.

\begin{prop}\label{rect:prop}
Let $c:[0,1]\rightarrow H$ be a  path transverse to $\mathcal
L$. Then $X=N^{-1}(c([0,1])\cap\Uu(1)$ is a locally rectifiable arc.

There exists a parametrization of $X$
\[
\widetilde c:I\rightarrow\widetilde \Uu(1)
\]
(where $I$ is an interval interval)  and a monotone
increasing function $s:I\rightarrow[0,1]$ such that
\[
     c(s(t))=N\circ\widetilde c(t)\,.
\]

The path $\widetilde c$ has finite length if and only if both $c(0)$
and $c(1)$ do not lie on $\partial H$.
\end{prop}
In order to prove this proposition, we need some technical lemmas.

\begin{lem}\label{proprd:lem}
Let $\hat\Uu(1)$ be the set of points $p$ in $\Uu(1)$ such that
$N(p)\in\mathring{H}$.
The map
\[
N:\hat\Uu(1)\rightarrow\mathring{H}
\]
is proper.
\end{lem}
\Dim
Given a divergent sequence $p_n\in\hat\Uu(1)$ we have to show that
$N(p_n)$ does not admit limit in $\mathring{H}$.

If $p_n\rightarrow p_\infty$ with
$p_\infty\in\Uu(1)\setminus\hat\Uu(1)$, then $N(p_n)\rightarrow
N(p_\infty)$, that, by definition of $\hat\Uu(1)$, is on $\partial H$.

Suppose that $p_n$ diverges in $\Uu(1)$. By contradiction, suppose
that there exists $x\in\mathring{H}$ such that $N(p_n)\rightarrow
x$ and let $p\in\Uu(1)$ such that $N(p)=x$. Consider the sequence of
segments $[p,p_n]$. Up to passing to a subsequence, it converges to a
spacelike infinite ray $R=p+\mr_{\geq 0} v$ for some $v\in T_x\mh^2$.
Since $[p,p_n]$ is contained in $\Uu$, so is $R$.

Since $x\in\ort v$, there exists some $y\in\mathring{H}$ such that
$\E{v}{y}>0$. Now, there exists some support plane of $\Uu$ orthogonal
to $y$, \emph{i.e.} there exists $C>0$ such that $\E{y}{q}\leq C$ for
every $q\in\Uu$. Since we have $\E{p+tv}{y}\rightarrow+\infty$ as
$t\rightarrow+\infty$ we get a contradiction.  \cvd

\begin{lem}\label{lemma-tecn-5.1}
Let $c:[0,1]\rightarrow \mathring H$ intersect each leaf 
of $\lambda$ at most in one point. Let $t_0=0\leq
t_1\leq\cdots\leq t_N=1$ be a partition of $[0,1]$ (notice that
$t_i=t_{i+1}$ is allowed) and for all $i\in\{0,\ldots, N\}$ let
$r_i\in\Sigma$ be such that $c(t_i)\in\mathcal F(r_i)$.\par 

Moreover
suppose that if $t_i=t_{i-1}$ then $r_i-r_{i-1}$ is a non-zero vector
pointing towards $c(1)$ (i.e. $\E{\dot c(t_i)}{r_i-r_{i-1}}>0$).\par
There exists a constant $C$ independent of the partition such that
\[
   \sum_{i=1}^N |r_i-r_{i-1}|\leq C\,.
\]
\end{lem}
\Dim
Let us set $x=r_N-r_0$. We have  $x=\sum_{i=1}^N r_i-r_{i-1}$ so
\begin{equation}\label{eq-tecn-5.1}
  \E{x}{x}=\sum_{i,j=1}^N\E{r_i-r_{i-1}}{r_j-r_{j-1}}\,.
\end{equation}
When $r_i-r_{i-1}$ is not zero, the orthogonal geodesic in $\mh^2$
separates $\mathcal F(r_i)$ from $\mathcal F(r_{i-1})$ (this descends
from~(\ref{fund-ineq})).  Since $c(t_i)\in\mathcal F(r_i)$ and $c$ is
geodesic, the geodesics of $\mh^2$ orthogonal to $r_i-r_{i-1}$ and
$r_j-r_{j-1}$ are disjoint or equal (if these vectors are not zero).
  Thus these vectors do not generate a
spacelike plane, that is, the reverse of the Schwartz inequality holds
\begin{equation}\label{eq-tecn-2-5.1}
\E{r_i-r_{i-1}}{r_j-r_{j-1}}^2\geq |r_i-r_{i-1}|^2|r_j-r_{j-1}|^2\,.
\end{equation}
Since $r_{i+1}-r_i$ and $r_{j+1}-r_{j}$ point in the same direction
(here we are using that $c$ intersects each leaf at most once) we
have $\E{r_i-r_{i-1}}{r_j-r_{j-1}}\geq |r_i-r_{i-1}||r_j-r_{j-1}|$.
We obtain
\[
\left(\sum_{i=1}^N \E{r_i-r_{i-1}}{r_i-r_{i-1}}^{1/2}\right)^2\leq\E{x}{x}\,.
\]
By Lemma~\ref{proprd:lem}, $N^{-1}(c([0,1]))$ is compact, so there
exists $C$ such that
\[
  |r-s|\leq C
\]
for all $r,s\in r\big(N^{-1}(c[0,1])\big)$.
This constant verifies the statement of the lemma.
\cvd
\begin{lem}\label{lemma-tecn2-5.1}
Let $c:[0,1]\rightarrow H$ be a  path transverse to the
lamination $\mathcal L$. Then for all $t\in [0,1]$ we have that $N^{-1}(c(t))$
is either a point or a segment in $\mx_0$. Moreover let $\mathcal S$ be the
subset of $[0,1]$ formed by points $t$ such that $N^{-1}(t)$ is not a single
point. Then $\mathcal S$ is numerable and
\[
   \sum_{t\in \mathcal S} l(t)<+\infty
\]
where $l(t)$ is the length of the segment $N^{-1}(t)$.
\end{lem}
\Dim The first statement is obvious. In order to prove the second part
of this lemma it is sufficient to prove that there exists a constant
$C$ such that for every finite subset $\mathcal S'$ of $\mathcal S$ we
have
\[
   \sum_{t\in \mathcal S'} l(t)\leq C\,.
\]
Since every transverse path is a finite composition of paths that
intersect each leaf at most once, we can assume $c$ to satisfy such a
property.

Now let us take $\mathcal S'=\{t_0<\cdots<t_{N}\}$ and consider the
partition of $[0,1]$ given by $t_1\leq t_1\leq t_2\leq
t_2\leq\cdots\leq t_N\leq t_N$ . For $i\leq N$ let $r_i,s_i\in\Sigma$
be the endpoints of $r(N^{-1}(t_i))$.  We have that
$l(t_i)=|r_i-s_i|$. Then by applying Lemma \ref{lemma-tecn-5.1} to the
family $r_1,s_1,r_2,s_2,\ldots,r_N,s_N$ we obtain that $\sum
l(t_i)\leq C$ where $C$ is the constant given by Lemma
\ref{lemma-tecn-5.1}.  \cvd

\emph{ Proof of Proposition~\ref{rect:prop}: } If $c(0)$ lies on
$\partial H$ then either $N^{-1}(c(0))$ is empty or it is a spacelike
geodesic ray. Thus in order to prove the first part of the proposition
we can assume that $c$ is contained in the interior of $H$.  Moreover
since $c$ is a finite composition of arcs that intersects each leaf at
most once, there is no loss of generality if we assume that it
satisfies such a property.

By using the above notation let us set $L=\sum_{t\in \mathcal S}
l(t)$. We want to construct an injective (and surjective) map from
$X=N^{-1}(c([0,1])\cap\widetilde \Uu(1))$ to $[0,L+1]$.

For $p\in X$ let us set $u(p)\in[0,1]$ such that $N(p)=c(u(p))$ and
$\mathcal S_p=\mathcal S\cap[0,u(p))$.\par
Now consider the map $\tau: X\rightarrow [0,L+1]$ so defined
\[
\tau(p)=\left\{\begin{array}{ll}
               u(p)+\sum_{t\in\mathcal S_p} l(t) &\textrm{ if
               }u(p)\notin\mathcal S\\
               u(p)+\sum_{t\in\mathcal S_p} l(t) + |p-q(p)| &
               \textrm{ otherwise}\,.
               \end{array}\right.
\]
where $q(p)$ is the endpoint of $N^{-1}N(p)$ such that $p-q(p)$ points
 towards $c(1)$.  It is easy to check that $\tau$ is a continuous
 map. Moreover it is injective. In fact if $u(p)<u(q)$ then we have
 $\tau(p)<\tau(q)$. On the other hand if $u(p)=u(q)$ then $u(p)\in
 \mathcal S$ and the first and second terms of the sum in the
 definition of $\tau$ are equal whereas the third ones are different.
 Thus $\tau(p)\neq \tau(q)$.

In order to prove that $X$ is a rectifiable arc consider the
parametrization $\widetilde c:[0,L+1]\rightarrow X$ given by
$\tau^{-1}$. Let us set $s(t)=u(\widetilde c(t))$. Notice that $s$ is
monotone increasing function and $N(\widetilde c(t))=c(s(t))$.

Let $t_0<t_1<\ldots<t_N$ be a partition of $[0,L+1]$. We have to show
that there exists a constant $K$ which does not depend on the
partition such that
\[
    \sum_{i=1}^N |\widetilde c(t_i)-\widetilde c(t_{i-1})|\leq K\,.
\]
Now let us set $r(t)=r(\widetilde c(t))$ and $N(t)=N(\widetilde c(t))=c(s(t))$.
We have  $\widetilde c(t)=N(t)+r(t)$, so
\[
  |\widetilde c(t_i)-\widetilde c(t_{i-1})|^2 
  =|N(t_i)-N(t_{i-1})|^2+ |r(t_i)-r(t_{i-1})|^2+
2\E{N(t_i)-N(t_{i-1})}{r(t_i)-r(t_{i-1})}\,.
\]
The last term in that sum is less than
$||N(t_i)-N(t_{i-1})||\,||r(t_i)-r(t_{i-1})||$ (where $||\cdot||$ is the
Euclidean norm). On the other hand since $N(t_i)$'s are contained in a compact
subset of $\mathbb H^2$ there exists a constant $A>1$ such that
\[
   ||N(t_i)-N(t_{i-1})||\leq A |N(t_i)-N(t_{i-1})|\,.
\]
Since the geodesic orthogonal to $r(t_i)-r(t_{i-1})$ intersects $c$
(indeed it separates $N(t_i)$ from $N(t_{i-1})$),
there exists a constant $A'>1$ such that
\[
   ||r(t_i)-r(t_{i-1})||\leq A' |r(t_i)-r(t_{i-1})|\,.
\]
There exists a constant $B$ such that
\[
   |\widetilde c(t_i)-\widetilde c(t_{i-1})|
    \leq B\left(|N(t_i)-N(t_{i-1})|+|r(t_i)-r(t_{i-1})|\right)\,.
\]
 It follows that
\[
    \sum_{i=1}^N |\widetilde c(t_i)-\widetilde c(t_{i-1})|\leq B\left(\ell(c)+
    \sum_{i=1}^N |r(t_i)-r(t_{i-1})|\right)\,.
\]
On the other by Lemma \ref{lemma-tecn-5.1} we have
$\sum_{i=1}^N |r(t_i)-r(t_{i-1})|\leq C$. So
$K=B(C+\ell(c))$ works.

To conclude the proof we have to show that if $c$ is an arc reaching
the boundary then the length of $\widetilde c$ is infinite.  Suppose
that $c(1)\in\partial H$. If $c(1)\in {\rm Im} N$ then $\widetilde c$
contains an infinite spacelike ray, so the claim is obvious.

On the other hand, if $c(1)\notin {\rm Im} N$, then $r(t)=r\circ
\widetilde c(t)$ escapes from compact sets of $\mx_0$ as $t\rightarrow
1$. In fact if $r(t_n)\rightarrow r_\infty$ for some subsequence then
$r_\infty$ is in $\partial\Uu$ and $r_\infty+\ort{c(1)}$ is a
spacelike support plane for $\Uu$, that is $c(1)\in {\rm Im} N$.  It follows
that for every $M$ there exists $\eps>0$ such that
\[
   \sum_{i=1}^{N} ||r(t_{i+1})-r(t_i)||>M
\]
for every subdivision $t_0\leq t_1\leq \ldots t_N=1-\eps$.  Since the
geodesic orthogonal to $r(t_{i+1})-r(t_i)$ converges to the boundary
component of $H$ as $t\rightarrow 1$, there exists
some constant $C$ such that $||r(t_{i+1})-r(t_i)||\leq
C|r(t_{i+1})-r(t_i)|$.  It easily follows that the length of
$c^{-1}(0,1-\eps)$ is bigger than $M$. Thus the length of $\widetilde
c$ is infinite.

\cvd

Now we can prove Theorem ~\ref{GEN_FLAT_CLASS}. 

\Dim
Let $\Ll$ be the geodesic lamination associated to $\Uu$ as in~\ref{lamrd}
We will construct a measure $\mu$ on  $\Ll$ such that $\Uu=\Uu^0_\lambda$.

Let us set $Y=\{x\in\mh^2|\#N^{-1}(x)\cap \Uu(1)>1\}$. Notice that if
$N(p)\in Y$ and $N$ is differentiable at $p$ then
$\mathrm{det}(dN)=0$. Since $N$ is Lipschitz it follows that $Y$ has
null Lebesgue measure (notice that $Y$ is a union of leaves of $\Ll$).

Let $c:[0,1]\rightarrow\mh^2$ be a geodesic segment transverse to
$\mathcal L$ with no end-point in $Y$.  Consider the inverse image
$\widetilde c$ of $c$ on the level surface $\Uu(1)$.  There exists an
arc-length parameterization $t\mapsto \widetilde c(t)$ of $\widetilde
c$. Since the intrinsic metric of $\Uu(1)$ is locally bi-Lipschitz with the
Euclidean one, the path $\widetilde c(t)$ turns to be Lipschitz. Let us
set $r(t)=r(\widetilde c(t))$ and $N(t)=N(\widetilde c(t))$. Since
$\widetilde c(t)=N(t)+r(t)$ and since $N(t)$ is Lipschitz it follows
that $r(t)$ is Lipschitz. So it has derivative almost everywhere and
\[
  r(t)=r(0)+\int_0^t\dot r(s)\mathrm d s\,.
\]
Since $r(t)=\widetilde c(t)- N(\widetilde c(t))$, $\dot r$ is a
spacelike vector almost everywhere and
\begin{equation}\label{trans-meas-eq1}
\dot r(t)\in
T_{\widetilde c(t)}\Uu(1)=T_{N(t)}\mh^2\,.
\end{equation}
On the other hand since
\[
\begin{array}{lr}
   \E{r(t+h)-r(t)}{x}\geq 0 & \E{r(t+h)-r(t)}{y}\leq 0
\end{array}
\]
for every $x\in \Ff(r(t+h))$ and $y\in \Ff(r(t))$
we obtain that $\dot r(t)=0$ except if $N(t)$ lies in some leaf of $\Ll$ and
in that case 
\begin{equation}\label{trans-meas-eq2}
\dot r(t)\in \ort{T_{N(t)}\Ff(r(t)}\,.
\end{equation}

Let us set $\mu_c=N_*(|\dot r|\mathrm dt)$.
Since $c$ is transverse to $\mathcal L$ it is easy to see that this set
is the closure of $\{t\in (0,1)|c(t)\in L\}=\{t\in (0,1)|\dot
c(t)\notin T_{c(t)}C(c(t))$ ($L$ is the support of $\mathcal L$).
By using (\ref{trans-meas-eq1})
and (\ref{trans-meas-eq2}) it is not hard to see that $\mu_c$ is a transverse
measure.

In this way we define  a transverse measure $\mu$ on $\mathcal L$.
Moreover if $c$ is a transverse path then the following identity holds
\[
 r(N^{-1}(c(1)))-r(N^{-1}(c(0))=\int_{\widetilde c}\dot r(t)\mathrm d t
\]
By (\ref{trans-meas-eq2}), $\dot r(t)=|\dot r(t)| v(N(t))$,
where $v(y)$ is the unit vector orthogonal to the leaf through $y$ if $y\in L$
 and is $0$ otherwise. It follows that
\[
r(N^{-1}(c(1)))-r(N^{-1}(c(0))=\int_{c}v(y)\mathrm d \mu_c(t)
\]
and so $\Uu=\Uu^0_{(\mathcal L, \mu)}$.
\cvd

\begin{remark}\label{dist:supliv}
\emph{ Given a geodesic segment, $c$, joining $x,y\in\mathring H$, let
us fix a Lipschitz parameterization $t\mapsto\widetilde c(t)$ of
$\widetilde c=N^{-1}(c)\cap\Uu(1)$.  Let us set $N(t)=N(\widetilde
c(t))$ and $r(t)=r(\widetilde c(t))$. Since $\widetilde
c(t)=N(t)+r(t)$ we have
\[
|\dot{\widetilde c}|\leq|\dot N|+|\dot r|\,.
\]
Thus $\ell(\widetilde c)\leq \ell(c)+\mu(c)$.  In particular given
$p,q\in\Uu(1)$ it follows that the distance between them (with respect
the intrinsic metric of $\Uu(1)$) is bounded by
$d_\mh(N(p)+N(q))+\mu([N(p)+N(q)]$ (where $d_\mh$ is the hyperbolic
distance and $\mu$ is the transverse measure associated to $\Uu$).  }
\end{remark}

\section{Initial singularities and $\R$-trees}\label{more:initial:sing} 
By construction the initial singularity $\Sigma$ of
a flat regular domain $\Uu = \Uu^0_\lambda$ is, in a sense, 
the dual object to $\lambda = (H,\Ll,\mu)$. We may go 
further in studying such a duality.

At the end of Section \ref{RD:sec} we have seen that the initial
singularity has natural length-space structure. We are going to see
that in fact it is a $\R$-tree.

Let us consider the measured geodesic lamination
$\hat\lambda$ on $\Uu(1)$ given by pulling back $\lambda$ (doing
like in Section \ref{ML_REGD}).

If $c$ is a geodesic arc in $H$ we
have seen that $\hat c=N^{-1}(c)\cap\Uu(1)$ is a locally rectifiable arc on
$\Uu(1)$. By (\ref{www}), if $u$ is an arc on $\Uu(1)$
joining two points $p,q\in\hat c$ then
\[
   \hat\mu(u)\geq\hat\mu(\hat c|_p^q)
\]
where $\hat c|_p^q$ is the sub-arc of $\hat c$ joining $p$ to $q$.
This remark is useful to characterize the geodesics of $\Sigma$.
\begin{prop}
If $c$ is a geodesic arc in $H$ joining a point in $\Ff(r_0)$ to a
point in $\Ff(r_1)$ then $k_c=rN^{-1}(c)$ is a geodesic arc of
$\Sigma$ passing through $r_0$ and $r_1$ (whenever $\Ff(r_i)$
are not weighted leaf, $r_0$ and $r_1$ are the endpoints of $k_c$).
\end{prop}
\Dim Denote by $\hat c$ the pre-image of $c$ on the level surface
$\Uu(1)$ and for every $p\in\Uu(1)$ let $r(a,p)$ be the intersection
of the integral line of the gradient of $T$ with the level surface
$\Uu(a)$.  For every $p,q\in\hat c$ let $u_a$ be the geodesic arc in
$\Uu(a)$ joining $r(a,p)$ to $r(a,q)$.  The length of $u_a$ of
$\Uu(a)$ is greater than the length of $r\circ u_a$, that, in turn,
is equal to the total mass of
$r(1,\cdot)\circ u_a$ that is a path joining $p$ to $q$.  Thus we have
\[
    \delta_a(r(a,p), r(a,q))\geq \mu_{\hat c}(\hat c|_p^q))=\ell(k_c|_p^q)
\]
where $k_c|_p^q$ is the sub-arc of $r\circ\hat c$ 
with end-points $r(p)$ and $r(q)$. 
As $\delta_a(r(a,p), r(a,q))$ tends to $\delta(r(p),r(q))$, we deduce
\[
   \delta(r(p), r(q))\geq \ell(k_c|_p^q)\,.
\]
Since $k_c$ is an arc containing $r(p)$ and $r(q)$, it is a
geodesic arc.  
\cvd 
Now we are able to prove that $(\Sigma,\delta)$ is
a $\R$-tree, dual to the lamination $\lambda$.\par Let us recall that
a $\R$-tree is a metric space such that for any pair of points $(x,y)$
there exists a \emph{unique} arc joining them. Moreover we require
that such an arc is isometric to the segment $[0, \delta(x,y)]$.
\begin{prop}\label{gen:tree:prop}
$(\Sigma,\delta)$ is a $\R$-tree.
\end{prop}

\Dim We have to show that given $x,y\in\Sigma$ a unique arc joining
them exists.\par It is not difficult to construct a numerable set of
geodesic arcs $c_n$ in $H_\Uu$ with a common end-point $p_0\notin
L_W$ such that every leaf of $\lambda$ is cut by at least one $c_n$.
It follows that $\Sigma$ is the union of all $k_n:=k_{c_n}$'s.  Now it
is not difficult to see that \emph{$k_i\cap k_j$ is a
sub-arc}. Moreover, the point $r_0=r(N^{-1}(p_0))$ lies on $k_i$ for
every $i$. Thus, by induction on $h$, we see that $k_{i_1}\cap
k_{i_2}\cap\ldots\cap k_{i_h}$ is a sub-arc for any $h$.  Let us set
\[
   \Sigma_i=k_1\cup\ldots\cup k_i\,.
\]
We have that $\Sigma_i$ is a simplicial tree ( non-compact if some
$c_j$ reach the boundary of $H_\Uu$). Moreover the inclusion of
$\Sigma_i$ into $\Sigma_j$ is isometric.

There exists $i$
sufficiently large such that $x,y\in\Sigma_i$. Denotes by $[x,y]$ the
(unique) arc joining them in $\Sigma_i$ and let
$\pi_i:\Sigma_i\rightarrow [x,y]$ the natural projection (that, in
particular, decreases the lengths).  Since the inclusion
$\Sigma_{i}\rightarrow\Sigma_{i+1}$ is isometric we have that
$\pi_{i+1}|_{\Sigma_i}=\pi_i$.  Thus the maps $\pi_i$ glue to a map
\[
\pi:\Sigma\rightarrow [x,y]
\]
that decreases the lengths.\par Let $\Bb$ be the set of points
$r\in\Sigma$ such that $\dim \Ff(r)=2$. We have that $\Bb$ is a
numerable set. Now we claim that every $z\in [x,y]\setminus\Bb$
disconnects $x,y$ in $\Sigma$. From the claim it follows that every
path $c$ joining $x,y$ in $\Sigma$ must contain
$[x,y]\setminus\Bb$. Since $\Bb$ is numerable
$\overline{[x,y]\setminus\Bb}=[x,y]$ and so $c$ contains $[x,y]$.\par
Now, the claim is proved by means of the projection
$\pi:\Sigma\rightarrow [x,y]$ we have defined.  In fact, since
$z\notin\Bb$ it is not difficult to see that $\pi^{-1}(z)=\{z\}$. Thus
the sets
\[
    \begin{array}{l}
    U_x=\pi^{-1}([x,z))\\
    U_y=\pi^{-1}((z,y])
    \end{array}
\]
are two disjoint open sets that cover $\Sigma\setminus\{z\}$.
\cvd

\section {Equivariant constructions}\label{flat:equiv}
Recall from Section ~\ref{FLAT:I} the definition of the set
$$\Rr^\Ee = \{(\Uu,\tilde{\Gamma})\}$$
and that $\ISO^+(\mx_0)$ naturally acts on it.
We are going to show that the map $\Uu^0$ extends to an
equivariant map
$$ \Uu^0: \Mm\Ll^\Ee \to \Rr^\Ee, \ \ \
(\lambda,\Gamma) \to (\Uu^0_\lambda,\Gamma^0_\lambda)$$ where
$\Gamma$ is the linear part of $\Gamma^0_\lambda$, 
that induces a map
$$ \Uu^0: \Mm\Ll^\Ee/SO(2,1) \to  \Rr^\Ee/\ISO^+(\mx_0) \ . $$ 
This eventually leads to the proof of Theorem ~\ref{FULLFLAT}, and Corollaries
 ~\ref{Teich-like} and  ~\ref{non-abel}, stated in Chapter ~\ref{INTRO}. 

\smallskip 

In fact, let $\Gamma\subset SO^+(2,1)$ be a discrete,
torsion free group of isometries of $\mh^2$, such that $(H,\lambda)$
is invariant under $\Gamma$.  We can construct a representation
\[
  h^0_\lambda:\Gamma\rightarrow\ISO^+_0(\mx_0)
\]
such that
\begin{enumerate}
\item
The linear part of $h^0_\lambda(\gamma)$ is $\gamma$.
\item
$\Uu^0_\lambda$ is $h^0_\lambda(\Gamma)$-invariant and the action of
$h^0_\lambda(\Gamma)$ on it is free and properly discontinuous.
\item
The Gauss map $N: \Uu^0_\lambda\rightarrow\mh^2$ 
is $h^0_\lambda$-equivariant:
\[
    N(h^0_\lambda(\gamma)p)=\gamma N(p).
\]
\end{enumerate}
In fact, we simply define
\[
  h^0_\lambda(\gamma)=\gamma\, +\, \rho(\gamma x_0) \ .
\]
Notice that $\tau(\gamma)=\rho(\gamma x_0)$ defines a {\it cocycle}
in $Z^1(\Gamma,\R^3)$; by changing the base point $x_0$,
that cocycle changes by coboundary, so we have a well defined
class in $H^1(\Gamma, \R^3)$ associated to $\lambda$.

Let us consider the hyperbolic surface $F= \mh^2/\Gamma$. $F$ is
homeomorphic to the quotient by $\Gamma$ of the image of the Gauss
map.  Then 
$$Y= Y(\lambda,\Gamma)=
\Uu^0_\lambda/h^0_\lambda(\Gamma)$$ 
is a flat maximal globally
hyperbolic, future complete spacetime homeomorphic to $F\times
\mr$. The natural projection $\Uu^0_\lambda \to Y$ is a locally
isometric universal covering map. The cosmological time $T$ of
$\Uu^0_\lambda$ descends onto the canonical cosmological time of
$Y$. Theorem ~\ref{FULLFLAT} and its corollaries easily follow. 

\paragraph{Cocompact $\Gamma$-invariant case}
Let us recall a few known facts that hold in this case (see \cite {M,
Bo, Bo1}\cite{BG}(1)). Any such a (future complete) spacetime is of
the form
$$ Y(\lambda,\Gamma)=\Uu(\mh^2,\lambda)/\tilde{\Gamma}$$ that is, with
the terminology introduced in Section \ref{ML:I}, it is an instance
of $\Mm\Ll(\mh^2)$-spacetime. $\tilde{\Gamma}$ and its linear part
$\Gamma$ are isomorphic groups, and $\Gamma$ is a Fuchsian group.
$F=\mh^2/\Gamma$ is a compact closed hyperbolic surface (of some genus
$g\geq 2$), homeomorphic to a Cauchy surface $S$ of
$Y(\lambda,\Gamma)$.

Hence, up to isometry homotopic to the identity, such flat spacetimes
are parametrized either by:
\smallskip
 
(a) $T_g\times \Mm\Ll_g$, where $T_g$ denotes the Teichm\"uller space
of hyperbolic structures on $S$, and $\Mm\Ll_g$ has been introduced in
Subsection \ref{Gamma:inv:lam}. 
\smallskip

\noindent or

\smallskip 
(b) the flat Lorentzian holonomy groups $h^0_\lambda(\Gamma)$'s,
up to conjugation by $\ISO_0(\mx_0)$. If we fix $\Gamma$, this
induces an identification between $\Mm\Ll_g$ and the cohomology
group $H^1(\Gamma, \mr^3)$ (where $\mr^3$ is identified with the
group of translations on $\mx_0$). 

\smallskip

\noindent Moreover, $Y(\lambda,\Gamma)$ is determined by the {\it
asymptotic states of its cosmological time}, that have in this case
the following clean description. For every $s>0$, denote by
$s\Uu^0_\lambda(a)$, the surface obtained by rescaling the metric on
the level surface $\Uu^0_\lambda(a) = T_\lambda^{-1}(a)$ by a constant
factor $s^2$. Clearly there is a natural isometric action of $\Gamma
\cong \tilde{\Gamma}$ on each $s\ \Uu^0_\lambda(a)$. Then

\smallskip
(i) when $a\to +\infty$, then the action of $\Gamma$ on
$(1/a)\Uu^0_\lambda(a)$ 
converges (in the sense of Gromov) to the action of $\Gamma$ on
$\mh^2$;
\smallskip

(ii) when $a\to 0$, then action of $\Gamma$ on $\Uu^0_\lambda(a)$
converges to the natural action on the initial singularity of
$\Uu^0_\lambda$. This is an action ``with small stabilizers'' on
such real tree that is dual to $\lambda$. Thanks to Skora's
duality Theorem it follows that these asymptotic states determine
the spacetime (for the notions of equivariant Gromov convergence,
and  Skora duality Theorem see e.g.  \cite{Ot}).


\chapter{Flat Lorentzian vs hyperbolic geometry}\label{HYPE}
Let $\Uu=\Uu^0_\lambda$ be a 
flat regular domain in $\mx_0$, according to Chapter~\ref{FGHST}.  
If $T$ denotes its cosmological time, recall
that $\Uu(a) = T^{-1}(a)$, $\Uu(\geq 1) = T^{-1}([1,+\infty[)$, and
so on.

The main aim of this chapter is to construct a local
$\mathrm C^1$-diffeomorphism
\[
   D=D_\lambda:\Uu(> 1)\rightarrow\mh^3
\]
such that the pull-back of the hyperbolic metric is obtained by a Wick
rotation of the standard flat Lorentzian metric, directed by the
gradient of the cosmological time of $\Uu$ (restricted to $\Uu(>1)$),
and with universal rescaling functions (in the sense of Section
~\ref{CANWR:I}).  Hence $D$ can be considered as a developing map of a
hyperbolic manifold.  By analyzing the asymptotic behaviour of $D$ at
the level surface $\Uu(1)$ as well as when $T\to +\infty$, we
recognize such a hyperbolic manifold as a ``$H$-hull" $M_\lambda$
described in Section ~\ref{ML:I} and in Section ~\ref{SPS}.  Eventually
we get a proof of the statement (1) of Theorem ~\ref{WR:I} of
Chapter ~\ref{INTRO}.

\section{Hyperbolic bending cocycles}\label{hyp:bend:cocy}
We fix once and for all an embedding of $\mh^2$ into $\mh^3$
as a totally geodesic hyperbolic plane.

In order to construct the map $D$ we have to recall the construction
of {\it bending} $\mh^2$ along $\lambda=(\Ll,\mu)$ (here we omit to
write $\mh^2$). This notion was first introduced by Thurston in
\cite{Thu}. We mostly refer to the Epstein-Marden paper \cite{Ep-M}
where bending has been carefully studied. In that paper a {\it
quake-bend} map is more generally associated to every {\it
complex-valued} transverse measure on a lamination $\Ll$. Bending maps
correspond to imaginary valued measures. So, given a measured geodesic
lamination $(\Ll,\mu)$ we will look at the quake-bend map
corresponding to the complex-valued measure $i\mu$.  In what follows
we will describe Epstein-Marden's construction referring to the cited
paper for rigorous proofs.
\smallskip

Given a measured geodesic lamination $\lambda$ on $\mh^2$, first let us
recall what the associated {\it bending cocycle} is. 
This is a map
\[
    B_\lambda:\mh^2\times\mh^2\rightarrow PSL(2,\mc)
\]
which satisfies the following properties:
\begin{enumerate}
\item
$B_\lambda(x,y)\circ B_\lambda(y,z)=B_\lambda(x,z)$ for every $x,y,z\in\mh^2$.
\item
$B_\lambda(x,x)=Id$ for every $x\in\mh^2$.
\item
$B_\lambda$ is constant on the strata of the stratification of $\mh^2$
determined by $\lambda$.
\item
If $\lambda_n\rightarrow\lambda$ on a $\eps$-neighbourhood of the
segment $[x,y]$ and $x,y \notin L_W$, then
$B_{\lambda_n}(x,y)\rightarrow B_{\lambda}(x,y)$ .
\end{enumerate}

\noindent
By definition, a $PSL(2,\mc)$-valued {\it cocycle} on an arbitrary set
$S$ is a map
\[
   b:S\times S\rightarrow PSL(2,\mc)
\]
satisfying the above conditions 1. and 2.
\smallskip

If $\lambda$ coincides with its simplicial part (this notion has been
introduced in Section \ref{laminations}), then there is an easy
description of $B_{\lambda}$.

If $l$ is an oriented geodesic of
$\mh^3$, let $X_l\in\sG\lG(2,\mc)$ denote the infinitesimal generator
of the positive rotation around $l$ such that $\exp(2\pi X_l)=Id$
(since $l$ is oriented the notion of \emph{positive} rotation is well
defined). We call $X_l$ the \emph{standard} generator of rotations around
$l$.

 Now take $x,y\in\mh^2$. If they lie in the same leaf of
$\lambda$ then put $B_\lambda(x,y)=Id$. If both $x$ and $y$ do not lie
on the support of $\lambda$, then let $l_1,\ldots,l_s$ be the
geodesics of $\lambda$ meeting the segment $[x,y]$ and $a_1,\ldots,
a_s$ be the respective weights.  Let us consider the orientation on
$l_i$ induced by the half plane bounded by $l_i$ containing $x$ and
non-containing $y$. Then put
\[
   B_\lambda(x,y)=\exp(a_1 X_1)\circ\exp(a_2 X_2)\circ\cdots\circ\exp(a_s X_s)
\, .
\]
If $x$ lies in $l_1$ use the same construction, but replace $a_1$
by $a_1/2$; if $y$ lies in $l_s$ replace $a_s$ by $a_s/2$.

If $\lambda$ is not simplicial, $B_\lambda(x,y)$ is defined as the limit of
$B_{\lambda_k}(x,y)$ where $\lambda_k$ is a standard approximation of
$\lambda$ in a box $B=[a,b]\times[c,d]$, such that $[x,y]=[a,b]\times\{*\}$. The fact that
$B_{\lambda_k}(x,y)$. converges is proved in~\cite{Ep-M}.

\begin{remark}
\emph{Even if $B_\lambda$ is defined in~\cite{Ep-M} only for measured
geodesic laminations of $\mh^2$, the same argument allows us to define
$B_\lambda$ for any $\lambda = (H,\Ll,\mu)\in \Mm\Ll$. In that case 
$B_\lambda(x,y)$ is defined only for $x,y\in \mathring{H}$.  }
\end{remark}

We will work in the general set-up indicated by the above remark.
The following estimate will play an important r\^ole in our study.  It is
a direct consequence of Lemma 3.4.4 (Bunch of geodesics) of
\cite{Ep-M}. We will use the operator norm on $PSL(2,\mc)$.
\begin{lem}\label{hyperbolic:bend:cont:lem}
For any compact set $K$ of $\mathbb H^2$ there exists a constant $C$
with the following property. For every $\lambda=(H,\Ll,\mu)\in \Mm\Ll$ 
such that $K\subset\mathring H$,  for every $x,y\in K$, and for
every leaf $l$ of $\Ll$ that cuts $[x,y]$,

\[
    || B_\lambda(x,y)- \exp(m X)||\leq C m d_\mh(x,y).
\]

where $X$ is the standard generator of the rotation along $l$, 
and $m$ is the total mass of the segment $[x,y]$.
\end{lem}
\cvd

\smallskip

The bending cocycle is not continuous on the whole definition set.  In
fact by Lemma~\ref{hyperbolic:bend:cont:lem} it is continuous on
$(\mathring{H}\setminus L_W)\times (\mathring{H}\setminus L_W)$ (recall
that $L_W$ is the support of the weighted part of $\lambda$).
Moreover, if we take $x$ on a weighted geodesic $(l,a)$ of $\lambda$
and sequences $x_n$ and $y_n$ converging to $x$ from the opposite
sides of $l$ then we have
\[
B_\lambda(x_n,y_n)\rightarrow \exp(aX)
\]
where $X$ is the generator of rotations around $l$.

Now we want to define a continuous ``pull-back'' of the bending
cocycle on the level surface $\Uu(1)$ of our spacetime.

\begin{prop}\label{lift:cocycle}
A determined construction produces a {\it continuous} cocycle
\[
    \hat B_\lambda:\Uu(1)\times\Uu(1)\rightarrow PSL(2,\mc)
\]  
such that 
\begin{equation}\label{hyperbolic:bend:ext:eq}
   \hat B_\lambda (p,q)=B_\lambda(N(p), N(q))
\end{equation}
for $p,q$ such that $N(p)$ and $N(q)$ do not lie on $L_W$.

The  map $\hat B_\lambda$ is locally Lipschitz.
For every compact subset of $\ \Uu(1)$, the Lipschitz constant
on $K$ depends only on $N(K)$ and on the diameter of $K$.
\end{prop}    
\Dim Clearly formula~(\ref{hyperbolic:bend:ext:eq}) 
defines $\hat B_\lambda$ on 
$\Uu \setminus N^{-1}(L_W)$. We claim
that this map is locally Lipschitz.
This follows from the following general lemma.

\begin{lem}\label{hyperbolic:bend:lip:lem}
Let $(E,d)$ be a bounded metric space, $b: E\times E \to PSL(2,\C)$
be a cocycle on $E$. Suppose there exists $C>0$ such that
\[
   ||b(x,y)-1||< C d(x,y).
\]
Then there exists a constant $H$, depending only on $C$ and on the
diameter $D$ of $E$, such that $b$ is $H$-Lipschitz.
\end{lem}

{\it Proof of Lemma \ref{hyperbolic:bend:lip:lem}:}  
For $x,x',y,y'\in E$ we have
\[
  || b(x,y)- b(x',y')||= ||b(x,y)-b(x',x)b(x,y)b(y,y')||. 
\]
It is not hard to show that, given three elements
$\alpha,\beta,\gamma\in PSL(2,C)$ such that $||\beta-1||<\eps$ and
$||\gamma - 1||<\eps$, there exists a constant $L_\eps>0$ such that
\[
   ||\alpha-\beta\alpha\gamma||<L_\eps ||\alpha||(||\beta-1||+||\gamma-1||).
\]
Since we have that $||b(x,y)-1||< CD$, we get
\[
  || b(x,y)-b(x',y')||\leq L_{CD} (D+1)C(d(x,x')+d(y,y')) \ .
\]
Thus $H=L_{CD} C(D+1)$ works.
\cvd
\smallskip

\noindent 
Fix a compact subset $K$
of $\Uu(1)$ and let $C$ be the constant given by
Lemma~\ref{hyperbolic:bend:cont:lem}. Then for 
$x,x'\in K'= K\setminus N^{-1}(L_W)$ we have
\[
   ||\hat B_\lambda (x,x')-1||\leq ||\exp\mu(c) X-1|| + 
C\mu(c)d_\mh(N(x),N(x'))
\]
where $X$ is the generator of the  infinitesimal rotation around a
geodesic of $\Ll$ cutting the segment $c=[N(x), N(x')]$ and
$\mu(c)$ is its total mass.

Recall that a measured geodesic lamination $(\hat\Ll,\hat\mu)$ has
 been defined in Section~\ref{lam:supliv} as the pull-back of
 $(\Ll,\mu)$.  By Lemma~\ref{geo:est:lem} $\mu(c)=\hat\mu(\hat c)$,
 where $\hat c$ is the geodesic path on $\Uu(1)$ joining $x$ to $x'$.

Thus, if $A$ is the maximum of the norm of
generators of rotations around geodesics cutting $K$ and $M$ is the
diameter of $N(K)$ we get
\begin{equation}\label{hyperbolic:bend:lip3:eq}
   ||\hat B_\lambda(x,x')-1||\leq (A+CM)\hat\mu(\hat c)\leq (A+CM) d(x,x') 
\end{equation}
where the last inequality is a consequence of Lemma~\ref{geo:est:lem}.

By Lemma~\ref{hyperbolic:bend:lip:lem} we have that $\hat B_\lambda$
is Lipschitz on $K'\times K'$. Moreover, since $A, C,
M$ depend only on $N(K)$, the Lipschitz constant depends
only on $N(K)$ and the diameter of $K$.

In particular $\hat
B_\lambda$ extends to a locally Lipschitz cocycle on the closure of
$\Uu(1)\setminus N^{-1}(L_W)$ in $\ \Uu(1)$.  Notice that this
closure is obtained  by removing from $\Uu(1)$ the interior part of the
bands corresponding to leaves of $\Ll_W$.  

Fix a band
$\Aa\subset\Uu(1)$ corresponding to a weighted leaf $l$.   For
$p,q\in\Aa$, let us set $u=r(p)$ and $v=r(q)$.  If $u=v$ then let us
put $\hat B_\Aa(p,q)=1$. Otherwise  $v-u$ is a spacelike vector
orthogonal to $l$.  Consider the
orientation on $l$ given by $v-u$. Let $X\in\sG\lG(2,\mc)$ be the standard
generator of positive rotation around $l$ . Then for $p,q\in\Aa$ let
us put
\[
   \hat B_{\Aa}(p,q)=\exp(|v-u|X)
\]
where $|v-u|=\E{v-u}{v-u}^{1/2}$.  Notice that $\hat B_{\Aa}$ is a
cocycle. 
Moreover, if $p,q\in\partial\Aa$, then
Lemma~\ref{hyperbolic:bend:cont:lem} implies that
\[
  \hat B_{\Aa}(p,q)= \hat B_\lambda (p,q) \ .
\]
Let us fix $p,q\in\Uu(1)$. If $p$ (resp. $q$) lies in a band $\Aa$
(resp. $\Aa'$) let us take a
point $p'\in\partial\Aa$  (resp. $q'\in\partial\Aa'$) 
otherwise put $p'=p$ ($q=q'$). Then let us define
\[
   \hat B_\lambda(p,q)=\hat B_{\Aa}(p,p')\hat B_\lambda(p',q')\hat
   B_{\Aa'}(q',q).
\]
By the above remarks it is easy to see that $\hat B(p, q)$ is
well-defined, that is it does not depend on the choice of $p'$ and $q'$.
Moreover it is continuous.  Now we can prove that there exists a
constant $C$ depending only on $N(K)$ and on the diameter of $K$ such that
\begin{equation}\label{hyperbolic:bend:lip:eq}
   ||\hat B_\lambda (p,q)-1||\leq C d(p,q).
\end{equation}
In fact we have found a constant $C'$ that works for
$p,q\in\Uu(1)\setminus N^{-1}(L_W)$.  On the other hand we have that
if $p,q$ lie in the same band $\Aa$ corresponding to a geodesic $l\in
\Ll_W$, then the maximum $A$ of norms of standard generators of rotations
around geodesics that cuts $K$ works.  If $p$ lies in $\Uu(1)\setminus
N^{-1}(L_W)$ and $q$ lies in a band $\Aa$, then consider the geodesic
arc $c$ between $p$ and $q$ and let $q'$ lie in the intersection of $c$
with the boundary of $\Aa$. Then we have
\begin{eqnarray*}
  ||\hat B_\lambda(p,q)-1||=
||\hat B_\lambda(p,q')\hat B_\lambda(q',q)-1|| \leq\\
\leq 
||\hat
    B_\lambda(p,q')-1||\,||\hat B_\lambda(q',q)||+||B_\lambda(q',q)-1||<
\\
< (C'Ad(q,q')+A) d(p,q).
\end{eqnarray*}

\noindent
Thus, if $D$ is the diameter of $K$, then the constant $C''=A(C'D+1)$
works.  In the same way we can find a constant $C'''$ working for
$p,q$ that lie in different bands. Thus the maximum $C$ between
$C',C'',C'''$ works. Notice that $C$ depends only on $N(K)$ and on the
diameter of $K$. Proposition \ref{lift:cocycle} is now
proved.
\cvd

\begin{remark}\label{pippo} \emph{
Lemma~\ref{hyperbolic:bend:cont:lem} implies that for a
fixed compact set $K$ in $\Uu(1)$ there exists a constant $C$ depending
only on $N(K)$ and on the diameter of $K$ such that for every
transverse arc $c$ contained in $K$ with end-points $p,q$ we have
\[
       ||\hat B_\lambda (p,q)-\exp(\hat\mu(c)X)||\leq C \hat\mu(c)
         d_\mh(N(p),N(q))
\]
where $X$ is the standard generator of a rotation around a geodesic of
$\lambda$ cutting the segment $[N(p),N(q)]$.
}\end{remark}

Let us extend now $\hat B_\lambda$ on the whole $\Uu\times\Uu$. If
$p\in \Uu$ we know that $p=r(p)+T(p)N(p)$. Let us denote by
$r(1,p)=r(p)+N(p)$ and put
\[
    \hat B_\lambda (p,q)=\hat B_\lambda(r(1,p),r(1,q)) \ .
\] 
Proposition~\ref{lift:cocycle} immediately extends to
the whole of $\Uu$.
\begin{cor}\label{hyperbolic:bend:lip:cor}
The map 
\[
  \hat B_\lambda:\Uu\times\Uu\rightarrow PSL(2,\mc)
\]
is locally Lipschitz (with respect to the Euclidean distance
on $\Uu$).  Moreover the Lipschitz constant on $K\times K$ depends
only on $N(K)$, on the diameter of $r(1,\cdot)(K)$ and on the maximum
$M$and minimum $m$ of $T$ on $K$.
\end{cor}
\Dim
It is sufficient to show that the map $p\mapsto r(1,p)$ is 
Lipschitz on $K$ for some constant depending only on $N(K)$, $m$, $M$. 

Take $p,q\in K$.
We have that $p=r(1,p)+(T(p)-1)N(p)$ and $q= r(1,q)+(T(q)-1)N(q)$.
Thus we have
\[
  r(1,p)-r(1,q)=p-q + (N(p)-N(q)) + T(q)N(q)-T(p)N(p) \ .
\]
Since $N(K)$ is compact there exists $C$ such that $||N(p)||<C$  and 
$||N(p)-N(q)||< C |N(p)-N(q)|$ for $p,q\in K$.
Now if we set $a=T(p)$ and $b=T(q)$ we have that
$|N(p)-N(q)|<1/b|p_b-q|$ where $p_b=r(p)+bN(p)$.
It follows that 
\[
   |N(p)-N(q)|<1/m(||p-q||+||p-p_b||)=1/m(||p-q||+|T(p)-T(q)|) \ .
\]
Hence we obtain the following inequality
\[
  || r(1,p)-r(1,q)||\leq ||p-q|| + C'||p-q|| + C''|T(q)-T(p)| \ .
\]
Since $N$ is the Lorentzian gradient of $T$ we have  that
\[
  |T(p)-T(q)|\leq C||p-q||
\]
and so the Lipschitz constant of $r(1,\cdot)$ is less than $1+C'+CC''$.
\cvd

In the last part of this subsection we will show that if
$\lambda_n\rightarrow \lambda$ on a $\eps$-neighbourhood $K_\eps$ of
a compact \emph{convex} set $K$, 
then $\hat B_{\lambda_n}$ tends to $\hat B_\lambda$ on
$N^{-1}(K)$.

More precisely, for $a<b$ let $U(K;a,b)$ denote the subset of
$\Uu^0_\lambda$ of the points in $N^{-1}(K)$ with cosmological time
greater than $a$ and less than $b$. By Prop.~\ref{piatto:conv:prop} we
know that $U(K;a,b)\subset\Uu_{\lambda_n}$, for $n$ sufficiently
large.  Then we can consider the map $\hat B_n$ given by the
restriction of $\hat B_{\lambda_n}$ on $U(K;a,b)$.
\begin{prop}\label{hyperbolic:bend:conv:prop}
The sequence $\{\hat B_n\}$ converges to the map $\hat B=\hat
B_{\lambda}$ uniformly on $U(K;a,b)$.
\end{prop}
\Dim For $n$ sufficiently large we have $N_n(p)\in K_\eps$ for $p\in
U(K;a,b)$.  Now let $C_n$ be the supremum of the total masses with
respect to $\lambda_n$ of geodesic arcs contained in $K$ (that is
compact).  By Remark~\ref{dist:supliv}, the diameter of
$N_n^{-1}(K_\eps)\cap\Uu_n(1)$ is bounded by ${\rm diam}(K)+C_n$.  On
the other hand, thanks to the compactness of $K$, $C_n$ is attained by
a geodesic arc in $K$ and converges to the supremum of the total
masses of geodesic arcs in $K$ with respect to $\lambda$.
   
Thus  there exists a constant $C$ such that every $\hat B_n$ is
$C$-Lipschitz on $U(K;a,b)$ for $n$ sufficiently large. It follows
that the family $\{\hat B_n\}$ is pre-compact in $\mathrm
C^0(U(K;a,b)\times U(K;a,b); PSL(2,\C))$.

So it is sufficient to
prove that if $\hat B_n$ converges to $\hat B_\infty$ then $\hat
B_\infty=\hat B$. Clearly  $\hat B_\infty$ is a cocycle
and it is sufficient to show that $\hat B(p_0,q)= \hat
B_\infty(p_0,q)$ for some $p_0\in K$. 
First suppose that $N(q)\notin L_W$. We can take
$q_n\in U(K;a,b)$ such that $q_n\rightarrow q$ and $N_n(q_n)$ are not
in $(L_W)_n$. Thus we have
\[
   \hat B_n(p_0, q_n)= B_n(N_n(p_0), N_n(q_n)) \ .
\]

By  Proposition~3.11.5 of \cite{Ep-M}, $B_n(N_n(p_0),
N_n(q_n))$ converges to $B(N(p_0), N(q))=\hat B(p_0,q)$.  Thus we have
that $\hat B(p,q)=\hat B_\infty(p,q)$ for $p,q$ lying in the closure
of $N^{-1}(\mh^2\setminus L_W)$.  Now take a point $q$ in 
$\Aa=N^{-1}(l)$ for some weighted leaf $l$. 
In order to conclude it is sufficient to show that $\hat
B_\infty(p,q)=\hat B(p,q)$ for $p\in\partial \Aa$ such that
$N(p)=N(q)$.  Notice that $r_n(p)$ is different from $r_n(q)$ for
$n$ sufficiently large so $[N_n(p), N_n(q)]$ intersects the lamination
$\lambda_n$. Choose for every $n$ a leaf $l_n$ intersecting $[N_n(p),
N_n(q)]$ and let $X_n$ be the standard generator of the rotation
around $l_n$.

Now consider the path $c_n(t)=r_n(1, tp+(1-t)q)$. It
is not hard to see that $c_n$ is a transverse arc in $\Uu(1)$ so that
a measure $\hat\mu_n$ is defined on it. Moreover, its total mass, $m_n$,
with respect to $\hat\mu_n$ is
\[
     m_n=\int_0^1|\dot r_n(t)|\d t \ .
\]

By Remark~\ref{pippo} there exists a constant $C$ such that
\[
   |\hat B_n(p,q)-\exp(m_nX_n)|<Cd_\mh(N_n(p), N_n(q))\ .
\]  
On the other hand since $N_n(p)$ and $N_n(q)$ converge to $N(p)=N(q)$,
$X_n$ tends to the generator of the
rotation around the leaf through $N(q)$.  In order to conclude it is
sufficient to show that $m_n$ converges to $|r(1,p)-r(1,q)|=|p-q|$.
We know that
\[
tp+(1-t)q=r_n(t)+T_n(t)N_n(t)\,,
\]
so deriving in $t$ we get
\begin{equation}\label{hyperbolic:bend:conv:eq}
p-q=\dot r_n(t)-\E{N_n(t)}{p-q}N_n(t)+T_n(t)\dot N_n(t) \ .
\end{equation}
Since $N_n(t)\rightarrow N(p)$, $\E{N_n(t)}{p-q}$ tends
to $0$.  We will prove that $\dot N_n(t)$ tends to $0$ in
$L^2([0,1];\mr^3)$ so $\dot r_n(t)$ tends to $p-q$ in
$L^2([0,1];\mr^3)$. From this fact it is easy to see that
$m_n\rightarrow |p-q|$.  

Since the images of $N_n$ are all contained in a compact set
$\overline K_\eps\subset\mh^2$, there exists $C$ such that
\[
    \int_0^1||\dot N_n(t)||^2\d t\leq C \int_0^1 |\dot N_n(t)|^2\d t \ .
\]
On the other hand by taking the scalar product of both sides of
equation~(\ref{hyperbolic:bend:conv:eq}) with $\dot N(t)$ we obtain
\[
 \E{p-q}{\dot N_n(t)}=\E{\dot N_n(t)}{\dot r_n(t)}+ T_n(t)|\dot N_n(t)|^2 \ .
\]
By inequality~(\ref{fund-ineq}) we can deduce $\E{\dot N_n(t)}{\dot
  r_n(t)}\geq 0$  so
\[
  \E{p-q}{\dot N_n(t)}\geq a|\dot N_n(t)|^2 \ .
\]
By integrating on $[0,1]$ we obtain that $\dot N_n$ tends to $0$ in
$L^2([0,1];\mr^3)\ $.  \cvd

\section{The Wick rotation}\label{theWR}
We are going to construct the local $\mathrm C^1$-diffeomorphism
\[
   D=D_\lambda:\Uu(>1)\rightarrow\mh^3
\]
with the properties outlined at the beginning of this Section.

Let $B=B_\lambda$ be the bending cocycle, and $\hat B = \hat
B_\lambda$ be the continuous cocycle defined on the whole $\Uu \times
\Uu$, as we have done above.

\noindent 
Fix a base point $x_0$ of $\mathring{H}$ ($x_0$ is supposed not to be
in $L_W$). The {\it bending} of $\mathring{H}$ along $\lambda$ is the
map
\[
    F=F_\lambda:\mathring{H}\ni x\mapsto B(x_0,x)x\in\mh^3\,.
\] 
It is not hard to see that $F$ satisfies the following properties:
\begin{enumerate}
\item
It does not depend on $x_0$ up to post-composition of elements of
$PSL(2,\mc)$.
\item
It is a $1$-Lipschitz map.
\item
If $\lambda_n\rightarrow\lambda$ then $F_{\lambda_n}\rightarrow F_\lambda$
with respect to the compact open topology.
\end{enumerate}
\begin{remark}\emph{
Roughly speaking, if we bend $\mathring{H}$ taking $x$ fixed then $B(x,y)$ is
the isometry of $\mh^3$ that takes $y$ to the corresponding point of
the bent surface.  }\end{remark}

\noindent
Since both $\mh^3$ and $\mh^2 \subset \mh^3$ 
are oriented, the normal bundle is oriented too. Let
$v$ denote the normal vector field on $\mh^2$ that is positive
oriented with respect to the orientation of the normal bundle.
Let us take $p_0 \in N^{-1}(x_0)$ and for
$p\in\Uu(>1)$ consider the geodesic ray $c_p$ of $\mh^3$
starting from $F(N(p))$ with speed vector equal to $w(p)=\hat
B(p_0,p)_*(v(N(p)))$.  Thus $D$ is defined in the
following way:
\begin{equation}\label{hyp:dev:map}
    D(p)=c_p(\arctgh(1/T(p)))=
\exp_{F(N(p))}\left(\arctgh\left(\frac{1}{T(p)}\right)w(p)\right) \ .
\end{equation}
\begin{teo}\label{hyperbolic:WR:teo}
The map $D$ is a local $\mathrm C^1$-diffeomorphism such that the
pull-back of the hyperbolic metric is equal to the Wick Rotation of
the flat Lorentz metric, directed by the gradient $X$ of the
cosmological time with universal rescaling functions:
\begin{equation}\label{wr:eq}
       \alpha =  \frac{1}{T^2-1} \ , \qquad\qquad \beta=\frac{1}{(T^2-1)^2} \ .
\end{equation}
\end{teo}

\begin{remark}\emph{
Before proving the theorem we want to give some heuristic motivations
for the rescaling functions we have found. Suppose $\lambda$ to be a
finite lamination on $\mh^2$.  If the weights of $\lambda$ are
sufficiently small then $F_\lambda$ is an embedding onto a bent
surface of $\mh^3$. In that case the map $D$ is a homeomorphism onto
the non-convex component, say $\Ee$, of $\mh^3\setminus
F_\lambda(\mh^2)$.  The distance $\delta$ from the boundary is a
$\mathrm C^1$-submersion. Thus the level surfaces $\Ee(a)$ give rise
to a foliation of $\Ee$.  The map $D$ satisfies the following
requirement:}
\begin{enumerate}
\item
\emph{
The foliation of $\Uu$ by $T$-level surfaces is sent to that foliation of
$\Ee$.
}
\item
\emph{
The restriction of $D$ on a level surface $\Uu(a)$ is a dilation by a factor
depending only on $a$.
}
\end{enumerate}
\emph{ The first requirement implies that $\delta(D(x))$ depends only
on $T(x)$ that means that there exists a function
$f:\mr\rightarrow\mr$ such that $\delta(D(x))=f(T(x))$.}\par

\noindent
\emph{ Denote by $\mh(\lambda)$ the surface obtained by replacing
every geodesic of $\lambda$ by an Euclidean band of width equal to the
weight of that geodesic.  We have that $\Uu(T)$ is isometric to the
surface $T\,\mh(\lambda/T)$, with the usual notation.
On the other hand the surface $\Ee(\delta)$ is
isometric to $\ch\delta\,\mh(\lambda\tgh \delta)$.  Now 
 $\mh(a\lambda)$ and $\mh(b\lambda)$ are related by a
dilation if and only if $a=b$.  By comparing $\Ee(\delta(T))$
with $\Uu(T)$ we can deduce that }
\[
    T=1/\tgh\delta(T)
\]
\emph{
so $t>1$. Moreover the dilation factor is}
\[
    (\alpha(t))^{1/2}=\frac{\ch\delta(T)}{T}=\frac{1}{(T^2-1)^{1/2}}.
\]
\emph{
In order to compute the vertical rescaling factor notice that the hyperbolic
gradient of $\delta$ is unitary. By requiring that $D$ induces a
Wick rotation directed by the gradient $X$  of $T$, we obtain that
$D_*X=f \mathrm{grad}\,\delta$ for some function $f$.
Thus we have}
\[
   \begin{array}{l}
   f=g(\mathrm{grad}\,\delta, D_*X)=X(D^*(\delta))=\mathrm
   d(\arctgh(1/T))[X]=\\
   =-\frac{1}{T^2-1}\d T(X)=\frac{1}{T^2-1}.
   \end{array}
\]
\end{remark}

\noindent
We will prove the theorem by analyzing progressively more complicated
cases.  First we will prove it in a very special case when $\Uu$ is
the future of a spacelike segment. Then, we will deduce the theorem
under the assumption that the lamination $\lambda$ consists of a finite
number of weighted geodesic lines. Finally, by using  standard
approximations (see Section \ref{laminations}), we will obtain the full
statement.

\paragraph{Wick rotation on  $\Uu_0$}
Let $\Uu_0$ be the future of the segment $I=[0,\alpha_0v_0]$, where
$v_0$ is a unitary spacelike vector and $0<\alpha_0<\pi$.  If $l_0$
denotes the geodesic of $\mh^2$ orthogonal to $v_0$, the
measured geodesic lamination corresponding to $\Uu_0$ is simply
$\lambda_0=(l_0,\alpha_0)$.   It is easy to see that in this case the map
$D_0:\Uu_0\rightarrow\mh^3$ is injective. We are going to point
out suitable $\mathrm C^{1,1}$-coordinates on $\Uu_0$ and on the image of
$D_0$ respectively, such that $D_0$ can be easily recognized with
respect to these coordinates.

\begin{figure}[h!]
\begin{center}
\input{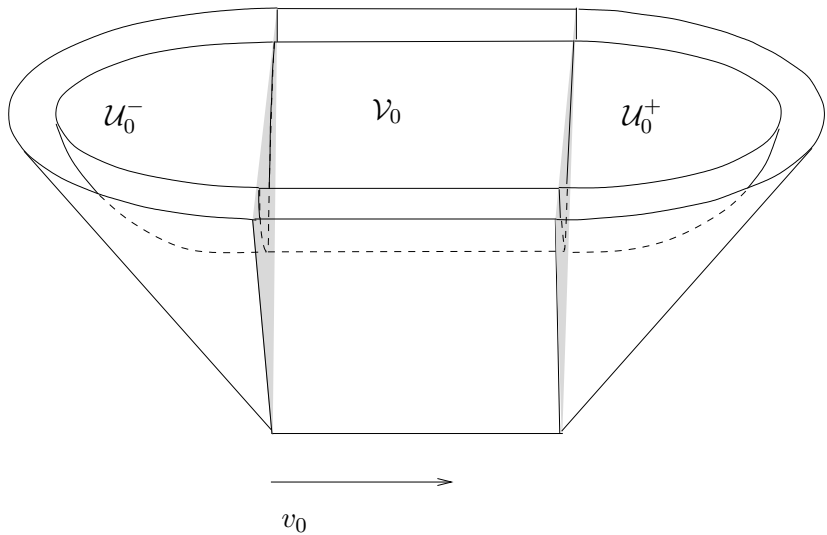}
\caption{{\small The domain $\Uu_0$ and its decomposition. Also a level
    surface $\Uu_0(a)$ is shown.}}
\end{center}
\end{figure}
We denote by $P_\pm$ the components
of $\mh^2 \setminus l_0$ in such a way that $v_0$ is outgoing from
$P_-$.

As usual, let $T$ be the cosmological time, $N$ denote the Gauss map
of the level surfaces of $T$ and $r$ denote the retraction on the
singularity $I$.  We have a decomposition of $\Uu_0$ in three pieces
$\Uu_0^-,\Uu_0^+,\Vv$ defined in the following way:
\[
\begin{array}{l}
\Uu_0^-=r^{-1}(0)=N^{-1}(P_-)\,;\\
\Vv=r^{-1}(0,\alpha_0 v_0)=N^{-1}(l_0)\,;\\
\Uu_0^+=r^{-1}(\alpha_0 v_0)=N^{-1}(P_+) \, .
\end{array}
\]

We denote by $\Uu_0^+(a),\Uu_0^-(a),\Vv(a)$ the intersections of
corresponding domains with the surface $\Uu_0(a)$. The Gauss map on
$\Uu_0^+(a)$ (resp. $\Uu_0^-(a)$) is a diffeomorphism onto $P^+$
(resp. $P^-$) that realizes a rescaling of the metric by a constant 
factor $a^2$.  On the other hand, the
parametrization of $\Vv$ given by
\[
  (0,\alpha_0)\times l_0\ni(t, y)\mapsto a y+t v_0\in \Vv(a)
\]
produces two orthogonal geodesic foliations on $\Vv$. The parametrization
restricted to horizontal leaves is an isometry, whereas on the
on vertical leaves it acts as a rescaling of factor $a$. Thus $\Vv(a)$ is a
Euclidean band of width $\alpha_0$.\\ 

Now we introduce $\mathrm C^{1,1}$
coordinates on $\Uu_0$.  Denote by $l_a$ the boundary of
$\Uu_0^{-}(a)$ and by $d_a$ the intrinsic distance of $\Uu_0(a)$. 
Fix a point $z_0$ on $l_0$ and denote by $\hat z_a\in l_a$
the point such that $N(\hat z_a)=z_0$.

For every $x\in\Uu_0(a)$ there is a unique point $\pi(x)\in l_a$ such
that $d_a(x,l_a)=d_a(x,\pi(x))$.  Then we consider coordinates $T,\zeta,
u$, where $T$ is again the cosmological time, and $\zeta, u$ are
defined in the following way
\begin{equation}\label{flatcoord}
\begin{array}{l}
\zeta(x)=\eps(x) d_{T(x)}(x,l_{T(x)})/T(x)\\
u(x)=\eps'(x)d_{T(x)}(\pi(x), \hat z_{T(x)})/T(x) 
\end{array}
\end{equation}
where $\eps(x)$ (resp. $\eps'(x)$ ) is $-1$ if $x\in\Uu_0^{-}$
(resp. $\pi(x)$ is on the left of $\hat z_{T(x)}$) and is $1$ otherwise.

Choose affine coordinates of Minkowski space $(y_0,y_1,y_2)$ such that
$v_0=(0,0,1)$ and $z_0=(1,0,0)$. Thus the parametrization induced by
coordinates $T,\zeta,u$ is given by
\[
 (T,u,\zeta)\mapsto\left\{\begin{array}{ll} T(\ch u\ch \zeta,\ \sh
                     u\ch\zeta,\ \sh\zeta) & \textrm{ if }\zeta<0\\
                     T(\ch u,\ \sh u,\ \zeta) & 
                     \textrm{ if }\zeta\in[0,\alpha_0/T]\\ 
                     T(\ch u\ch\zeta',\ \sh
                     u\ch\zeta',\ \sh\zeta'+\alpha_0/T) &\textrm{otherwise}
                     \end{array}\right.
\]
where we have put $\zeta'=\zeta-\alpha_0/T$.\par
With respect to these coordinates the metric take the following form: 
\begin{equation}\label{flatmetric}
 h_0(T,\zeta, u)=\left\{\begin{array}{ll}
                           -\d T^2+T^2(\d \zeta^2+\ch^2\zeta \d u^2) &
                           \textrm{ if }\zeta<0\,,\\
                           -\d T^2 + T^2(\d \zeta^2+\d u^2) & 
                           \textrm{ if }\zeta\in[0,\alpha_0/T]\,,\\
                           -d T^2+ T^2(\d \zeta^2+\ch^2(\zeta')\d u^2) &
                           \textrm{ otherwise.}
                            \end{array} \right.
\end{equation}
Notice that the gradient of $T$ is just the coordinate field
$\frac{\partial\,}{\partial T}$.

The Gauss map takes the following form
\begin{equation}\label{flatGM}
   N(T,\zeta,u)=\left\{\begin{array}{ll} (\ch u\ch \zeta,\ \sh
                     u\ch\zeta,\ \sh\zeta) & \textrm{ if }\zeta<0\\
       (\ch u,\ \sh u,\ 0) &\textrm{if }\zeta\in[0,\alpha_0/T]\\ 
       (\ch u\ch\zeta',\ \sh u\ch\zeta',\ \sh\zeta') & \textrm{otherwise}
                     \end{array}\right.
\end{equation}
and the bending cocycle $\hat B_0(p_0, (T,\zeta,u))$ is the rotation
around $l_0$ of angle equal to $0$ if $\zeta<0$, $\zeta$ if
$\zeta\in[0,\alpha_0/T]$, $\alpha_0/T$ otherwise.

Let $\mh^3$ be identified with the set of timelike unit vectors in the
$3+1$-Minkowski space $\mm^4$. We can choose affine coordinates on
$\mm^4$ in such a way the inclusion $\mh^3\subset\mh^4$ is induced by
the inclusion $\mx_0\rightarrow\mm^4$ given by
$(x_0,x_1,x_2)\mapsto(x_0,x_1,x_2,0)$.  Thus the general rotation
around $l_0$ of angle $\alpha$ is represented by the linear
transformation $T_\alpha$, such that
\begin{equation}\label{flatbend}
T_\alpha(e_0)=e_0,\ T_\alpha(e_1)= e_1, \ T(e_2)=\cos\alpha\
e_2+\sin\alpha\ e_3,\ T_\alpha(e_3)=-\sin\alpha\ e_2+\cos\alpha\
e_3\,.
\end{equation}
where $(e_0, e_1, e_2, e_3)$ is the canonical basis of $\mr^{4}$.

Thus by means of (\ref{flatGM}) and (\ref{flatbend}) we can write
(\ref{hyp:dev:map}) in local coordinates

\begin{figure}
\begin{center}
\input{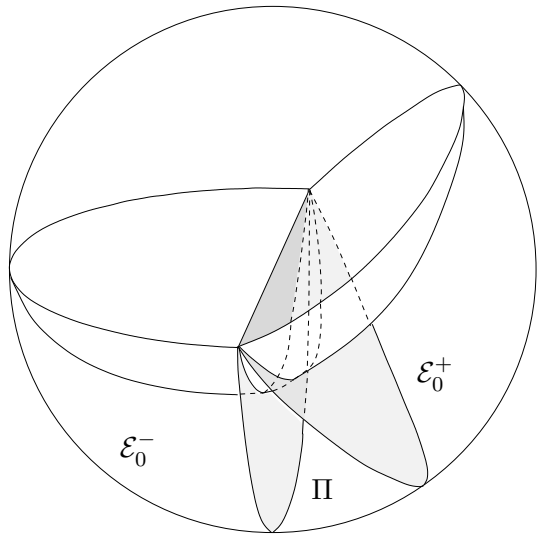}
\caption{{\small The image $\Ee_0$ of $D_0$ and its decomposition.}}
\end{center}
\end{figure}
\[
D_0(T,u, \zeta)\mapsto\left\{\begin{array}{ll}
                          \ch\delta\left(\ch\zeta\ch u,\ \ch\zeta\sh u,\
                          \sh\zeta,\ 0\right)\ +\ \sh\delta(0,0,0,1) \\
                          \textrm{if }\zeta\leq 0\ ;\\ \ch\delta\left(\ch
                          u,\ \sh u,\ 0,0\right)\ + \
                          \sh\delta\left(0,\ 0,\
                          -\sin\frac{\zeta}{\tgh\delta},\,
                          \cos\frac{\zeta}{\tgh\delta}\right) \\
                          \textrm{if }\zeta\in[0,\alpha_0/T]\ ;\\

                          \ch\delta\left(\ch\zeta'\ch u,\ \ch\zeta'\sh
                          u,\ \sh\zeta'\cos\alpha_0,\
                          \sh\zeta'\sin\alpha_0\right)\ + & \\
                          \sh\delta(0,0,-\sin\alpha_0,\
                          \cos\alpha_0)\\ \textrm{otherwise}
                           \end{array}\right.
\] 
where $\delta=\arctgh(1/T)$ and  $\zeta'=\eta-\alpha_0/T$.\par

This map is clearly smooth for $\zeta\neq 0,\alpha_0/T$.
Since the derivatives of $D_0$ with respect the coordinates fields glue 
along $\zeta=0$ and $\zeta=\alpha_0 T$ the map $D$ is $\mathrm C^1$.

By computing the pull-back of the hyperbolic metric we have 
\begin{equation}\label{hypmetric}
D_0^*(g)(T,\zeta, u)=\left\{\begin{array}{ll}
                           \d\delta^2+
                    \ch^2\delta(\d\zeta^2+\ch^2\zeta\d u^2) &
                           \textrm{ if }\zeta<0\\
                           \d\delta^2 + \ch^2\delta(\d\zeta^2+\d u^2) & 
                           \textrm{ if }\zeta\in[0,\alpha_0/T]\\
                   \d\delta^2+\ch^2\delta(\d\zeta^2+\ch^2(\zeta')\d u^2) &
                           \textrm{ otherwise}.
                            \end{array}\right.
\end{equation}
Since $\d\delta=\frac{1}{T^2-1}\d T$ and
$\ch^2\delta=\frac{T^2}{T^2-1}$, comparing (\ref{hypmetric}) and
(\ref{flatmetric}) shows that $D_0^*(g)$ is obtained by the Wick
rotation along the gradient of $T$ with rescaling functions given
in~(\ref{wr:eq}).

\begin{remark}
\emph{ The map $D_0$ is not $\mathrm C^2$ along $\zeta=0$ and
$\zeta=\alpha_0/T$.  On the other hand it is not hard to see that the
derivatives are locally Lipschitz.  }\end{remark}

\paragraph{Finite laminations} 
Suppose that $\lambda$ is a finite lamination on $\mh^2$.  We want to
reduce this case to the previous one. In fact, we will show that for
any $p\in\Uu^0_\lambda$ there exist a small neighbourhood $U$ in
$\Uu^0_\lambda$, an isometry $\gamma$ of $\mx_0$ and an isometry
$\sigma$ of $\mh^3$ such that
\begin{enumerate}
\item
$\gamma(U)\subset\Uu_0$.
\item
   $\gamma$ preserves the cosmological time, that is
  \[
     T(\gamma(p)) = T_\lambda(p)
  \] 
  for every $p\in U$.
\item 
 We have
\begin{equation}\label{1}
    \sigma\circ D_\lambda(p)=D_0\circ \gamma(p)
\end{equation}
for every $p\in U$.
\end{enumerate}

\noindent
First suppose that $p$ does not lie in any Euclidean band.  Fix 
$\eps>0$ so that the disk $B_\epsilon(x)$ in $\mh^2$, with center at
$x=N(p)$ and ray equal to $\eps$, does not intersect any geodesic of
$\lambda$.
Thus, we can choose an isometry $\gamma_0$ of $\mh^2$ such that the
distance between $z=\gamma_0(x)$ and $z_0$ is less than $2\eps$ ($z_0$ being
the base point fixed in the previous Subsection).  Now
let us set $U=N^{-1}(B_\eps(x))$, $\gamma =\gamma_0 - r(p)$
(where $r$ denotes here the retraction of $\Uu^0_\lambda$) and 
 $\sigma=\hat B_\lambda(p_0,p)$.  In fact we have
\[
  D_\lambda(\xi)= \hat B_\lambda(p_0, p)\circ D_0\circ \gamma(\xi)
\]
for $\xi\in U$.\\

\noindent
If $p$ lies in the interior of a band $\Aa=r^{-1}[r_-,r_+]$
corresponding to a weighted leaf $l$, let $I$ be an interval $[s_-,
s_+]$ contained in $[r_-,r_+]$ centered in $r(p)$ of length less than
$\alpha_0$ and set $U=r^{-1}(I)$.  Let $q$ be a point in $U$ such
that $r(q)=s_-$ and $N(q)=N(p)$.  Let $\gamma_0$ be an isometry of
$\mh^2$ which sends $N(q)$ onto $z_0$ and $l$ onto $l_0$; set
$\gamma =\gamma_0-r(q)$.  Now for $\xi\in U$ we have
\[
  D_\lambda(\xi)= \hat B_\lambda(p_0,q)\circ D_0\circ\gamma (\xi).
\]

\noindent
Finally suppose that $p$ lies on the boundary of a band
$\Aa=r^{-1}[r_-,r_+]$ corresponding to the weighted leaf $l$.  Without
loss of generality we can suppose that $r(p)=r_-$.  Now let us fix a
neighbourhood $U$ of $p$ that does not intersect any other Euclidean
band and such that $r(U)\cap [r_-, r_+]$ is a proper interval of
length less than $\alpha_0$.  Let $\gamma_0\in PSL(2,\mr)$ send $N(p)$
onto $z_0$ and $l$ onto $l_0$; set  $\gamma = \gamma_0-r(p)$.  Also in this
case we have
\[
  D_\lambda(\xi)= \hat B_\lambda(p_0,p)\circ D_0\circ\gamma (\xi).
\]
\paragraph{General case}
Before proving the theorem in the general case we need some remarks
about the regularity of $D_\lambda$, when $\lambda$ is finite.  We use
the notations of the proof of Proposition
\ref{hyperbolic:bend:conv:prop}.
\begin{lem}\label{hyperbolic:WR:lem}
Fix a bounded domain $K$ of $\mh^2$, a constant $A>0$,
and $1<a<b$. For every finite lamination $\lambda$ 
let us set $U_\lambda=U_{\lambda}(K;a,b)$.
Then there exists a constant $C$ depending only on $K$, $A$ and
$a,b$ such that for every finite lamination $\lambda$ 
such that the sum of the weights is less than $A$, 
the first and the second derivatives of $D_\lambda$ on
$U_\lambda$ are bounded by $C$.
\end{lem}
\Dim For every point $p$ of $U_\lambda$, the above construction
gives  a neighbourhood $W$, an isometry $\gamma_W$ of $\mx_0$,
and an isometry $\sigma_W$ of $\mh^3$ such that
\[
   D_\lambda=\sigma_W\circ D_0\circ \gamma_W \ .
\]
Moreover, we can choose $W$ so small in such a way that $\gamma_W$ is
contained in $U_0(B_{2\alpha_0}(z_0);a,b)$. Let us fix a constant $C'$
such that first and second derivatives of $D_0$ are bounded by $C'$ on
this set.  By construction the linear part of $\gamma_W$ is an
isometry of $\mh^2$ sending a point in $K$ close to $z_0$. So the
linear parts of $\gamma_W$ form a bounded family in $SO(2,1)$. On the
other hand the Euclidean norm of the translation part of $\gamma_W$ is
bounded by some constant depending only on $K$ and $A$.  Finally
$\sigma_W=B_\lambda(p_0, p)$ for some $p\in K$. By
Lemma~\ref{hyperbolic:bend:cont:lem}, its norm is bounded by some
constant depending only on $K$ and $A$.  Eventually the family
$\{\gamma_W\}$ (resp. $\{\sigma_W\}$) is contained in some compact
subsets of $\ISO(\mx_0)$ (resp. $PSL(2,\mc)$) depending only on $K$
and $A$.

Hence there exists a constant
$C''$ such that first and second derivatives of both $\gamma_W$ and
$\sigma_W$ are bounded by $C''$.  Thus first and second
derivatives of $D_\lambda$ are bounded by $C=27(C'')^2C'$.  
\cvd

\noindent We can finally prove Theorem \ref{hyperbolic:WR:teo} in the
general case.  Take a point $p\in\Uu^0_\lambda$ and consider a
sequence of standard approximations $\lambda_n$ of $\lambda$ on a
neighbourhood $K$ of the segment $[N(p_0), N(p)]$.  By
Propositions~\ref{piatto:conv:prop},
and~\ref{hyperbolic:bend:conv:prop} $D_{\lambda_n}$ converges to
$D_\lambda$ on $U(K;a,b)$.  On the other hand by
Lemma~\ref{hyperbolic:WR:lem} we have that $D_{\lambda_n}$ is a
pre-compact family in $\mathrm C^1(U(K;a,b); \mh^3)$.  Thus it follows
that the limit of $D_{\lambda_n}$ is a $\mathrm C^1$-function.
Finally, as $D_{\lambda_n}$ ${\rm C}^1$-converges to $D_\lambda$, the
cosmological time of $\Uu^0_{\lambda_n}$ ${\rm C}^1$-converges on
$U(K;a,b)$ to the one of $\Uu^0_\lambda$ (see
Proposition~\ref{piatto:conv:prop}), and the pull-back on
$\Uu^0_{\lambda_n}$ of the hyperbolic metric is obtained via the
determined Wick rotation, the same conclusion holds on $\Uu$.  \cvd

\section{On the geometry of $M_\lambda$}\label{hyp:geom:sec}
The hyperbolic $3$-manifold $M_\lambda$ arising by performing
the Wick rotation described in Theorem~\ref{hyperbolic:WR:teo}
consists of the domain $\Uu^0_\lambda(>1)$ endowed
with a determined hyperbolic metric, say $g_\lambda$. 

We are going to study some geometric properties of $M_\lambda$. As
usual, $T$ denotes the cosmological time of the spacetime
$\Uu^0_\lambda$, and $N$ its Gauss map.\\

\subsection{Completion of $M_\lambda$}
Let $\delta$ denote the length-space-distance on $M_\lambda$
associated to $g_\lambda$. 
The following theorem summarizes the main features
of the geometry of  the completion, $\overline M_\lambda$, of $M_\lambda$. 
The remaining part of the Section is devoted to prove it. 

\begin{teo}\label{hyperbolic:compl:teo}
(1) The completion of $M_\lambda$ is 
$\overline M_\lambda=M_\lambda\cup H$ ($H$ being the straight convex set 
on which $\lambda$ is defined)
 endowed with the distance $\overline\delta$ 
\[
\begin{array}{ll}
   \overline\delta (p,q)=\delta(p,q) & \textrm{ if }p,q\in M_\lambda\,,\\
   \overline\delta (p,q)=d_\mh(p,q)  & \textrm{ if }p,q\in H\,,\\
   \overline\delta (p,q)=\lim_{n\rightarrow+\infty}\delta(p,q_n) 
& \textrm{ if } p\in M_\lambda\textrm{and }q\in H
\end{array}
\]
where $(q_n)$ is any sequence in $\Uu^0_\lambda$ such that $T(q_n)=n$ and
$N(q_n)=q$. The copy of $H$ embedded into $\overline M_\lambda$
is called the {\rm hyperbolic boundary}
$\partial_h M_\lambda$ of $M_\lambda$. 
\smallskip

(2) The developing map $D_\lambda$ continuously extends to a map
defined on $M_\lambda \cup \mathring H$. Moreover, the restriction of
$D_\lambda$ to the hyperbolic boundary $\partial_h 
M_\lambda$ coincides with the bending map $F_\lambda$.
\smallskip

(3) Each level surface of the cosmological time $T$ restricted to
$\Uu(>1)$ is also a level surface in $\overline M_\lambda$
of the distance function $\Delta$ from its hyperbolic boundary 
$\partial_h M_\lambda$. Hence the inverse Wick rotation is directed
by the gradient of $\Delta$.
\end{teo}

For simplicity, in what follows we denote by $\delta$ both the
distance on $M_\lambda$ and the distance on $M_\lambda\cup H$.

We are going to establish some auxiliary results.
\begin{lem}\label{hyperbolic:compl:lem}
The map $N:M_\lambda\rightarrow\mh^2$  is $1$-Lipschitz.\\
\end{lem}
\Dim 
Let $p(t)$ be a $\mathrm C^1$-path in $M_\lambda$. We have to
show that the length of $N(t)=N(p(t))$ is less than the length of
$p(t)$ with respect to $g_\lambda$.  (Since $N$ is locally
Lipschitz with respect to the Euclidean topology $N(t)$ is a
Lipschitz path in $\mh^2$.)

By deriving the identity
\[
   p(t)= r(t)+ T(t)N(t)
\]
we get
\[
  \dot p(t) =\dot r(t) + \dot T(t) N(t) + T(t)\dot N(t)
\]
As $\dot r$ and $\dot N$ are orthogonal to $N$ (that up to the sign is the
gradient of $T$) we have
\begin{equation}\label{hyperbolic:compl:eq}
  g_\lambda(\dot p(t), \dot p(t))= 
\frac{\dot T(t)^2}{(T(t)^2-1)^2}+\frac{1}{T(t)^2-1}\E{\dot r(t) + T(t)\dot
  N(t)}{\dot r(t) + T(t)\dot N(t)}.
\end{equation}
By inequality ~(\ref{fund-ineq}), we have
$\E{\dot r(t)}{\dot N(t)}\geq 0$, so
\[
   g_\lambda(\dot p(t), \dot p(t))\geq 
\frac{T(t)^2}{T(t)^2-1}\E{\dot N(t)}{\dot
   N(t)}\geq\E{\dot N(t)}{\dot N(t)} \ .
\]
\cvd

\begin{lem}\label{hyperbolic:compl2:lem}
The map $\arctgh(1/T)$ is $1$-Lipschitz on $M_\lambda$.
Moreover, the following inequality holds
\begin{equation}\label{hyp:ineq}
   \delta(p,q)\leq \arctgh(1/T(p))+\arctgh(1/T(q))+ d_\mh(N(p),N(q))\,.
\end{equation}
\end{lem}
\Dim
By using equation~(\ref{hyperbolic:compl:eq}) we can easily see that
$\arctgh(1/T)$ is $1$-Lipschitz function.

Let us take $p,q\in M_\lambda$ and for $a>\max(T(p), T(q))$ let us set
$p_a=r(p)+a N(a)$ and $q_a=r(q)+aN(q)$. Finally let $c_a$  be the geodesic on
$\Uu^0_\lambda(a)$ joining $p_a$ to $q_a$. Clearly the distance between $p$ and
$q$ is less than the length of the path $[p,p_a]*c_a*[q_a,q]$ (with respect to
the hyperbolic metric $g_\lambda$). By an explicit computation we get
\begin{equation}\label{hyp:ineq2}
   \begin{array}{l}
   \delta(p,q)\leq
   \arctgh(1/T(p))-\arctgh(1/T(p_a))+\frac{d_a(p_a,q_a)}{\sqrt{a^2-1}}+\\
   +\arctgh(1/T(q))-\arctgh(1/T(q_a))
   \end{array}
\end{equation}
where $d_a$ is the distance of $\Uu^0_\lambda(a)$ as slice of the Lorentzian
manifold $\Uu^0_\lambda$. 
In~\cite{Bo} it has been shown that
\[
 \frac{1}{a} d_a(p_a,q_a)\rightarrow d_\mh(N(p),N(q))
\]
as $a\rightarrow+\infty$. So, by letting $a$ go to $+\infty$
in~(\ref{hyp:ineq2}) we get
\[ 
   \delta(p,q)\leq\arctgh(1/T(p))+\arctgh(1/T(q))+ d_\mh(N(p),N(q))\,.
\]
\cvd

\emph{Proof of statements (1) and (2) and of
Theorem~\ref{hyperbolic:compl:teo}:} By
Lemmas~\ref{hyperbolic:compl:lem} and~\ref{hyperbolic:compl2:lem} both
$N$ and $\arctgh(1/T)$ extend to continuous functions of $\overline
M_\lambda$ and if $(q_n)$ is a Cauchy sequence in $M_\lambda$ then
$N(q_n)$ and $\arctgh(1/T(q_n))$ are Cauchy sequences.  In particular
either $T(q_n)$ converges to $a>1$ or to $+\infty$. In the former case
the sequence $r(1, q_n)=r(q_n)+N(q_n)$ is a Cauchy sequence of
$\Uu_\lambda^0(1)$: in fact by~(\ref{hyp:ineq}) the map
$r(1,\cdot):M_\lambda\rightarrow \Uu^0_\lambda(1)$ is Lipschitz on
$\delta^{-1}(\alpha,\beta)$ for $1<\alpha<\beta<+\infty$. Thus $q_n$
converges in $\Uu^0_\lambda$.

Now suppose $(q_n)$ is a sequence such that $N(q_n)\rightarrow
x_\infty$ and $T(q_n)\rightarrow +\infty$. Inequality~(\ref{hyp:ineq})
shows that $(q_n)$ is a Cauchy sequence.  Thus the map
\[
  N:\overline M_\lambda\rightarrow\mh^2
\]
is injective on $\partial M_\lambda=\overline M_\lambda\setminus M_\lambda$
and $N(\partial M_\lambda)=H$.

Finally we have to prove that $N:\partial M_\lambda\rightarrow\mh^2$
is an isometry. Since $N$ is $1$-Lipschitz it is sufficient to show
that $N$ does not decrease the distance on $\partial M_\lambda$. By the
previous description of non convergent Cauchy sequences of $M_\lambda$
we see that $\arctgh(1/T(p))=0$ for every $p\in\partial M_\lambda$. So,
inequality~(\ref{hyp:ineq}) gives the estimate we need.  
\cvd

We are going to prove statement (3) of  Theorem~\ref{hyperbolic:compl:teo}.
\begin{cor}
The function $\Delta$ is $\mathrm C^1$. Moreover the following formula
holds
\[
   \Delta(p)=\arctgh(1/T(p)).
\]
For every point $p\in M_\lambda$ the unique point realizing $\Delta$
on the boundary is $N(p)$ and the geodesic joining $p$ to $N(p)$ is
parametrized by the path
\[
    c:[T(p),+\infty)\ni t\mapsto r(p)+tN(p)\in M_\lambda.
\]
\end{cor}
\Dim
If $p(t)$ is a $\mathrm C^1$-path, by (\ref{hyperbolic:compl:eq}) we have
\[
   g_\lambda(\dot p(t),\dot p(t))\geq (\dot T(t))^2/(T^2-1)^2
\]
and the equality holds if and only if $\dot r(t)=0$ and $\dot N(t)=0$.
Thus we obtain $\Delta(p)\geq\arctgh(1/T(p))$. The 
hyperbolic length of $c$ is equal to $\arctgh(1/T(p))$ so $\Delta(p)=
\arctgh(1/T(p))$.  Moreover, if $p(t)$ is a geodesic realizing the
distance $\Delta$ we have that $\dot r=0$ and $\dot N=0$ so $p$ is a
parametrization of $c$.  \cvd
\smallskip
When $H=\mh^2$ the topology of the completion is described in the
following proposition. Later we will get information in the general
case, together with the study of the AdS rescaling.
\begin{prop}
Suppose $\lambda$ to be a measured geodesic lamination of the whole
 $\mh^2$. Then $\overline M_\lambda$ is a topological manifold with
 boundary, homeomorphic to $\mr^2\times[0,+\infty)$. Moreover,
 $\overline M_\lambda(\Delta\leq\eps)$ is a collar of
 $\mh^2=\partial_h M_\lambda$.
\end{prop}
  
\Dim
It is sufficient to show that the for every
$\eps>0$ the set $\overline M_\lambda(\Delta\leq\eps)$ is homeomorphic
to $\mh^2\times [0,\eps]$.  Unfortunately the map
\[
  \overline M_\lambda(\Delta\leq\eps)\ni x\mapsto (N(x),
  \Delta(x))\in\mh^2\times[0,\eps]
\]
works only if $L_W$ is empty. Otherwise it is not injective.
Now the idea to avoid this problem is the following.
Take a point $z_0\in M_\lambda$ and consider the surface
\[ 
\mh(z_0)=\{x\in\fut(r(z_0))|\E{x-r(z_0)}{x-r(z_0)}= -T(z_0)^2\}.
\]
It is a spacelike surface of $\Uu^0_\lambda(>1)$ (in fact $\mh(z_0)$ is
contained in $\Uu^0_\lambda(>a)$ for every $a<T(z_0)$).  Denote by $v$
the Gauss map of the surface $\mh(z_0)$.  It sends the metric of
$\mh(z_0)$ to the hyperbolic metric multiplied by a factor $1/T(z_0)$.
Now we have an embedding
\[
  \varphi:\mh(z_0)\times [0,+\infty)\ni (p,t)\mapsto
  p+tv(p)\in\Uu^0_\lambda(>1)
\]
that parameterizes the future of $\mh(z_0)$.  Clearly if we cut the
future of $\mh(z_0)$ from $M_\lambda$ we obtain a manifold
homeomorphic to $\mr^2\times (0,+1]$.  Thus in order to prove that
$\overline M_\lambda$ is homeomorphic to $\mr^2\times[0,1]$ is
sufficient to prove the following claim.

The map $\varphi$ extends to a map
\[
   \varphi:\mh(z_0)\times [0,+\infty]\mapsto\overline M_\lambda
\]
that is an embedding onto a neighbourhood of $\partial M_\lambda=\mh^2$ in
$\overline M_\lambda$ such
that
\[
   \varphi(p,+\infty)=v(p) \ .
\]

The remaining part of the proof will be devoted to prove the claim.

A fundamental family of neighbourhoods
of a point $v_0\in\mh^2=\partial M_\lambda$ in $\overline M_\lambda$
is given by
\[
    V(v_0;\eps,a)=\{x\in\Uu^0_\lambda| d_\mh(N(x),v_0)\leq\eps\textrm{ and
    }T(x)\geq a\}\cup \{v\in\mh^2|d_\mh(v,v_0)\leq\eps\} \ .
\]
To prove the claim it is sufficient to see that for any compact sets $H\subset\Uu^0_\lambda$,
$K\subset\mh^2$,  and $\eps,a>0$ there exists $M>0$ such that
\[
     p_0 + t v_0\in V(v_0;\eps,a)
\]
for every $p_0\in H$, $v_0\in K$ and $t\geq M$.

Let us set $p(t)=p_0+t v_0$ and denote by
$r(t)$, $N(t)$, $T(t)$ the retraction, the Gauss map and the
cosmological time computed at $p(t)$.  Notice that $T(t)>t+T(0)$ and
since $p_0$ runs in a compact set there exists $m$ that does not
depend on $p_0$ and $v_0$ such that $ T(t)>t+m$.

On the other hand
by deriving the identity
\[
   p(t)=r(t)+T(t)N(t)
\]
we obtain 
\begin{equation}\label{hyperbolic:compl:top:eq1}
   v_0=\dot p(t)=\dot r(t)+T(t)\dot N(t)+\dot T(t)N(t) \ .
\end{equation}
By taking the scalar product with $\dot N$ we obtain
\[
  \E{v_0}{\dot N}=\E{\dot r}{\dot N}+T\E{\dot N}{\dot N}>0\ .
\]
Since  $\ch(d_\mh(v_0, N(t))=-\E{v_0}{N(t)}$,
the function 
\[
  t\mapsto d_\mh(v_0, N(t))
\]
is decreasing.  Thus there exists a compact set $L\subset\mh^2$ such that
$N(t)\in L$ for every $p_0\in H$ , $v_0\in K$, $t>0$.  By Lemma~\ref{proprd:lem}
there exists a compact set $S$ in $\mx_0$ such that $r(t)\in S$ for every
$t>0$, $p_0\in H$ and $v_0\in K$.  We can choose a point $q\in\mx_0$
such that $S\subset\fut(q)$.  Notice that
\[
   T(t)=|p(t)-r(t)|\leq\sqrt{-\E{p(t)-q}{p(t)-q}} \ .
\]
By using this inequality it is easy to find a constant $M$ (that depends
only on $H$ and $K$) such that
\[
   T(t)\leq t + M.
\]
This inequality can be written in the following way:
\[
   \int_0^{t} (\dot T(s)-1)\mathrm ds \leq M \ .
\]
On the other hand, by ~(\ref{hyperbolic:compl:top:eq1}) we have
\[
  \ch\d_\mh(v_0, N(t))=-\E{v_0}{N(t)}=\dot T(t)
\]
so $\dot T>1$. It follows that the measure of the set
\[
    I_\eps=\{s| \dot T(s)-1>\eps\}
\]
is less than $M/\eps$. Since $T$ is concave, $I_\eps$ is an
interval (if non-empty) of $[0,+\infty)$ with an endpoint at $0$.
Thus $I_\eps$ is contained in $[0, M/\eps]$.

Eventually we have proved
that for $t>\max(M/\eps, a)$ we have $T(t)>a$ and
\[
   \ch\d_\mh(v_0, N(t))=-\E{v_0}{N(t)}=\dot T(t)\leq 1+\eps \ .
\]
Thus the claim is proved.\\

In order to conclude the proof we have to show that $\varphi$ is
proper, and the image is a neighbourhood of $\mh^2=\partial M_\lambda$
in $\overline M_\lambda$.  For the last statement we will show that
for every $v_0$ and $\eps>0$ there exists $a>0$ such that
$V(v_0;\eps,a)\subset\fut(\mh(z_0))$.  In fact, since $N$ is
proper on level surfaces, there exists a compact set $S$ such
that $r(V(v_0;\eps,a))$ is contained in $S$ for every $a>0$. In
particular it is easy to see that there exist constants $c,d$ (depending on
$v_0$ and $\eps$) such
that if we take $p\in V(v_0;\eps,a)$ we have
\[
  |p-r(z_0)|\sim -T(p)^2+cT(p)+d \ .
\]
We can choose $a_0$ sufficiently large such that if $T(p)>a_0$ then
$\E{p-r(z_0)}{p-r(z_0)}<-T(z_0)$. So $V(v_0;\eps,a)\subset\fut(\mh(z_0))$
for $a>a_0$.\\

Finally we have to prove that if $q_n=\varphi(p_n, t_n)$ converges to
a point then $p_n$ is bounded in $\mh(z_0)$.\\ Let us set
$p_n(t)=p_n+tv(p_n)=r_n(t)+T_n(t)N_n(t)$.  Since $N_n(t_n)$ is compact
there exists a compact $S$ such that $r_n(t_n)\in S$.  Thus by arguing
as above we  can find $M>0$ such that
\[
    T(q_n)\leq t_n+M \, .
\]
On the other hand we have
\[
    T(q_n)-t_n>\int_0^{t_n}(-\E{N_n(t)}{v(p_n)}-1)\mathrm ds \ .
\]
Since $-\E{N_n(t)}{v(p_n)}-1=\dot T_n-1$ is a decreasing positive
function  for every $\eps>0$
\[
   0<-\E{N_n(t)}{v(p_n)}<1+\eps
\]
for $t_n\geq t\geq M/\eps$.
In particular, for $n$ sufficiently large
we have that $\E{N_n(t_n)}{v(p_n)}<2$. Thus $v(p_n)$ runs is a compact
set. Since $p_n=v(p_n)+r_0$, the  conclusion follows.
\cvd
\smallskip

\subsection{Projective boundary of $M_\lambda$}\label{proj-bound}
Let us define $\hat M_\lambda=\overline
M_\lambda\cup\Uu^0_\lambda(1)$. In this section we will
prove that the map $D_\lambda:\overline M_\lambda\rightarrow\mh^3$ can
be extended to a map
\[
D_\lambda:\hat M_\lambda\rightarrow\overline\mh^3
\]
in such a way that the restriction of $D_\lambda$ on $\Uu^0_\lambda(1)$
takes value on $S^2_\infty=\partial\mh^3$ and is a $\mathrm C^1$-developing
map for a projective structure on $\Uu^0_\lambda(1)$.

In fact, for a point $p\in\Uu^0_\lambda(1)$ we know that the image via
$D_\lambda$ of the integral line $(p_t)_{t>1}$ of the gradient of
$T$ is a geodesic ray in $\mh^3$ of infinite length starting
at $F_\lambda(N(p))$.  Thus, we define $D_\lambda(p)$ to be the
end-point in $S^2_\infty$ of such a ray.  
In the following statement we use the lamination $\hat \lambda$ on the
level surface $\Uu^0_\lambda(1)$ defined in Section
\ref{ML_REGD}. Moreover we will widely refer to Section \ref{SPS}.
\begin{teo}\label{hyperbolic:proj:teo}
The map $D_\lambda:\Uu^0_\lambda(1)\rightarrow S^2_\infty$ is a local
$\mathrm C^1$-conformal map. In particular it is a developing map for
a projective structure on $\Uu^0_\lambda(1)$.\\ The {\rm canonical
stratification} associated to this projective structure coincides with
the stratification induced by the lamination $\hat\lambda$ and its
{\rm Thurston metric} coincides with the intrinsic spacelike surface
metric $k_\lambda$ on $\Uu^0_\lambda(1)$ .
\end{teo}
\Dim 
The first part of this theorem is proved just as
Theorem~\ref{hyperbolic:WR:teo}.  In fact an explicit computation
shows that $D_\lambda:\Uu^0_\lambda(1)\rightarrow S^2_\infty$ is a $\mathrm
C^1$-conformal map if $\lambda$ is a weighted geodesic. Thus it
follows that $D_\lambda$ is a $\mathrm C^1$-conformal map if $\lambda$
is a simplicial lamination.  Then by using standard approximations we
can prove that $D_\lambda$ is a $\mathrm C^1$-conformal map for any
$\lambda$. Indeed $\Uu^0_\lambda(1)$ can be regarded as the graph of a
$\mathrm C^1$-function $\varphi_\lambda$ defined on the horizontal
plane $P=\{x_0=0\}$.  Moreover if $\lambda_n\rightarrow\lambda$ on a
compact set $K$, then $\varphi_{\lambda_n}$ converges to
$\varphi_\lambda$ on $P(K)=\{x|N(\varphi(x),x)\in K\}$ in $\mathrm
C^1$-topology.  Thus, by using parameterizations of
$\Uu^0_{\lambda_n}(1)$ given by
\[
\sigma_{\lambda_n}(x)=(\varphi_{\lambda_n}(x),x)\,,
\]
we obtain maps
\[
   d_{\lambda_n}:P(K) \rightarrow S^2_\infty \ .
\]
The same argument used in Theorem~\ref{hyperbolic:WR:teo} shows that
$d_{\lambda_n}$ converges to $d_\lambda$ on $P(K)$ in $\mathrm
C^1$-topology.  Finally if $k_n$ is the pull-back of $k_{\lambda_n}$
on $P(K)$, we have that $k_n$ converges on $P(K)$ to the pull-back of
$k_\lambda$.  Since $d_{\lambda_n}:(P,g_n)\rightarrow S^2_\infty$ is a
conformal map by taking the limit we obtain that $d_\lambda$ is 
conformal on $P(K)$.\\

The proof of the second part of the statement is more difficult.
Consider the round disk $\md_0$ in $S^2_\infty$ such that
$\partial\md_0$ is the infinite boundary of the
right half-space bounded by $\mh^2$ in $\mh^3$.  Notice that the
retraction $\md_0\rightarrow\mh^2$ is a conformal map (an isometry if
we endow $\md_0$ with its hyperbolic metric). We denote by
$\sigma:\mh^2\rightarrow\md_0$ the inverse map.  With this
notation the map $D_\lambda:\Uu_\lambda^0(1)\rightarrow S^2_\infty$ can 
be expressed in the following way:
\[
    D_\lambda(p)=\hat B(p_0,p)\sigma(N(p)) \ .
\]
Now for every point $p\in\Uu^0_\lambda(1)$ let us consider the round circle
$\md_p=\hat B(p_0,p)(\md_0)$ and define $\Delta_p$ to be the connected
component of $D_\lambda^{-1}(\md_p)$ containing $p$.

To conclude the proof we will use the following estimate whose proof is
postponed.

\begin{lem}\label{proj:est:lem}
Let $g_{\md_p}$ be the hyperbolic metric on $\md_p$. For
$q\in\Delta_p$ we have
\[
    D_\lambda^*(g_{\md_p})(q)=\eta k_\lambda(q)
\]
where $k_\lambda$ is as usual the intrinsic metric of 
$\Uu^0_\lambda(1)$ and $\eta$ is a positive number such that
\begin{equation}\label{hyperbolic:proj:eq}
    \log\eta>\int_{[N(p),N(q)]}\delta(t)\d\mu(t)+a(p,q) \ .
\end{equation}
where $\delta(t)$ is the distance of $N(q)$ from the stratum of $\lambda$
containing $t$ (that is a point on the geodesic segment $[N(p),N(q)]$) 
and $a(p,q)$ is equal to $|r(p)-r(q)$ if $N(p)=N(q)$ and $0$ otherwise.
\end{lem}

By Lemma~\ref{proj:est:lem}, $D_\lambda:\Delta_p\rightarrow\md_p$
increases the lengths.  Thus a classical argument shows that it is a
homeomorphism.  Since $\Delta_p$ contains $F_p=\Uu^0_\lambda(1)\cap
r^{-1}r(p)$ it is a maximal round ball.  Notice that on $F_p$ we have
$D_\lambda|_{F_p}=B(p_0,p)\circ N$. Thus the image of $F_p$ in $\md_p$
is an ideal convex set. Moreover if $\Delta'_p$ is the stratum
corresponding to $\Delta_p$ the same argument shows that
$F_p\subset\Delta_p'$ and in particular $\Delta'_p$ is the stratum
through $p$.

The map
\[
(D_\lambda)_{*,p}:T_p\Uu^0_\lambda(1)\rightarrow T_{D_\lambda(p)}\md_p
\] 
is a conformal map, moreover its restriction on $T_p F_p$ is an
isometry (with respect to the hyperbolic metric of $\md_p$). Thus it
is an isometry and this shows that $k_\lambda$ coincides with the
Thurston metric.

Finally we have to show that $F_p=\Delta'_p$. If
$q\notin F_p$, formula~(\ref{hyperbolic:proj:eq}) implies that
\[
(D_\lambda)_{*,q}:T_q\Uu^0_\lambda(1)\rightarrow T_{D_\lambda(q)}\md_p
\] 
is not an isometry. Thus $\Delta_p$ is different from
$\Delta_q$ and  $q\notin\Delta'_p$.
\cvd

\begin{figure}
\begin{center}
\input{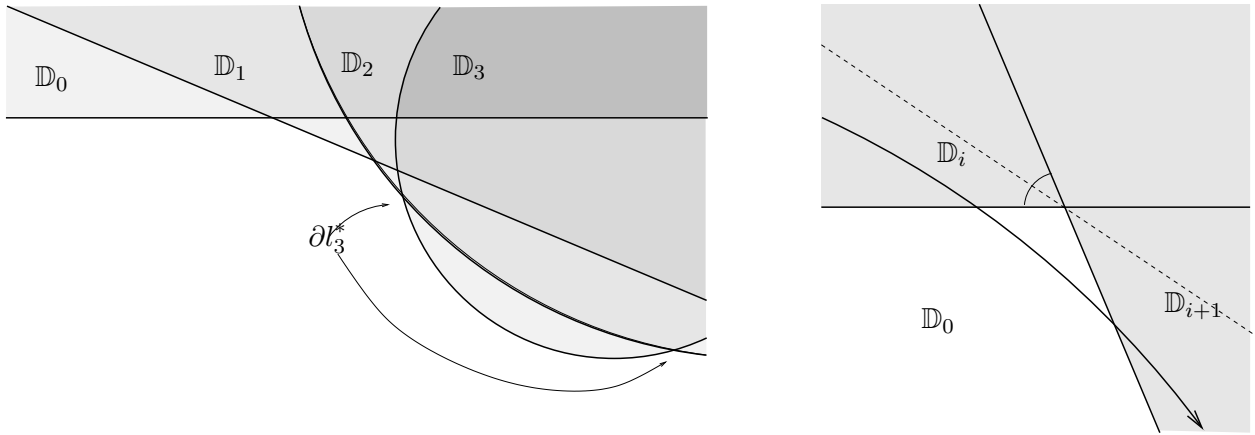}
\caption{{\small On the left it is shown how disks $\md_i$ intersect
    each other. On the right picture shows that if
    (\ref{hyperbolic:proj:eq3}) is not verified then
    (\ref{hyperbolic:proj:eq2}) does not hold}.}\label{proj:fig}
\end{center}
\end{figure}
Now we have to prove Lemma~\ref{proj:est:lem}. For the proof we need the
following estimate.

\begin{lem}\label{hyperbolic:proj:lem}
Let $l$ be an oriented geodesic of $\md_0$.  Denote by $R\in
PSL(2,\mc)$ the positive rotation around the corresponding geodesic of $\mh^3$
with angle $\alpha$.  Take a point $p\in\md_0$ in the left half-space
bounded by $l$ and suppose $R(p)\in\md_0$. Then if $g$ is the
hyperbolic metric on $\md_0$ 
\[
    R^*(g)(p)= \eta(p) g(p)
\]
where $\eta(p)=(\cos\alpha-\sh(d)\sin\alpha)^{-1}>0$ where $d$
is the distance from $p$ and $l$.
\end{lem}
\Dim Up to isometries we can identify $\md_0$ with the half-plane
$\{(x,y)|y>0\}$ in such a way that $l=\{x=0\}$ is oriented from $0$
towards $\infty$.  In these coordinates we have
\[
    R(x,y)=(x \cos\alpha - y\sin\alpha , x \sin\alpha + y \cos\alpha )\ .
\]
Since $p=(x,y)$ is in the left half-plane bounded by $l$ then $x<0$. 
Moreover as $R(p)\in\md_0$ we
have that $|y/x|>\tan\alpha$.  By an explicit computation we have
\[
   R^*(g)(p)=  \frac{1}{y\cos\alpha + x\sin\alpha }(\d x^2+\d y^2)
\]
then we see that $\eta(p)=(\cos\alpha-u\sin\alpha)^{-1}$ where
$u=|x/y|$.  On the other hand a classical hyperbolic formula shows that
$|x/y|=\sh d$ where $d$ is the distance of $p$ from $l$.  
\cvd

\emph{Proof of Lemma~\ref{proj:est:lem}:}
It is sufficient to consider the case $\lambda$ simplicial. The general case
will follows by an approximation argument.

Up to post-composition by an element of
$PSL(2,\mc)$ we can suppose that the point $p$ is the base point so
$\md_p$ is $\md_0$.  Take $q\in\Delta_p$ and consider a path $c$ in
$\Delta_p$ containing $p$ and $q$. The
intersection of every stratum of $\hat\lambda$ with $\Delta_p$ is
convex. Thus we can suppose that $c$ intersects every leaf at most
once.\\ Denote by $l_0,l_1,\ldots,l_n$ the leaves intersecting $N(c)$
and let $a_0,\ldots, a_n$ be the respective weights with the following
modifications.  If $q$ lies in a Euclidean band denote by $a_n$ the
distance from the component of the boundary of the band that meets
$c$.  In the same way if $p$ lies in a Euclidean band $a_1$ is the
distance of $p_1$ from the component of the boundary hit by $c$.
Finally if $p$ and $q$ lie in the same Euclidean boundary (that is the
case when $N(p)=N(q)$) then $n=1$ and $a_1$ is the distance between
$r(p)$ and $r(q)$ (that is by definition $a(p,q)$).\par Let us set
$B_i=\exp(a_{1}X_{1})\circ\cdots\circ\exp(a_iX_i)$ where $X_i$ is the
standard generator of the rotation around $l_i$.  Notice that
$B_n=\hat B(p,q)$.

We want to prove that $D_\lambda(q)\in\md_i=B_i(\md_0)$.  In
fact we will prove that $\md_i\cap\md_0$ is a decreasing sequence of
sets (with respect to the inclusion).  By the hypothesis on $c$ we
have that $\md_i\cap\md_0\neq\varnothing$. Moreover if we denote by
$X^*_{i+1}$ the standard generator of rotation around the geodesic
$l_{i+1}^*=B_i(l_{i+1})$ then
\begin{equation}\label{hyperbolic:proj:eq2}
   \exp(tX^*_{i+1})\md_i\cap\md_0\neq\varnothing
\end{equation}
for $0<t<a_{i+1}$
(in fact there exists a point $q'\in c$ lying on the Euclidean band of
$\Uu^0_\lambda(1)$ corresponding to
$l_{i+1}$ with distance from the left side equal to $t$ and $D_\lambda(q')$
lies in the intersection (\ref{hyperbolic:proj:eq2})).
Now by induction we can show that
\begin{equation}\label{hyperbolic:proj:eq3}
\left\{
\begin{array}{l}
  \md_0\cap\md_{i+1}\subset\md_0\cap\md_{i}\\
   \textrm{ the component of }\partial\md_i-\partial l^*_{i}\textrm{
   containing }l^*_{i+1}\textrm{ does not meet }\md_0\ .
\end{array}\right.
\end{equation}
Suppose $\md_0\cap\md_{i+1}$ is not  contained in
$\md_0\cap\md_i$.  Since $\md_{i+1}$ is obtained by the rotation along
$l_{i+1}^*$ whose end-points are outside $\md_0$ it is easy to see
that there should exist $t_0<a_{i+1}$ such that
$\exp(t_0X^*_{i+1})\md_i$ does not intersect $\md_0$ (see
Fig.~\ref{proj:fig}).\\

Let $g_i$ denote the hyperbolic metric on $\md_i$. We have that
$D_\lambda^*(g_n)$ is the intrinsic metric on $\Uu^0_\lambda(1)$ at
$q$. Moreover we have that
\[
 g_{i}(D_\lambda(q))=\eta_i g_{i+1}(D_\lambda(q))
\]
with $\eta_i^{-1}=\cos a_i-u_i\sin a_i$ where $u_i=\sh d_i$ where $d_i$ is the
distance of $N(q)$ from $l_i$.
Since $\eta=\prod_{i=0}^{n-1}\eta_i$ we obtain
\[
   -\log\eta=\sum\log(\cos a_i)+\sum\log(1-u_i\tan a_i) \ .
\]
Now $\log\cos(a_i)<-a_i^2/2$ and $\log(1-u_i\tan a_i)>d_ia_i$ so we get
\[
   \log\eta>\sum d_ia_i + a(p,q) \ .
\]
\cvd.
\begin{cor}
Every level surface $\Uu^0_\lambda(a)$ is equipped with a projective
structure.  Moreover, the corresponding Thurston distance is equal to
the intrinsic distance up to a scale factor.
\end{cor}

\Dim
The map
\[
   \Uu^0_\lambda\ni x\mapsto tx\in\Uu^0_{t\lambda}
\]
rescales the metric by a factor $t^2$. Moreover, it takes
$\Uu^0_\lambda(1/t)$ onto $\Uu^0_{t\lambda}(1)$.  \cvd

\section{Equivariant theory}\label{hyp:eq:sec}
Assume that the lamination $\lambda$ is invariant under the action of
a discrete group $\Gamma$, that is $\lambda$ is the lifting of a
measured geodesic lamination defined on some straight convex set $H$
of the hyperbolic surface $F=\mh^2/\Gamma$.  The following lemma,
proved in~\cite{Ep-M}, determines the behaviour of the cocycle
$B_\lambda$ under the action of the group $\Gamma$.
\begin{lem}\label{hyperbolic:group:lem1}
Let $\lambda$ be a measured geodesic lamination on $H$ invariant
under the action of $\Gamma$. Then if
$B_\lambda:\mathring H\times\mathring H\rightarrow PSL(2,\mc)$ is the cocycle
associated to $\lambda$ we have
\[
   B_\lambda(\gamma x,\gamma y)=\gamma\circ B(x,y)\circ\gamma^{-1}
\]
for every $\gamma\in\Gamma$.
\end{lem}
\cvd
Now let us fix base point $x_0\in H$ and  consider the bending map 
\[
   F_\lambda:\mathring H\rightarrow\mh^3\ .
\]

For $\gamma\in\Gamma$ let us define
\[
   h^1_\lambda(\gamma)= B_\lambda(x_0,\gamma x_0)\circ\gamma\in PSL(2,\mc)
\]
Lemma~\ref{hyperbolic:group:lem1} implies that $h^1_\lambda:\Gamma\rightarrow
PSL(2,\mc)$ is a homomorphism.
Moreover by definition it follows that $F_\lambda$ is
$h^1_\lambda$-equivariant.\\

On the other hand in~\ref{flat:equiv} we have seen that there exists a
homomorphism
\[
  h^0_\lambda:\Gamma\rightarrow\ISO_0(\mx_0)
\]
such that $\Uu^0_\lambda$ is $h^0_\lambda$-invariant and the
Gauss map is $h^0_\lambda$-equivariant that is
\[
    N(h^0_\lambda(\gamma)(p))=\gamma(N(p)) \ .
\]
By using this fact it is easy to see that
\[
   \hat B_\lambda(h^0_\lambda(\gamma)p,h^0_\lambda(\gamma)q)=
   \gamma\hat B_\lambda(p,q)\gamma^{-1}.
\]
hence that
\[
    \hat D_\lambda(h^0_\lambda(\gamma)p)=h^1_\lambda(\gamma)(D_\lambda(p)) \ .
\]
In particular we have that the map $\hat D_\lambda$ is a developing
map for a hyperbolic structure on $M_\lambda/h^0_\lambda(\Gamma)$.  The
completion of such a structure is a manifold with boundary
homeomorphic to $F\times [0,+\infty)$. The boundary is
isometric to $H$.\\

The map $D_\lambda:\Uu^0_\lambda(1)\rightarrow S^2_\infty$ is
$h^1_\lambda$-equivariant so it is a developing map for a projective
structure on $\Uu^0_\lambda(1)/\Gamma$.

Notice that given a marking $F\rightarrow\Uu^0_\lambda(1)/h^0_\lambda(\Gamma)$, by using the flow of the gradient of the
cosmological time, we obtain a marking
$F\rightarrow\Uu^0_\lambda(a)/h^0_\lambda(\Gamma)$. Thus we
obtain a path in the Teichm\"uller-like space of projective structures
on $F$ and clearly an underlying path of conformal structures in the
Teich\"uller space of $F$.

The following corollary is a consequence of
Theorems~\ref{hyperbolic:WR:teo}, ~\ref{hyperbolic:proj:teo}.
\begin{cor}
Let $Y$ be a maximal globally hyperbolic flat spacetime such that
$\tilde Y$ is a regular domain. If $T$ denotes the cosmological time
and $X$ is the gradient of $T$ then the Wick rotation on $Y(>1)$,
directed by $X$, with rescaling functions
\[  
   \alpha=\frac{1}{T^2-1} \qquad \beta=\frac{1}{(T^2-1)^2}
\]
is a hyperbolic metric. 

A projective structure of hyperbolic type  is defined on $Y(1)$  
by extending the developing map
on $\tilde Y(\geq 1)$. The intrinsic metric on $Y(1)$ coincides with the
Thurston metric associated to such a structure and $Y(>1)$ equipped 
with the hyperbolic metric given by the Wick
rotation coincides with the $H$-hull of the projective surface $Y(1)$.

The canonical stratification associated to the projective structure on 
$Y(1)$ coincides with the partition given by the fibers of the retraction
on $\tilde Y$.\par
The measured lamination corresponding to $\tilde Y$ (according to
Theorem~\ref{FULLFLAT}) coincides with the measured lamination
associated to the projective structure on $Y(1)$ (according to~\cite{Ku}).
\end{cor}
\cvd

Let $S$ be an orientable surface with non-Abelian fundamental group.
By~\cite{Ba} the universal coverings of maximal globally hyperbolic
flat spacetimes with a complete spacelike surface homeomorphic to $S$
are regular domains.  On the other hand, by~\cite{Ku}, projective
structures on $S$ are of hyperbolic type.

\begin{cor}
Let $S$ denote an orientable surface such that $\pi_1(S)$ is not
Abelian. Then the set of maximal globally hyperbolic flat
spacetimes containing a complete Cauchy surface homeomorphic to $S$ is
non-empty and, up to isometries, bijectively corresponds to the set of
projective structures on $S$.
\end{cor}
\cvd

\paragraph{Cocompact $\Gamma$-invariant case}
If the group $\Gamma$ is cocompact, we can relate this construction
with the Thurston parametrization of projective structures on a base
compact surface $F$ of genus $g\geq 2$. In fact, it is not hard to see
that the projective structure on
$\Uu^0_\lambda(1)/h^0_\lambda(\Gamma)$ is simply the structure
associated to $(\Gamma,\lambda)$ in Thurston parametrization.  We have
that the conformal structure on $\Uu^0_\lambda(1)/h^0_\lambda(\Gamma)$
is the {\it grafting} of $\mh^2/\Gamma$ along $\lambda$ (see
\cite{McM,Sc}).  It follows that the surface
$\Uu^0_\lambda(a)/h^0_\lambda(\Gamma)$ corresponds to
$gr_{\lambda/a}(F)$; $a\mapsto[\Uu^0_\lambda(a)/h^0_\lambda(\Gamma)]$
is a real analytic path in the Teichm\"uller space $\Tt_g$.  Such a
path has an endpoint in $\Tt_g$ at $F$ as $a\rightarrow +\infty$ and
an end-point in Thurston boundary $\partial\Tt_g$ corresponding to the
lamination $\lambda$ (or equivalently to the dual tree $\Sigma$).


\chapter{Flat vs de Sitter Lorentzian geometry}\label{dS}
In this chapter we will construct a map
\[
   \hat D:\Uu_\lambda(<1)\rightarrow\mx_1
\]
where $\mx_1$ is the de Sitter space. Hence $\hat D$ can be considered as
a developing map of spacetime $\Uu^1_\lambda$ of constant curvature
$\kappa = 1$.  The pull-back of the de Sitter metric is obtained by a
rescaling of the standard flat Lorentzian metric, directed by the
gradient of the cosmological time and with universal rescaling
functions.  The map $\hat D$ is the semi-analytic continuation of the
hyperbolic developing map $D$ constructed in the previous chapter,
regarding $\mh^3$ and $\mx_1$ as open sets of the real projective
space (Klein models), separated by the quadric $S^2_\infty$.  By
studying $\Uu^1_\lambda$ we will eventually achieve
Theorem ~\ref{DSRESC:I} and the statement (3) of of Theorem ~\ref{WR:I}
of Chapter ~\ref{INTRO}.

Finally, a suitable equivariant version of all constructions (together
with the results of Chapters ~\ref{FGHST} and ~\ref{HYPE}) will lead
us to the classification Theorem ~\ref{FULL_CLASS}, in the cases of
constant curvature $\kappa = 0,1$.
\begin{remark}\emph{ We will widely refer to \cite{Sc} in which maximal 
globally hyperbolic de Sitter spacetimes with {\it compact} Cauchy
surface are classified in terms of projective structures. Anyway we
have checked that essentially all constructions work as well by simply
letting the Cauchy surface be complete.}
\end{remark}
\section{Standard de Sitter spacetimes}\label{standarddS}
The main idea of \cite{Sc} is to associate to any projective structure
on a surface a so called \emph{standard} de Sitter spacetime.  It
turns out that the canonical de Sitter rescaling on $\Uu(<1)$ produces
the standard spacetime associated to the projective structure on
$\Uu(1)$ previously obtained thanks to the Wick rotation.  In this way
we will eventually see that, apart from a few exceptions, maximal
globally hyperbolic \emph{flat} spacetimes containing a complete
Cauchy surface homeomorphic to a given surface $F$, bijectively
corresponds to maximal globally hyperbolic \emph{de Sitter} spacetimes
with a complete Cauchy surface homeomorphic to the same surface $F$
(the canonical rescaling giving the bijection).
\smallskip

We start by recalling the construction of standard de Sitter
spacetimes corresponding to projective structure of hyperbolic type
(also called ``standard spacetimes of hyperbolic type''). This
construction is, in fact, dual to the construction of the
$H$-hulls.\par Given a projective structure of hyperbolic type on a
surface $S$ with developing map
\[
     d:\tilde S\rightarrow S^2_\infty
\]
recall the canonical stratification of $\tilde S$ (see Section
\ref{SPS}).  For every $p\in\tilde S$ let $U(p)$ denote the stratum
passing through $p$ and $U^*(p)$ be the maximal ball containing
$U(p)$.  Now $d(U^*(p))$ is a ball in $S^2_\infty$ so it determines a
hyperbolic plane in $\mh^3$. Let $\rho(p)$ denote the point in $\mx_1$
corresponding to this plane: the map $\rho:\tilde S\rightarrow\mx_1$
turns out to be continuous. There exists a unique timelike geodesic
$c_p$ in $\mx_1$ joining $\rho(p)$ to $d(p)$ so we can define the map
\[
    \hat d:\Delta\times(0,+\infty)\ni (p,t)\mapsto c_p(t)\in\mx_1
\]
This map is a developing map for a de
Sitter structure on $S\times (0,+\infty)$ that is called the {\it
standard spacetime} corresponding to the given projective structure.

\begin{remark}\emph{
In~\cite{Sc} a
standard spacetime is associated to every complex projective surface (also of
parabolic or elliptic type). However we will deal only with standard
spacetimes of hyperbolic type.
}\end{remark}
\section{The rescaling}\label{rescdS}
We are going to prove
\begin{teo}\label{memgen:dsmain:teo}
Let $\Uu$ be a regular domain . The spacetime, say $\Uu^1$, 
obtained by rescaling $\Uu(<1)$ along the gradient of the cosmological
time $T$ and rescaling functions
\begin{equation}\label{dsresc:eq}
   \alpha =\frac{1}{1-T^2} \qquad \beta= \frac{1}{(1-T^2)^2}\ .
\end{equation}
is a standard \emph{de Sitter} spacetime of hyperbolic type
corresponding to the projective structure on $\Uu(1)$ produced by the
Wick rotation.
\end{teo}
\Dim

Let $\lambda$ be the measured geodesic lamination defined on some straight
convex set $H$ such that $\Uu=\Uu^0_\lambda$.
We construct a map
\[
      \hat D:\Uu^0_\lambda(<1)\rightarrow\mx_1
\]
that, in a sense, is the map dual to the map $D$ constructed in the
previous chapter.  We prove that such a map is $\mathrm C^1$ and
the pull-back of the de Sitter metric is a rescaling of the flat
metric of $\Uu^0_\lambda$.  

The construction of $\hat D$ is very simple. In fact if $s$ is a geodesic
integral line of the gradient of cosmological time,
$s_{>1}=s\cap\Uu^0_\lambda(>1)$ is sent by $D$  onto a geodesic ray of
$\mh^3$. We define $\hat D$ on $s_{<1}$ in such a way that it
parameterizes the timelike geodesic ray in $\mx_1$ contained in the
projective line (in the Klein model) determined by $D(s_{>1})$
(that is the continuation of $D(s_{<1})$, see Section~\ref{dSMod}).\\

Let us be more precise. Consider the standard inclusion
$\mh^2\subset\mh^3$. Since $\mh^2$ is oriented there is a well-defined
dual point $v_0\in\hat\mx_1$ (that is the positive vector of the
normal bundle).

Now let us take the base point $x_0\in\mathring H$ for the
bending map and a corresponding point $p_0\in\Uu^0_\lambda(1)$. For
$p\in\Uu^0_\lambda$ let us define
\[
\begin{array}{l}
  v(p)=\hat B_\lambda(p_0,p)v_0\in\hat\mx_1\\  
  x(p)=\hat B_\lambda(p_0,p)N(p)=F_\lambda(N(p))\ .
\end{array}
\]

Thus let us set
\[
    \hat D(p)=[\ch\tau(p) v(p) +\sh\tau(p) x(p)]
\]
where we have put $\tau(p)=\arctgh T(p)$ (notice that for $p\in\Uu^0_\lambda(>1)$ we have $D(p)=[\ch\delta(p) x(p)+\sh\delta(p) v(p)]$ where $\delta(p)=\arctgh 1/T(p)$).

We claim that the map
\[
   \hat D:\Uu^0_\lambda(<1)\rightarrow\mx_1
\]
is $\mathrm C^1$-local diffeomorphism. The pull-back of the metric of
$\mx_1$ is the rescaling of the metric of $\Uu^0_\lambda(<1)$ along the
gradient of $T$ with rescaling functions
\begin{equation}\label{resc:ds:eq}
\begin{array}{ll}
   \alpha = \frac{1}{1-T^2}\,\qquad & \beta =\frac{1}{(1-T^2)^2}\ .
\end{array}
\end{equation}

Clearly the claim proves Theorem~\ref{memgen:dsmain:teo}.

 The proof of the claim is quite similar to the proof of
Theorem~\ref{hyperbolic:WR:teo}.  In fact by an explicit computation
we get the result in the case when $\lambda$ is a weighted geodesic.
Thus the statement of the theorem holds when $\lambda$ is a simplicial
lamination. Moreover by proving the analogous of
Lemma~\ref{hyperbolic:WR:lem} and using standard approximations we
obtain the proof of the general case.

Then the same argument of
Theorem~\ref{hyperbolic:WR:lem} works in the same way and we omit
details.

\begin{figure}
\begin{center}
\input{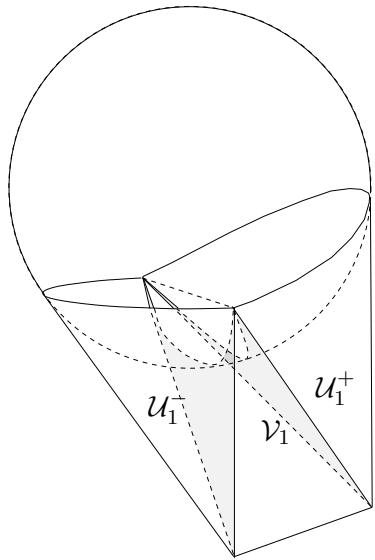}
\caption{{\small The domain image of $\hat D_0$.}}
\end{center}
\end{figure}

To make the computation for $\Uu_0=\Uu^0_{\lambda_0}$, let us use the
same notation as in Section~\ref{theWR}. In particular let us identify
$\mx_{-1}$ as the set of spacelike lines through $0$ of $\mm^4$ and
consider coordinates $(u,\zeta,T)$ on $\Uu_0$ given
in~(\ref{flatcoord}).

In these coordinates we have
\[
\hat D_0(T,u, \zeta)\mapsto\left\{\begin{array}{l}
                          \sh\tau\left(\ch\zeta\ch u,\ \ch\zeta\sh u,\
                          \sh\zeta,\ 0\right)\ +\ \ch\tau(0,0,0,1)  \\
                          \textrm{if }\eta\leq 0\ ;\\ 
                          \sh\tau\left(\ch  u,\ \sh u,\ 0,0\right)\ + \
                          \ch\tau\left(0,\ 0,\
                          -\sin (\zeta\tgh\tau),\,
                          \cos(\zeta\tgh\tau)\right)  \\
                           \textrm{if }\eta\in [0,\alpha_0/T] \\
                          \sh\tau\left(\ch\zeta'\ch u,\ \ch\zeta'\sh
                          u,\ \sh\zeta'\cos\alpha_0,\
                          \sh\zeta'\sin\alpha_0\right)\ +  \\
                           \ch\tau(0,0,-\sin\alpha_0,\
                          \cos\alpha_0)\\ \textrm{otherwise} 
                           \end{array}\right.
\] 
where $\tau=\tgh T$ and  $\zeta'=\eta-\alpha_0/T$.

This map is clearly smooth for $\zeta\neq 0,\alpha_0/T$.
Since the derivatives of $\hat D_0$ with respect the coordinates fields glue 
along $\zeta=0$ and $\zeta=\alpha_0 T$ the map $D$ is $\mathrm C^1$.

By computing the pull-back of the de Sitter metric we have 
\begin{equation}\label{dsmetric}
D_0^*(g)(T,\zeta, u)=\left\{\begin{array}{ll}
                           -\d\tau^2+
                    \sh^2\tau(\d\zeta^2+\ch^2\zeta\d u^2) &
                           \textrm{ if }\zeta<0\\
                           -\d\tau^2 + \sh^2\tau(\d\zeta^2+\d u^2) & 
                           \textrm{ if }\zeta\in[0,\alpha_0/T]\\
                   -\d\tau^2+\sh^2\tau(\d\zeta^2+\ch^2(\zeta')\d u^2) &
                           \textrm{ otherwise}.
                            \end{array}\right.
\end{equation}
Since $\d\tau=\frac{1}{1-T^2}\d T$ and $\sh^2\tau=\frac{T^2}{1-T^2}$,
comparing (\ref{dsmetric}) and (\ref{flatmetric}) shows that $\hat
D_0^*(g)$ is obtained by the rescaling along the gradient of $T$ with
rescaling functions given in~(\ref{dsresc:eq}).  \cvd

\begin{cor}
The rescaling of Theorem~\ref{memgen:dsmain:teo} gives rise to a
bijective correspondence between the set of regular domains and the set of
standard de Sitter spacetimes of hyperbolic type.
\end{cor}
\cvd

\begin{cor}
Let $\Sigma$ be the initial singularity of $\Uu^0_\lambda$ defined in
Section~\ref{RD:sec}.  The map $\hat D$ extends to a continuous
map
\[
   \Uu^0_\lambda(\leq 1)\cup\Sigma\rightarrow \mx_1\cup S^2_\infty \ .
\]
Moreover, $\hat D$ restricted to $\Uu^0_\lambda(1)$ coincides with $D$.
\end{cor}
The extension of $\hat D$ on $\Uu^0_\lambda(1)$ follows by construction. On
the other hand, we see that the cocycle $\hat B$ is induced by a
cocycle
\[
    \overline B:\Sigma\times\Sigma\rightarrow\SOO^+(3,1)
\]
and $\hat D$ can be extended to $\Sigma$ by putting
\[
  \hat D(r)= \overline B(r(p_0),r)v_0 \ .
\]
\cvd
\begin{remark}\emph{
By ~(\ref{hyperbolic:bend:lip3:eq}), the cocycle
\[
   \overline B:\Sigma\times\Sigma\rightarrow\SOO^+(3,1)
\]
can be shown to be continuous with respect to the intrinsic distance $d_\Sigma$.
}\end{remark}

\begin{remark}\emph{
The above construction allows to identify $\Sigma$ with
the space of maximal round balls of $\Uu^0_\lambda(1)$.  }
\end{remark}

In what follows, we denote by $\Uu^1_\lambda$ the domain
$\Uu^0_\lambda(<1)$ endowed with the de Sitter metric induced by
$\hat D$.
\begin{prop}\label{desitter:ct:prop}
The cosmological time of $\Uu^1_\lambda$ is
\[
   \tau=\arctgh(T).
\]
Every level surface $\Uu^1_\lambda(\tau=a)$ is a Cauchy surface.
\end{prop}
\Dim Let $\gamma:[0,a]\rightarrow\Uu^1_\lambda$ denote a timelike-path
with future end-point $p$ parametrized in Lorentzian arc-length.  Now
as path in $\Uu^0_\lambda$ we have a decomposition of $\gamma$
\[
   \gamma(t)=r(t)+T(t)N(t) \ .
\]
By computing derivatives we obtain
\[
  \dot\gamma =\dot r+ T\dot N + \dot T N
\]
so the square of de Sitter norm is
\begin{equation}\label{desitter:ct:eq}
   -1=-\frac{\dot T^2}{(1-T^2)^2}+ \frac{|\dot r + T\dot N|^2}{1-T^2} 
\end{equation}
where $|\cdot|$ is the Lorentzian flat norm.
It follows that
\[
   1< \frac{\dot T}{1-T^2}
\]
and integrating we obtain
\[
   \arctgh T(p)-\arctgh T(0)> a
\]
\emph{i.e.} the de Sitter proper time of $\gamma$ is less than
$\arctgh T(p)$.  On the other hand the path $\gamma(t)=r(p)+tN(p)$ for
$t\in[0,T(p)]$ has proper time $\arctgh T(p)$ so the cosmological time
of $\Uu_\lambda^1$ is
\[
   \tau=\arctgh T \ .
\]

Now let $\gamma:(a,b)\rightarrow\Uu^1_\lambda$ be an inextensible
timelike-curve parametrized in Lorentz arc-length such that
$\gamma(0)=p$. We want to show that the range of $T(t)=T(\gamma(t))$
is $(0,1)$.\par Suppose $\beta=\sup T(t)<1$. Since $T(t)$ is
increasing then $\beta=\lim_{t\rightarrow b}T(t)$.  Then the path
\[
   c(t)=r(t)+N(t)
\]
should be inextensible (otherwise we could extend $\gamma$ in
$\Uu^1_\lambda$).
Now we have
\[
   \dot c =\dot r + \dot N \ .
\]
For $t>0$ we have $T(t)>T(p)=T_0$ so
\[
  T_0 |\dot c| <|\dot r+T\dot N|
\]
Multiplying by the horizontal rescaling factor we have
\[
     \frac{T_0}{\sqrt{1-T^2}}|\dot c|\leq\frac{|\dot r+T\dot
     N|}{\sqrt{1-T^2}}\,.
\]
Since $T(t)<\beta<1$ it results
\[
 \frac{T_0}{\sqrt{1-\beta^2}}
  |\dot c|\leq\frac{|\dot r+T\dot
     N|}{\sqrt{1-T^2}}.
\]
By looking at equation~(\ref{desitter:ct:eq}) we deduce
\[
 \frac{T_0}{\sqrt{1-\beta^2}}|\dot c|\leq  \frac{\dot T}{1-T^2} \ .
\]
Thus the length of $c$ is bounded. On the other hand since
$\Uu^0_\lambda(1)$ is complete it follows that $c$ is extensible.  Thus
we have proved that $\sup T(t)=1$.  The same computation applied to
$\gamma(a,0)$ shows that $\inf T(t)=0$.  \cvd

\begin{cor}
Any standard de Sitter spacetime of hyperbolic type 
contains a complete Cauchy surface. Moreover its cosmological time is regular.
\end{cor}
\cvd

\section{Equivariant theory}

Suppose $F=\mh^2/\Gamma$ be a complete hyperbolic surface and $\lambda$ be a
measured geodesic lamination defined on some straight convex set $H$ of $F$.

We have seen that
there exists an affine deformation of $\Gamma$
\[
    h^0_\lambda:\Gamma\rightarrow\ISO_0(\mx_0)
\]
such that $\Uu^0_\lambda$ is $h^0_\lambda(\Gamma)$-invariant
and the Gauss map is $h^0_\lambda$-equivariant.  Moreover in
the previous section we have constructed a representation
\[
    h^1_\lambda:\Gamma\rightarrow PSL(2,\C)=\SOO^+(3,1)
\]
such that
\[
   D\circ h^0_\lambda(\gamma)=h^1_\lambda(\gamma)\circ D 
\qquad\textrm{ for }\gamma\in\Gamma.
\]
Now it is straightforward to see that the same holds changing $D$ by $\hat D$.

Thus $\hat D$ is a developing map for a maximal globally hyperbolic
structure, say $Y^1_\lambda$ on $F\times\mr$.  By construction,
$Y^1_\lambda$ is the standard de Sitter space-time associated to
$\Uu^0_\lambda(1)/h^0_\lambda(\Gamma)$ (that carries a natural
projective structure by~\ref{proj-bound}).  Notice that $Y^1_\lambda$
is obtained by a canonical rescaling on
$Y^0_\lambda=\Uu^0_\lambda/h^0_\lambda(\Gamma)$ along the
gradient of the cosmological time with rescaling functions as
in~\ref{memgen:dsmain:teo}.

Let $S$ be a compact closed surface. In~\cite{Sc} 
Scannell proved that  the universal covering of any
maximal globally hyperbolic de Sitter spacetime $\cong S\times\mr$ 
is a standard spacetime. In fact, the same argument proves the following a bit
more general fact.

\begin{prop}
Let $S$ be any surface and $M$ be a maximal globally hyperbolic de
 Sitter spacetimes containing a \emph{complete} Cauchy surface $\cong
 S$. Then it is a standard spacetime.
\end{prop}

As a consequence of this proposition we get the following
classification theorem.

\begin{teo}\label{desitter:class:teo}
The  correspondence
\[
   Y_\lambda\rightarrow Y^1_\lambda
\]
induces a bijection between flat and de Sitter maximal globally hyperbolic 
spacetimes  admitting regular cosmological time and a complete Cauchy surface.
\end{teo}
\cvd

\begin{cor}
Let $S$ be a surface with non-Abelian fundamental group.  The
rescaling given in Theorem~\ref{DSRESC:I} establishes a
bijection between flat and de Sitter maximal hyperbolic spacetimes
admitting a complete Cauchy surface homeomorphic to $S$.
\end{cor}
\cvd 


\chapter{Flat vs AdS Lorentzian geometry}
\label{AdS}
First we perform a canonical rescaling on any given flat regular
domain $\Uu^0_\lambda$, obtaining a globally hyperbolic AdS spacetime
$\Pp_\lambda$ containing complete Cauchy surfaces. The AdS
$\Mm\Ll$-spacetime $\Uu^{-1}_\lambda$ is by definition the maximal
globally extension of $\Pp_\lambda$. This AdS canonical rescaling runs
parallel to the Wick rotation of Chapter \ref{HYPE}. Every
spacelike plane $P$ is a copy of $\mh^2$ into the Anti de Sitter space
$\mx_{-1}$. So the core of the construction consists in a suitable
{\it bending procedure} of $P$ along any given $\lambda \in \Mm\Ll$.
However, in details there are important differences. Both spacetime
and time orientation will play a subtle r\^ole.

Then we characterize the class of our favourite simply connected
maximal globally hyperbolic AdS spacetimes as to coincide with the
class of so called {\it standard} AdS spacetimes ({\it i.e} the {\it
Cauchy developments of achronal curves} on the boundary of
$\mx_{-1}$), and we study their geometry. In particular this allows us
to recognize $\Pp_\lambda$ as the {\it past part} of
$\Uu^{-1}_\lambda$ (that is {\it the past of the future boundary of
its convex core}).  Finally we show that
$$\lambda \to \Uu^{-1}_\lambda$$ actually establishes a bijection onto
the set of maximal globally hyperbolic AdS spacetimes containing a
complete Cauchy surface. By combining all these results (including
their equivariant version) we will eventually prove Theorem
~\ref{ADSRESC:I}, (2) of Theorem ~\ref{WR:I}, and the case $\kappa =
-1$ of Theorem ~\ref{FULL_CLASS}, stated in Chapter ~\ref{INTRO}.
At the end of the Chapter we discuss (broken) $T$-symmetry
and relations with the theory of generalized earthquakes.

\section{Bending in AdS space}\label{bend-ads}
The original idea of bending a spacelike plane in $\mx_{-1}$ was
already sketched in \cite{M}.  We go deeply in studying this
notion and we relate it to the bending cocycle notion of Epstein and
Marden.

First let us describe a {\it rotation around a spacelike geodesic
$l$}. By definition such a rotation is simply an isometry $T$ which
point-wise fixes $l$.  Up to isometries $l$ can be supposed to lie on
$P_0=P(Id)$, that is the plane dual to the identity in $PSL(2,\mr)$
(see Section~\ref{AdSMod}). The dual geodesic $l^*$ is a hyperbolic
$1$-parameter subgroup, as we have remarked in Section~\ref{AdSMod}.

\begin{lem}\label{adesitt:bend:lem}
Let $l$ be a geodesic contained in $P_0$ and $l^*$ denote its dual
line.  For $x\in l^*$, the pair $(x,
x^{-1})\in PSL(2,\R)\times PSL(2,\R)$ represents a rotation around
$l$. The map
\[   R:l^*\ni x\mapsto (x,x^{-1})\in PSL(2,\R)\times PSL(2,\R)\]
is an isomorphism onto the subgroup of rotations around $l$. 
\end{lem}
\Dim First of all, let us show that the map
\[  \mx_{-1}\ni y\mapsto xyx\in\mx_{-1}\]
fixes point-wise $l$ (clearly $l$ is invariant by this transformation
because so is $l^*$). If $c$ is the axis of $x$ considered as an
isometry of $\mh^2$, then 
 $l$ is the set of rotations by $\pi$ around points in
$c$. Thus it is enough to show that if $p$ is the fixed point of $y\in l$,
then \[ xyx(p)=p \ .\] If we orient $c$ from the repulsive fixed point of
$x$ towards the attractive one, $x(p)$ is obtained by translating $p$
along $c$ in the positive direction, in such a way that $d(p,x(p))$ is
the translation length of $x$.  Since $y$ is a rotation by $\pi$ around
$p$, we have that $yx(p)$ is obtained by translating $p$ along $c$ in the
\emph{negative} direction, in such a way $d(p,yx(p))=d(p,x(p))$.  Thus we get
$xyx(p)=p$.

Now $R$ is clearly injective. On the other hand, 
the group of rotations around a geodesic has
dimension at most $1$ (for the differential of a rotation at $p\in l$
fixes the vector tangent to $l$ at $p$).  Thus $R$ is surjective onto
the set of rotations around $l$. \cvd

\begin{cor}
Rotations around a geodesic $l$ act freely and transitively on the
dual geodesic $l^*$. Such action induces an isomorphism between the
set of rotations around $l$ and the set of translations of $l^*$. 

By duality, rotations around $l$ act freely
and transitively on the set of spacelike planes containing $l$.  
Given two spacelike planes $P_1,\ P_2$ such that $l\subset P_i$, then
there exists a unique rotation $T_{1,2}$ around $l$ such that
$T_{1,2}(P_1)=P_2$.
\end{cor}  
\cvd

Given two spacelike planes $P_1,\ P_2$ meeting each other along a
geodesic $l$, the dual points $x_i=x(P_i)$ lie on the
geodesic $l^*$ dual to $l$.  Then we define the {\it angle between
$P_1$ and $P_2$} as the distance between $x_1$ and $x_2$ along $l^*$.
Notice that:
\smallskip

{\it By varying the couple of distinct spacelike 
planes, the angles between them are well defined numbers that
span the whole of the interval $(0,+\infty)$}.
\smallskip

\noindent This is a difference with respect to the hyperbolic case,
that will have important consequences for the result of the bending
procedure.
\begin{figure}
\begin{center}
\input{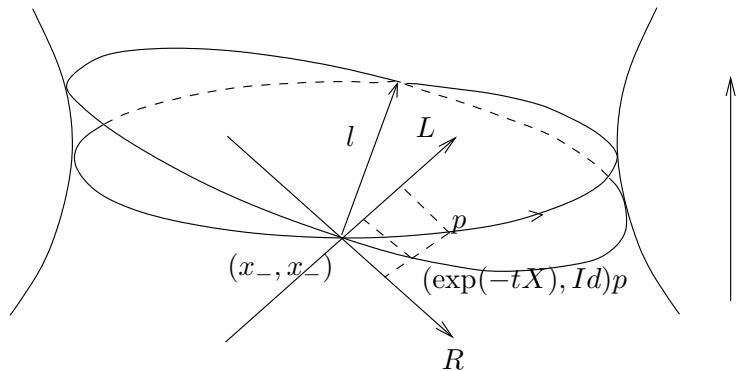}
\caption{{\small $(\exp(-tX),Id)$ rotates planes around $l$ in the
positive sense.}}\label{ADS:rot:fig}
\end{center}
\end{figure}
\begin{cor} 
An isometry $T$ of $\mx_{-1}$ is a rotation around a geodesic if and
only if it is represented by a pair $(x,y)$ such that $x$ and $y$ are
isometries of $\mh^2$ of hyperbolic type with the same translation length. 

Given two spacelike planes $P_1,P_2$ meeting along a geodesic $l$,
let $(x,y)$ be the rotation taking $P_1$ to $P_2$. Then the  
translation length $\tau$ of $x$ coincides with the angle between $P_1$ 
and $P_2$.
\end{cor}
\Dim Suppose that $(x,y)$ is a pair of hyperbolic transformations with
the same translation length. Then there exists $z\in PSL(2,\R)$ such
that $zyz^{-1}=x^{-1}$.  Hence $(1,z)$ conjugates $(x,y)$ into
$(x,x^{-1})$. Thus $(x,y)$ is the rotation along the geodesic
$(1,z)^{-1}(l)$ where $l$ is the axis of $(x,x^{-1})$.  

Conversely, if $(x,y)$ is a rotation, it is conjugated to a
transformation $(z,z^{-1})$ with $z$ a hyperbolic element of
$PSL(2,\R)$. Thus, $x$ and $y$ are hyperbolic
transformations with the same translation length.

In order to make the
last check, notice that, up to isometry, we can suppose
$P_1=P(Id)$. Thus, if $(x,x^{-1})$ is the isometry taking $P_1$
onto $P_2$, then the dual points of $P_1$ and $P_2$ are $Id$ and $x^2$
respectively. If $d$ is the distance of $x^2$ from $Id$ there exists a
unitary spacelike element $X\in\sG\lG(2,\mr)$ such that
\[
    x^2=\ch d I +\sh d X
\]
Thus we obtain that $\tr x^2=2\ch d$. By a classical identity,  
$\tr x^2=2\ch u/2$ where $u$ is the translation length of $x^2$. Since
$u=2\tau$ the conclusion follows. \cvd
\smallskip

There is a natural definition of positive rotation around an oriented
spacelike geodesic $l$ (depending only on the orientations of $l$ and
$\mx_{-1}$).  Thus, an orientation on the dual line $l^*$ is induced
by requiring that positive rotations act by positive translations on
$l^*$.

In particular, if we take an oriented geodesic
$l$ in $P(Id)$, and denote by $X$ the infinitesimal generator
of positive translations along $l$ then it is not difficult to show that the
positive rotations around $l$ are of the form $(\exp(-tX),\exp(tX))$ for
$t>0$. Actually, by looking at the action on the boundary we  deduce
that both the maps $(\exp(-tX),Id)$ and
$(Id,\exp(tX))$ rotate planes through $l$ in the positive direction
(see Fig.~\ref{ADS:rot:fig}).

\subsection{AdS bending cocycle}\label{ads:bend:coc:sec}
We can finally define the bending along a measured geodesic lamination.
First, take a finite measured geodesic lamination $\lambda$ of
$\mh^2$.  Take a pair of points $x,y \in\mh^2$ and enumerate the
geodesics in $\lambda$ that cut the segment $[x,y]$ in the natural way
$l_1,\ldots,l_n$.  Moreover, we can orient $l_i$ as the boundary of
the half-plane containing $x$.  With a little abuse, denote by $l_i$ 
also the geodesic in $P(Id)$
corresponding to $l_i$, then let $\beta(x,y)$ be the isometry of
$\mx_{-1}$ obtained by composition of positive rotations around $l_i$
of angle $a_i$ equal to the weight of $l_i$.  In particular, if $X_i$
denotes the unit positive generator of the hyperbolic transformations with
axis equal to $l_i$, then we have
\[
\begin{array}{l}
\beta_\lambda(x,y)=
(\beta_-(x,y),\beta_+(x,y))\in PSL(2,\R)\times PSL(2,\R)\qquad\textrm
{where}\\
\beta_-(x,y)=\exp(-a_1 X_1/2)\circ\exp(-a_2
X_2/2)\circ\ldots\circ \exp(-a_n X_n/2)\\ 
\beta_+(x,y)=\exp(a_1
X_1/2)\circ\exp(a_2 X_2/2) \circ\ldots\circ\exp(a_n
X_n/2)
\end{array}
\] 
with the following possible modifications: $a_1$
is replaced by $a_1/2$ when $x$ lies on $l_1$ and $a_n$ is replaced by
$a_n/2$ when $y$ lies on $l_1$ The factor $1/2$ in the definition of
$\beta_\pm$ arises because the  translation length of $\exp tX$ is
$2t$. 

Notice that $\beta_-$ and $\beta_+$ are the Epstein-Marden cocycles
corresponding to the {\it real-valued} measured laminations $-\lambda$
and $\lambda$.  Thus, for a general lamination $\lambda=(H,\Ll,\mu)$,
we can define $\beta_\lambda(x,y)$ for $x,y\in\mathring H$, just by
taking the limit of $\beta_{\lambda_n}(x,y)$, where $(\lambda_n)$ is a
standard approximation of $\lambda$ in a box containing the segment
$[x,y]$. The convergence of $\beta_{\lambda_n}(x,y)$ is proved
in~\cite{Ep-M}.

\begin{remark}\label{E-M-lambda}
{\rm We stress that the above $-\lambda$ is just obtained from
$\lambda = (\Ll,\mu)$ by taking the {\it negative-valued} measure
$-\mu$. Although this is no longer a measured lamination in the sense
of Section \ref{laminations}, the construction of \cite{Ep-M} does
apply. In Section \ref{der} (and in the Introduction) we  use the
notation $-\lambda$ in a different context and with a different
meaning.}
\end{remark} 

Let us enlist some properties of the bending cocycle that will be
useful in this work.

\begin{enumerate}
\item $\beta_\lambda(x,y)\circ\beta_\lambda(y,z)=\beta_\lambda(x,z)$ 
for every $x,y,z\in\mathring H$ (this means that $\beta_\lambda$ is a 
$PSL(2,\R)\times PSL(2,\R)$-valued cocycle);
\item $\beta_\lambda(x,x)=Id$;
\item $\beta_\lambda$ is constant on the strata of the stratification
determined by $\lambda$.
\item If $x,y$ lie in different strata then $\beta_+(x,y)$
(resp. $\beta_-(x,y)$) is a non-trivial hyperbolic transformation
whose axis separates the stratum through $x$ and the stratum through
$y$. Moreover the translation length is bigger than the total mass of
$[x,y]$.
\item If $\lambda_n\rightarrow \lambda$ on a $\eps$-neighbourhood of
the segment $[x,y]$ and $x,y\notin L_W$, then
$\beta_{\lambda_n}(x,y)\rightarrow\beta_{\lambda}(x,y)$.
\end{enumerate}

Properties $1),2),3)$ follow by definition, property $5)$ is proved
in~\cite{Ep-M}.  Finally property $4)$ follows from the following lemma.
\begin{lem}\label{ads:bend:lem}
If $g,h\in PSL(2,\mr)$ are hyperbolic transformations whose axes are
disjoint and point in the same direction then the composition $g\circ
h$ is hyperbolic, its axis separates the axis of $g$ from the axis of
$h$ and its translation length is bigger than the sum of the
translation lengths of $g$ and $h$.
\end{lem}
\Dim 
Notice that 
\[
\begin{array}{l}
g=\ch d_1 I + \sh d_1 X\\
h=\ch d_2 I + \sh d_2 Y
\end{array}
\]
with $X,Y\in\sG\lG(2,\mr)$ unitary elements and $d_1$ and $d_2$ equal
respectively to the half of the translation lengths of $g$ and $h$.
The plane generated by $X,Y$ is not spacelike: otherwise its dual
point (in $P(id)$) would be the intersection of the axis of $g$ and of
$h$.  Thus the reverse of the Schwarz inequality holds
\[
  1/2\ \tr XY=\eta(X,Y)\geq 1
\]
that, in turn, implies
\[
  1/2\ \tr (g\circ h)=\ch d_1\ch d_2 +1/2\sh d_1\sh d_2 \tr(XY)\geq
  \ch(d_1+d_2)\,.
\]
Since the interval $I_+$ (resp. $I_-$) in $\partial\mh^2$ whose
end-points are respectively the attractive (resp. repulsive)
end-points of $g$ and $h$ is sent into itself by $g\circ h$ (resp
$(g\circ h)^{-1}$) it contains the attractive (resp. repulsive) fixed
point of $g\circ h$. Thus the axis of $g\circ h$ separates the axis of
$g$ from the axis of $h$.  \cvd

Another important property of $\beta_\cdot$ is that for close points
$x,y$ the map $\beta_\lambda(x,y)$ is approximatively equal to a
hyperbolic transformation whose axis is a leaf of $\lambda$ cutting
$[x,y]$ and whose translation length is the total mass $m$ of
$[x,y]$. The following lemma gives a more precise estimate. Similarly
to Lemma~\ref{hyperbolic:bend:cont:lem}, it is an immediate
consequence of Lemma~3.4.4~(Bunch of geodesics) of~\cite{Ep-M} applied
to real-valued measured geodesic laminations.

\begin{lem}\label{adesitt:bend:cont:lem}
For any compact set $K$ in $\mh^2$ and any $M>0$, there exists a
constant $C>0$ with the following property. Let $\lambda=(H,\Ll,\mu)$
be a measured geodesic lamination on a straight convex set $H$ such
that $K\subset\mathring H$ and the total mass of any geodesic segment
joining points in $K$ is bounded by $M$.  For every $x,y\in K$ and
every geodesic line $l$ of $\Ll$ that cuts $[x,y]$, let $X$ be the
unit infinitesimal positive generator of the hyperbolic group with
axis $l$ and $m$ be the total mass of $[x,y]$. Then we have
\[     
||\beta_\lambda(x,y)-(\exp(-m/2X),\exp(m/2X))||<Cmd_\mh(x,y)\,.
\]
(On $\ PSL(2,\R)\times PSL(2,\R)$ the product norm of
the norm of $\ PSL(2,\R)$ is considered.)
\end{lem} \cvd

\subsection{AdS bending map}\label{ads-bend-map}

Take a base point $x_0$ in $\mathring H$. The {\it bending map} of $H$
with base point $x_0$ is simply
\[\varphi_\lambda: \mathring H\ni x\mapsto\beta_\lambda(x_0,x)x.\]
 
\begin{prop}\label{ads:bend:descr:prop}
The map $\varphi_\lambda$ is a local isometric $\mathrm C^0$ embedding
of $\mathring H$ into $\mx_{-1}$.
\end{prop}
\Dim By local isometric $\mathrm C^0$ embedding, we mean a
Lipschitzian map that preserves the length of the rectifiable paths
(and in particular sends rectifiable paths to spacelike paths of
$\mx_{-1}$).

Lemma~\ref{adesitt:bend:cont:lem} implies that $\varphi_\lambda$ is locally
Lipschitzian. 
Take a rectifiable arc  $k$ of $\mathring H$ parameterized in arc length
and let $\mu$ be the transverse measure on $k$.  
We claim that
\begin{equation}\label{pp}
   \frac{\d}{\d t}\left(\varphi_\lambda\circ k\right)(t)=\beta(x_0, k(t))\dot k(t)
\end{equation}
for almost every $t$ (with respect the Lebesgue measure $d t$).  From
(\ref{pp}) it follows that the length of rectifiable arcs is preserved
by $\varphi_\lambda$, so it turns to be a local $\mathrm
C^0$-embedding.

To prove the claim, take a point $t_0$.  We can suppose that $t_0$
lies in a bending line $l$ (the other case being obvious). Moreover we
can suppose that $l$ is not weighted (in fact there are at most
numerable many $t$ such that $k(t)\in L_W$)

By Lemma~\ref{adesitt:bend:cont:lem} there exists a constant $K$ such
that
\begin{equation}\label{compl2:eq}
||\beta(k(t_0), k(t_0+\eps))-(\exp(-m_\eps X_l/2), \exp(m_\eps X_\eps/2))||\leq
  K m_\eps\eps
\end{equation}
where we have put $m_\eps=\mu([k(t_0), k(t_0+\eps)])$.  Now we have
\[
  \varphi_\lambda(k(t_0+\eps))=\beta(x_0, k(t_0))\circ \beta(k(t_0),
  k(t_0+\eps)) k(t_0+\eps)
\]
so we can write
\begin{eqnarray*}
  \frac{1}{\eps}\left(\beta(k(t_0), k(t_0+\eps))k(t_0+\eps)- k(t_0)\right) =\\
  \frac{1}{\eps}\left(\beta(k(t_0), k(t_0+\eps))-(\exp(-m_\eps X_l/2),
  \exp(m_\eps X_l/2))\right) k(t_0+\eps) +\\
   \frac{1}{\eps} \left((\exp(-m_\eps X_l/2),\exp(m_\eps
  X_l/2))-Id\right)k(t_0+\eps) +\\
   \frac{1}{\eps} (k(t_0+\eps)-k(t_0))\,.
\end{eqnarray*}
The first term on the right hand tends to $0$ because of (\ref{compl2:eq})
(and the assumption that $k(t_0)\notin L_W$).

The last term converges to $\dot k(t_0)$.

Finally, for almost every $t_0$ the second term converges to
\[
   -\frac{h(t_0)}{2}(X_l k(t_0)+ k(t_0)X_l)
\]
where $h(t_0)$ is the derivative of $\mu$ with respect the Lebesgue
measure.  Deriving the identity $\exp(-t X_l)k(t_0)\exp(t X_l)=k(t_0)$
shows that the last quantity is $0$.  This proves that (\ref{pp})
holds for almost every point.  \cvd

\begin{remarks}\label{ads:bend:inj:rem}
{\rm 
(1) It turns out that the map $\varphi_\lambda$ is always
injective onto an achronal set of $\mx_{-1}$.  This fact will follow
as a corollary of the rescaling theory will describe in the next
sections (see Proposition~\ref{Plambda-past}). 
However, the reader could directly check it by proving that given
$x,y\in\mathring H$ the transformation
$\beta_+(x,y)I(y)\beta_-(x,y)^{-1}I(x)$ is a non trivial hyperbolic
element of $PSL(2,\mr)$ (that means that $\varphi_\lambda(x),
\varphi_\lambda(y)$ are joint by a non trivial spacelike geodesic
segment).
\smallskip

(2) Suppose $H\neq\mh^2$ and consider the behaviour of
$\varphi_\lambda$ in a neighbourhood of a boundary component, say $l$,
of $H$.  If $l$ is a weighted leaf then $\beta_\lambda(x_0,y)$
converges in $PSL(2,\mr)$ as $y$ goes towards $l$.  Thus
$\varphi_\lambda$ extends to $l$ and the image of $l$ is a spacelike
geodesic of $\mx_{-1}$. If $l$ is not weighted then property
4. implies that $\beta_\lambda(x_0,y)$ is not convergent and
$\varphi_\lambda$ does not extend on $l$.  On the other hand take a
sequence of leaves $l_n$ converging to $l$. There are three
possibilities: either $\varphi_\lambda(l_n)$ converges to a spacelike
geodesic (in the Hausdorff topology of $\overline\mx_{-1}$), or it
converges to a segment (left or right) leaf of $\partial\mx_{-1}$, or
it converges to a point of $\partial\mx_{-1}$. In fact if
$p,q\in\partial\mh^2$ are the end-point of $l$, the limit of
$\varphi_\lambda(l_n)$ has end-points (in $PSL(2,\mr)\times
PSL(2,\mr)=\partial\mx_{-1}$)
$\lim_{n\rightarrow+\infty}(\beta_+(x_0,y_n)p, \beta_-(x_0, y_n)p)$
and $\lim_{n\rightarrow+\infty}(\beta_+(x_0,y_n)q,\beta_-(x_0, y_n)q)$
where $y_n$ is any point on $l_n$.  Roughly speaking the difference
among these cases depends on how fast the measure goes to infinity
along a geodesic segment joining a point in $\mathring H$ to a point
on the boundary. In~\cite{BSK} some computations in this sense are
given.}
\end{remarks}

\subsection{AdS bending cocycle on $\Uu^0_\lambda$}\label{bend-cocy-U0}

Let $\Uu=\Uu^0_\lambda$ be the flat Lorentzian spacetime 
corresponding to $\lambda$, as in Section \ref{ML_REGD}.
Just as in the hyperbolic case we want to ``pull-back'' the bending cocycle
$\beta_\lambda$ to a continuous bending cocycle
\[   \hat\beta_\lambda:\Uu \times\Uu \rightarrow PSL(2,\R)\times 
PSL(2,\R) \ .\]

By using Lemma~\ref{adesitt:bend:cont:lem} we can prove the analogous of
Proposition ~\ref{lift:cocycle}. The proof is similar 
so we omit the details.

\begin{prop}\label{adesitt:bend:ext:prop}
A determined construction produces a continuous cocycle
\[\hat\beta_\lambda:\Uu(1)\times \Uu(1)\rightarrow PSL(2,\R)
\times PSL(2,\R)\]
such that
\[    \hat\beta_\lambda(p,q)=\beta_\lambda(N(p),N(q))\]
for $p,q$ such that $N(p)$ and $N(q)$ do not lie on $L_W$.  Moreover,
the map $\hat\beta_\lambda$ is locally Lipschitzian.  For every
compact subset $K$ on $\Uu(1)$, the Lipschitz constant on $K$ depends
only on the diameter of $N(K)$ and the maximum of the total masses of
geodesic path in $H$ joining points in $N(K)$.
\end{prop}  \cvd

\noindent
Finally we can extend the cocycle $\hat\beta$ on the whole $\Uu$ by
requiring that it be constant along the integral geodesics of the
gradient of the cosmological time $T$. In particular, if we set
$r(1,p)=r(p)+N(q)\in\Uu(1)$, $\hat\beta$ satisfies
\[    \hat\beta(p,q)=\hat\beta(r(1,p),r(1,q)).\]

\begin{cor}\label{adesitt:bend:ext:cor}
The map
\[    \hat\beta_\lambda:\Uu\times \Uu\rightarrow PSL(2,\R)\times PSL(2,\R)\]
is locally Lipschitzian (with respect to the Euclidean distance on
$\Uu$).  Moreover, the Lipschitz constant on $K\times K$ depends only
on $N(K)$, the maximum of the total masses of geodesic paths of $H$
joining  points in $N(K)$ and the maximum and the minimum of $T$ on
$K$.  

If $\lambda_n\rightarrow\lambda$ on a $\eps$-neighbourhood
of a compact set $K$ of $\mathring H$, then $\hat \beta_{\lambda_n}$ converges
uniformly to $\hat \beta_\lambda$ on $U(H;a,b)$ (that is the set of points
in $\Uu$ sent by Gauss map on $H$ and with cosmological time in the
interval $[a,b]$).
\end{cor}\cvd

\section{Canonical AdS rescaling}\label{can-ads-resc}
In this section we define a map
\[\Delta_\lambda:\Uu^0_\lambda\rightarrow\mx_{-1}\]
such that the pull-back of the Anti de Sitter metric is a rescaling of
the flat metric directed by the gradient of the cosmological time,
with universal rescaling functions. We obtain in this way a globally
hyperbolic spacetime $\Pp_\lambda$.  A main difference with respect to
the Wick rotation map of Chapter \ref{HYPE} will be that the
developing map $\Delta_\lambda$ is always an isometric {\it embedding}
onto a convex domain of $\mx_{-1}$. By definition, the maximal
globally extension $\Uu^{-1}_\lambda$ of $\Pp_\lambda$ will be the
corresponding AdS $\Mm\Ll$-{\it spacetime} (as in Section
~\ref{ML:I}).
\smallskip

Recall that to construct the hyperbolic manifold $M_\lambda$ in
Chapter \ref{HYPE}, we have constructed the bending map $f_\lambda:
\mathring H\rightarrow\mh^3$, noticed that $f_\lambda$ is a locally convex
embedding in $\mh^3$, then $M_\lambda$ has been obtained by following
the normal flow, that is the flow on $\mh^3$ obtained by following the
geodesic rays normal to $f_\lambda(\mathring H)$ in the {\it non-convex} side bounded by $f_\lambda(\mathring H)$ (the flow on the convex side would produce
singularities). Eventually the developing map $D_\lambda$ has been
obtained by requiring that the integral lines of the cosmological
times would be sent to the integral lines of the normal flow.

In the same way $\varphi_\lambda:\mathring H\rightarrow\mx_{-1}$ is a locally
convex embedding (in fact an embedding), so the map $\Delta_\lambda$
can be constructed by requiring that the integral lines of the
cosmological time of $\Uu^0_\lambda$ are sent to the integral line of
the normal flow. An important difference with respect to the
hyperbolic case is that the normal flow is followed now in the {\it
convex} side bounded by $\varphi_\lambda(\mathring H)$ (otherwise singularities
would be reached).
 
For every $p\in\Uu^0_\lambda$, we define $x_-(p)$ as the dual point of
the plane $\hat\beta_\lambda(p_0,p)(P(Id))$, and
$x_+(p)=\hat\beta_\lambda(p_0,p)(N(p))$.  Thus let us choose
representatives $\hat x_-(p)$ and $\hat x_+(p)$ in $\SL{2}{R}$ such
that  the geodesic segment between $\hat x_-(p)$ and
$\hat x_+(p)$, is future directed.  Let us set
\begin{equation}\label{ads:dev:eq}
    \Delta_\lambda(p)=[\cos\tau(p)\hat x_-(p)+\sin\tau(p)\hat x_+(p)]
\end{equation}
where $\tau(p)= \arctan T(p)$ .

\begin{prop}\label{adesitt:rescaling:teo} 
The map
\[    \Delta_\lambda:\Uu^0_\lambda \rightarrow \mx_{-1}\]
is a local $\mathrm C^1$-diffeomorphism.  Moreover, the
pull-back of the Anti de Sitter metric is equal to the rescaling of
the flat Lorentzian metric, directed by the gradient of the
cosmological time $T$, with universal rescaling functions:
\begin{equation}\label{resc:ads:eq}
\begin{array}{ll}   
\alpha =\frac{1}{1+T^2} \ , \qquad &   
\beta=\frac{1}{(1+T^2)^2} \ .
\end{array}
\end{equation}
\end{prop}
\Dim 
To prove the theorem it is sufficient to analyze the map
$\Delta_\lambda$ in the case when $\lambda$ is a single weighted
geodesic. In fact, if we prove  the theorem in that
case, the same result will be proved when $\lambda$ is a finite
lamination.  The proof is completed in the general case 
by using standard approximations as in the final part of the proof of 
Theorem~\ref{hyperbolic:WR:teo}.
 
\begin{figure}
\begin{center}
\input{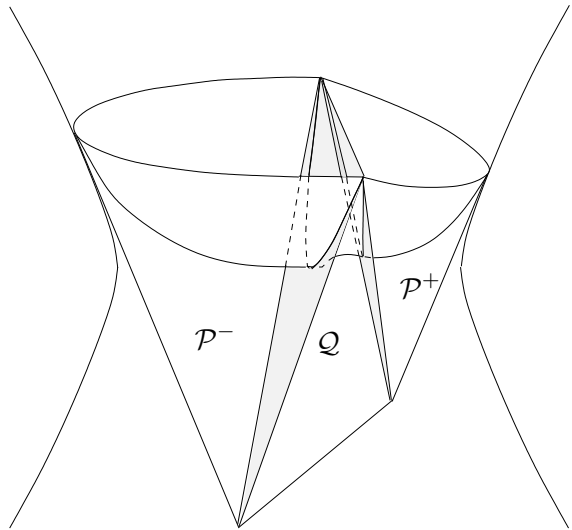}
\caption{{\small The domain $\Pp$ with its decomposition. Also the surface
    $\Pp(a)$ is shown.}}
\end{center}
\end{figure}

Let us set $\lambda_0=(l_0,a_0)$ and choose a base point
$p_0\in\mh^2-l_0$.  The surface $P=\varphi_\lambda(\mh^2)$ is simply
the union of two half-planes $P_-$and $P_+$ meeting each other along a
geodesic (that, with a little abuse of notation, is denoted by $l_0$).
We can suppose that $p_0$ is in $P_-$, and $l_0$ is oriented as the
boundary of $P_-$.  If $v_\pm$ denote the dual points of the planes
containing $P_\pm$ we have $v_-=Id$ and $v_+=\exp - a_0X_0$, $X_0$
being the standard generator of translations along $l_0$. By
Remark~\ref{mod:ort:rm} the vector $X_0$ is tangent to $P(id)$ along
$l_0$, orthogonal to it, and points towards $p_0$.

By definition , the image, say $\Pp$, of $\Delta_0=\Delta_{\lambda_0}$
is the union of three pieces: the cone with vertex at $v_-$ and basis
$P_-$, say $\Pp_-$, the cone with vertex at $v_+$ and basis $\Pp_+$,
and the join of the geodesic $l_0$ and the segment $[v_-,v_+]$, say
$\Qq$.

Fix a point in $l_0$, say $p_0$, and denote by $v_0$ the unit tangent
vector of $l_0$ at $p_0$ (that we will identify with a matrix in
$M(2,\mr)$).  Consider the coordinates on $\Uu_0$, say $(T,u,\zeta)$
introduced in Section~\ref{theWR}. With respect to these coordinates
we have
\[
 \Delta_0(T,u,\zeta)=\left\{\begin{array}{ll} \sin\tau\big
  (\ch\zeta(\ch u\ \hat p_0 + \sh u\ v_0) - \sh\zeta\ X_0 \big)\ +\
  \cos\tau\ \hat v_- & \textrm{ if } \zeta<0\\ \sin\tau( \ch u\ \hat
  p_0+\sh u\ v_0)\ +\ \cos\tau\exp(-\zeta \tan\tau\ X_0) & \textrm{ if
  } \zeta\in[0, a_0/T]\\ \sin\tau \big(\ch \zeta'(\ch u\ \hat p_0 +
  \sh u\ v_0) - \sh\zeta' X_0 \big)\ +\ \cos\tau\ \hat v_+& \textrm{
  otherwise}
\end{array}\right.
\]
where $\zeta'=\zeta-a_0/T$, $\tau=\arctan T$ and $\hat p_0, \hat v_+,
\hat v_-\in SL(2,\mr)$ are chosen as in~(\ref{ads:dev:eq}).

Clearly $\Delta_0$ is $\mathrm C^{\infty}$ for $\zeta\neq 0, a_0/T$. A
direct computation shows that the derivatives along the coordinate
fields glue on $\zeta=0$ and $\zeta= a_0/T$ and this proves that
$\Delta_0$ is $\mathrm C^1$.

By a direct computation we have
\begin{equation}\label{adsmetric}
 \Delta_0^*(\eta)=\left\{\begin{array}{ll}
  -\d\tau^2+\sin^2\tau(\d\zeta^2+\ch^2\zeta\d u^2)
 & {\rm if}\ \zeta<0\\
  -\d\tau^2+\sin^2\tau(\d\zeta^2+\d u^2) & {\rm if}\ \zeta\in[0, a_0/T]\\
 -\d\tau^2+\sin^2\tau(\d\zeta^2+\ch^2\zeta'\d u^2)
 & {\rm otherwise.}
\end{array}\right.
\end{equation}

Since $\d\tau^2=\frac{1}{(1+T^2)^2}$ and $\sin^2
\tau=\frac{T^2}{1+T^2}$, comparing (\ref{adsmetric}) with the
expression of the flat metric given in (\ref{flatmetric}) proves that
$\Delta_0$ is obtained by a rescaling along the gradient of $T$ with
rescaling functions given in (\ref{resc:ads:eq}).  \cvd

We denote by $\Pp_\lambda$ the domain $\Uu^0_\lambda$ endowed with the
Anti de Sitter metric induced by $\Delta_\lambda$. In other words, we
are considering $\Delta_\lambda$ as a developing map for such an AdS
structure; note that the smooth structure of $\Pp_\lambda$ is only
$\mathrm C^1$-diffeomorphic to the original one on $\Uu^0_\lambda$.

\begin{prop}\label{memgen:adsresc2:prop}
Let $T$ be the cosmological time of the flat regular domain
 $\ \Uu^0_\lambda$. Then $T$-level surfaces are Cauchy surfaces of
 $\Pp_\lambda$ (hence it is globally hyperbolic). Moreover
 the function
\[
    \tau=\arctan T\,.
\]
is the cosmological time of $\Pp_\lambda$. Every integral line of the
gradient of $\tau$ is a geodesic realizing the cosmological time.
\end{prop}
\Dim The causal-cone distribution on $\Pp_\lambda$ is contained in the
flat one of $\Uu$.  Thus, causal curves of $\Pp_\lambda$ are causal
also with respect to the flat metric. Since $\Uu(a)$ is a Cauchy
surface with respect to the flat metric, it is a Cauchy surface also
with respect to the Anti de Sitter metric.

Let $\gamma(t)$ denote a causal curve of $\Pp_\lambda$. Consider
the orthogonal decomposition of its tangent vector
\[
   \gamma(t)=\alpha(t)X(t)+h(t)
\]
where $X$ is the gradient of $T$.
The equalities
\[
\begin{array}{l}
  \E{\gamma(t)}{\gamma(t)}_{AdS}=
-\frac{\alpha(t)^2}{(1+T(t)^2)}+\frac{\E{h(t)}{h(t)}_F}{1+T(t)^2}\\
  \tau'(t)=\frac{\alpha(t)}{(1+T(t))}
\end{array}
\]
imply
\[
    \tau'(t)>\sqrt{-\E{\dot\gamma}{\dot\gamma}_{AdS}}
\]
so $\tau$ is greater or equal to the cosmological time.  On the other
hand the length of a integral line of the gradient of $\tau$ (that up
to re-parameterization coincides with the integral line of the (flat)
gradient of $T$) is just $\tau$. So $\tau$ is equal to the
cosmological time and the integral line of the gradient of $\tau$
realizes it.

The same computation shows that, up to re-parametrization, it is
the \emph{unique} curve realizing the cosmological time.  The fact
that it is a geodesic can be shown either by means of general facts
(see~\cite{A}) or by looking at the developing map.  \cvd

\section{Maximal globally hyperbolic AdS spacetimes}\label{MGHADS}
We introduce here the so called ``standard'' AdS spacetimes and we
study their geometry. Finally we will show that the class of simply
connected maximal globally hyperbolic AdS spacetimes that contain
complete Cauchy surfaces does coincide with the class of standard
ones. Proposition ~\ref{memgen:adscompl:prop} will be a key point
of our discussion.

\begin{prop}\label{memgen:adsstand:prop}
Let $Y$ be an Anti de Sitter {\rm simply connected} spacetime, and
$F\subset Y$ be a spacelike Cauchy surface. Suppose the induced
Riemannian metric on $F$ is complete. Then the developing map $
Y\rightarrow\mx_{-1}$ is an embedding onto a convex subset of
$\mx_{-1}$.

 The closure of $F$ in $\overline{\mx}_0$ is a closed
disk and its boundary $\partial F$ is a nowhere timelike curve of
$\partial\mx_{-1}$.

If $Y$ is the maximal globally hyperbolic Anti
de Sitter spacetime containing $F$ then $Y$ coincides with the Cauchy
development of $F$ in $\mx_{-1}$.  The curve $\partial F$ determines
$Y$, namely $p\in Y$ iff the dual plane $P(p)$ does not meet $\partial
F$.

Conversely $\partial F$ is determined by $Y$, in fact
$\partial F$ is the set of accumulation points of $Y$ on
$\partial\mx_{-1}$.  If $F'$ is another complete spacelike Cauchy
surface of $Y$ then $\partial F'=\partial F$.
\end{prop}
The proof of this proposition can be found in Section 7 of~\cite{M}.  
Notice however that by
``Cauchy development'' of a surface we mean the interior part of its
``domain of dependence'' in the sense of Mess.
\begin{remark}\label{ads:bend:inj:rem:2}
{\rm When $H=\mh^2$, then the first claim in Remarks ~\ref{ads:bend:inj:rem}
also follows from the above proposition.}
\end{remark}

Since the $T$-level surfaces in $\Uu^0_\lambda$ are complete, then
also the $\tau$-level surfaces in $\Pp_\lambda$ are complete, so  
the above proposition does apply. Hence we have:
\begin{cor}\label{AdS-risc-embed}
$\Delta_\lambda$ is a $\mathrm C^1$-homeomorphism of $\ \Uu^0_\lambda$, hence
an isometry of $\ \Pp_\lambda$, onto a open convex domain in $\mx_{-1}$.
\end{cor}
For simplicity we will confuse $\Pp_\lambda$ with its isometric image
in $\mx_{-1}$ via $\Delta_\lambda$. 

\subsection{Standard AdS spacetimes }\label{standard-ads} 
Given a nowhere timelike simple closed curve $C$ embedded in
$\partial\mx_{-1}$ its \emph{Cauchy development} is defined as
\[
   \Yy(C)=\{p| \partial P(p)\cap C=\varnothing\}
\]
When $C$ is different from a (left or right) leaf of the natural
 double foliation of $\partial\mx_{-1}$ (that it is a so called {\it
 admissible achronal curve}), then $\Yy(C)$ is also called a (simply
 connected) {\it standard AdS spacetime}, and $C$ is its {\it curve at
 infinity}.  We collect some easy facts about standard spacetimes.

\emph{1. }As $C$ is a closed nowhere timelike curve in the boundary of
 $\mx_{-1}$ different from a leaf of the double
 foliation of $\partial\mx_{-1}$, then it is homotopic to the meridian
 of $\partial\mx_{-1}$ with respect to $\mx_{-1}$.
\medskip\par\noindent 

\emph{2. } There exists a spacelike plane $P$ not intersecting it (see
Lemma 5 of~\cite{M}). In the Klein model we can cut $\Pm^3$ along the
projective plane $\hat P$ containing $P$ and we have that $\Yy(C)$ is
contained in $\mr^3=\Pm^3\setminus\hat P$. Since $C$ is nowhere
timelike then for every point $p\in C$ the plane $P(p)$ tangent to
$\partial\mx_{-1}$ at $p$ (that cuts $\mx_{-1}$ in a null totally
geodesic plane) does not separate $C$. It follows that the convex hull
$\Kk(C)$ of $C$ in $\mr^3$ is actually contained in $\mx_{-1}$.
\medskip\par\noindent 

\emph{3. } Support planes of $\Kk(C)$ are non-timelike and the closure
$\overline{\Yy(C)}$ of $\Yy(C)$ in $\mx_{-1}$ coincides with the set
of dual points of spacelike support planes of $\Kk(C)$ whereas the set
of points dual to null support planes of $\Kk(C)$ coincides with $C$.
\medskip\par\noindent

\emph{4. } $\overline{\Yy(C)}$ is convex and the closure of $\Yy(C)$ in
 $\overline\mx_{-1}$ is $\overline{\Yy(C)}\cup C$.
It follows that $\Kk(C)\subset\overline{\Yy(C)}$.
A point $p\in\partial\Kk(C)$ lies in $\Yy(C)$ if and only if it is touched
 only by spacelike support planes. 
\medskip\par\noindent 

We call $\Kk(C)$ the \emph{convex core} of $\Yy(C)$. 

\subsection{The boundary of the convex core}\label{boundary-K}
By general facts about convex sets, $\partial\Kk(C)\cup C$ is
homeomorphic to a sphere (in fact it is the boundary of a convex set
in $\mr^3$).  In particular $\partial\Kk(C)$ (that is the boundary of
$\Kk(C)$ in $\mx_{-1}$) is obtained by removing a circle from a
sphere, so it is the union of two disks. These components will be
called {\it the past and the future boundary} of $\Kk(C)$, and denoted
$\partial_-\Kk(C)$ and $\partial_+\Kk(C)$ respectively.  In fact given
any timelike ray contained in $\Kk(C)$, its future end-point lies on
the future boundary, and the past end-point lies on the past boundary.

By property \emph{4.} $\partial_+\Kk(C)\cap\Yy(C)$ is
obtained by removing from $\partial_+\Kk(C)$ the set of points that
admits a null support plane.  Now suppose that a null support plane
$P$ passes through $x\in\partial_+\Kk(C)$. The set $P\cap\Kk(C)$ turns
to be the convex hull of $\overline P\cap C$.  Now, $\overline
P\cap\partial\mx_{-1}$ is the union of the left and the right leaves
passing through the dual point $x(P)$. Since $P$ does not separate
$C$, $C$ is not a leaf of the double foliation of $\partial\mx_{-1}$,
$C$ is achronal, the only possibility is that $C\cap P$ is the union
of a segment on the left leaf through $x(P)$ and a segment on the
right leaf. Thus $P\cap\Kk(C)$ is a triangle with a vertices at
$x(p)$, two ideal edges (that are segments on the leaves of the double
foliation of $\partial\mx_{-1}$ and a complete geodesic of $\Kk(C)$
(that is a bending line).

It follows that the set $\partial_+\Kk(C)\cap\Yy(C)$ is obtained by
removing from $\partial_+\Kk(C)$ (at most) numerable many ideal
triangles, so  it is homeomorphic to a
disk.  Moreover the only case for $\partial_+\Kk(C)\cap\Yy(C)$ to be
empty is that it is formed by two null triangles, that is the case
when the curve $C$ is obtained by joining the end-points of a
spacelike geodesic $l$ with the end-points of its dual geodesic
$l^*$. Notice that in that case $\Yy(C)=\Kk(C)$, we call it the {\it
degenerate standard spacetime}, $\Pi_{-1}$. It could be regarded as
the analogous of the future of a line in Minkowski context.
We postpone the analysis of this case to  Chapter~\ref{QD}.

So, from now on, \emph{standard spacetimes are assumed to be not
degenerate}. Moreover, since we will be mainly interested in
$\partial_+\Kk(C)\cap\Yy(C)$, from now on we will denote
$\partial_+\Kk(C)$ that set by (with a bit abuse of notation).  When
there is no ambiguity on the curve $C$, we will denote only by $\Kk$
and $\Yy$ the convex core and the Cauchy development of $C$.

\begin{prop}\label{ads:boundary:prop} 
$\partial_+\Kk$ is locally $\mathrm C^0$-isometric to $\mh^2$. 
\end{prop}
\Dim 
Let $c$ be a small arc in
$\partial_+\Kk$ containing $x$ that intersects every bending line in at
most one point and every face in a sub-arc. Take a dense sequence $x_n$
on $c$ such that $x_1=x$, and for every $n$ choose a support plane
$P_n$ of $\partial_+K$ at $x_n$.  For every $n$ let $S_n$ be the
future boundary of the domain obtained intersecting the past of $P_i$
in $\mx_{-1}\setminus P(x)$ for $i=1\ldots n$. It is easy to see that
$S_n$ is a finitely bent surface.  In fact if $x_i, x_j, x_k$ is an
ordered set of points on $c$ then it is easy to see that $P_i\cap P_j$
and $P_j\cap P_k$ are geodesics of $P_j$ (if non-empty) that are
separated by the face (or bending line) $P_j\cap\partial_+ \Kk$.

Now for any $n$ we can choose an isometry $J_n:S_n\rightarrow\mh^2$
that is constant on $\partial_+\Kk\cap P_1$. Moreover the bending locus
on $S_n$ produces a measured geodesic lamination
$\lambda_n=(\Ll_n,\mu_n)$ on $\mh^2$. 

We claim that for every
geodesic arc $k$ the total mass $\alpha_n$ of $k$ with respect to
$\mu_n$ is a decreasing function.  The claim follows from Lemma~\ref{3planes}
that is the strictly analogous of Lemma 1.10.1 (Three planes)
of~\cite{Ep-M}.  We postpone the proof of the Lemma (and the claim) to
the end of this proof.

By the claim we see that $\lambda_n$ converges to a measured geodesic
lamination $\lambda$ on $\mh^2$. In particular, the bending map
$\varphi_{\lambda_n}$ converges to $\varphi_\lambda$.

Since $\varphi_{\lambda_n}$ is the inverse of $J_n$ we see that
$\varphi_{\lambda}(\mh^2)$ contains all the faces and the bending
lines of $\partial_+\Kk$ passing through $c$. The union of all this
strata is a neighbourhood $U$ of $\partial_+\Kk$, and we have proved
there exists an open set $V$ on $\mh^2$ such that
$\varphi_\lambda(V)=U$.  By Proposition~\ref{ads:bend:descr:prop} we
know that $\varphi_\lambda$ is an isometry.  \cvd

\begin{lem}\label{3planes}
Let $P_1,P_2,P_3$ be three spacelike planes in $\mx_{-1}$ without a
common point of intersection, such that any two intersect
transversely. Suppose that there exists a spacelike plane between
$P_2$ and $P_3$ that does not intersect $P_1$.  Then the sum of the
angles between $P_1$ and $P_2$ and $P_1$ and $P_3$ is less than the
angle between $P_2$ and $P_3$. \end{lem}

\emph{Proof of Lemma \ref{3planes}:} Denote by $x_i$ the dual point of
$P_i$. The hypothesis implies that segment between $x_i$ and $x_j$ is
spacelike, but the plane containing $x_1$, $x_2$ and $x_3$ is
non-spacelike.  The existence of a plane between $P_2$ and$P_3$ implies
that every non-spacelike  geodesic starting at $x_1$ meets the segment
$[x_2,x_3]$. Thus there exists a point $u\in[x_2,x_3]$ and a unit
timelike vector $v\in T_u\mx_{-1}$ orthogonal to $[x_2,x_3]$, such
that $x_1=\exp_u(tv)$. Choose a lift of $[x_2,x_3]$, say $[\hat
x_2,\hat x_3]$, on $SL(2,\mr)$ and denote by $\hat u$ the lift of
$u$ on that segment.  Denote by $l_1$ (resp. $l_2$, $l_3$) the length
of the segment $[x_2,x_3]$ (resp. $[x_1,x_3]$, $[x_1,x_2]$). We have
that\[\ch l_i=|\E{\hat x_j}{\hat x_k}|\] where $\{i,j,k\}=\{1,2,3\}$.
Now we have that $\hat x_1=\cos t\, \hat u+\sin t\, \hat v$ so that

$|\E{\hat x_1}{\hat x_i}|=\cos t |\E{\hat u}{\hat x_i}|$. 
Hence
\[ \begin{array}{ll} \ch l_2+<\ch l_2' & \ch l_3<\ch l_3' \end{array}\]
where $l_2'$ and $l_3'$ are the lengths of $[x_2,u]$ and $[x_3,u]$.
Finally we have $l_2+l_3<l_2'+l_3'=l_1$. 
\cvd

\begin{remark}\label{ads:bend:inj:rem:3}
\emph{
If $\partial_+ \Kk$ is complete then it is isometric to $\mh^2$ and the
bending lamination gives rise to a bending lamination of $\mh^2$, say
$\lambda$, such that $\partial_+ \Kk$ coincides with the image of
$\varphi_\lambda$.}

\emph{ In general \emph{$\partial_+ \Kk$ is not complete} even if the
curve $C$ does not contain any segment on a leaf ($C$ can be chosen to
be the graph of a homeomorphism of $S^1$ onto itself),
see~\ref{earthex} for an example.}

\emph{ We will show that $\partial_+ \Kk$ is isometric to a straight
convex set, $H$, of $\mh^2$ and the bending lamination on it gives
rise to a bending lamination on $H$, say $\lambda$, such that
$\partial_+ \Kk$ coincides with the image of $\varphi_\lambda$. The proof
is based on the rescaling of Theorem~\ref{adesitt:rescaling:teo} and
we are not able to prove it by a direct argument.}
\end{remark}

\subsection{The past part of a standard spacetime}\label{past-part}

The \emph{ past part } $\Pp = \Pp(C)$ of a standard AdS spacetime
$\Yy(C)$ is the past in $\Yy(C)$ of the future boundary $\partial_+
\Kk$ of its convex core.  The complement of $\partial_+ \Kk$ in the
frontier of $\Pp(C)$ in $\mx_{-1}$ is called the {\it past boundary}
of $\Yy(C)$, denoted by $\partial_- \Pp$.
\begin{prop}\label{ads:ct:prop}
Let $\Pp$ be the past part of some $\Yy(C)$.  Then $\Pp$ has
cosmological time $\tau$ and this takes values on $(0,\pi/2)$. For
every point $p\in \Pp$ there exist only one point $\rho_-(p)\in
\partial_- \Pp$, and only one point $\rho_+(p)\in \partial_+ \Kk$ such
that
\medskip\par\noindent
\emph{1. }$p$ is on the timelike segment joining $\rho_-(p)$ to $\rho_+(p)$.
\medskip\par\noindent
\emph{2. }$\tau(p)$ is equal to the length of the segment $[\rho_-(p),p]$.
\medskip\par\noindent
\emph{3. }the length of $[\rho_-(p),\rho_+(p)]$ is $\pi/2$. 
\medskip\par\noindent 
\emph{4. }$P(\rho_-(p))$ is a support plane for
$\Pp$ passing through $\rho_+(p)$ and $P(\rho_+(p))$ is a support
plane for $\Pp$ passing through $\rho_-(p)$.
\medskip\par\noindent \emph{5. }The map $p\mapsto \rho_-(p)$ is
continuous. The function $\tau$ is $\mathrm C^1$ and its gradient at
$p$ is the unit timelike tangent vector $\mathrm{grad}\,\tau(p)$ such
that
\[
    \exp_p\left(\tau(p)\mathrm{grad}\,\tau(p)\right)=\rho_-(p)\,.
\]
\end{prop}
\Dim For $p\in\mx_{-1}$ denote by $\Gg^{+}(p)$ (resp. $\Gg^-(p)$) 
the set
of points related to $p$ by a future-pointing (resp. past-pointing)
timelike-geodesic of length less than $\pi/2$.  Given $p\in\Pp$ it is
not hard to see that
\[
  \fut_{\Pp}(p)=\Gg^+(p)\cap\Pp\qquad\pass_{\Pp}(p)=\Gg^-(p)\cap\Pp\,.
\]
For every $q\in\fut_{\Pp}(p)$ the Lorentzian distance between $p$ and
$q$ in $\Pp$ is realized by the unique geodesic segment joining $p$ to
$q$ in $\Pp$.

Given $p\in\Pp$, as the dual plane $P(p)$ is disjoint from $\Kk=\Kk(C)$, it
is not hard to see that $\Gg^+(p)\cap\partial_+\Kk$ is a non-empty
pre-compact set. So there exists a point $\rho_+(p)$ on $\partial_+\Kk$
which maximizes the distance from $p$.  For each $a\in (0,\pi/2)$
consider the surface
\[
    \mh_p(a)=\{\exp_p(av)|v\textrm{ future-directed unitary vector in
    }T_p\mx_{-1}\}
\]
It is strictly convex in the future and the
tangent plane at $q\in\mh_p(a)$ is the plane orthogonal to the segment
$[p,q]$ contained in $\Gg^+(p)$ (these facts can be proved directly or
by means of the Lorentzian version of the Gauss Lemma in Riemannian
geometry).  Since the set of points in $\Gg^+(p)$ with assigned
Lorentzian distance from $p$ is a strictly convex in the future surface and
$\partial_+\Kk$ is convex in the past it follows that $\rho_+(p)$ is
unique and the plane passing through $\rho_+(p)$ orthogonal to the segment
$[p,\rho_+(p)]$ is a support plane for $\partial\Kk$.

The point $\rho_-(p)$, dual point of this plane,
is contained in the past boundary of $\Yy$. The dual plane of
$\rho_+(p)$ is a support plane of $\Yy$ passing through $\rho_-(p)$: in fact
for $q\in P(\rho_+(p))$ we have that $\rho_+(p)\in P(q)$, $\partial
P(q)\cap C\neq\varnothing$, so $q\notin\Yy$. It follows that the
cosmological time at $\rho_+(p)$ is exactly $\pi/2$. Thus 
the cosmological time of $p$ is the length of $[\rho_-(p),p]$.

If $p_n$ is a sequence converging to $p_\infty$ in $\Pp$ we have that
the sequence $(\rho_+(p_n))$ runs in a compact set of $\partial_+\Kk$:
indeed if we choose $q\in\pass_\Pp(p_\infty)$ then
$\rho_+(p_n)\in\Gg^+(q)\cap\partial\Kk$ that is a compact
set. Since the limit of any converging sub-sequence 
is $\rho_+(p_\infty)$, $\rho_+(p_n)$ converges to $\rho_+(p_\infty)$. Thus,
$\rho_+$ is continuous and so is $\rho_-$.

Finally given a point $p\in\Pp$, there exists a
neighbourhood $U$ of $p$ that is contained in $\Gg^+(\rho_-(q))$ and in
$\Yy(\partial P(\rho_+(p)))$. Denote by $\tau_1$ the Lorentzian distance
from $\rho_-(q)$ and by $\tau_2$ the Lorentzian distance from $P(\rho_+(p))$:
they are smooth functions defined on $U$.  Moreover we have
\[
   \begin{array}{l}
   \tau_1(q)\leq\tau(q)\leq\tau_2(q)\quad\textrm{ for all }q\in U\\
   \tau_1(p)=\tau(p)=\tau_2(p)\\
   \mathrm{grad}\,\tau_1(p)=\mathrm{grad}\,\tau_2(p)=v_0
  \end{array}
\]
where $v_0$ is the unit timelike vector at $p$ such that
 $exp_p \tau(p)v_0=\rho_-(p)$. It follows that $\tau$ is differentiable at
$p$ and $\nabla\tau(p)=v_0$.
\cvd

Summing up, given the past part $\Pp$ of a standard AdS spacetime
$\Yy(C)$, we have constructed
\medskip\par\noindent
The cosmological time $\tau:\Pp\rightarrow (0,\pi/2)$.
\medskip\par\noindent
The future retraction  $\rho_+:\Pp\rightarrow\partial_+\Kk$.
\medskip\par\noindent
The past retraction $\rho_-:\Pp\rightarrow\partial_-\Pp$.

\begin{cor}\label{ads:ct:cor}
\emph{1.} Given $r$ in the past boundary of $\Yy$, $\rho_-^{-1}(r)$ is
the set of points $p$ such that the ray starting from $r$ towards $p$
meets at time $\pi/2$ the future boundary of $\Kk$.
\medskip\par\noindent \emph{2.} The image of $\rho_-$ is the set of
points of $\partial_-\Pp$ whose dual plane meets $C$ at least in two
points.
\medskip\par\noindent
\emph{3.} The image of $\rho_+$ is the whole $\partial_+\Kk$.
\end{cor}
\Dim
Point 1. follows from points 4. and 5. of Proposition~\ref{ads:ct:prop}.
It, in turn, implies point 3.
Moreover the image of $\rho_-$ turns to be the set of points of $\partial_-\Pp$, whose dual plane is a support plane of $\Kk$ touching $\partial_+\Kk$. 
Given $p\in\partial_-\Pp$, its dual plane, $P(p)$, meets $C$ (otherwise $p$ would lie in $\Yy$).
On the other hand, since $p$ is limit of points in $\Yy$, 
$P(p)$ does not intersect the interior of $\Kk$. Thus, $P(p)$ is a support plane of $\Kk$ and $P(p)\cap\Kk$ is the convex hull of $P(p)\cap C$. 
Since $p$ lies on the past boundary, $P(p)$ does not intersect $\partial_-\Kk$. Thus, it contains points of $\partial_+\Kk$ iff $P(p)\cap C$ contains at least two points. 
\cvd

The image of the past retraction  is called the
\emph{initial singularity} of $\Yy(C)$ .

Proposition~\ref{memgen:adsstand:prop} ensures that a simply connected
maximal globally hyperbolic Anti de Sitter spacetime containing a
complete Cauchy surface is a standard spacetime.  We are going to show
that also the converse is true, that is every spacetime $\Yy(C)$
contains a complete Cauchy surface.

\begin{prop}\label{memgen:adscompl:prop}
If $\Pp$ is the past part of $\Yy(C)$
then every level surface $\Pp(a)$ of the cosmological
time is complete.
\end{prop}

\begin{remark}\emph{
The same result has been recently achieved by Barbot~\cite{Ba}(2)
with a different approach.}
\end{remark}
Since the proof of Proposition ~\ref{memgen:adscompl:prop} is quite
technical we prefer to give first the scheme. Let us fix
$p_0\in\Pp(a)$: we have to prove that the balls centered at $p_0$ are
compact.  Given a point $p\in\Pp(a)$ there exists a unique spacelike
geodesic in $\mx_{-1}$ joining $p_0$ to $p$.  Denote by $\xi(p)$ the
length of such a geodesic. We will prove the following facts
\medskip\par\noindent \emph{Step 1.} $\xi$ is proper and
$\xi(p)\rightarrow+\infty$ for $p\rightarrow\infty$;
\medskip\par\noindent \emph{Step 2.} If $c$ is a path in $\Pp(a)$
joining $p_0$ to $p$ then the length of $c$ is bigger than $M \xi(p)$
where $M$ is a constant depending only on $a$.
\medskip\par\noindent The proof of \emph{ Step 1.} is based on the
remark that the dual plane of $p_0$ is disjoint from the closure of
$\Pp(a)$ in $\overline\mx_{-1}$ so the direction of the geodesic
joining $p_0$ to $p$ cannot degenerate to a null direction.\\
\emph{Step 2.} is more difficult. We prove that the second fundamental
form of $\Pp(a)$ is uniformly bounded by the first fundamental
form. By using this fact we will be able to to conclude the proof.\\

Before proving the proposition let us just recall how the second
fundamental form is defined:

Given a spacelike surface $S$ in a
Lorentzian manifold $M$ denote by $N$ the future-pointing unit vector
on $S$. Then the second fundamental form on $S$ is a symmetric
bilinear form defined by
\[
     II(x,y)=\E{\nabla_x N}{y}
\]
where $\nabla$ is the Levi-Civita connection on $M$. It is symmetric, and
$\nabla_x N\in T_pS$ for every
$x\in T_pS$.

\begin{lem}
Let $II$ denote the second fundamental form on $\Pp(a)$ then we have
\[
     II(x,x)\leq\frac{1}{\tan a}\E{x}{x}
\]
for every $x\in T_p \Pp(a)$ and $p\in \Pp(a)$.
\end{lem}
\Dim Let us fix $p_0\in\Pp(a)$ and set $r_0=\rho_-(p_0)$. The set
$U:=\Gg^+(r_0)\cap\Pp(a)$ is a neighbourhood of $p_0$ in $\Pp(a)$ and for every
$p\in U$ there exists a unique timelike geodesic contained in
$\Gg^+(r_0)$ joining $r_0$ to $p$. Denote by $\sigma(p)$ the length of
such a geodesic. By definition we have that $0<\sigma(p)\leq\tau(p)=a$
and $\sigma(p_0)=a$ so $\sigma$ takes a maximum at $p_0$.  Equivalently the
function
\[
   h(p)=\cos\sigma(p)
\]   
takes a minimum at $p_0$ so $\grad
h(p_0)=0$ and the symmetric form
\[
    \omega:T_{p_0}\Pp(a)\times T_{p_0}\Pp(a)\ni (u,v)\mapsto
    \E{\nabla_u \grad h}{v}
\]
is positive semi-definite.
On the other hand by looking at the exponential map in $\mx_{-1}$,
it is not difficult to see that
\[
   h(p)=-\E{p}{r_0}
\]
where $\E{\cdot}{\cdot}$ is the form on $M(2,\mr)$ inducing the Anti
de Sitter metric on $\mx_{-1}$.  It follows that the gradient of $h$
on $\Pp(a)$ is given by
\[
    \grad h(p)=- r_0 - \E{r_0}{p} p - \E{r_0}{N}{N}
\]
so we have
\[
   \omega(u,u)= -\E{r_0}{p_0}\E{u}{v} - \E{r_0}{N(p_0)}II(u,v)
\]
Since $N(p_0)$ is the tangent vector at $p_0$
to the geodesic joining $r_0$ to $p_0$ we have
\[
    \begin{array}{l}
     p_0=\cos a\ r_0 +\sin a\ n_0\\
     N(p_0)=-\sin a\ r_0+\cos a\ n_0.
    \end{array}
\]
By using these equalities we get
\[
   \omega(u,u)=\cos a\E{u}{u} -\sin a\ II(u,u)\,.
\]
\cvd
\begin{remark}\emph{
The surface $\Pp(a)$ is $\mathrm C^{1,1}$ so its second fundamental
form is defined almost every-where. Anyway the inequality proved in
the Lemma holds in each point on which $II$ is defined and this will
be sufficient for our computation.  }
\end{remark} 
\emph{Proof of Proposition~\ref{memgen:adscompl:prop}} : First let us
prove step $1$. Suppose by contradiction that we can find a divergent
sequence $p_n\in \Pp(a)$ such that $\xi(p_n)$ is bounded by $A$.  Since
$p_n=\exp_{p_0}\xi(p_n)v_n$ with $|v_n|=1$ we obtain that $v_n$
diverges. So up to a subsequence the direction of $v_n$ tends to a
null direction. It follows that the geodesic $c_n$ joining $p_0$ to
$p_n$ converges to a null direction with end-point $p_\infty=\lim
p_n$. We have that $p_\infty\in\overline\Pp(a)\cap P(p_0)$ and this
is a contradiction because $\overline\Pp(a)=C\cup\Pp(a)$ and by
definition $C \cap P(p_0)=\varnothing$.\\

On $\Pp(a)$ the function $g(p)=-\ch\xi(p)=\E{p}{p_0}$ 
is $\mathrm C^{1,1}$, proper
and has a unique maximum at $p_0$.  It follows that if $c(t)$ is a
maximal integral line of $\grad g$ defined on the interval $(a,b)$
then
\[
    \lim_{t\rightarrow a} c(t)=p_0\,.
\]
Now we claim that there exists $K$ such that
\[
  \E{\grad g}{\grad g}< K (g^2-1)\,.
\]
Let us first show how the proof of \emph{Step 2.} follows from the claim. 
If $c(t)$ is any arc in $\Pp(a)$ joining $p_0$ to $p$ we have
\[
   \xi(p)=\int_c\frac{\E{ \grad g}{\dot x}}{\sqrt{g^2-1}}
\]
so from the claim we get
\[
  \xi(p)\leq\int_c K^{1/2}|\dot x|=K^{1/2}\ell(c)\,.
\]
Finally let us prove the claim.
By an explicit computation we have that
\[
   \grad g=  -(p_0+g(p)p+\E{N}{p_0}N)
\]
and
\[
   \E{\grad g}{\grad g}= g^2-1+\E{N}{p_0}^2\,.
\]
We see that it is sufficient to show that the function
\[
   f(p)=\E{N(p)}{p_0} 
\]
is less that $H(g(p)-1)$ for some $H>0$.  The function
$g(p)=-\ch\xi(p)$ is Lipschitz, proper and has a unique maximum at
$p_0$.  It follows that if $c:(t_-,t_+)\rightarrow\Pp(a)$ is a maximal
integral line of $\grad g$ passing through $p$ then
\[
    \lim_{t\rightarrow t_-} c(t)=p_0\,.
\] 
Now consider the integral line $c$ passing through $p$ and compare the
functions $f(t)=f(c(t))$ and $h(t)=g(c(t))-1$. We have that
\begin{equation}\label{memgen:lim:eq}
    \lim_{t\rightarrow t_-}f(t)=\lim_{t\rightarrow t_-}g(t)=0\,.
\end{equation}
On the other hand we have 
\[
   \begin{array}{l}
   \dot f =\E{\nabla_{\dot c}N}{p_0}=\E{\nabla_{\grad g} N}{\grad
   g}\leq\frac{1}{\tan a}\E{\grad
   g}{\grad g}\\
   \dot g =\E{ \grad g}{\dot c}=\E{\grad g}{\grad g}\,.
   \end{array}
\]
Since $\dot f(t)\leq\frac{1}{\tan a}\dot g(t)$
by~(\ref{memgen:lim:eq}) we can argue that
\[
   f(p)\leq\frac{1}{\tan a}(g(p)-1)\,.
\]
\cvd    

\begin{cor}\label{memgen:adscompl:cor}
 For every level surface $\Pp(a)$ of the past part $\Pp$ of a standard
AdS spacetime $\Yy(C)$:

(1) $\Pp(a)$ is a complete Cauchy surface of $\Yy(C)$ and this is
 the maximal globally hyperbolic AdS spacetime that extends $\Pp$;

(2) $\tau$ extends to the cosmological time of $\Yy(C)$, that takes
values on some interval $(0, a_0(C))$, for some well defined 
$\pi/2 < a_0(C) < \pi$.
\end{cor}
\Dim For (1) it is sufficient to show that every inextensible null ray
contained in $\Yy(C)$ intersects $\Pp(a)$. Let $l$ be a null ray passing
through $x\in\Yy(C)$ that does not intersect $\Pp(a)$. Since $\overline
\Pp(a)$ is a compression disk of $\overline\mx_{-1}$, either $l$
intersects $\Pp(a)$ or the dual point of $l$ lies on $C$. Since the
dual point of $l$ lies on the dual plane of $x$ the last possibility
cannot happen.

For (2), since there exists a plane that does not intersect $\Yy(C)$,
its cosmological function is a finite-valued function (actually it
takes values in $(0,\pi)$).  Moreover, by (1) every inextensible
causal curve intersects $\Pp$. So the cosmological function converges
to $0$ along the past side of any inextensible causal curve. Thus it
is the cosmological time (see Section \ref{CT}).  The
 value $\pi/2$ is taken on the future boundary of the convex core
of $\Yy(C)$ (that is also the future boundary of $\Pp$). It follows
that $\tau$ is $\mathrm C^{1,1}$ on $\Pp$ but not everywhere.  \cvd

The following corollary is a consequence of
Propositions~\ref{memgen:adsstand:prop}
and~\ref{memgen:adscompl:prop}.
\begin{cor}
The correspondence
\[
    C\mapsto\Yy(C)
\]
induces a bijection between the set of admissible achronal closed curves of
$\partial\mx_{-1}$ (up to the action of $PSL(2,\mr)\times PSL(2,\mr)$)
and the set of simply connected maximal globally hyperbolic Anti de
Sitter spacetimes containing a complete Cauchy surface (up to
isometries).
\end{cor}

Let us go back to the AdS $\Mm\Ll$-spacetimes $\Uu^{-1}_\lambda$.
We can now clarify the geometry of the embedding of $\Pp_\lambda$.
\begin{prop}\label{Plambda-past}
Let $\Uu=\Uu^0_\lambda$ be a flat regular domain, $\Pp_\lambda$ be the
Anti de Sitter spacetime obtained by the canonical rescaling of $\Uu$
(see Theorem~\ref{adesitt:rescaling:teo}), and $\Yy =
\Uu^{-1}_\lambda$ be the $\Mm\Ll$ AdS spacetime that extends
$\Pp_\lambda$. Then $\Pp_\lambda$ coincides with the past part $\Pp$
of $\Yy$.
\end{prop}
\Dim Notice that the developing map $\Delta$ of $\Pp_\lambda$ can be
written in the following form
\[ 
   \begin{array}{ll}
   \Delta:\Pp_\lambda\ni p\rightarrow [\cos\tau(p)\hat x_-(p)+\sin\tau(p)\hat
   x_+(p)]\in\mx_{-1}  & \textrm{where}\\
   \tau(p)=\arctan T(p)& \\
   \hat x_-(p)=\hat\beta(p_0,p)(Id)&\\
   \hat x_+(p)=\hat\beta(p_0,p)(N(p))&
   \end{array}
\]
(we have considered the standard identification $P(Id)=\mh^2$ described in
Chapter~\ref{MOD}).

Since $\hat\beta(p,q)=Id$ if and only if $r(p)=r(q)$, a cocycle
is induced on $\Sigma$ that  will be denoted  by $\hat\beta$ as well. In
particular the map $\Delta$ can be extended on $\Sigma$ by setting
$\Delta(s)=\hat\beta(r_0,s)(Id)$.

We claim that $\Delta(s)$ lies on the initial singularity of $\Yy$ for
every $s\in\Sigma$. To show that $\Delta(s)$ does not lie in
$\Yy$, it is sufficient to check that
$\Gg^+(\Delta(s))\cap\Delta(\Pp_\lambda(1))$ is not pre-compact. On
the other hand this set is bigger than
\begin{equation}\label{memgen:sing:eq}
\begin{array}{l}
   \Gg^+(\Delta(s))\,\cap\,\Delta(\,\Pp_\lambda(1)\cap r^{-1}(s)\,)=\\
    =\hat\beta(r_0,s)(\{[\frac{\sqrt{2}
   Id+\sqrt{2} x}{2}]|x\in\Ff(s)\})
\end{array}
\end{equation}
(where $\Ff(s)=N(r^{-1}(s))$ is the stratum of $H$ corresponding to
$s$, see Section \ref{RD:sec}) that is not pre-compact in $\mx_{-1}$.  Since
$\Delta(\Pp_\lambda)\subset\Yy$ the image of $\Sigma$ turns to be
contained in the past boundary of $\Yy$.

Denote by $C$ the curve at infinity of $\Yy$ (that coincides with the
set of accumulation points of the image of $\Delta$ on the boundary).
From equation~(\ref{memgen:sing:eq}) we may deduce that $C\cap\partial
P(\Delta(s))$ contains $\hat\beta(r_0,s)(\partial_\infty\Ff(s))$ so,
by Corollary~\ref{ads:ct:cor}, $\Delta(s)$ lies on the initial
singularity of $\Yy$ and $P(\Delta(s))\cap\partial_+\Kk$ contains
$\hat\beta(r_0,s)(\Ff(s))=\varphi_\lambda(\Ff(s))$ (where
$\varphi_\lambda:\mathring H\rightarrow\mx_{-1}$ is the bending map).

For $p\in\Pp_\lambda$, the integral line of the gradient of $T$ is
sent by $\Delta$ onto the geodesic segment $c$ joining $\Delta(r(p))$
to $\varphi_\lambda(N(p))$. From Corollary~\ref{ads:ct:cor} the
cosmological time of $\Delta(p)$ (as point in $\Yy$) coincides with
$\tau$ and we have
\[
\begin{array}{l}
  p\in\Pp\\
  \rho_-(\Delta(p))=\Delta(r(p))\\
  \rho_+(\Delta(p))= \varphi_\lambda(N(p))\,.
\end{array}
\]
 
So for every $a\in(0,\pi/2)$ the developing map induces a local isometry
\[
  \Delta_a:\Pp_\lambda(a)\rightarrow\Pp(a)\,.
\]
Since $\Pp_\lambda(a)$ is complete $\Delta_a$ is an isometry (in
particular injective and surjective). Thus $\Delta$ is an isometry on
$\Pp$.  \cvd

\begin{cor}
Let $\lambda=(H,\Ll,\mu)\in \Mm\Ll$. Then the bending map
$\varphi_\lambda:\mathring H\rightarrow \mx_{-1}$ is an isometry onto
the future boundary  of the convex core of 
$\Yy=\Uu^{-1}_\lambda$.
\end{cor}
\cvd

\section{Classification via AdS rescaling}\label{resc-class}
We have a map that associates to every flat regular
domain (hence to every $\lambda=(H,\Ll,\mu)\in \Mm\Ll$) a simply connected
maximal globally hyperbolic AdS spacetime (that is a standard one):
$$ \lambda \leftrightarrow \Uu^0_\lambda \to \Yy_\lambda =
\Uu^{-1}_\lambda \ .$$ We are
going to show that such a correspondence is bijective. In particular
we will show that given the past part $\Pp$ of any standard Anti de
Sitter spacetime $\Yy$, the rescaling along the gradient of the
cosmological time $\tau$ of $\Pp$ with rescaling functions
\begin{equation}\label{memgen:adsfact:eq}
    \hat\alpha=\frac{1}{\cos^2\tau} \qquad\qquad\hat\beta=\frac{1}{\cos^4\tau}
\end{equation}
produces a regular domain and this makes the inverse of the previous one.
\smallskip

In fact, it is immediate that such a rescaling performed on 
$\Pp_\lambda$, actually recovers the original regular domain
$\Uu^0_\lambda$. In other words, we see that the above map is {\it
injective}. Moreover, the future boundary of the convex core of the
standard spacetime $\Yy_\lambda$ is isometric to the image of the
Gauss map of $\Uu^0_\lambda$ via an isometry that sends the bending
locus to the measured geodesic lamination $\lambda$.
\smallskip

Conversely, by inverting the construction, we have that the past part
$\Pp$ of any standard spacetime $\Yy$ is obtained by the canonical
rescaling on a regular domain iff the future boundary of its convex
core is isometric to a straight convex set pleated along a measured
geodesic lamination. Hence, in order to prove that our favourite map
is also {\it surjective} we have to show that this fact always happens.
\smallskip

\noindent{\bf The $\Mm\Ll(\mh^2)$ case.}
This is easy to achieve in the particular case such that the future
boundary is {\it complete}.  In fact, the following Proposition is a
consequence of Remarks ~\ref{ads:bend:inj:rem},
~\ref{ads:bend:inj:rem:2} and ~\ref{ads:bend:inj:rem:3} (and all the
already established facts). This also establishes the characterization
of AdS $\Mm\Ll(\mh^2)$-{\it spacetimes} given in Proposition
~\ref{-1mlh2}.
\begin{prop}\label{mh2-case} A standard AdS spacetime $\Yy$ 
is obtained by the canonical rescaling of a regular domain with {\rm
surjective Gauss map} if and only if the future boundary of the convex
core of $\Yy$ is complete.
\end{prop}
Before going on, we will check that the class of standard AdS
spacetimes verifying this property is large.

\begin{prop}\label{memgen:cpt:prop}
Let $C$ be a no-where timelike curve in the boundary of Anti de Sitter
space and $\Kk$ denote its convex hull in $\mx_{-1}$. Suppose the set
of spacelike support planes touching $\partial_+\Kk$ to be compact, then
$\partial_+\Kk$ is isometric to $\mh^2$.
\end{prop}
In order to prove this proposition we need the following technical lemma:
\begin{lem}\label{memgen:cpt:lem}
There exists a timelike vector field $X$ on $\mx_{-1}$ such that
\medskip\par\noindent
\emph{1.} It extends to a timelike vector field on $\overline\mx_{-1}$.
\smallskip\par\noindent
\emph{2.}  The metric $g_X$ obtained by  Wick rotation along $X$ is complete.
\end{lem}
\Dim
Consider the covering
\[
  \mr^3\ni(x,y,\lambda)\mapsto\left(\begin{array}{ll}
                                    \xi\sin\lambda+x &
                                    \,y\,+\,\xi\cos\lambda\\
                                    y\,-\,\xi\cos\lambda &
                                    \,\xi\sin\lambda\,-\,x
                                    \end{array}\right)\in PSL(2,\mr)
\]
where we have set $\xi=\sqrt{1+x^2+y^2}$.  
Let us put
$X=\frac{\partial\,}{\partial\lambda}$.
It is not difficult to see that
\[
   X(A)=\frac{AX_0+X_0A}{2}
\]
where 
\begin{equation}\label{memgen:field:eq}
   X_0=\left(\begin{array}{ll} 0 & -1 \\1 & 0\end{array}\right)\ .
\end{equation}
In particular $X(A)$ extends to a timelike vector field on the whole
$\overline\mx_{-1}$.\par By an explicit computation it turns out that
the Anti de Sitter metric takes the form
\[
  -\xi^2\d\lambda^2+\frac{(1+y^2)\d x^2+ xy\d x\d y+ (1+x^2)d y^2}{\xi^2}\,.
\]
The horizontal metric is independent of $\lambda$; moreover, it is a
complete hyperbolic metric.  This can be shown either by noticing that
the surface $\{\lambda=\theta_0\}$ is the dual plane to the point with
coordinates $x=0,y=0,\lambda=\pi/2-\theta_0$ or by noticing that with
respect to the parametrization of $\mh^2$ given by
\[
   \mr^2\ni(x,y)\mapsto (x,y,\xi(x,y))\in\mx_0
\]
the hyperbolic metric takes exactly the form written above.\\

Eventually $g_X$ take the form
\[
    \xi^2(x,y)\d^2\lambda+ g_\mh(x,y)\,.
\]
Now if $c(t)=(x(t),y(t),\lambda(t))$ is an arc-length path defined on $[0,a)$
we have to show that $c(t)$ can be extended. 
Since $g_\mh$ is complete there exist
\[
   x(a)=\lim_{t\rightarrow a}x(t)\qquad y(a)=\lim_{t\rightarrow a}y(t)\,.
\]
On the other hand since $|\dot\lambda|<1$ there exists
\[
   \lambda(a)=\lim_{t\rightarrow a}\lambda(t)\,.
\]
\cvd
\emph{Proof of Proposition~\ref{memgen:cpt:prop}}
Let $X$ be the field given by Lemma~\ref{memgen:cpt:lem} and for every $v\in
T\mx_{-1}$ denote by $v_V$ the projection of $v$ along $X$ and $v_H$ the
projection on $X^\perp$. If $v$ is spacelike let us define
\[
   \eta(v)=\frac{|v_V|}{|v_H|}<1\,.
\]
Now we claim that given a compact set of spacelike planes, say $\Qq$, 
there exists a constant $M<1$ (depending only on $\Qq$ and $X$) such that
\[
   \eta(v)\leq M
\]
for every $v\in TP$ for $P\in\Qq$.\\

We can quickly conclude the proof of the proposition from the claim.
Indeed let $\Qq$ denote the family of support planes for
$\partial_+\Kk$ and $M$ be the constant given by the claim.  Denote by
$g_X$ the Riemannian metric obtained by the Wick Rotation along $X$.
If $c$ is a rectifiable arc contained in $\partial_+\Kk$ we have that
$\E{\dot c}{\dot c}> (1-M^2)/(1+M^2) g_X(\dot c,\dot c)$. This
inequality implies that the pseudo-distance on $\partial_+\Kk$ is
a complete distance.

So in order to conclude it is sufficient to prove the claim.
First let us fix a plane $P$: we are going to prove that
\[
    M(P)=\sup\{\eta(v)|v\in TP\}
\]
is less than $1$.
Since $\eta(\lambda v)=\eta(v)$ we can suppose $\E{v}{v}=1$. If we set
$\xi(v)=|v_V|$ we have $\eta(v)=\xi(v)/\sqrt{1+\xi(v)^2}$. So we have to show
that the there exists $K=K(P)$ such that
\[
  \xi(v)\leq K
\]
for every $v\in TP$ with $|v|=1$.\par
Consider the following function
\[
  H:(M(2\times 2,\mr)\setminus\{0\})\times \sG\lG(2,\mr)\ni (A,Y)\mapsto
  \frac{| \E{AY}{AX_0+X_0A}|}{-\E{AX_0+X_0A}{AX_0+X_0A}^{1/2}}\in\mr_{\geq 0}
\]
where $X_0$ is as in~\ref{memgen:field:eq}.
The following are easy remarks:
\medskip\par\noindent
- It is well-defined because $\E{AX_0+X_0A}{AX_0+X_0A}$ is equal to
$-(a^2+b^2+c^2+d^2)$;
\medskip\par\noindent
- it is homogenous in $A$ so it induces a map $\hat
H:\Pm^3\times\sG\lG(2,\mr)\rightarrow\mr_{\geq 0}$;
\medskip\par\noindent
- if $A\in P$ and $X\in\sG\lG(2,\mr)$ then $\xi(AX)=\hat H(A, X)$.
\medskip\par\noindent
Now it is not hard to see that there exists a spacelike subspace $T_0$ of
$\sG\lG(2,\mr)$ such that
\[
    T_A P=A T_0
\]
for every $A\in P$ so we can set
\[
    K=\sup\{\hat H(A,X)|A\in\overline P, X\in UT_0\}<+\infty\,.
\]
where $UT_0$ is the set of unit vectors in $T_0$.

Since $M(P)$ continuously depends on
$P$ the claim is proved.
\cvd 

\noindent{\bf General case.}  Let us turn now to a general standard
spacetime $\Yy$. We will deal with it by the following steps:
\smallskip

(i) we first prove that the rescaling on $\Pp$ directed by the
gradient of $\tau$ with rescaling functions given in
(\ref{memgen:adsfact:eq}) produce a flat spacetime; 

(ii) then we show the so obtained spacetimes is a maximal globally
hyperbolic flat spacetime with cosmological time. Hence, it is a
regular domain thanks to the results of Chapter ~\ref{FGHST}.
\smallskip

Given two standard AdS spacetimes $\Yy$ and $\Yy'$, with past
parts $\Pp$ and $\Pp'$, we say that $\Yy,\Yy'$ are
{\it ``locally equivalent around $p\in \Pp$ and $p'\in \Pp'$''}
if there exists a neighbourhood $U$ of $p$ in $\Pp$ and an
isometric embedding of $U$ onto a neighbourhood $U'$ of $p'$ in $\Pp'$ 
that preserves the cosmological time.

Since the first step (i) is a local check, thanks to Proposition
\ref{mh2-case}, in order to get it in general it is enough
to show the following proposition.

\begin{prop}\label{memgen:equiv:prop}
Let $\Yy$ be a simply connected Anti de Sitter standard spacetime and
$p$ be in the past part $\Pp$ of $\Yy$.  Then there exists a standard
Anti de Sitter spacetime $\Yy'$ and $p'$ in its past part $\Pp'$, such
that such that $\Yy$ and $\Yy'$ are equivalent around $p$ and $p'$,
and the future boundary of the convex core of $\Yy'$ is complete.
\end{prop}
\Dim We can distinguish two cases: either there exists a null support
plane passing through $\rho_+(p)$ or there exists a neighbourhood $U$
of $\rho_+(p)$ such that support planes touching $U$ are spacelike and
form a compact set $\Qq(U)$.

First consider the latter case.
Consider a small compact arc 
$k\subset\Uu$ passing through $p$ and intersecting transversally every
bending line of $\partial_+\Kk$ at most once. Denote by $\Qq$ the set of
support planes touching $k$ and define 
$$\Cc=\bigcap_{Q\in\Qq} \pass_{\mx_{-1}\setminus P(\rho_+(p))}(Q).$$
It is not hard to see that the boundary of $\Cc$ in $\mx_{-1}\setminus
P(\rho_+(p))$  is a connected achronal surface $\partial_+\Cc$
satisfying the following properties.
\medskip\par\noindent \emph{1.} It does not contain vertex (which
means that coincides with the future boundary of the convex hull of
$C_\infty:
=\overline{\partial_+\Cc}^{\overline\mx_{-1}}\setminus\partial_+\Cc$).
\medskip\par\noindent
\emph{2.} The support planes of $\Cc$ touching $\partial_+\Cc$ are in
$\Qq$.
\medskip\par\noindent
\emph{3.} Every face or bending line intersecting $k$ is contained in
$\partial_+\Cc$. In particular there exists a neighbourhood $U'$ of
$\rho_+(p)$ in $\partial_+\Kk$ contained in $\partial_+\Cc$.
\medskip\par
From Proposition~\ref{memgen:cpt:lem} we have that $\partial_+\Cc$ is complete
so $C_\infty$ is a closed no-where timelike curve. Denote by $\Pp'$ the past
part of the Cauchy development of $C_\infty$. The future boundary of $\Pp'$ is
$\partial_+\Cc$. If we denote by $\rho_+'$ the future retraction on $\Pp'$ we
easily see that $\rho_+^{-1}(U')=(\rho_+')^{-1}(U')=V$ and $\rho_+=\rho_+'$ on
$V$. It follows that $\Pp$ and $\Pp'$ are locally equivalent around $p$.\\
 
Suppose now that only one null support plane $P_0$ passes through
$\rho_+(p)$. Denote by $P_1$ the spacelike plane orthogonal to the segment
$[p,\rho_+(p)]$ at $\rho_+(p)$. Since $\partial_+\Kk$ does not have vertices we
have that $l=P_0\cap P_1$ is contained in $\partial_+\Kk$.
On the other hand, the point dual to $P_0$ is contained in $C$ too. 
Then $\partial_+K$ contains a \emph{null triangle} bounded by $l$. 
Take a spacelike plane $Q$ between
$P_0$ and $P_1$ and consider the surface $\Ss$ obtained by replacing
$P_0\cap\partial_+K$ with the half-plane bounded by $l$ in
$Q$. $\partial\Ss$ is a no-where timelike curve in the boundary and
consider its domain $\Uu'$. It follows that $\Ss$ is the future boundary of
its convex hull.
So  $p$ is contained in $\Uu'$ and $(p,\Uu)$ and $(p,\Uu')$ are locally
equivalent. Moreover, no null support plane passes through $\rho_+'(p)$.
 
Finally suppose that two null support planes $P_1,P_2$ pass through
$\rho_+(p)$ then we see that $\partial_+\Kk$ is the union of two null
triangles respectively lying on $P_1$ and $P_2$ and bounded by
$P_1\cap P_2$. So $\Pp$ is $\Pi_{-1}$ (see Chapter~\ref{QD}).  \cvd

So we know now that the rescaling of the past part $\Pp$ of any
standard AdS spacetime $\Yy$, directed by the gradient of its
cosmological time with rescaling functions~(\ref{memgen:adsfact:eq})
yields a flat spacetime. It remains to prove that this is in fact a
regular domain.  The key point to prove this fact is the completeness
of the level surfaces $\Pp(a)$ of the cosmological time.

\begin{teo}\label{ads:flat:teo} For every standard AdS domain $\Yy=\Yy(C)$,
the rescaling of its past part $\Pp$, directed by the gradient of the
cosmological time $\tau$, with universal rescaling functions
\begin{equation}\label{ads:flat:resc:eq}
   \alpha=\frac{1}{\cos^2\tau} \qquad \beta=\frac{1}{\cos^4\tau}
\end{equation}
produces a regular domain, whose cosmological time is given by the formula
\[
    T=\tan\tau\,.
\] 
\end{teo} 

\Dim Denote by $\Pp^0$ the flat spacetime produced by such a rescaling
of $\Pp$ \par By a computation like that one in
Proposition~\ref{desitter:ct:prop} we can prove that 
the cosmological time of $\Pp^0$ is
given by
\[
   T(p)=\tan\tau(p)\,.
\]

and the level surfaces $\Pp^0(a)$ are Cauchy surfaces.  Consider the
developing map
\[
   D:\Pp^0\rightarrow\mx_{0}\,.
\]
Since $\Pp^0(a)=\frac{1}{\cos^2(\arctan a)}\Pp(\arctan a)$, Proposition~\ref{memgen:adscompl:prop} implies that $\Pp^0(a)$ is a complete surface. Then  $D$
is an embedding.  Denote by $\Uu$ the domain of dependence of
$\Pp^0(1)$ in $\mx_{0}$. 

For every $p$ in the initial singularity of $\Pp$ let us consider
the set $U=\rho_-^{-1}(p)$. Up isometries of $\mx_{0}$ we have that
$D(U\cap\Pp(\pi/4))$ is a straight convex set in $\mh^2$ (with respect
to the standard embedding of $\mh^2$ into $\mx_{-1}$).  So
$0\notin\Uu$ (indeed $\fut(0)\cap\Pp^0(\pi/4)$ is not bounded) and
there exist at least $2$ null rays starting from $0$ disjoint from
$\Pp^0(1)$.

\begin{figure}
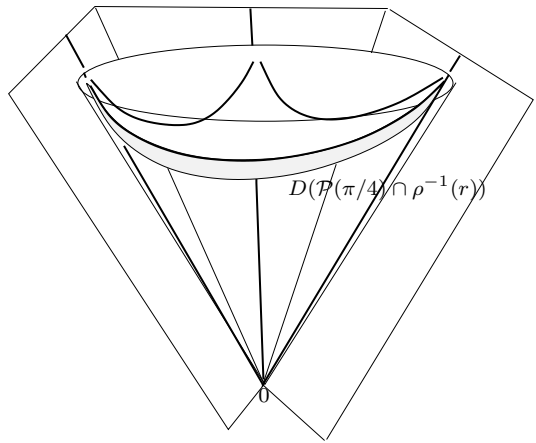

\begin{center}
\input GR2005_fig_adstofl.pstex_t
\caption{{\small The image through $D$ of
$U\cap\Pp(\pi/4)$.}}\label{adstofl:fig}
\end{center}
\end{figure}

It follows that $\Uu$ is a regular domain containing
$\Pp^0$.  Moreover, the image of $D(U)$ via the retraction of
$\Uu$ is $0$. Thus $T$
coincides with the restriction of the cosmological time $t$ on $\Uu$.
The inclusion $\Pp^0\rightarrow\Uu$ gives rise to a locally
isometry
\[
  D_a:\Pp^0(a)\rightarrow\Uu(a)
\]
since $\Pp^0(a)$ is complete $D_a$ is an isometry. So the map
\[
   D:\Pp\rightarrow\Uu
\]
is surjective.  \cvd
 
\section{Equivariant rescaling}\label{equiv-resc}
Let $Y$ be a maximal globally hyperbolic flat spacetime containing a
complete Cauchy surface and equipped with cosmological time. By
Theorem~\ref{FULLFLAT}, there exists a discrete group
$\Gamma<PSL(2,\mr)$ and an invariant measured geodesic lamination
$\lambda$ defined on some straight convex set $H$ such that
$Y=Y(\lambda, \Gamma)$. By Theorem~\ref{adesitt:rescaling:teo}, the
space obtained by a rescaling along the gradient of the cosmological
time of $Y$ with rescaling functions given in~(\ref{resc:ads:eq}) is a
maximal globally hyperbolic AdS spacetime, that is denoted
by $Y^{-1}=Y^{-1}(\lambda, \Gamma)$.

As in Section~\ref{hyp:eq:sec} we can easily compute the holonomy
$h^{(-1)}:\Gamma\rightarrow PSL(2,\mr)\times PSL(2,\mr)$ of $Y^{-1}$.
Indeed if $x_0\in\mathring H$ is a base point then for every
$\gamma\in\Gamma$ then
\[
    h^{(-1)}(\gamma)=\beta_\lambda(x_0,\gamma x_0)\circ(\gamma,\gamma)\,.
\]

Conversely if $U$ is a maximal globally hyperbolic Anti de Sitter
spacetime containing a complete Cauchy surface, then by
Proposition~\ref{memgen:adscompl:prop} its universal covering is a
standard spacetime.

Thus $U$ has cosmological time, $\tau$, taking values on $(0,a)$ with
$\pi/2< a\leq \pi$. Moreover $\tau$ is $\mathrm C^1$ on
$P=\tau^{-1}(0,\pi/2)$, and by Theorem~\ref{ads:flat:teo} the space
obtained by rescaling $P$ along the gradient of $\tau$ with functions
given in~\ref{ads:flat:resc:eq} is a maximal globally hyperbolic flat
spacetime.
\smallskip

Theorem~\ref{FULL_CLASS} is now completely proved.

\paragraph{Cocompact case}

Mess \cite{M} proved that the holonomy of any globally hyperbolic Anti
de Sitter spacetime containing a closed surface $\Sigma$ of genus
$g\geq 2$ is given by a pair of Fuchsian representations of
$\pi_1(\Sigma)$.

Conversely, given a pair of Fuchsian representations of
$\pi_1(\Sigma)$, say $(h_-,h_+)$, there exists an orientation
preserving homeomorphism of $\overline{\mh}^2 = \mh^2 \cup S^1_\infty$
which conjugates the action of $h_-$ on $\overline{\mh}^2$ with the
one of $h_+$.  In fact its restriction, $u$, to $S^1_\infty$ is
determined by $(h_-, h_+)$.  The graph of $u$ is a curve, say $C$, of
$S^1_\infty\times S^1_\infty=\partial\mx_{-1}$ that is $(h_-,h_+)$
invariant.  In fact it can be easily shown that it is the unique
$(h_-,h_+)$-invariant curve.  Since it is the graph of an
\emph{orientation preserving} homeomorphism it is no-where
timelike. Its Cauchy development $\Yy$ is invariant for $(h_-,h_+)$
and the action of $\Gamma$ on it is free and properly discontinuous
and the quotient $\Yy/(h_-,h_+)$ is a globally hyperbolic spacetimes
$\cong\Sigma\times\mr$ (all these results are discussed in~\cite{M}).

So there are two natural parameterizations of the set of maximal
globally hyperbolic AdS spacetimes with closed Cauchy surface of genus
$g\geq 2$.

The first one by looking at the future boundary of the convex core that
we have discussed in previous Sections. In this case the parameter
space is $\Tt_g\times\Mm\Ll_g$.

The second one by considering the holonomy. In this case the above
remarks show that the parameter space is $\Tt_g\times\Tt_g$.

The induced map $\Tt_g\times\Mm\Ll_g\rightarrow\Tt_g\times\Tt_g$ can
be explicitly described in terms of earthquakes as we are going to
explain in the next section, in a more general framework.

In~\cite{Ba}(2,3) Barbot has studied the holonomies of AdS spacetimes
containing a non-compact spacelike surface.  A generalization of these
results of Mess has been achieved in those papers.
\begin{prop}\cite{Ba}(2,3)\label{discrete-rep}
Given a standard Anti de Sitter spacetime $Y$ let
$h=(h_-,h_+):\pi_1(Y)\rightarrow PSL(2,\mr)\times PSL(2,\mr)$ be the
holonomy.  If $\tilde Y$ is different from $\Pi_{-1}$, then $h_-$ and
$h_+$ are discrete representations such that
$\mh^2/h_-\cong\mh^2/h_+$. Moreover, $Y\cong\mh^2/h_-\times\mr$.

Conversely given a pair of discrete representations $h=(h_-,h_+)$ of
the fundamental group of a surface $F$, such that
$\mh^2/h_-\cong\mh^2/h_+\cong F$, then there exists a standard
spacetime $Y\cong F\times\mr$ whose holonomy is $h$.
\end{prop}

However, if $F$ is not compact, it is not true that globally
hyperbolic AdS structures on $F\times\mr$ are determined by the
holonomy. A counterexample will be given in Section~\ref{3cusp}.

\section{AdS rescaling and generalized earthquakes}\label{ge:quake}
Given a measured geodesic lamination $\lambda$ on a straight convex set 
$H$ the \emph{generalized left earthquake} along $\lambda$ is the map
\[
   \Ee_L:\mathring H\ni x\mapsto\beta_+ (x_0,x)x\in\mh^2\,.
\]
where $\beta_+$ is the Epstein-Marden cocycle corresponding to
$\lambda$ (see Section~\ref{ads:bend:coc:sec}).  The generalized right
earthquake is defined by replacing $\beta_+$ by $\beta_-$
(that is the cocycle corresponding to the negative-valued measure
$-\lambda$).

If $H=\mh^2$ and $\Ee_L$ is surjective, then $\Ee_L$ is a ``true'' left
earthquake, according to the definition given by Thurston
in~\cite{Thu2}. In that case $\Ee_L$ extends in a natural way to a map
$\overline\mh^2\rightarrow\overline\mh^2$ and the restriction to
$\partial\mh^2$ is a homeomorphism.  Conversely, any homeomorphism of
$\partial\mh^2$ (up to post-composition by elements of $PSL(2,\mr)$)
is the trace on the boundary of a unique left earthquake of $\mh^2$
\cite{Thu2}.

The following interesting relation between earthquakes and Anti de
Sitter geometry was pointed out in~\cite{M}.
\begin{prop}\label{earth:class:prop}\cite{M}
Let $Y$ be an Anti de Sitter spacetime with compact Cauchy surface of
genus $g\geq 2$ and denote by $h=(h_-,h_+):\pi_1(Y)\rightarrow
PSL(2,\mr)\times PSL(2,\mr)$. If $F$ denotes the future boundary of
the convex core of $Y$ and $\lambda$ is its bending lamination then
$h_+$ (resp. $h_-$) is the holonomy of the hyperbolic surface obtained
by a left earthquake (resp. right earthquake) on $F$ along $\lambda$.
\end{prop}

\noindent We stress that Proposition~\ref{earth:class:prop} actually
gives a new ``AdS'' proof of the ``classical'' Earthquake Theorem in
the cocompact case. For, given $F,F'$ two hyperbolic structures on a
compact surface $\Sigma$, there exists a unique spacetime $Y$ whose
holonomy is $(h,h')$ where $h$ is the holonomy of $F$ and $h'$ is the
holonomy of $F'$. If $\lambda$ is the bending lamination of the
future boundary of the convex core of $Y$, then the earthquake along
$2\lambda$ transforms $F$ into $F'$.
\smallskip

On the other hand, there is in ~\cite{Thu2} a formulation of the
Earthquake Theorem that strictly generalizes the cocompact case. In
this section we study the relations between generalized earthquakes
defined on straight convex sets of $\mh^2$ and standard Anti de Sitter
spacetimes. As a corollary, we will point out an ``AdS'' proof of such
a general formulation. 

\begin{prop}\label{earth:gen:prop}
The map $\Ee_L$ is injective and the image is a straight convex set.
Moreover, the image of the lamination $\lambda$ is a lamination
$\lambda'$ on $\Ee_L(H)$.
\end{prop}
\Dim By Lemma~\ref{ads:bend:lem}, $\beta_+(x,y)$ is a hyperbolic
transformation, whose axis separates the stratum through $x$ from the
stratum through $y$ and whose translation distance is bigger than the
total mass of $[x_,y]$.  This fact easily implies that $\Ee_L$ is
injective.

For every unit tangent vector $v$ at $x_0$ let $u(v)$ be the end-point
of the intersection of the geodesic ray $c(t)=\exp(tv)$ with the
boundary of $H$ (notice that $u(v)$ can lie on $\partial\mh^2$).
First suppose that $u(v)$ is an accumulation point for $L$ (that is
the support of $\lambda$).  Then for every $t\in [x_0, u(v)]\cap L$
let $P_v(t)$ be the half-plane of $\mh^2$ bounded by $\beta_L(x_0, t)
l_t$ (where $l_t$ is the leaf through $t$) and containing $x_0$. By
Lemma~\ref{ads:bend:lem} we have that
\[
   P_v(t)\subset P_v(s) \qquad\textrm{ if } t<s\,.
\]
Thus, $P_v(s)$ converges to either the whole $\mh^2$ or to a
half-plane for $t\rightarrow u(v)$. Let us denote by $P_v$ such a
limit (that is none but the union of all $P_v(s)$).

If $u(v)$ is not an accumulation point of $L$ then let us put $P_v=\mh^2$. 
It is not difficult to see that if $v\neq v'$ then either $P_v=P_{v'}$ or
\[
  \partial P_v\cap\partial P_{v'}=\varnothing
\]
Hence, the intersection of all $P_v$'s is a straight convex set.
Now we claim that
\[
   \Ee_L(H)=\bigcap_{v\in T^1_{x_0}\mh^2} P_v\,.
\]
The inclusion $(\subset)$ is quite evident.  So let us prove the other
inclusion. First let us prove that the image of $\Ee_L$ is convex.
Given $x,y\in\mathring H$ there exists a rectifiable
arc $c$ in $H$ of length equal to
$d_{\mh}(\Ee_L(x),\Ee_L(y))+\mu([x,y])$ such that $\Ee_L(c)=[\Ee_L(x),
\Ee_L(y)]$.  If $[x,y]$ meets only finite leaves then it is clear how
to construct $c$. The general case follows by using standard
approximations.

Now suppose there exists $x\notin\Ee_L(H)$. Then there exists a point
$y$ on $[x_0,x]$ such that $[x_0,y)\subset\Ee_L(H)$ and $(y,x]$ does
not intersect $\Ee_L(H)$. Thus we see that there exists a locally
rectifiable transverse arc $k$ in $\mathring H$ such that
$\Ee_L(k)=[x_0,y)$.  Since $k$ intersects each stratum in a convex set, $k$ has limit, $u_\infty$, lying on the boundary of $H$ in
$\overline\mh^2$. Moreover, the segment $[x_0,u_\infty]$ is homotopic
to $k$ through a family of transverse arcs.  This implies that there
exists $t_n$ on $[x_0, u_\infty)\cap L$ such that $y$ is an
accumulation point of $\beta_L(x_0,t_n)l_{t_n}$ (where $l_{t_n}$ is
the leaf through $t_n$). Thus we get $x\notin P_{u_\infty}$.

The image of the leaves of $\lambda$ form a geodesic lamination
$\hat\Ll$ on $\hat H=\Ee_L(H)$.  Now let us take a geodesic arc $k$ in
the interior of $\Ee_L(H)$ We have seen that there exists a
rectifiable arc $k'$ in $H$ such that $\Ee_L(k')=k$.  Let us set
$\hat\mu_k$ the image of the measure $\mu_{k'}$.  It is easy to show
that $\hat\mu$ satisfies points $1.$ and $2.$ in the definition of
transverse measure given in
Section~\ref{laminazioni:meas:sec}. Moreover, if $k_1$ and $k_2$ are
geodesic arcs homotopic through a family of transverse arcs then
$\mu_{k_1}(k_1)=\mu_{k_2}(k_2)$. In order to conclude we should see
that the total mass of an arc reaching the boundary of $\hat H$ in
$\mh^2$ is infinite.  Before proving this fact, notice that, however,
we can define the Epstein-Marden right cocycle, say $\hat\beta_R$ on the interior of $\hat H$.  Now we want to prove that the map
\[
   \mathring{\hat H}\ni x\mapsto \hat\beta_R(x_0,x)x\in\mh^2
\]
is the inverse of $\Ee_L$.
In fact it is sufficient to prove that
\begin{equation}\label{memgen:inveart:eq}
   \hat\beta(\Ee_L(x),\Ee_L(y))\circ\beta(x,y)=Id
\end{equation}
for every $x,y\in\mathring H$.  Choose a standard approximation
$\lambda_n$ of $\lambda$ between $x$ and $y$.  The image of
$\lambda_n$ through $\Ee_L$ is a standard approximation, say
$\hat\lambda_n$, of $\hat\lambda$ between $\Ee_L(x)$ and $\Ee_L(y)$
Denote by $\beta_L^n$ and $\hat\beta_R^n$ the left and right cocycle
associated to $\lambda_n$ and $\hat\lambda_n$ respectively.  For a
fixed $n$ denote by $x_1=x,\ldots,x_n=y$ the intersection points of
$\lambda_n$ with the segment $[x,y]$ and let $g_i\in PSL(2,\mr)$ be
the translation along the leaf of $\lambda_n$ through $x_i$ with
translation length equal to the mass of $x_i$.  It turns out that
$\beta_L^n(x,y)=g_1\circ\cdots\circ g_n$, and
$\hat\beta_R^n(\Ee_L(x),\Ee_L(y))= g_1^{-1}\beta_1
g_2^{-1}\beta_1^{-1}\cdots \beta_{n-1} g_n\beta_{n-1}^{-1}$ where
$\beta_i=\beta_L(x,x_{i+1})$.  By Lemma~\ref{adesitt:bend:cont:lem}
$||\beta_{n-1}- g_1\circ\ldots\circ g_n||\leq C/n$, whereas
$||\beta_{i-1}^{-1}\beta_i- g_i||\leq C m_i/n$ where $m_i$ is the mass
of the segment $[x_i,x_{i+1}]$. It follows that
\[
   ||\hat\beta_R^n(\Ee_L(x), \Ee_L(y))\circ\beta_L^n(x,y)-Id||\leq C'/n\,.
\]
Passing to the limit shows the identity~(\ref{memgen:inveart:eq}).\\

Now we can prove that if $k$ is an arc reaching the boundary of $\hat
H$ in $\mh^2$ then its total mass is infinite. Let $x_n$ be a sequence
of points on $k$ converging to a point on the boundary of $\hat H$. We
have that $\hat\beta_R(x_0,x_n)x_n$ either converge to a boundary leaf
of $\lambda$ or converges to a point on $\partial\mh^2$. In the former
case, it follows that the measure of $k$ is $+\infty$ by the
hypothesis on $\lambda$. In the latter case, we have that the measure
of $k$ must be infinite, otherwise the estimates of
Lemma~\ref{adesitt:bend:cont:lem} should imply that
$\hat\beta_R(x_0,x_n)$ is a precompact family in $PSL(2,\mr)$.  \cvd

\begin{remark}
\emph{ Let $\tilde Y$ (different from $\Pi_{-1}$) be the universal
covering of a spacetime $Y$, with holonomy representation
$h=(h_-,h_+):\pi_1(Y)\rightarrow PSL(2,\mr)\times PSL(2,\mr)$.  The
future boundary of the convex core of $\tilde Y$ is isometric to a
straight convex set, $H$, of $\mh^2$ bent along a measured geodesic
lamination $\lambda$.  There exists a discrete representation
$h_0:\pi_1(Y)\rightarrow PSL(2,\mr)$ such that $H$ and $\lambda$ are
$h_0$-invariant and the bending map $\varphi_\lambda:\mathring
H\rightarrow\partial_+\Kk$ is $\pi_1(Y)$-equivariant.  Since
$\varphi_\lambda(x)=(\beta_-(x_0,x),\beta_+(x_0,x))I(x)$, the
generalized right earthquake along $\lambda$, say $\Ee_R$, conjugates
$h_0$ with $h_-$. Thus one can  see that $h_-$ is a discrete
representation and $\mh^2/h_-$ is homeomorphic to $\mathring H/h_0$. This
furnishes another proof of the first part of
Proposition~\ref{discrete-rep}.  }
\end{remark}

By means of the generalized earthquakes we can characterize the
measured laminations that produce standard AdS spacetimes whose
boundary at infinity is the graph of a homeomorphism.

\begin{prop}\label{earth:homeo:prop}
Let $\lambda$ be a measured geodesic lamination on a straight convex
set. Then the following statements are equivalent.
\medskip\par\noindent
1) The left and right earthquakes along $\lambda$ are surjective maps on
$\mh^2$.
\medskip\par\noindent
2) The boundary curve of $\Yy_\lambda$ is the graph of a homeomorphism.
\end{prop}

First suppose that right and left earthquakes are surjective
maps. Denote by $\hat\lambda$ the image of $\lambda$ via the
(generalized) right earthquake $\Ee_R$ along $\lambda$. By hypothesis
the left earthquake along $2\hat\lambda$ is a true-earthquake (that
means that it is surjective). Thus, it extends to a homeomorphism
$u_{\hat\lambda}$ of $\partial\mh^2=S^1_\infty$.  We  show that
the curve at infinity of $\Uu^{-1}_\lambda$ (that is a subset of
$\partial\mx_{-1}=S^1_\infty\times S^1_\infty$) is the graph,
$C_\lambda$, of $u_\lambda$.  In fact, it is sufficient to show that
such a curve contains the pair $(x,u_\lambda(x))$ for any vertex $x$ of
any stratum of $\hat\lambda$.  Now for such an $x$ there exists a
vertex $y$ of a stratum $T$ of $\lambda$ such that
\[
    x=\beta_-(x_0, z)y
\]
where $z$ is any point in $T$. Thus the image of $x$ via the left
earthquake along $2\hat\lambda$ is $\beta_+(x_0,z)y$. On the other
hand, since the image through the bending map
$\varphi_\lambda:\mathring H\rightarrow\mx_{-1}$ of $T$ is
$(\beta_-(x_0,z), \beta_+(x_0,z))(T)$, the point $(x,
u_\lambda(x))=(\beta_-(x_0,z), \beta_+(x_0,z))(y)$ lies on the
boundary curve of $\Yy_\lambda$.\\

Conversely suppose that the boundary curve $C_\lambda$ is the graph of
a homeomorphism.  Take a point $(x,y)$ on $C_\lambda$ and let $p_n$ be
a sequence of points in $\partial_+\Kk_\lambda$ such that
$p_n\rightarrow (x,y)$ in $\overline\mx_{-1}$. If $T_n$ is a face or a
bending line through $p_n$ we have three cases:
\medskip\par
1) $T_n$ converges to a stratum $T_\infty$ of $\partial_+\Kk_\lambda$.
\medskip\par
2) $T_n$ converges to  $(x,y)$.
\medskip\par
3) $T_n$ converges to a segment on a leaf of $\partial\mx_{-1}$.
\medskip\par\noindent Since $C_\lambda$ is the graph of a
homeomorphism, we can discard the last case. In the other cases, it is
easy to construct a sequence of vertices $q_n$ of faces or bending
lines $T_n$ converging to $(x,y)$.

If $S_n$ is the stratum of $\lambda$ corresponding to $T_n$ via the
bending map, we can find a sequence of end-points of $S_n$, say $z_n$,
such that
\[
\begin{array}{l}
  x=\lim_{n\rightarrow +\infty} \beta_-(x_0, u_n)z_n\\
  y=\lim_{n\rightarrow+\infty}\beta_+(x_0,u_n)z_n
\end{array}
\]
where $u_n$ is any point of $S_n$.  Notice that $\beta_-(x_0,u_n)z_n$
lies in the closure of the image of the right earthquake $\Ee_R$ along
$\lambda$. Thus $x\in\overline{\Ee_R( H)}$. In an analogous way we see
that $y\in\overline{\Ee_L( H)}$.

Summarizing, we have proved that if $(x,y)\in C_\lambda$ then
$x\in\overline{\Ee_R( H)}$ and $y\in\overline{\Ee_L( H)}$.
Since $C_\lambda$ is the graph of a homeomorphism we have that
$S^1_\infty$ is contained in the closure of the image of $\Ee_L$ (and
$\Ee_R$).
Since such images are convex, both $\Ee_L$ and $\Ee_R$ are surjective.
\cvd

We stress that Proposition~\ref{earth:homeo:prop} does {\it not} imply
that if the boundary curve at infinity is the graph of a homeomorphism
then the future boundary is complete.  We show a counterexample.

\begin{figure}
\begin{center}
\input{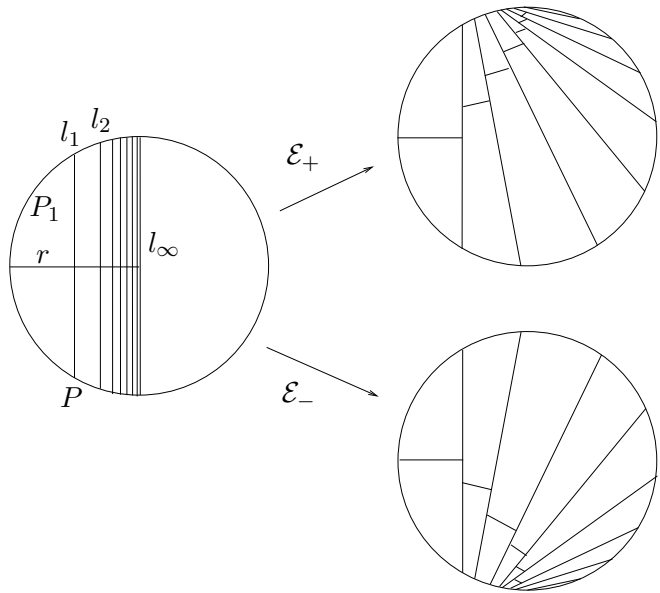}
\caption{{\small The maps $\Ee_\pm$ are injective. Since the geodesic
    $l_\infty$ escapes to infinity then they are also surjective.}}
\end{center}
\end{figure}

\begin{exa}\label{earthex}
{\rm Let $H$ be a half-plane bounded by a geodesic $l$ and $r$ be a
geodesic ray contained in $H$ starting from some point $x_\infty$ of
$l$ and orthogonal to $l$.  Denote by $x_n$ the point on the ray such
that $d_\mh(x_\infty,x_n)=1/n$.  Let $l_n$ (resp. $l_\infty$) be the
geodesic through $x_n$ (resp. $x_\infty)$ orthogonal to $r$. Then
$\Ll=\{l_n|n\in \mn_+\}\cup{l}$ is a geodesic lamination on
$H$. Putting the weight $1$ on each $l_n$ equips $\Ll$ with a
transverse measure, say $\mu$. We claim the the left and right
earthquakes along $\lambda=(H,\lambda,\mu)$ are surjective maps}
\[
   \Ee_L,\Ee_R:\mathring H\rightarrow\mh^2\,.
\]
{\rm Denote by $H_n$ the half-plane bounded by $l_n$ and contained in
$H$ and take for each $n$ a point $u_n\in H_n\setminus H_{n-1}$.  By
Proposition~\ref{earth:gen:prop} we have to prove that the sequence of
geodesics $\beta_\pm(u_0,u_n)(l_n)$ (where $\beta_\pm$ is the
Epstein-Marden cocycle associated to $\pm\lambda$) is divergent.  The
transformation $\beta_n(u_0,u_n)$, that for simplicity will be denoted
by $g_n$, is the composition of hyperbolic transformations with axes
$l_1,\ldots,l_{n-1}$ and translation lengths equal to $1$.  By
Lemma~\ref{ads:bend:lem} the following facts hold}\par
\smallskip
{\rm 1) $g_n$ is a hyperbolic transformation and its axis is contained
in $H-H_{n-1}$.}\par
\smallskip
{\rm 2) The translation length of $g_n$, say $\alpha_n$, is greater
than $n-1$.}\par
\smallskip
{\rm 3) The distance of $l_n$ from the axis of $g_n$, say $\eps_n$, is
greater then $1/n(n-1)$.}\par

{\rm Let $z_n$ denote the point on the axis of $g_n$ that realizes the
  distance from $l_n$. The distance of $z_n$ from $g_n(l_n)$ is then
  equal to $\sh\eps_n\,\ch\alpha_n$. Because of points 2) and 3)  this
  number is bigger than
  $$\sh(\frac{1}{n(n-1)})\ch(n-1)>\frac{1}{n(n-1)}\ch(n-1)\,.$$ Thus
  the distance of $z_n$ from $g_n(l_n)$ tends to $+\infty$.  Because
  of point 1)  $z_n$ runs in a compact set and this imply that
  $g_n(l_\infty)$ is divergent.}\par

{\rm The image of $\lambda$ via $\Ee_R$, say $\hat\lambda$, is a
measured geodesic lamination of $\mh^2$. Thus we see that the left
generalized earthquake along $\hat\lambda$ is not a true earthquake,
whereas the  one along $2\hat\lambda$ is.}
\end{exa}

In Section~\ref{3cusp} we will prove that also the converse is
false. That is, we will show examples of spacetimes whose boundary
curve is not a homeomorphism and such that the future boundary of the
convex core is complete (this is related to the well-known fact that
there are measured geodesic laminations of $\mh^2$ that do not give
rise to true earthquakes).

\section{$T$-symmetry}\label{T-symm}
Let $Y$ be any maximal globally hyperbolic AdS spacetime containing a
complete Cauchy surface (that is the quotient of some standard domain
$\Yy=\Yy(C)$ in $\mx_{-1}$). It is evident from our previous
discussion that by reversing the time orientation we get a spacetime
$Y^*$, quotient of $\Yy^* = \Yy(C^*)$ where $C^*$ is the image of the
curve $C$ under the involution of $\partial\mx_{-1}=S^1_\infty\times
S^1_\infty$
$$ (x,y)\mapsto(y,x)\ .$$ Moreover, the holonomy of $Y^*$ is
obtained by exchanging the components of the holonomy of $Y$
$$ (h_-, h_+) \leftrightarrow (h_+, h_-) \ .$$
Thanks to the classification, this induces an involution on the
set of AdS $\Mm\Ll$-spacetimes, hence on $\Mm\Ll^\Ee$, called
$T$-symmetry. That is Proposition \ref{TSYM:PROP} is now proved.

In Section~\ref{3cusp} we will illustrate by some examples
the broken $T$-symmetry on general $\Mm\Ll(\mh^2)$-spacetimes
(according to Section \ref{BROKEN:I}).

\section{Examples}\label{3cusp}
Let us summarize some nice properties satisfied by any maximal
globally hyperbolic spacetime $Y$ of constant curvature $\kappa$ that
contains a closed Cauchy surface of genus $g\geq 2$:
\smallskip

(a) It is maximal in the strong sense. In fact it cannot be embedded
    in a bigger constant curvature globally hyperbolic spacetime
    (usually, maximal means that there exists no isometric embedding
    sending in a bigger spacetime $Y'$ sending a Cauchy surface of $Y$
    onto a Cauchy surface of $Y'$, see Section~\ref{glob:hyp:ST}).
\smallskip

(b) If $\kappa=0,-1$, $Y$ is determined by its holonomy.
\smallskip

(c) If $\kappa=-1$, both the future and past boundaries of its convex
core are complete (with respect to the intrinsic metric).
\smallskip

(d) If $\kappa=-1$, the boundary curve of the universal covering of
    $\Yy$ is the graph of a homeomorphism of $S^1_\infty$ into itself.
\smallskip

(e) If $\kappa=-1$, the set of such spacetimes is closed for the
$T$-symmetry.

\medskip

In these section we will illustrate examples that show that these
properties fail for general $\Mm\Ll(\mh^2)$-spacetimes, even if the
surface is of finite type with negative Euler characteristic.

The elements in $\Mm\Ll^\Ee$ corresponding to spacetimes we are going
to construct will be of the form $(\lambda, \Gamma)$ where $F=
\mh^2/\Gamma$ is the finite area hyperbolic surface homeomorphic to the
{\it thrice-punctured sphere}.

It is well known that $F$ is {\it rigid}, that is the corresponding
Teichm\"uller space is reduced to one point. Hence, $F$ has no
measured geodesic laminations with compact support.

$F$ can be obtained by gluing two geodesic ideal triangles along their
edges as follows.  In any ideal triangle there exists a unique point
(the "barycenter") that is equidistant from the edges. In any edge
there exists a unique point that realizes the distance of the edge
from the barycenter. Such a point is called the mid-point of the edge.
Now the isometric gluing is fixed by requiring that mid-points of
glued edges match (and that the so obtained surface is topologically
a three-punctured sphere - by a different pattern of identifications
we can obtain a $1$-punctured torus).  It is easy to see that the
resulting hyperbolic structure is complete, with a cusp for any
puncture, and equipped by construction with an {\it ideal
triangulation}.

The three edges of this triangulation form a geodesic lamination
$\Ll_F$ of $F$. A transverse measure $\mu_F = \mu_F(a_1,a_2,a_3)$ on
such a lamination consists of giving each edge a positive weight
$a_i$.  The ideal triangles in $F$ lift to a tessellation of the
universal cover $\mh^2$ by ideal triangles.  The $1$-skeleton $\Ll$ of
such a tessellation is the pull-back of $\Ll_F$; a measure $\mu_F$
lifts to a $\Gamma$-invariant measure $\mu$ on $\Ll$. So, we will
consider the $\Gamma$-invariant measured laminations $\lambda =
(\Ll_F,\mu_F)$ on $\mh^2$ that arise in this way. From now on in the
present Section we will refer to such a family of $\Mm\Ll$-spacetimes.

\paragraph {Blind flat Lorentzian holonomy}
That is we show that in general the above property (b) fails.  Varying
the weights $a_i$, we get a $3$-parameters family of flat spacetimes
$\Uu^0_{\lambda}= \Uu^0_{\lambda(a_1,a_2,a_3)}$, with associated
quotient spacetimes $\hat \Uu^0_{\lambda} =
\Uu^0_{\lambda}/h^0_\lambda(\Gamma)$, where $h^0_\lambda(\Gamma)$
denotes the flat Lorentzian holonomy.

The spacetimes $\Uu^0_{\lambda}$ have {\it homeomorphic} initial
singularities $\Sigma_\lambda$.  A topological model for them is given
by the {\it simplicial} tree (with $3$-valent vertices) which forms
the $1$-skeleton of the cell decomposition of $\mh^2$ dual to the
above tessellation by ideal triangles. The length-space metric of each
$\Sigma_\lambda$ is determined by the fact that each edge of the tree
is a geodesic arc of length equal to the weight of its dual edge of
the triangular tessellation. In fact, every $\Sigma_\lambda$ is
realized as a spacelike tree embedded into the frontier of
$\Uu^0_{\lambda}$ in the Minkowski space, and 
$h^0_\lambda(\Gamma)$ acts on it by isometries.

The behaviour of the {\it asymptotic states} of the cosmological time
of each $\Uu^0_{\lambda}$, is formally the same as in the cocompact
case. In particular, when $a\to 0$, then action of 
$h^0_\lambda(\Gamma)$
on the level surface $\Uu^0_\lambda(a)$ converges to action on the
initial singularity $\Sigma_\lambda$.  The marked length spectrum of
$\hat \Uu^0_\lambda(a)$ (which coincides with the minimal displacement
marked spectrum of the action of $h^0_\lambda(\Gamma)$ on
$\Uu^0_\lambda(a)$), converges to the minimal displacement marked
spectrum of the isometric action on the initial singularity.  If
$\gamma_i$, $i=1,2,3$, are (the conjugacy classes of) the parabolic
elements of $\Gamma$ corresponding to the three cusps of $F$, the last
spectrum takes values $\gamma_i \to a_i+a_{i+1}$, where we are
assuming that the edges of the ideal triangulation of $F$ with weights
$a_i$ and $a_{i+1}$ enter the $i$th-cusp, $a_4=a_1$.  By the way, this
implies that these spacetimes are not isometric to each other.

However, it follows from \cite{Ba} that:
\medskip

{\it 
(1) The flat Lorentzian holonomies $h^0_\lambda(\Gamma)$ are 
all conjugated (by $\ISO_0(\mx_0)$) to their common 
{\rm linear part} $\Gamma$.
\smallskip

(2) All non trivial classes in $H^1(\Gamma, \mr^3)$ are not realized
by any flat spacetime having $\Gamma$ as linear holonomy.}
\medskip

\noindent Hence, the flat Lorentzian holonomy is completely 
``blind'' in this case, and, on the other hand, the non trivial algebraic
$\ISO_0(\mx_0)$-extensions of $\Gamma$ are not correlated to the
geometry of any spacetime.

\begin{remarks}
{\rm (1) The above facts would indicate that the currently accepted
  equivalence between the classical formulation of 3D gravity in terms
  of Einstein action on metrics, and the formulation via Chern-Simons
  actions on connections (see \cite{W} and also Section \ref{END} of
  Chapter \ref{INTRO}), should be managed instead very carefully
  outside the cocompact $\Gamma$-invariant range (see \cite{Mat} for
  similar considerations about flat spacetimes with particles).
\smallskip

(2) Every $\Uu^0_{\lambda}$ in the present family of examples can be
embedded in a $\Gamma$-invariant way in the static spacetime $\fut(0)$
as the following construction shows (this is the geometric meaning of
point (1) above, and shows by the way that (a) above fails).
\smallskip

Up to conjugating $h^0_\lambda$ by an isometry of $\ISO(\mx_{0})$ we
can suppose $h^0_\lambda(\gamma)=\gamma$ for every $\gamma\in\Gamma$.
We want to prove that $\Uu^0_\lambda$ is contained in $\fut(0)$. By
contradiction suppose there exists $x\in\Uu^0_\lambda$ outside the
future of $0$. Since $\Uu^0_\lambda$ is future complete
$\fut(x)\cap\partial\fut(0)$ is contained in $\Uu^0_\lambda$. But we
know there is no open set of $\partial\fut(0)$ such that the action of
$\Gamma$ on it is free. 

There is a geometric way to recognize such domains inside $\fut(0)$.
Take a $\Gamma$-invariant set of horocycles $\{B_n\}$ in $\mh^2$
centered to points corresponding to cusps.  We know that every
horocircle $B_n$ is the intersection of $\mh^2$ with an affine null
plane orthogonal to the null-direction corresponding to the center of
$B_n$.  Now it is not difficult to see that the set
\[
    \Omega= \bigcap\fut(P_n)
\]
is a regular domain invariant by $\Gamma$.  This is clear if $B_n$ are
sufficiently small (in that case we have that
$\Omega\cap\mh^2\neq\varnothing$). For the general case denote by
$\hat B_n(a)$ the intersection of $P_n$ with the surface
$a\mh^2$. Then the map
\[
    f_a: a\mh^2\ni x\mapsto x/a\in\mh^2
\]
sends $\hat B_n(a)$ to a horocircle $B_n(a)$ smaller and smaller as
$a$ increases.  It follows that $\Omega\cap a\mh^2\neq\varnothing$ for
$a>>0$.  Since a regular domain is the intersection of the future of
its null-support planes it follows that every $\Uu^0_\lambda$ can be
obtained in this way.
\smallskip
}
\end{remarks}

\paragraph{Earthquake ``failure'' and broken $T$-symmetry}
It was already remarked in
\cite{Thu2} that such a $\lambda$ produces neither left nor right true
earthquake (in particular the boundary curve of the universal covering
of $Y(\Gamma,\lambda)$ is not a homeomorphism by
Proposition~\ref{earth:homeo:prop} - this could be checked also
directly). On the other hand $\partial_+ \Kk$ is complete (they are
$\Mm\Ll(\mh^2)$-spacetimes). Hence this is the example promised at
the end of Section \ref{ge:quake}.

Nevertheless, we can consider the Epstein-Marden cocycles
$\beta_-=\beta_{-\lambda}$ and $\beta_+=\beta_{+\lambda}$.  The
associated representations
\[
\begin{array}{l}
 h_L(\gamma)=\beta_{-}(x_0,\gamma(x_0))\gamma\\
 h_R(\gamma)=\beta_+(x_0,\gamma(x_0))\gamma
\end{array}
\]
are faithful and discrete.  Both the limit set $\Lambda_L$ of $h_L$
and $\Lambda_R$ of $h_R$ are Cantor sets such that both quotients of
the respective convex hulls are isometric to the same pair of
hyperbolic pants with totally geodesic boundary $\Pi_\lambda$. If
$\gamma_i$ is as above, then we have that $h_L(\gamma_i)$ and
$h_R(\gamma_i)$ are hyperbolic transformations corresponding to the
holonomy of a boundary component of $\Pi_\lambda$, with translation
length equal to $a_i+a_{i+1}$.

Thus $h_L$ and $h_R$ are conjugated in $PSL(2,\R)$.  The
above remarks imply that there exists a spacelike plane $P$ in
$\mx_{-1}$ that is invariant for the representation $h=(h_L, h_R)$.
The limit set of the action of $h$ on $P$ is the Cantor set
\[
\Lambda=\overline{\{(x^+_L(\gamma)), x^+_R(\gamma))\}}\, ,
\]
$x^+_L(\gamma)$ (resp. $x^+_R(\gamma)$) denoting the attractive fixed
point of $h_L(\gamma)$ (resp. $h_R(\gamma)$).

It is not hard to see that every closed $\Gamma$-invariant subset of $\partial\mx_{-1}$ must contain $\Lambda$.

\smallskip

{\it Broken $T$-symmetry.} Consider now the bending of
$\mh^2$ along $\lambda$. The key point is to describe the curve
$C_\lambda$, that is the curve at infinity of the universal covering
of $Y^{-1}(\Gamma,\lambda)$.

\begin{figure}
\begin{center}
\input{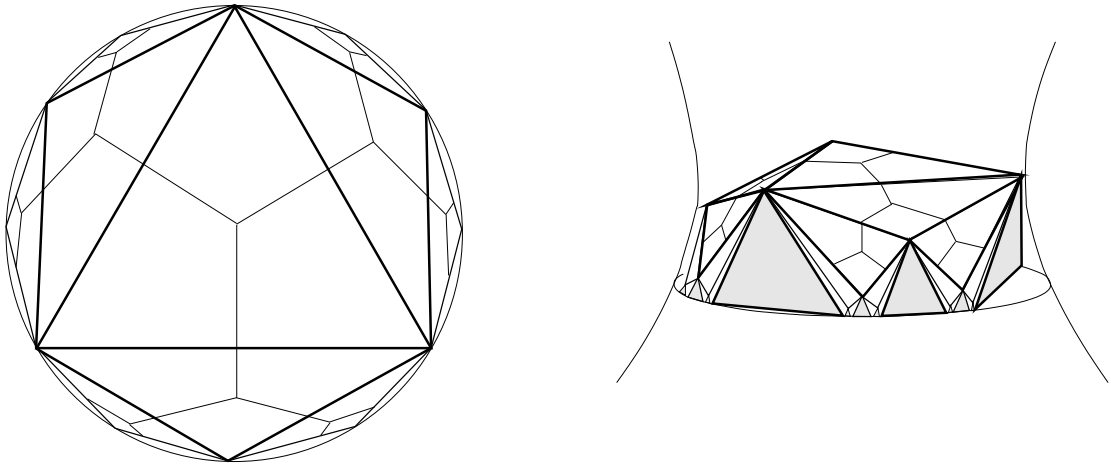}
\caption{{\small On the left the lamination $\Ll$ with its dual spine. On the
    right the bending of $\mh^2$ along $\lambda$ in $\mx_{-1}$. Grey regions
    are null components of the past boundary of $\Kk_\lambda$.}} 
\end{center}
\end{figure}

Take a point $x\in\partial \mh^2$ that is a vertex of a triangle $T$ of
$\lambda$. The point $x$ corresponds to a puncture of $F$ so there is
a parabolic transformation $\zeta\in\Gamma$ conjugated to one
$\gamma_i$ that fixes $x$.  Moreover, we can choose $\zeta$ in such a
way that it is conjugated to a translation $z\mapsto z+a$ with $a>0$
in $PSL(2,\R)$.  If we take $z\in T$, the point
$\beta_-(x_0,z)x$ is the repulsive fixed point of $h_L(\zeta)$
whereas $\beta_+(x_0,z)x$ is the attractive fixed point of
$h_R(\zeta)$.  
Since $T$ is sent via the bending map to $\beta(x_0,z)(T)$, we see that
$C_\lambda$ contains $u=(x^-_L(\zeta), x^+_R(\zeta))$.

Since $C_\lambda$ is $\Gamma$-invariant it also contains
\[
\begin{array}{ll}
u_\infty=(x^+_L(\zeta),x^+_R(\zeta)) & v_\infty=(x^-_L(\zeta),x^-_R(\zeta))\,.
\end{array}
\]

Since $C_\lambda$ is achronal and $u$ and $v_\infty$
are in the same right leaf, it follows that the future directed segment
on the right leaf from $v_\infty$ towards $u$ is contained in
$C_\lambda$.  In the same way we have that the future directed segment
on the left leaf from $u_\infty$ towards $u$ is contained in
$C_\lambda$.  \par

Notice that $u_\infty$ and $v_\infty$ are the
vertices of a boundary component of the convex hull $H$ of
$\Lambda$ in $P$.  

Now take a boundary component $l$
of $H$ oriented in the natural way and let $u_-$ and $u_+$ be its
vertices. Then the right leaf through $u_-$ and the left leaf through
$u_+$ meet each other at a point $u_l$.  The future directed segments
from $u_\pm$ towards $u$ in the respective leaves is contained in
$C_\lambda$. Denote by $V_l$ the union of such segments.  We have that
the union of $V_l$, for $l$ varying in the boundary components of $H$,
is contained in $C_\lambda$. But  its
closure is a closed path so that $C_\lambda$ coincides with it.

Since the description of the curve $C_\lambda$ is quite simple we can
describe also the past boundary $\partial_- \Kk$ of the convex hull
of $C_\lambda$, \emph{i.e.} of the AdS convex core of $\Uu^{-1}_\lambda$. 
Notice that $P$ is the unique spacelike support plane
touching $\partial_- \Kk$.  Then for every component $l$ of $H$, there
exists a unique null support plane $P_l$ with dual point at $u_l$. Thus,
$\partial_- \Kk$ is the union of $H$ and an infinite number of
null half-planes, each attached to a boundary component of $H$.  
It follows that $\partial_- \Kk$ is {\it not} complete. This shows
that properties (c) and (e) stated at the beginning of this Section
fail.


\newpage{\pagestyle{empty}\cleardoublepage} 
\chapter {$\Qq\Dd$-spacetimes}\label{QD}

In this chapter we treat the Wick rotation-rescaling theory on
``degenerate'' regular domains in $\mx_0$, \emph{i.e.} the future
$I^+(r)$ of spacelike geodesic lines $r$ of $\mx_0$. In fact we will
extend the theory on (flat) spacetimes modeled on $I^+(r)$, that are
governed by quadratic differentials on Riemann surfaces.

\section{Quadratic differentials}\label{basic_QD}
We recall a few general facts about quadratic differentials on
Riemann surfaces. We refer for instance to \cite{Ab, MS, Ke} for details.
\smallskip

Let $F$ be a Riemann surface. A {\it meromorphic quadratic differential}
$\omega$ on $F$ is a meromorphic field of quadratic forms on $F$. 
In local coordinates $\omega$ looks like
$$\omega = \phi dz^2$$ for some meromorphic function $\phi$. We will
limit ourselves to consider quadratic differentials that have poles of
order $2$ at most. If $\omega \neq 0$, then its {\it singular locus}
$X(\omega)= X'(\omega)\cup X''(\omega)$
is the discrete subset of $F$ where $X'$ is the union of zeros and simple poles
of $\omega$, while $X''$ is made by the order $2$ poles. 
Set $F' = F\setminus X(\omega)$. If $\omega$ is {\it holomorphic} at $p\in F$,
then there are local {\it normal coordinates} $z=u+iy$ such that, $p=0$
and
$$\omega = z^mdz^2, \ m\in \N \ , $$ where $m\geq 1$ iff $p$ is a zero
of order $m$.  The local normal form at a simple pole is
$$\omega = \frac{1}{z}dz^2 \ . $$
The local normal forms at poles of order $2$ depend on one complex 
{\it modulus} as they are
$$ \omega = \frac{a}{z^2}dz^2,\ \ a\in \C \ . $$ 

On $F'$ we have the $\omega$-{\it metric}
$$ds^2=ds_\omega^2  = (\omega \overline{\omega})^{\frac{1}{2}}$$
\emph{i.e.}, in local coordinates
$$ ds^2=|\phi||dz|^2 \ .$$ 
By using the local normal forms, it is easy to see that $ds_\omega^2$ 
satisfies the following properties:

(1) it is {\it flat} on $F'$;

(2) it extends to a conical singularity with cone angle
$(m+2)\pi$ at each zero of order $m$; it extends to a conical
singularity with cone angle $\pi$ at each simple pole;

(3) every pole of order $2$ (with local normal form as above)
gives rise to a cylindrical ``end'' of $F'$ isometric to 
$S^1\times \R^+$ endowed with the metric 
$$|a|^2d\theta^2 + dt^2  \ .$$

\smallskip

A vector $v\in TF'_p$ is said to be $\omega$-horizontal ( $\omega$-vertical)
if $\omega(v)$ is real and strictly positive (negative). This induces
on $F'$ two fields of directions, that are orthogonal each other
with respect to $ds^2_\omega$.  The integral lines of these fields give rise to
the {\it $\omega$-horizontal} and {\it $\omega$-vertical foliations}
respectively, denoted $\Ff_*=\Ff_*(\omega)$, $*=h,v$.  
By using local normal forms, we see that:

(i) at non-singular points they are given by the $y$-constant and
$u$-constant lines, respectively; 

(ii) at any zero $p$ of order $m$, both foliations extend to singular
foliations having a saddle singularity, with $(m+2)$ germs of
singular leaves emanating from $p$. In a similar way, they extend at 
any simple pole $p$, with one germ of singular leaf emanating from $p$; 

(iii) at a pole of order $2$, the foliations are induced by a pair of
{\it constant} orthogonal vector fields on the end $S^1\times \R^+$
(with coordinate $(\theta, t)$ as above). In fact if the parameter of
the pole is $a=|a|e^{i\alpha}$, then the horizontal direction is given
by rotating $\frac{\partial}{\partial t}$ as follows
$$e^{-\frac{i\alpha}{2}}\frac{\partial}{\partial t} \ .$$ In
particular, there is a closed $\omega$-horizontal leaf iff $a$ is real
and $a<0$. In such a case every closed curve $S^1\times \{*\}$ is in
fact a closed horizontal leaf.
\smallskip

We will say that a leaf of these singular foliations is {\it
non-singular} if it does not end at any zero or simple pole of
$\omega$. 

The foliations $\Ff_*$  are endowed with {\it transverse
measures} $\mu_*$ which in a normal coordinate at a
non-singular point are given by $|dy|$ and $|du|$, respectively.
\smallskip

In fact, either the data $(ds^2,\ (\Ff_*,\mu_*), \ *=h,v)$ or
$(F,\omega)$ determine each other. Note, in particular, that
$(F,-\omega)$ corresponds to $(ds^2,\ (\Ff_v,\mu_v),\ (\Ff_h,\mu_h))$
\emph{i.e.} the flat metric is the same, but the measured foliations exchange
each other.

This suggests a very convenient geometric way to deal with quadratic
differentials, by using the $(\mx,\Gg)$-structure machinery (see
section \ref{XG}).
\medskip

{\it Local model.} Consider the complex plane $\C$ with coordinate
$z=u+iy$, and the non singular quadratic differential $\omega_0 =
dz^2$. Then $ds^2_0 = |dz|^2 $ is the $\omega_0$-metric and we have
the associated $\omega_0$-measured foliations $(\Ff_*^0,\mu_*^0)$.

The subgroup Aut$(\C,\omega_0)$ of Aut$(\C)$ preserving $\omega_0$
(equivalently, the group of direct isometries of $ds^2_0$ preserving
$(\Ff_*^0,\mu_*^0)$) is generated by translations
$$\sigma_{v}(u+iy)= (u+iy + v), \ v=p+iq $$ and by the
rotation $R_\pi$ of angle $\pi$. This last inverts the orientation of
the $\omega_0$-foliations.

The {\it Teichm\"uller ray} based on $\omega_0$ is given by
the $1$-parameter family of structures
$$(ds^2_\tau, (\Ff_h, \tau\mu_h), (\Ff_v,\mu_v)),\ \ \tau\geq 1 $$
where
$$ ds_{\tau}^2 = \tau^2du^2 + dy^2, \tau \geq 1  \ .$$
By taking $\tau >0$ we have the associated Teichm\"uller {\it line}.  
\smallskip

Let $S$ be an oriented surface, $X=X'\cup X''$ be a discrete subset of
$S$.  Set $S' = S\setminus X$. Consider any $(\C, {\rm
Aut}(\C,\omega_0))$-{\it structure} on $S'$. This is equivalent to
giving $S'$ a Riemann surface structure $F'$ equipped with a
non-singular quadratic differential $\omega'$. In fact a (maximal)
$(\C, {\rm Aut}(\C,\omega_0))$-atlas coincides with a maximal family of
local normal coordinates for $(F',\omega')$.  We say that this
$(\C, {\rm Aut}(\C,\omega_0))$-structure is of {\it meromorphic type}
(with singular set equal to $X$) if $(F',\omega')$ extends to a
Riemann surface structure $F$ on the whole of $S$, equipped with a
meromorphic differential $\omega$, in such a way that $X'=X'(\omega)$
and $X''=X''(\omega)$. This is  determined by
the behaviour of $(ds^2_{\omega'},(\Ff_*(\omega'),\mu_*(\omega'),
*=h,v)$, around each point $p\in X$ (see the above local models).

The Teichm\"uller line based on $\omega_0$ lifts to the Teichm\"uller
line based on $\omega'$. This eventually leads to a $1$-parameter
family of structures $(F_\tau,\omega_\tau)$ on $S$. Each $\omega_\tau$
is a meromorphic quadratic differential on the Riemann surface
$F_\tau$. When $\tau$ varies several objects remain constant:
$X'(\omega_\tau)=X'$ and $X''(\omega_\tau)=X''$; every $p\in X'$ has
constant (pole or zero) order; the (unmeasured) foliations
$\Ff_*(\omega_\tau) = \Ff_*(\omega)$; $\mu_h(\omega_\tau) =
\mu_h(\omega)$.  On the other hand, the $\omega_\tau$-metric, the
measure $\mu_v(\omega_\tau)$ and the moduli of the poles $p\in X''$
vary with $\tau$.
\smallskip

As usual, for every Riemann surface $F$ as above, we often prefer to
consider a conformal universal covering $\Omega \to F$, so that
$F=\Omega/\Gamma$ for a group of conformal transformations of
$\Omega$. Then we can develop the above theory on $\Omega$, possibly
in a $\Gamma$-invariant way. By the uniformization theorem, $\Omega =
\C, \ \D^2, \ \mP^1(\C)$. A simple application of the Gauss-Bonnet formula
for flat metrics with conical singularities, shows that $\Omega$ is
not compact if it carries a {\it holomorphic} quadratic differential.

\section{Flat $\Qq\Dd$-spacetimes}\label{flatQD}

\subsection{Basic facts about $I^+(r)$}\label{basicPI}
We consider the Minkowski space $\mx_0$ with coordinates $(x,y,t)$ and
metric $k_0=dx^2+dy^2-dt^2$. As usual, it is oriented in such a way
that the standard basis is positive; time-oriented in such a way that
$\partial/\partial t$ is future directed.  We can assume that the
spacelike line $r$ is the $y$-coordinate line, so that $I^+(r)$
coincides with
$$\Ii = \{t>0, \ x^2 - t^2 <0\} \ . $$ Clearly, the function $\tau =
(t^2-x^2)^{1/2}$ is the {\it cosmological time} of $\Ii$, and $r$ is
its {\it initial singularity}.  The image of the Gauss map of $\Ii$
consists just of one geodesic line of $\mh^2$, i.e. $\{ x^2 - t^2 = -1,
\ y=0 \}$.
\smallskip

It is useful to perform the following change of coordinates that
establishes an immediate relationship with quadratic differentials:
$$D_0: \Pi_0 \to \mx_0$$
$$x= \tau\sh(u)\ , \ y= y \ , \ t= \tau\ch(u) \ .$$ This is an
isometry onto $\Ii$ of the open upper half-space $\Pi_0 = \{\tau >0\}$
of $\mr^3$, endowed with the metric $$g_0=\tau^2du^2 + dy^2 - d\tau^2
\ .$$ In fact $D_0$ can be considered as a developing map for a 
$(\mx_0,{\rm Isom}^+(\mx_0))$-structure on $\Pi_0$.  

\smallskip

The coordinate $\tau$ is the cosmological time of $\Pi_0$.
Note that this is real analytic, and it is also a CMC time: each
$\tau$-level surface $\Pi_0(a)$ has mean curvature equal to $1/2a$.
We will use the following notations. For every $X\subset (0,+\infty)$,
$\Pi_0(X)= \tau^{-1}(X)$. Sometimes we also use $\Pi_0(>a)$ instead of
$\Pi_0((a,+\infty))$ and so on.
\smallskip

Let us identify $\Pi_0(1)$ with the complex plane $\C$, by setting $z=
u+iy$.  We see immediately that: 

{\it The $1$-parameter family of flat
metrics on the level surfaces $\Pi(a)$, $a>0$, coincides with the one
of the Teichm\"uller line associated to the quadratic differential
$\omega_0$.}
\smallskip

Together with the image of the Gauss map, a ``dual'' object to the
initial singularity of $\Pi_0$ is just the $\omega_0$- vertical measured
foliation $(\Ff_v^0,\mu_v^0)$ on $\Pi_0(1)$.  In fact, if
$\gamma$ is an arc transverse to $\Ff_v^0$, then $\mu_v^0(\gamma)$
coincides with the length of its image via the retraction onto the
initial singularity.

\begin{remarks}\label{vari_QD}
  {\rm (1) Let us change the time orientation of $\Pi_0$ (by keeping the
spacetime one). This corresponds to the change of coordinates $u'=y,
y'=u, \tau' = -\tau$.  The new spacetime has now a {\it final
singularity}. If we change again the time orientation, but reversing the
spacetime one, we get the metric $du^2 + \tau^2dy^2 - d\tau^2$, and
the above description holds by replacing the quadratic differential
$\omega_0$ with $-\omega_0$.

(2) Let us consider now the $1$-parameter family of quadratic
differentials $s^2dz^2$, $s>0$. Then, the corresponding family of
Teichm\"uller lines of flat metrics is given by $s^2(\tau^2du^2 +
dy^2)-d\tau^2$, $\tau>0$; these metrics are given by the the pullback
of $h_0$ via the maps $\gG_s: \Pi_0 \to \Pi_0$, $\gG_s(u,y,\tau)=
(su,sy,\tau)$.}
\end{remarks}

The group ${\rm Isom}^+(\Pi_0)$ of isometries preserving both
spacetime and time orientations is generated by ``horizontal''
translations and the rotation $R_\pi$ of angle $\pi$, around the
vertical $\tau$-axis.  Hence it is canonically isomorphic to
Aut$(\C,\omega_0)$. Note that, for every $a>0$, the isometry group of
$\Pi_0(>a)$ coincides with the whole of ${\rm Isom}^+(\Pi_0)$ (just
via the restriction map).

We have the faithful representation:
$$h_0:{\rm Isom}^+(\Pi_0) \to {\rm Isom}^+(\mx_0)$$
given by

$$h_0(\sigma_v)=\left(\begin{array}{ccc}
\ch(p)& \sh(p)& 0\\
\sh(p)& \ch(p)& 0\\
0& 0& 1
\end{array}\right)+ \left(\begin{array}{c}
0\\
0\\
q
\end{array}\right)$$
$$ h_0(R_\pi)(x,y,t)=(-x,y,t) \ .$$ The representation $h_0$ can be
considered as a {\it compatible universal holonomy} for the above
developing map $D_0$. In fact for every $\xi \in {\rm Isom}^+(\Pi_0)$,
for every $p\in \Pi_0$, we have that
$$ D_0(\xi(p))= h_0(\xi)(D_0(p)) \ .$$
Moreover, the image of $h_0$ coincides with ${\rm Isom}^+(\Ii)$,
and $\Ii = D_0(\Pi_0)$. 

\subsection {Globally hyperbolic flat $\Qq\Dd$-spacetimes}\label{glob_hypQD}
These spacetimes had been already pointed out in \cite{BG}(3).  With
the notations of Section \ref{basic_QD}, let us take $S$, $X=X'\cup
X''$, $S'= S\setminus X$, related to a $(\C, {\rm
Aut}(\C,\omega_0))$-structure of meromorphic type, corresponding to
the structure $(F,\omega)$ on $S$.  This immediately induces a
$(\Pi_0, {\rm Isom}^+(\Pi_0))$-structure on $S'\times (0,+\infty)$.
The resulting flat spacetime is denoted by $Y_0(F,\omega)$.  The
cosmological time $\tau$ of $\Pi_0$ lifts to a submersion (still
denoted) $$\tau: Y_0(F,\omega)\to (0,+\infty) \ .$$ 

Set
$$\overline{S}=
S' \cup X'$$ 
$$L=L'\cup L'' = X\times (0,+\infty)\subset S\times
(0,+\infty) \ . $$ Then $Y_0(F,\omega)$ extends to a {\it cone} spacetime
$$\overline{Y}_0(F,\omega)$$ supported by $\overline{S}\times
(0,+\infty)$. Every component $\{p\}\times (0,+\infty)$, $p\in X'$, of
$L'$ corresponds to the world line of a conical singularity (a
``particle''); the cone angle coincides with the one of $p$ as a zero
or a simple pole of $\omega$. Every pole of order $2$ corresponds to a
so-called {\it peripheral end} of $Y_0(F,\omega)$ homeomorphic to
$(S^1\times \R^+)\times (0,+\infty)$.

The function $\tau$ extends to $\overline{Y}_0(F,\omega)$ and the
$1$-parameter family of spacelike metrics (with conical singularities)
on the level surfaces coincides with the Teichm\"uller line based on
$\omega$.

By extending in a natural way the notion of causal curve on
$\overline{Y}_0(F,\omega)$ (by allowing that such a curve intersects
the particle world lines), we realize that in fact $\tau$ is the {\it
cosmological time} of $\overline{Y}_0(F,\omega)$, every $\tau$-level
surface is a Cauchy surface and that $\overline{Y}_0(F,\omega)$ is
maximal globally hyperbolic.
\smallskip

\noindent Such a $\overline{Y}_0(F,\omega)$ is said to be a {\it globally
hyperbolic flat $\Qq\Dd$-spacetime}. $Y_0(F,\omega)$ is its associated
{\it non-singular} flat $\Qq\Dd$-spacetime. $(S\times (0,+\infty),
L=L'\cup L'')$ is said to be  its {\it support}.

Here is the {\it local models} for any globally hyperbolic flat
$\Qq\Dd$-spacetime:

$$\Pi_0 = \overline{Y}_0(\C,\omega_0)=Y_0(\C,\omega_0) \ ,$$

$$\overline{Y}_0(\C,z^mdz^2) \ ,$$

$$\overline{Y}_0(\C,\frac{1}{z}dz^2) \ ,$$ 
and
$$\overline{Y}_0(\C,\frac{a}{z^2}dz^2) = Y_0(\C,\frac{a}{z^2}dz^2), \ a\in
\C \ .$$
\smallskip

{\bf On the initial singularities.} As usual let us assume that
$S=\Omega$ is simply connected, that $\omega$ is defined on $\Omega$,
possibly in a $\Gamma$-invariant way, for a suitable group of
conformal automorphisms.  We want to define and describe the initial
singularity of $\overline{Y}_0(\Omega,\omega)$.  This is not so simple
to figure out as in the case of flat $\Mm\Ll$-spacetimes because 
now the developing map is not in general an embedding.
However, their intrinsic descriptions are in fact very close to each other
(see Section \ref{more:initial:sing}).    
Consider the level surface
$\SG=\overline{Y}_0(\Omega,\omega)(1)$ of the cosmological time $\tau$.
By construction, it can be naturally identified with
$(\Omega,\omega)$.  For every $p\in \SG$, let $\gamma_p$ be the past
directed ray of the integral line of the gradient of $\tau$ that
starts at $p$. For every $s<1$, set $p(s)=\gamma_p \cap
\overline{Y}_0(\Omega,\omega)(s)$. Let us say that $\gamma_p$ and
$\gamma_q$ are {\it asymptotically equivalent} if $d_s(p(s),q(s))\to
0$ when $s\to 0$, where $d_s$ denotes the distance on
$\overline{Y}_0(\Omega,\omega)(s)$ associated to its spacelike
structure.  This is in fact an equivalence relation ``$\cong$'' on
$\SG$. We define the initial singularity (as a set) to be the quotient
set
$$\Sigma = \Sigma(\overline{Y}_0(\Omega,\omega))= \SG/\cong \ .$$  
We describe now another
equivalence relation on $\SG$ via its natural identification with
$(\Omega,\omega)$. Define 
$$\delta_h(p,q)= \inf \{\mu_h(c)\}$$ where $c$ varies among the arcs
connecting $p$ and $q$ and that are piecewise contained in leaves of
the $\omega$-horizontal foliation or are transverse to it. The set
$p\cong_h q$ iff $\delta_h(p,q)=0$. The quotient set $\SG/\cong_h$ is
in fact a metric space with distance (induced by) $\delta_h$.  Note
that every non-singular leaf of $\Ff_h$ determines one equivalence
class, and the same fact holds for every connected component of the
union of the singular ones. It is immediate that the two equivalence
relations coincide on the above local models and this eventually
holds on the whole $\overline{Y}_0(\Omega,\omega)$.  Summarizing,
$(\Sigma ,\delta_h)$ is a metric space, the natural {\it
retraction} 
$$r: \SG\to \Sigma$$ is continuous; in the
$\Gamma$-invariant case, $\Gamma$ acts on $(\Sigma, \delta_h)$ by
isometries, and $r$ is $\Gamma$-equivariant.
\begin{remark}\label{QDcompact}
{\rm In the special case when $\Omega=\mh^2$, $S=\mh^2/\Gamma$ is
compact, and $\omega$ is {\it holomorphic}, $(\Sigma, \delta_h)$ is
just the {\it $\R$-tree} associated to the $\Gamma$-invariant
quadratic differential $\omega$ by Skora theory. On the other hand
there is a natural $\Gamma$-invariant measured geodesic lamination
$\lambda=(\Ll,\mu)$ on $\mh^2$ associated to $\omega$. Roughly speaking,
each non singular leaf of $\Ff_h$ has two limit points at
$S^1_\infty$, hence it determines a geodesic in $\mh^2$. The support
of $\lambda$ is given by the closure on the union of the geodesics
obtained in this way; moreover $\mu_h$ induces a transverse measure
$\mu$ so that $\mu_h$ and $\mu$ share the same dual $\R$-tree
$(\Sigma, \delta_h)$, equipped with an isometric action of $\pi_1(S)$
with small stabilizers. In fact this establishes a bijection between
holomorphic quadratic differentials and measured geodesic laminations
on $S$. Hence the $\Qq\Dd$-spacetime $\overline{Y}_0(\Omega,\omega)$
and the flat $\Mm\Ll(\mh^2)$-spacetime $\Uu^0_\lambda$ share the same
initial singularity.}
\end{remark}
  
\paragraph{ $\Pi_0$-quotients.}
Quotient spacetimes of $\Pi_0$ are the simplest but important examples
of globally hyperbolic flat $\Qq\Dd$-spacetimes.

Let $\Lambda$ be a non-trivial group of isometries acting freely and
properly discontinuously on $\Pi_0$. Then either it is generated by
one or two $\mr$-independent translations: either $\Lambda =
<\sigma_v>$, $v= p+iq$, or $\Lambda = <\sigma_{v_1}, \sigma_{v_2}>$, $v_j =
p_j +iq_j$, where we assume that $v_1,v_2$ make a positive $\R$-basis
of $\C$. Clearly $\Pi_0/\Lambda$ is homeomorphic to $S\times
(0,+\infty)$, where either $S= S^1\times \R$ or $S= S^1\times
S^1$. In fact, it is easy to see that either: 
$$\Pi_0/\Lambda = Y_0(\C,\frac{a}{z^2}dz^2)$$ where, 

$$v =|v|e^{\beta i}, \ a= |v|e^{2\beta i}$$  
or
$$\Pi_0/\Lambda = Y_0(F_\Lambda,\omega_\Lambda)=
\overline{Y}_0(F_\Lambda,\omega_\Lambda)$$ 
where $F_\Lambda$ is the
complex torus $F_\Lambda = \C/\Lambda$, and $\omega_0$ descends to the
non-singular differential $\omega_\Lambda$ on $F_\Lambda$.
Note that  the rotation $R_{\pi}$ conjugates $\Lambda$ and 
$-\Lambda$, where this last group is obtained by replacing
each generator $w$ by $-w$. 

It is well-known \cite{M, Mon} that all {\it non static} maximal
globally hyperbolic flat spacetimes with {\it toric Cauchy surfaces}
arise in this way, and that the corresponding Teichm\"uller-like space
is parametrized by the pairs $(F,\omega)=(F_\Lambda, \omega_\Lambda)$
as above (equivalently, by the groups $\Lambda$ up to conjugation by
${\rm Isom}^+(\Pi_0)$).

\subsection { General flat $\Qq\Dd$-spacetimes}\label{GenQD}
Let $(M,L=L'\cup L'')$ be formed by an oriented $3$-manifold $M$ and a
$1$-dimensional submanifold $L$ of $M$, such that every component of
$L$ either belongs to $L'$ or $L''$ and is diffeomorphic to $\R$.
\begin{defi}{\rm A {\it flat $\Qq\Dd$-spacetime with support 
$(M,L=L'\cup L'')$}
is given by a $(\Pi_0, {\rm Isom}^+(\Pi_0))$-structure on $M' =
M\setminus L$ such that the so obtained spacetime $Y_0$ satisfies the
following properties. Note that the cosmological time $\tau$ on
$\Pi_0$ lifts to a function (still denoted) $$\tau: Y_0 \to ]0,+\infty[$$
without critical points. The $\tau$-level surfaces give a foliation of
$Y_0$ by spacelike surfaces which is in fact a foliation by Riemann
surfaces endowed with non-singular quadratic differentials. We require
that $Y_0$ extends to a {\it cone} spacetime $\overline{Y}_0$ supported by
$\overline{M}= M' \cup L'$, which locally looks like a globally
hyperbolic flat $\Qq\Dd$-spacetime, respecting the $\tau$-functions.
This means in particular that for every $p\in L''$ there is a
neighbourhood $U$ of $p$ in $M$, such that $U\cap M'$ considered as an
open set of $Y$, isometrically embeds in some
$Y_0(\C,\frac{a}{z^2}dz^2)$ respecting the $\tau$-functions. Similarly
for $p\in L'$. Note that the $\tau$-function continuously extends to
the whole of $\overline{Y}$ and the same fact holds for its gradient
vector field.}
\end{defi}

The following proposition will show that flat $\Qq\Dd$-spacetimes can
realize {\it arbitrarily complicated topology and global causal
structure}. Consider any $(N,S_+,S_-,v)$, where $N$ is a compact
oriented $3$-manifold; $\partial N = S_+\cup S_-$, and each of 
$S_-$ and $S_+$ is
union of connected components of $\partial N$; $v$ is a nowhere
vanishing vector field on $N$ such that $v$ is transverse to $\partial
N$ and is ingoing (outgoing) at $S_+$ ($S_-$). The field $v$ is
considered up to homotopy through non singular vector fields
transverse to $\partial N$ (sometimes $v$ is called a {\it combing} of
$(N,S_+,S_-)$). Possibly $\partial N = \emptyset$.  Recall that such a
combing exists iff $\chi(N)-\chi(S_+)=0$ (equivalently
$\chi(N)-\chi(S_-)=0$), where $\chi(.)$ denotes the Euler
characteristic. We say that a closed $2$-sphere $\Sigma$ embedded in
the interior of $N$ {\it splits $v$} (in ``traversing'' pieces) if:

- $\Sigma$ bounds a $3$-ball $\Bb$ embedded in the interior of $N$;

- the restriction of $v$ on a neighbourhood of $\Bb$ looks like the
field $\frac{\partial}{\partial x_3}$ at a standard round ball in
$\R^3$;

- each integral line of the restriction of $v$ on 
$\hat{N}= N\setminus {\rm Int}(\Bb)$
is homeomorphic to $[0,1]$ with end-points on $\partial \hat{N}$.
\smallskip

\noindent Note that $\hat{v}=v|_{\hat{N}}$, has just one simple closed
curve $C$ of simple tangency on the $2$-sphere $\Sigma$. Hence
$\hat{v}$ is ingoing (outgoing) at one component $\Sigma_+$
($\Sigma_-$) of $\Sigma \setminus C$. Set $\hat{S}_\pm = S_\pm\cup
\Sigma_\pm$.  $\hat{v}$ is considered up homotopy through vector fields
with the same qualitative properties. A {\it tangle of orbits} for
$(\hat{N},\hat{S}_+,\hat{S}_-, \hat{v})$ is a finite union $E=E'\cup
E''$ of (generic) integral lines of $\hat{v}$ that are not tangent to
$\Sigma$.

Finally we can state
\begin{prop}\label{complicati_QDST} Let  $(N,S_+,S_-,v)$ be as above.
Then there exist a sphere $\Sigma$ that splits $v$, and a tangle of
orbits $E = E'\cup E''$ for $(\hat{N},\hat{S}_+,\hat{S}_-, \hat{v})$
such that:

(1) $(\hat{N},E=E'\cup E'')$ embeds into $(M,L=L'\cup L'')$ which
is the support of a flat $\Qq\Dd$-spacetime $\overline{Y}_0$;

(2) the gradient of the $\tau$-function of $\overline{Y}_0$ is 
transverse to (the image of) $(\hat{S}_+ \cup \hat{S}_-)\setminus E''$
and the restriction of $\hat{v}$ to $\hat{N}\setminus E''$ coincides
with the restriction of the gradient of the $\tau$-function.
\end{prop} 

We limit ourselves to a rough sketch of a proof based on the
theory of {\it branched spines} of $3$-manifolds (see \cite{BP2}). In
fact we use a variation of the construction made in Section $5$ of
\cite{BP2}(3).  For the existence of the splitting sphere $\Sigma$,
see \cite{BP2}(2).  In fact, there is a branched standard spine $P$ of
$\hat{N}$, smoothly embedded in the interior of $\hat{N}$, which
carries $\hat{v}$ as traversing normal field. We can prove that $P$
admits a {\it branched} Riemann surface structure $F$ equipped with a
meromorphic quadratic differential $\omega$, with singular set
contained the union of the open $2$-regions of $P$. In fact, we can
consider the {\it Euler cochain} $e$ (see \cite{BP2}(1)), constructed
by means of the {\it maw} vector field on $P$, which gives each
oriented $2$-region $R$ an integer value $\leq 0$. Then we can
construct $\omega$ in such a way that: $R$ contains one zero of order
$2|e(R)|-2$ if $e(R)<1$; $R$ contains one pole of order $2$ with
positive integral modulus if $e(R)=1$; there no other singularities of
$\omega$. By using $(F,\omega)$, as for the globally hyperbolic
spacetimes, we can construct a {\it branched} flat $\Qq\Dd$-spacetime
structure supported by $(P\times ]0,+\infty[, X'(\omega)\times
]0,+\infty[, X''(\omega)\times ]0,+\infty[)$. Finally, as in Section 5
of \cite{BP2}(3), we can {\it smoothly} embed $\hat{N}$ in $P\times
]0,+\infty[$ in such a way that we get, by restriction, an ordinary
spacetime structure satisfying the statement of the Proposition.

\section{$\Qq\Dd$ Wick rotation-rescaling theory}\label{WR_QD}
We are going to construct the canonical Wick rotation-rescaling theory on
$\Pi_0$. More precisely, for every $\mx = \mh^3,\mx_{\pm 1}$, for every
$*=\mh,\pm 1$ respectively, we are going to construct canonical pairs
$(D_*,h_*)$, where
$$D_*: \Pi_0(X_*) \to \mx$$
is a {\it real analytic} developing map defined on a suitable open
subset $\Pi_0(X_*)$ diffeomorphic to $\Pi_0$,
$$h_*: {\rm Isom}^+(\Pi_0) \to {\rm Isom}^+(\mx)$$ is a faithful {\it
universal holonomy} representation, such that $(D_*,h_*)$ formally
satisfies the same properties stated above for $(D_0,h_0)$. Moreover,
the pull back $g_*$ on $\Pi_0(X_*)$ via $D_*$ of the canonical metrics
$k_*$ of constant curvature on $\mx$ will be related to $g_0$ via
canonical Wick rotation-rescaling directed by the constant field
$\frac{\partial}{\partial \tau}$.  Later this shall be generalized on
arbitrary flat $\Qq\Dd$-spacetimes.

\subsection{Canonical Wick rotation on $\Pi_0$}\label{PIWR}
We use the half-space model for the hyperbolic space, $\mh^3= \mc
\times ]0,+\infty[$, with coordinates $(w=a+ib,c)$.

The {\it canonical hyperbolic developing map} is defined by
$$D_{\mh}:\Pi_0(>1) \to \mh^3$$ 
$$D_{\mh}(u,y,\tau) = (\frac{1}{\tau} \exp(u+iy),\ 
\frac{(\tau^2-1)^{1/2}}{\tau}\exp (u)) \ .$$
The {\it compatible universal holonomy} is defined by:
$$ h_{\mh}(\sigma_{(p+iq)})(w,c) = (\exp(p+iq)w,\  \exp(p)c) $$
$$ h_{\mh}(R_{\pi})(w,c)= \frac{1}{c^2+w\overline{w}}(\overline{w},c) \ .$$

It is straightforward that $D_{\mh}$ verifies the following
properties:
\begin{enumerate}
\item  It is a developing map for a hyperbolic structure on $\Pi_0(>1)$
(\emph{i.e.} it is a local diffeomorphism).

\item Its image is the open set $\mh^3 \setminus \gamma$ of $\mh^3$, where 
$\gamma = \{a=b=0\}$.

\item It sends any $\tau$-level surface $\Pi_0(a)$ onto a level surface of
the distance function, say $\Delta$, from the geodesic line $\gamma$. 

\item It analytically extends to the exponential map defined
on $\Pi_0(1)$.

\item It admits a compatible universal holonomy, say $h_{\mh}$.

\item The tautological hyperbolic metric $g_\mh$ on $\Pi_0(>1)$ which makes
$D_{\mh}$ a local isometry is obtained from the flat Lorentzian metric
$g_0|\Pi_0(>1)$ via a Wick rotation directed by the gradient of $\tau$,
with rescaling functions which only depend on the value of $\tau$.
\end{enumerate}

We claim that these properties {\it completely determine}
$(D_{\mh},h_{\mh})$.  A posteriori, we can verify that the rescaling
functions of the Wick rotation are $\alpha = \frac{1}{\tau^2-1}$, $\ \beta =
\alpha^2$, in agreement with the ones obtained for the flat
$\Mm\Ll(\mh^2)$-spacetimes, so that
$$g_{\mh} = \frac{\tau^2}{\tau^2-1}du^2+ \frac{1}{\tau^2-1}dy^2 +
 \frac{1}{(\tau^2-1)^2}d\tau^2 \ .$$ This claim can be proved via
 straightforward computations, based on the following
 considerations:
\medskip

(a) The choice of $\Pi_0(>1)$ as domain of definition of $D_{\mh}$ is
a first normalization. In fact, for every $r>0$, let us consider the
map $f_r:\Pi_0 \to \Pi_0$, $f_r(u,y,\tau)=(u,ry,r\tau)$. Clearly, it
maps $\Pi_0(\geq 1/r)$ onto $\Pi_0(\geq 1)$, and $f_r^*(h_0) =
r^2ds^2$. So the composition $D_{\mh}\circ f_r$ is a developing map
defined on $\Pi_0(\geq 1/r)$, which satisfies the above
requirements. The rescaling functions are now $\alpha =
\frac{1}{r^2(r^2\tau^2 - 1)}$, $\beta = \alpha^2$. On  $\Pi_0(>1)$ we have
agreement with the $\Mm\Ll$-spacetimes. 
\smallskip

(b) Also condition (4) is just a normalization, in order to fix
$D_{\mh}$ among all the other developing maps obtained via
post-composition with directed hyperbolic isometries that share the
geodesic $\gamma$ as axis.  Moreover, this has a clear geometric
meaning, because the exponential map is a local isometry of $\Pi_0(1)$
onto $\C \setminus \{0\}\subset \partial \mh^3$, endowed with its
Thurston's metric. 
\smallskip

(c) In order to satisfy conditions (1), (2), (3) and (5), such a map
must be of the form 
$$F:(u,y,\tau)\to (a(\tau) e^{(u+iy)},\  b(\tau) e^u)$$ 
where $a$,$b$ are positive functions. Moreover, as the map must
be defined on $\Pi_0(>1)$, we have to require that  $a(1)=1$, $b(1)=0$,
$\frac{b}{a}$ is an increasing function which tends to $+\infty$ when 
$\tau \to +\infty$.
\smallskip

(d) In order to satisfy (6), we realize that the functions
$a$ and $b$ must satisfy the condition $a'a+b'b=0$; by integrating it,
we get  $1-a^2(\tau)=b^2(\tau)$.
\smallskip

(e) Let us denote by $g$ the pull-back of the hyperbolic metric, and
by $h=h_0|\Pi_0(>1)$ the Lorentzian metric on $\Pi_0(>1)$. We have
\[
\begin{array}{lr} 
g(\frac{\partial\,}{\partial u}, \frac{\partial\,}{\partial
u})=\frac{1}{b^2} &\qquad
 h(\frac{\partial\,}{\partial u}, \frac{\partial\,}{\partial u})=\tau^2\\
\hspace{5pt} & \\
 g(\frac{\partial\,}{\partial y}, \frac{\partial\,}{\partial
y})=\frac{a^2}{b^2} &\qquad
 h(\frac{\partial\,}{\partial y}, \frac{\partial\,}{\partial y})=1 \ .
\end{array}
\]
In order to get (6), for some rescaling functions $\alpha$, $\beta$,
we finally obtain
$$a=\frac{1}{\tau} \ , \    b=\frac{(\tau^2-1)^{1/2}}{\tau}$$
$$\alpha=\frac{1}{\tau^2-1}$$
$$\beta=g(\frac{d}{d\tau}, \frac{d}{d\tau} )=
\frac{a'^2+b'^2}{b^2}=\frac{1}{(\tau^2-1)^2} \ .$$

\paragraph{Wick rotation along a ray of quadratic differentials.}
Let us consider as above the family of quadratic differentials
$s^2dz^2$, $s>0$. Set
$$\widetilde{D}_s = D_{\mh}\circ \gG_s: \Pi_0(>1)\to \mh^3 $$ 
where $\gG_s$ has been defined in Remark \ref{vari_QD}. Hence
$$\widetilde{D}_s(u+iy,\tau)=(\frac{\exp(s(u+iy))}{\tau},\ 
\frac{(\tau^2-1)^{1/2}}{\tau}\exp (su))\ .$$ So we have the developing maps
for a family of hyperbolic structures defined on $\Pi_0(>1)$. 
The associated compatible representation is given by
$$\widetilde{h}_s(\sigma_v)(w,c)=(\exp(sv)w,\ |\exp(sv)|c) \ .$$ More
precisely, we have a family of {\it marked spaces} $g_s:\Pi_0(>1)\to
\Pi_0(>1)$, where the target $\Pi_0(>1)$ is considered as fixed {\it
base space}. On the target space we have the hyperbolic structure
specified by the canonical developing map $D_{\mh}$, and we use the
marking $\gG_s$ to pull-back the structure on the source
space. Similarly, we consider on the target $\Pi_0(>1)$ the usual
Lorentzian structure, and we pull it back via the marking. We consider
these marked spaces up to {\it Teichm\"uller-like equivalence}.
Moreover, we can modify the developing map up to post-composition with
direct isometry of $\mh^3$ (or $\Pi_0$).

We want to study the limit behavior of these structures when $s\to 0$.
Let us consider the new marked spaces
$$\gG_s\circ \phi_s : \Pi_0 \to \Pi_0(>1)$$ where
$$\phi_s: \Pi_0 \to \Pi_0(>1)$$
$$ \phi_s(u+iy,\tau) = (u+iy + \frac{1}{s}\log(\frac{(s^2\tau^2 +
1)^{1/2}}{s}),\ (s^2\tau^2 + 1)^{1/2}) \ .$$ Then, we have the family of
hyperbolic developing maps on $\Pi_0$
$$ D_s = \rho_s \circ \widetilde{D}_s \circ \phi_s$$ where
$$ \rho_s(w,c)=(w-\frac{1}{s},c) \ .$$
Hence
$$ D_s(u+iy,\tau) = (\frac{\exp(s(u+iy)-1)}{s},\  \exp(su)\tau) \ .$$
The corresponding compatible universal holonomy representations
are given by
$$ h_s(\sigma_{v})(w,c)= (\exp(sv)w + \frac{\exp(sv)-1}{s},\  |\exp(sv)|c) \
.$$ 

Moreover, by a direct computation, we see that the pull back of the
hyperbolic metric $k_{\mh}$ equals:
$$D^*_s(k_{\mh})=
\frac{1}{\tau^2}((1+s^2\tau^2)du^2+dy^2)+\frac{1}{\tau^2}d\tau^2
+\frac{2s}{\tau}dudt . $$ So we can easily conclude
\begin{prop}\label{limitWR}
When $s\to 0$, then:
\smallskip 

(1) the hyperbolic developing maps $D_s$ tend (uniformly
on the compact sets) to the ``identity'' map 
$(w,c)= D^0(u+iy,\tau)= (u+iy,\tau)$;
\smallskip

(2) The holonomy representations $h_s$ tends to
$$h_0(\sigma(v))(w,c)=(w+v,c)$$
which is clearly compatible with $D^0$.
\smallskip

(3) $D_s^*(k_{\mh})\to k_{\mh}$ (uniformly on the compact sets,
and up to trivial renaming of the variables).
\end{prop}

Consider now $\Pi_0(>1)$ as a subset of a Minkowski space with
coordinates $(u,y,\tau)$ ($\tau \in ]-\infty,+\infty[$), and metric
$du^2 + dy^2 - d\tau^2$. There is an evident Wick rotation on
$\Pi_0(>1)$ directed by the vector field $\partial_\tau =
\frac{\partial}{\partial \tau}$ that converts it into the hyperbolic
space, and is equivariant for the action of ${\rm Isom}^+(\Pi_0)$
(note that it still acts by isometries on both sides of the Wick
rotation). The above discussion would suggest that also the
$1$-parameter family of Wick rotations tends to this limit one, when
$s\to 0$. We are going to see that this is the case up to suitable
re-parametrization, i.e. at the Teichm\"uller-like space level.  Let us
try with the sequence of markings used above:
$$ \Delta_s= \gG_s\circ \phi_s: \Pi_0 \to \Pi_0(>1) \ .$$
As usual denote by $g_0$ the standard Lorentzian metric on $\Pi_0$.
An easy computation shows that
$$\Delta_s^*(g_0)= 
s^2((s^2\tau^2+1)du^2 + dy^2) + 2s^3\tau^2dud\tau $$
$$\Delta_s^*(\partial_\tau)= (-(\frac{1}{s(s^2\tau^2+1)})^{1/2},\  0,\ 
\frac{(s^2\tau^2+1)^{1/2}}{s^2\tau}) \ .$$ Hence, we see that these
Lorentzian metrics degenerate when $s\to 0$, while the slope of
$\Delta_s^*(\frac{\partial}{\partial_\tau}) \to -\infty$.

However, let us consider the further markings
$$ \hat{\Delta}_s = \Delta_s \circ \gG_{\frac{1}{s}}$$
$$ \hat{\Delta}_s(u+iy,\tau)= 
(u+iy + \frac{\log((s^2\tau^2 + 1))^{1/2}}{s},\ (s^2\tau^2 + 1)^{1/2})\ .$$
Finally  we can easily check:
\begin{prop}\label{limitLOR} When $s\to 0$
$$\hat{\Delta}_s^*(g_0)= (s^2\tau^2+1)du^2 + dy^2 +(s^2\tau^2+
\frac{s^2-1}{s^2\tau^2+1})d\tau^2 + 2s^2\tau^2dud\tau$$ converges to the
usual Minkowski metric $du^2+dy^2 -d\tau^2 $.  Moreover, the
action of the group ${\rm Isom}^+(\Pi_0)$ on $\Pi_0$ is isometric for every
$\hat{\Delta}_s^*(g_0)$, $s\geq 0$.
\end{prop}

\subsection{Canonical dS rescaling on $\Pi_0$}\label{dSPI}
Like for $\Mm\Ll(\mh^2)$ spacetimes, the construction of $D_{1}$ is
somewhat dual to the one of $D_{\mh}$. Let us identify the half-space
model of $\mh^3$ used before with the usual hyperboloid model 
embedded in the Minkowski space $\mm^4$ (with metric $-dx_0^2 +
\sum_{j=1,\ldots,3} dx_j^2$) in such a way that the geodesic $\gamma$
becomes the intersection of $\mh^3$ with $\{x_2=x_3=0\}$, and the
point $(0,0,1)$ becomes the point $(1,0,0,0)\in \mh^3$. Thus the
developing map $D_\mh: \Pi_0(>1)\to \mh^3$ becomes
$$ D_\mh(u+iy,\tau) = \ch(\Delta(\tau))\left(\begin{array}{c}
(c(\tau)+1)/2c\\
(c(\tau)-1)/2c\\
0\\
0
\end{array}\right)+ \sh(\Delta(\tau)\left(\begin{array}{c}
0\\
0\\
u/||(u,y)||\\
y/||(u,y)||
\end{array}\right)$$
where 
$$\Delta(\tau)={\rm arctgh}(\frac{1}{\tau})$$
$$c=\frac{(\tau^2-1)^{1/2}}{\tau}$$
$$||(u,y)|| = (u^2+y^2)^{1/2} \ . $$
By taking $[D_\mh]$ we immediately get such a developing map with values
in the projective model of $\mh^3$. Note that every vertical line in
$\Pi_0(>1)$ parameterizes a geodesic ray $\gamma_{u+iy}(\tau)$ in $\mh^3$.

The universal holonomy representation transforms in
$$h_\mh(\sigma_v)= \left(\begin{array}{cccc} \ch(p)& \sh(p)& 0&0\\
\sh(p)& \ch(p)& 0&0\\ 0& 0& \cos(q)&\sin(q)\\ 0&0&-\sin(q)&\cos(q)
\end{array}\right) \ .$$

Thus $$D_{1}:\Pi_0(<1)\to \mx_{1}$$
is defined in such a way that every vertical line in $\Pi_0(<1)$
parameterizes the geodesic arc  $\gamma^*_{u+iy}(\tau)$ in $\mx_1$
that is dual to  $\gamma_{u+iy}(\tau)$ and shares the same pair of points
at infinity. Precisely,
$$D_1(u+iy,\tau)=[\ch(\Delta^*(\tau))\left(\begin{array}{c}
0\\
0\\
u/||(u,y)||\\
y/||(u,y)||
\end{array}\right)+ \sh(\Delta^*(\tau))\left(\begin{array}{c}
(c(\tau)+1)/2c\\
(c(\tau)-1)/2c\\
0\\
0
\end{array}\right)]$$
$$  \Delta^*(\tau)= {\rm arctgh(\tau)} \ .$$
In fact, this last is the cosmological time of the so obtained spacetime
of constant curvature $\kappa = 1$. 

We have the coincidence of the universal holonomies $h_1=h_\mh $.

It is easy to verify that
$$g_1=D_1^*(k_1)= 
\frac{\tau^2}{1-\tau^2}du^2 + \frac{1}{1-\tau^2}dy^2 - 
\frac{1}{(1-\tau^2)^2}d\tau^2 \ . $$

The behaviour of $(D_\mh, h_\mh)$ straightforwardly dualizes to 
$(D_1, h_1)$.

\subsection{Canonical AdS rescaling on $\Pi_0$}\label{AdSPI}
Let us describe now the canonical AdS pair $(D_{-1},h_{-1})$.
Let us fix an ordered pair $(l,l')$ of spacelike geodesic lines in
$\mx_{-1}$ that are dual each other. Let us denote by $x_{\pm}$ and
$x'_{\pm}$ the respective endpoints. This determines four segments on
$\partial \mx_{-1}$ with endpoints $x_*,x'_*$.  These segments belong
to suitable null lines for the natural conformal Lorentzian structure
on the boundary of $\mx_{-1}$. Denote by $C$ the curve given by the
union of these four segments, and $\Kk(C)$ its convex hull in the
projective space.  Abstractly, we can look at $\Kk(C)$ as a closed {\it
oriented} positively embedded into the closure of $\mx_{-1}$ in the
projective space. The curve $C= \Kk(C)\cap \partial \mx_{-1}$, $\partial
\Kk(C)$ is made by $4$ triangular faces contained in
four distinct null-planes. Let us order its four vertices as $
x_-,x_+,x'_-,x'_+$. This induces another orientation of $\Kk(C)$.  We
stipulate that the two orientations coincide, and that $l'$ is in the
future of $l$ with respect to the time orientation of $\Kk(C)$. Fix two
interior points $x_0\in l$ and $x_0'\in l'$ respectively. Clearly the
dual plane $P(x_0)$ contains $l'$, while $l\subset P(x'_0)$. The
time-like geodesic orthogonal to $P(x'_0)$ at $x_0$ is orthogonal also
to $P(x_0)$ at $x'_0$. It is easy to see that given two patterns of
data $D(i)=(l_i,l'_i,x(i)_\pm,x'(i)_\pm,x(i)_0,x'(i)_0)$, $i=1,2$, as
above, there is a unique isometry $f\in PSL(2,\R)\times PSL(2,\R)$ of
$\mx_{-1}$ such that $f(D(1))=D(2)$. Thus we fix a choice by setting:
$$ l(t)= [\ch (t)\left(\begin{array}{cc}
1 & 0\\
0 & 1
\end{array}\right)+ \sh (t) \left(\begin{array}{cc}
1 &  0\\
0 & -1
\end{array}\right)]= [\left(\begin{array}{cc}
\exp(t) & 0\\
0 & \exp(-t)
\end{array}\right)] \ ,$$
\smallskip

$$ l'(t)= [\left(\begin{array}{cc}
0&\exp(t)\\
-\exp(-t)&0
\end{array}\right)],\  t\in \R $$
\smallskip

$$x_0 = [\left(\begin{array}{cc}
1 & 0\\
0 & 1
\end{array}\right)]\ , \ \ x'_0 = [\left(\begin{array}{cc}
0& 1\\
-1 & 0
\end{array}\right)]$$
where we have also specified a lifting in $\hat{\mx}_{-1}=SL(2,\R)$.
Denote $\Pi_{-1}$ the interior of $\Kk(C)$.
\smallskip

\noindent  
$D_{-1}$ is the embedding of $\Pi_0$ onto $\Pi_{-1}$, defined 
as follows
$$ D_{-1}(u,y,\tau)= [\cos(t) \left(\begin{array}{cc}
\exp(u) & 0\\
0 & \exp(-u)
\end{array}\right)+
\sin(t)\left(\begin{array}{cc}
0& \exp(y)\\
-\exp(-y) & 0
\end{array}\right)]$$

$$t=  {\rm arctan}(\tau) \ .$$
We can straightforwardly verify that:
\smallskip

(1) $t= {\rm arctan}(\tau)$ is the  cosmological time  of
$\Pi_{-1}$. The spacelike geodesic $l$ is its initial singularity.

\smallskip

(2) The Lorentzian metric $g_{-1}$ on $\Pi_0$, of constant curvature $-1$,
which makes $D_{-1}$ an isometry is a rescaling of the flat metric $g_0$,
directed by the gradient of $\tau$, with rescaling functions $\alpha =
\frac{1}{1+\tau^2}$, $\beta = \alpha^2$, that is
$$g_{-1} = \frac{\tau^2}{1+\tau^2}du^2+ \frac{1}{1+\tau^2}dy^2 +
 \frac{1}{(1+\tau^2)^2}d\tau^2 \ .$$ 
\medskip 
  
As universal  holonomy we simply have
$$ h_{-1}(\sigma_{p+iq})=$$ 
$$= \left(\left(\begin{array}{cc}
\exp(p-q) & 0\\
0 & \exp(q-p)
\end{array}\right),\left(\begin{array}{cc}
\exp(p+q) &  0\\
0 & \exp(-(p+q))
\end{array}\right)\right)\ ,$$
\smallskip

$$ h_{-1}(R_{\pi}) = \left(\left(\begin{array}{cc}
0& 1\\
-1 & 0
\end{array}\right),\left(\begin{array}{cc}
0&  1\\
-1 & 0
\end{array}\right)\right)\\ $$
where we have again specified a lifting in $SL(2,\R)\times SL(2,\R)$.

\begin{remark}\label{tutto_Pi-1}
{\rm The behaviour of this AdS rescaling is quite ``degenerate''
with respect to the $\Mm\Ll$-spacetimes. For in that case
the convex hull of the boundary line at infinity of the level surfaces
of the cosmological time produced a sort of AdS-convex core, {\it
strictly} contained in the AdS spacetimes. Here the convex hull coincides
with the whole of the spacetime $\Pi_{-1}$.
}
\end{remark}
\paragraph{AdS rescaling along a ray of quadratic differentials.}
Similarly to the Wick rotation case treated above, let us consider the
family of developing maps
$$ D_s=D_{-1}\circ \gG_s \ .$$
For the corresponding holonomy representations we have, in particular
$$ h_s(\sigma_{p+iq})=$$ 
$$= \left(\left(\begin{array}{cc}
\exp(s(p-q)) & 0\\
0 & \exp(s(q-p))
\end{array}\right),\left(\begin{array}{cc}
\exp(s(p+q)) &  0\\
0 & \exp(-s(p+q)
\end{array}\right)\right)\ .$$
Let us conjugate them by $\rho_s=(A_s,A_s)\in PSL(2,\R)\times
PSL(2,\R)$, where
$$A_s = \left(\begin{array}{cc}
1 & 1/2s\\
0 & 1
\end{array}\right) \ ,$$
obtaining the family of compatible pairs 
$$(\rho_sD_s,\ \rho_s h_s\rho_s^{-1}) \ .$$
As in the Wick rotation case, it is easy to verify that, for $s\to 0$,
we have a nice convergence to the parabolic representation:
$$ \rho_s h_s(\sigma_{p+iq})\rho_s^{-1}\to
\left(\left(\begin{array}{cc}
1& p-q\\
0 & 1
\end{array}\right),\left(\begin{array}{cc}
1&  p+q\\
0 &  1
\end{array}\right)\right)\ .$$
On the other, in the present case, the {\it images} of the developing
maps $\rho_sD_s$ degenerate to the null plane passing through the
common fixed point of those parabolic transformations.
 
\subsection{Wick rotation-rescaling on $\Qq\Dd$-spacetimes}\label{WR_on QD}
We have constructed above the canonical developing maps $D_*$ for the
Wick rotation-rescaling theory on $\Pi_0$. Given any flat $\Qq\Dd$-spacetime
$\overline{Y}_0$, by using the maps $D_*$, we can immediately
construct the corresponding Wick rotation-rescaling on the associated non
singular spacetime $Y_0$. In this way we obtain a hyperbolic
$3$-manifold $M_Y$, and spacetimes $Y_\kappa$, of constant curvature
$\kappa = \pm 1$, all diffeomorphic to $Y$.  It is easy to see that
also $Y_\kappa$ completes to a cone spacetime $\overline{Y}_\kappa$
homeomorphic to $\overline{Y}_0$, with the same cone angles along the
singular world lines. They are all said $\Qq\Dd$-{\it spacetimes}.  We are
going to show some instances of these constructions.
\medskip

\paragraph{Wick rotation on spacetimes with toric Cauchy surfaces.}  
Consider a group of isometries of $\Pi_0$ of the form $\Lambda =
<\sigma_{v_1}, \sigma_{v_2}>$, $v_j = p_j +iq_j$, where we assume that
$v_1,v_2$ make a positive $\R$-basis of $\C$. Apply the Wick rotation
on $\Pi_0(>1)/\Lambda$. This gives a {\it non-complete} hyperbolic
structure on $(S^1\times S^1)\times \R^+$. In order to study its
completion we can apply verbatim the discussion occurring in the proof
of Thurston's {\it hyperbolic Dehn filling theorem} (see \cite{Thu} or
\cite{BP}).  Recall that the feature of the completion is determined
by the unique real solutions, say $(r,s)$, of the equation
$$ r(p_1+iq_1) + s(p_2+iq_2)= 2\pi i \ .$$ In particular, if $r/s \in
\mq$, then the completion is homeomorphic to the interior of a tube
$\D^2\times S^1$. Notice that it happens exactly when the leaves of
the vertical foliation of the quadratic differential $\omega_\Lambda$
are parallel {\it simple closed curves}. If $l$ is any such a leaf,
the completion of the image, via $D_{\mh}$, of the annulus $l \times
\{\tau >1\}$ is a meridian disk of the tube. The core of the tube ({\it i.e.}
the line made by the centers of those disks, when $l$ varies) is a
simple closed geodesic. Along this core, there is a conical
singularity of cone angle $2\pi/\theta$, where $(r,s)=\theta (r',s')$,
$r'/s'\in \mq$, and $r'$, $s'$ are co-prime.  

If $\Lambda$ corresponds to the structure $(F,\omega)$ on 
the torus, then the $1$-parameter family of spacetimes $\Pi/r\Lambda$, $r>0$,
corresponds to the family of pairs $(F, r^2\omega)$.  Hence, we get a
family of hyperbolic structures on $(S^1\times S^1)\times
]0,+\infty[$, and by completion a family of hyperbolic cone manifold
structures on the tube $\D^2\times S^1$.  The cone angles are
$r2\pi/\theta$, so they tend to $0$ when $r\to 0$, {\it opening in a
cusp}.

\paragraph{BTZ black-holes.} Note that for every quotient spacetime
$\Pi_0/\Lambda$, the corresponding AdS $\Qq\Dd$-spacetime is just the
quotient $\Pi_{-1}/h_{-1}(\Lambda)$. Recall that $\Lambda$ is
conjugate to $-\Lambda$. By reversing the time orientation we get the
spacetime $\Pi_{-1}/h_{-1}(\widetilde{\Lambda})$, where
$\widetilde{\Lambda}$ is obtained by replacing every generator $w=
a+ib$ with $\widetilde{w}= b+ia$. If we allow also isometries that
change the spacetime orientation we can get
$\Pi_{-1}/h_{-1}(\overline{\Lambda})$, obtained by replacing every
generator $w= a+ib$ with $\overline{w}= a-ib$.  Let us consider the
simplest case $\Lambda = <\sigma_{v}>$, $v=p+iq$. Then, up to
conjugation and orientation reversing we can (and we will) assume that
$p\geq q \geq 0$.

Recall that $\Pi_0/\Lambda$ is equal to $Y_0(\C,\frac{a}{z^2}dz^2)$,
$a=|v|e^{2\beta i},\ v=|v|e^{\beta i}$. By varying $p+iq$, we get the
family of so-called {\it regions of type II} of the BTZ {\it black
holes}. We refer to \cite{BTZ}, \cite{Ca} for a detailed description
of these spacetimes. Here we limit ourselves to a few remarks.
\smallskip

{\it The Kerr-like black hole metric on $\Pi_{-1}/h_{-1}(\Lambda)$.}
The BTZ solution of Einstein field equations in three spacetime
dimensions shares many characteristics of the classical $(3+1)$ Kerr
rotating black hole and this is its main reason of interest.  Recall
that the (2+1) Kerr-like metric in coordinates $(r,v,\phi)$ is of the
form
$$ds^2 = -fdv^2+f^{-1}dr^2 +r^2(d\phi - \frac{J}{2r^2}dv)^2 $$
where, with the usual notations of ADM approach to gravity, we have:
$$f=(\ort{N})^2= -M+r^2 + \frac{J^2}{4r^2}, \ N^{\phi} = - \frac{J}{2r^2}$$
where $M$ and $J$ are constant. Here we assume that the ``mass'' 
$M=r_+^2+r_-^2>0$, and that the ``angular momentum''
$J=2r_+r_- $. Note that $M-|J|\geq 0$; moreover $M$ and $J$
determine $r_{\pm}$ up to simultaneous change of sign, and we
stipulate that $r_{+} \geq r_-$. From now on we also assume we are in
the {\it generic} case so that $r_+>r_-\geq 0$.  In fact the
coordinates $(r,\phi)$, $r>0$, should be considered as ``polar
coordinates'' (ie $\phi$ is periodic) on the $v$-level surfaces, so
that the topological support of the metric should be homeomorphic to
$S^1\times \R^2$.  Note that this metric is singular at $r=0$ and
$r=r_\pm$; otherwise it is non-singular, and a direct computation
shows that it is of constant curvature $\kappa = -1$. Note that we can
rewrite the Kerr metric in the form
$$ds^2 = (M-r^2)dv^2 + f^{-1}dr^2 + r^2 d\phi^2 -Jdvd\phi$$
hence $\frac{\partial}{\partial v}$ is timelike for $r> M^{1/2}$.
Under our assumptions, $r_+<M^{1/2}$ and $r$ is timelike on $]r_-,r_+[$.

Let us go back to a quotient spacetime $\Pi_0/\Lambda$ homeomorphic to
$S^1\times \R^2$, and we assume we are in the {\it generic} case $p>q$; set
$p=r_+$, $q=r_-$.  We want to point out a rather simple
re-parametrization of $\Pi_0$ (depending on $r_\pm$)
$$(u,y,\tau)= F(r,\phi,v), \ F=F_{r_\pm}$$  
$$r\in ]r_-,r_+[=]q,p[,\ \phi\in \R, v\in \R$$ such that the pull back
$(D_{-1}\circ F)^*(k_{-1})$ of the AdS metric on $\Pi_{-1}$ gives us
the above Kerr-like metric on the slab $\{r_- < r <r_+\}$ and this
passes to the quotient spacetime $\Pi_0/\Lambda$.

We know that
$$g_{-1}=D_{-1}^*(k_{-1})= \alpha(\tau)(\tau^2du^2 +
dy^2)-\alpha^2(\tau)d\tau^2$$
$$\alpha(\tau)= \frac{1}{1+\tau^2} \ .$$ Recall that $t={\rm
arctan}(\tau)$ is the cosmological time of $\Pi_{-1}$.
The constant spacelike vector field on $(\Pi_0,g_{-1})$ 
$$\xi(u,y,\tau)= r_+\frac{\partial}{\partial u} + 
r_-\frac{\partial}{\partial y}$$
is the infinitesimal generator of $\Lambda$. The
AdS quadratic form $q$ on this field produces the function

$$r^2=q(\xi)=\frac{\tau^2r_+^2+r_-^2}{1+\tau^2}= \sin^2(t)r^2_+ +
\cos^2(t)r_-^2 = $$

$$= \frac{r^2 - r^2_-}{r^2_+ - r^2_-}r^2_+ + \frac{r^2_+ - r^2}{r^2_+ -
r^2_-}r^2_- , \ $$

$$ \tau^2 = \frac{r^2-r_-^2}{r^2_+ - r^2} \ .$$ Note that
$$ \frac{d}{dt}r^2 = \sin(t)\cos(t)(r^2_+ - r^2_-)>0 \ .$$ Hence we
can consider the coordinate transformation $(u,y,\tau)=F(r,\phi,v)$
given by
$$r=(q(\xi)(\tau))^{1/2}>0, \ \ \tau = (\frac{r^2-r^2_-}{r^2_+ - r^2})^{1/2}$$ 
$$ u+iy = (r_+ + ir_-)(\phi + iv) \ .$$ Note that this includes a
positive re-parametrization of the canonical time of the spacetime
$(\Pi_0,g_{-1})$, i.e. of $\Pi_{-1}$. This function is defined
also on the closure of $\Pi_{-1}$ in $\mx_{-1}$; precisely it takes
the value $r_+$ on the two null faces that contain the line $l'$; the
value $r_-$ at the line $l$. Finally, a simple computation shows that
\begin{lem}
$(D_{-1}\circ F)^*(k_{-1})= -fdv^2+f^{-1}dr^2 +r^2(d\phi -
\frac{J}{2r^2}dv)^2 $.
\end{lem}
Consider the other canonical pairs $(D_*,h_*)$ of the Wick rotation-rescaling
theory on $\Pi_0$. We can consider the spaces
$$(\Pi_0,(D_*\circ F)^*(k_{*}))\ .$$ It is evident by construction that
these metrics are related to Kerr-like one either by natural rescaling
or Wick rotation, directed by the vector field
$\frac{\partial}{\partial r}$. In fact the rescaling functions are the
usual ones, once we consider them as functions of $r$, via
$\tau=\tau(r)$.  Moreover all of this is $\Lambda$-invariant.
In particular
$$(D_\mh\circ F)^*(k_\mh)= 
W_{(\frac{\partial}{\partial r},\gamma^2,\gamma)}((D_{-1}\circ F)^*(k_{-1}))$$
where
$$\gamma(r) = \frac{1+\tau^2}{\tau^2-1} = \frac{r_+^2-r_-^2}{2r^2-M}\ .$$

{\it On the extreme cases.}  Let us consider the {\it critical} case
$r_+ = r_- = r_0>0$. In this case $q(\xi)$ is the constant
$1$-function. The above maps $(u,y,\tau)= F_{r_\pm}(r,\phi,v)$
degenerate when $r_+\to r_-$. Moreover, the action of $\Lambda$ on
$\Pi_{-1}$, even considered up to diffeomorphism, has a very different
dynamic with respect to the generic cases, for the line $l'$ is {\it
point-wise} fixed.  On the other hand, both the Kerr-like metric
and $(D_{-1}\circ \hat{F}_{r_\pm})^*(k_{-1})$
are well defined also for $r_+=r_-=r_0$ where
$$(u,y,\tau)=\hat{F}_{r_\pm}(\tau,\phi,v)=((r_+ + ir_-)(\phi +
iv),\tau) \ .$$ In fact there are rather complicated coordinate
transformations (defined on different patches of $\Pi_0$) that
transform each metric in the other.  These can be effectively computed
by using $D_{-1}$, the so called ``Poincar\'e coordinates'' on $\Pi_{-1}$
and (3.11), (3.34)-(3.37) of \cite{BTZ}.

The Kerr-like metric makes sense also when $M=0$ (hence $r_\pm=0$) as it
becomes:
$$ -r^2dv^2 + r^{-2}dr^2 + r^2d\phi^2 \ , r\neq 0 \ .$$
By setting $r=z^{-1}$ we get
$$\frac{dz^2+d\phi^2-dv^2}{z^2}$$
that is the ``Poincar\'e coordinates on the future of a suitable null plane.
This agrees with the above discussion on the behaviour of the AdS rescaling
along ray of quadratic differentials.

\smallskip 

{\it On the maximal spacetime containing the BTZ black hole.}
$\Pi_{-1}/h_{-1}(\Lambda)$ is only a region of a bigger (non globally
hyperbolic) spacetime $B(r_\pm)/h_{-1}(\Lambda)$ of constant curvature
$\kappa = -1$ that actually contains the BTZ black hole. We want to
briefly describe $B(r_\pm)$. Recall that a lifted copy of $\mx_{-1}$
in $\hat{\mx}_{-1}=SL(2,\R)$ is given by the matrices of the form
$$X=\left(\begin{array}{cc}
T_1+X_1& T_2+X_2\\
-T_2+X_2& T_1-X_1
\end{array}\right)$$ 
such that ${\rm det}(X)=1$, $0<T_1^2 -X_1^2 <1$, $X_1, T_1$ have a
definite sign. In fact, in defining $D_{-1}$ we have also specified
such a lifting over $\Pi_{-1}$.  The group $h_{-1}(\Lambda)$ acts on
the whole of $\mx_{-1}$, again with the constant vector field $\xi$ as
infinitesimal generator.  Hence the function
$$q(\xi) =(T_2^2-X_2^2)r_+ +  (T_1^2 - X^2_1)r_-$$ makes sense on the
whole AdS spacetime. Roughly speaking, $B(r_\pm)$ is the {\it maximal}
region of $\mx_{-1}$ such that:
\smallskip

(1) $\Pi_{-1} \subset B(r_\pm)$;
\smallskip

(2) $q(\xi)>0$ on $B(r_\pm)$, so that we can take the function
    $r=q(\xi)^{1/2}>0$;
\smallskip

(3)$B(r_\pm)$ is $h_{-1}(\Lambda)$-invariant, the group acts
nicely on $B(r_\pm)$ and the quotient spacetime does not contain
closed timelike curves ({\it causality} condition).
\medskip

In fact (see \cite{BTZ}) such a $B(r_\pm)$ admits a
$h_{-1}(\Lambda)$-invariant ``tiling'' by regions of three types I,
II, III contained in $\{r>r_+\}$, $\{r_+>r>r_-\}$, $\{r_->r\}$
respectively.  Each region supports the above Kerr-like metric, is
bounded by portions of null-planes at which $r=r_\pm$; in particular
$\Pi_{-1}$ itself is a region. Moreover, at $\{r=r_\pm\}$ there are
only ``coordinate singularities''.

\begin{remark} \label{multiBH}
{\rm Every globally hyperbolic AdS $\Qq\Dd$ spacetime
$\overline{Y}_{-1}(F,\omega)$ contains a {\it peripheral end}
corresponding to each pole of $\omega$ of order 2. Every peripheral
end is homeomorphic to $(S^1\times \R)\times \R$ and can be
isometrically embedded in a suitable BTZ region of type II as above,
by respecting the canonical time. We say that such an end is {\it
static} if the leaves of the horizontal $\omega$-foliation at the pole
are simple closed curves; in such a case $J=0$. Otherwise it is {\it
``rotating''}. Note that portions of static BTZ regions of type II
also occur in AdS $\Mm\Ll$-spacetime associated to some
$(\lambda,\Gamma)$ such that $H/\Gamma$ has at least one closed
boundary component that is also an isolated ($+\infty$ weighted) leaf
of the lamination (for example, the $T$-symmetric spacetimes to the
$\Mm\Ll(\mh^2)$-spacetimes breaking the $T$-symmetry in Section
\ref{3cusp} were of this kind). If such a component is not isolated we
can say that it is a {\it rotating $\Mm\Ll$ end}.  Similarly to
$B(r_\pm)/h_{-1}(\Lambda)$ above with respect to a BTZ region of type
II, it should be interesting to investigate maximal causal (non
globally hyperbolic) AdS extension of any such a space with static or
rotating ends, that would contain a {\it ``multi black hole''}.  This
kind of situations are studied for example in \cite{Ba}(3) and
\cite{BSK}.  }
\end{remark} 

\paragraph{$T$-symmetry.}
It follows from remark \ref{tutto_Pi-1} that the rescaling of any
globally hyperbolic flat spacetime  $\overline{Y}_{0}(S,\omega)$ produces
the whole of the associated AdS one  $\overline{Y}_{-1}(S,\omega)$.
The level surface $\overline{Y}_0(S,\omega)(1)$ transforms in the
``middle'' surface $\overline{Y}_{-1}(S,\omega)(\pi/4)$. If $S$ is compact,
this is the one of largest area.

These AdS spacetimes are closed for the $T$-symmetry. In fact,
by inverting the time orientation we simply get 
$$\overline{Y}_{-1}(S,\omega)^* = \overline{Y}_{-1}(S,-\omega) \ .$$


\chapter{Complements}\label{COMP}
\section{Moving along a ray of laminations}\label{der}
Let us fix $(\lambda, \Gamma)\in\Mm\Ll^\Ee$ and put $F=\mh^2/\Gamma$.
The ray of ($\Gamma$-invariant) measured laminations determined by
$\lambda$ is given by $t\lambda = (\Ll,t\mu)$, $t\geq 0$.  So we have
$1$-parameter families of spacetimes $\hat\Uu^\kappa_{t\lambda}$, of
constant curvature $\kappa\in\{0,1,-1\}$, diffeomorphic to
$F\times\mr_+$, having as universal covering
$\Uu^\kappa_{t\lambda}$. We have also a family of hyperbolic
$3$-manifolds $M_{t\lambda}$, obtained via the canonical Wick
rotation.  $\hat\Pp_{t\lambda}$ is contained in
$\hat\Uu^{-1}_{t\lambda}$ and is the image of the canonical rescaling
of $\hat\Uu^0_{t\lambda}$. Its universal covering is
$\Pp_{t\lambda}\subset \Uu^{-1}_{t\lambda}$.

First, we want to (give a sense and) study the ``derivatives'' at $t=0$
of the spacetimes $\hat\Uu^\kappa_{t\lambda}$, of their holonomies
and ``spectra'' (see below).

\subsection{Derivatives of spacetimes at $t=0$}\label{derivateST}
Let $$\frac{1}{t}\hat\Uu^\kappa_{t\lambda}$$ be the spacetime of
constant curvature $t^2\kappa$ obtained by rescaling the Lorentzian
metric of $\hat\Uu^\kappa_{t\lambda}$ by the constant factor
$1/t^2$. We want to study the limit when $t\to 0$. For the present
discussion it is important to recall that all these spacetimes are
well defined only up a Teichm\"uller-like equivalence relation. So we
have to give a bit of precision on this point. Fix a base copy of
$F\times \mr_+$ and let
\[
\varphi:F\times\mr_+\rightarrow\hat\Uu_\lambda^0
\]
be a 
{\it marked spacetime} representing the equivalence class of 
$\hat\Uu_\lambda^0$. Denote by $k_0$ the flat Lorentzian
metric lifted on $F\times \mr_+$ via $\varphi$.
A developing map with respect to such a metric is a diffeomorphism
\[
   D:\tilde F\times\mr_+\rightarrow\Uu_\lambda^0\subset\mx_{0} \ .
\] 
Up to translation, we can suppose $0\in\Uu^0_\lambda$.
Notice that, for every $s>0$, the map
\[
   g_s:\Uu_\lambda^0\ni z\mapsto sz\in\mx_{0}
\]
is a diffeomorphism onto $\Uu_{s\lambda}^0$. Moreover, it is
$\Gamma$-equivariant, where $\Gamma$ is supposed to act on
$\Uu_\lambda^0$ (resp. $\Uu_{s\lambda}^0$) via $h^0_\lambda$
(resp. $h^0_{s\lambda}$) as established in Chapter \ref{FGHST}.  
Thus $g_s$ induces to the quotient a
diffeomorphism
\[
   \hat g_s:\hat \Uu_\lambda^0\rightarrow\hat \Uu_{s\lambda}^0
\]
such that the pull-back of the metric is simply obtained by
multiplying the metric on $\hat\Uu_\lambda^0$ by a factor $s^2$.

Thus the metric $k_s=s^2k$ makes $F\times\mr_+$ isometric to
$\hat\Uu^0_{s\lambda}$.
We want to prove now a similar result for $\kappa=\pm 1$.

The cosmological time of $(F\times\mr_+,k_s)$ is
$T_s=s\tau$, where $T$ is the cosmological time of
$(F\times\mr_+,k_0)$. It follows that the gradient with respect to $k_s$
of $\tau_s$ does not depend on $s$ and we denote by $X$ this field.
Now suppose $\kappa=-1$ and denote by $h_s$ the metric obtained by
rescaling $k_s$ around $X$ with rescaling functions
\[
   \begin{array}{ll}
   \alpha=\frac{1}{1+T_s^2} \ ,\ \  & \beta=\frac{1}{(1+T_s^2)^2} \ .
   \end{array}
\]  
We know that $(F\times\mr_+,h_s)$ is isometric to $\Pp_{s\lambda}$.
Moreover the metric $h_s/s^2$ is obtained by a rescaling of the metric $k_0$
along $X$ by rescaling functions
\[
 \begin{array}{ll}
    \alpha=\frac{1}{1+s^2 T^2} \  \ , \ \ & \beta=\frac{1}{(1+s^2 T^2)^2}.
   \end{array}
\]  
Thus, we obtain $\lim_{s\rightarrow 0} h_s/s^2=k$.

Finally suppose $\kappa=1$, then we can define a metric $h'_s$ on the
subset $\Omega_s$ of $F\times\mr_+$ of points
$\{x|T_s(x)<1\}=\{x|T<1/s\}$ such that $(\Omega_s,
h'_s)=\hat\Uu^1_\lambda$. In fact we can set $h'_s$ to be the metric
obtained by rescaling $k_s$ by rescaling functions
\[
   \begin{array}{ll}
   \beta=\frac{1}{(1-T_s^2)^2}, & \alpha=\frac{1}{1-T_s^2}.
   \end{array}
\]
Choose a continuous family of embeddings $u_s:F\times\mr_+\rightarrow
F\times\mr_+$ such that
\begin{enumerate}
\item
   $u_s(F\times\mr_+)=\Omega_s$;
\item
 $  u_s(x)=x \textrm{ if } T(x)<\frac{1}{2s}$.
\end{enumerate}
Then the family of metrics  $h_s=u_s^*(h'_s)$ works.
\medskip

We can summarize the so obtained results as follows:

\begin{prop}\label{lim:space}
For every $\kappa=0,\pm 1$, 

\[ \lim_{t\rightarrow 0}\ \frac{1}{t} \Uu^\kappa_{t\lambda}=  \Uu^0_\lambda
\,.\]
\end{prop}
\cvd

\noindent 
For $\kappa = 0$ we have indeed the strongest fact that
for every $t>0$
\[\frac{1}{t}\Uu^0_{t\lambda} =  \Uu^0_\lambda \ . \]
Note that this convergence is in fact like a convergence
of pointed-spaces; for example, the convergence of 
spacetimes $\frac{1}{t}\Uu^{-1}_{t\lambda}$ only concerns the
past side of them, while the future sides simply disappear. 

\subsection{Derivatives of representations}\label{derivateREP}
For any $\kappa \in\{0,1,-1\}$ the set of holonomies
of $\ \hat\Uu_{t\lambda}^\kappa$ gives rise to continuous families of
representations
\[
    h_t^\kappa:\Gamma\rightarrow\ISO(\mx_{\kappa})
\]

We compute the derivative of such families at $t=0$.  The
following lemma contains the formula we need.  In fact this lemma is
proved in~\cite{Ep-M, McM}, and we limit ourselves to a sketch
of proof.

\begin{lem}\label{var:der:lem}
Let $\lambda$ be a {\rm complex-valued} measured geodesic lamination 
on a straight convex set, and denote by
$E_\lambda$ the Epstein-Marden {\rm bending-quake} cocycle.
Fix two points $x,y\in\mh^2$ then the function
\[
    u_\lambda:\mc\ni z\mapsto E_{z\lambda}(x,y)\in PSL(2,\mc)
\]
is holomorphic. Moreover, if $\lambda_n\rightarrow\lambda$ on a
neighbourhood of $[x,y]$, then $u_{\lambda_n}\rightarrow u_\lambda$ in
the space of holomorphic functions of $\mc$ with values in $
PSL(2,\mc)$.
\end{lem}
\Dim The statement is obvious when $\lambda$ is a finite lamination.
On the other hand, for every $\lambda$ there exists a sequence of
standard approximations $\lambda_n$. Now it is not hard to see that
$u_{\lambda_n}$ converges to $u_\lambda$ in the compact-open topology
of $\mathrm C^0(\mc; PSL(2,\mc))$. Since uniform limit of holomorphic
functions is holomorphic the first part of the lemma is achieved. In
fact, the same argument proves also the last part.  \cvd
\smallskip

The computation of the derivative of $u_\lambda$ at $0$ follows easily
from Lemma~\ref{var:der:lem}.

Notice that $\sG\lG(2,\mc)$ is the complexification of
$\sG\lG(2,\mr)$ that is
\[
    \sG\lG(2,\mc)=\sG\lG(2,\mr)\oplus i\sG\lG(2,\mr) \ .
\]
Now if $l$ is an oriented geodesic denote by $X(l)\in\sG\lG(2,\mr)$
the unitary generator of positive translations along $l$.  The element
$iX(l)/2$ is the standard generator of positive rotation around $l$
(see Section~\ref{hyp:bend:cocy}).  Thus if $\lambda$ is a finite
lamination and $l_1,\ldots,l_n$ are the geodesics between $x$ and $y$
with respective weights $a_1,\ldots,a_n\in\mc$ we have that
\[
    \frac{\d E_{z\lambda}(x,y)}{\,z}|_0=\frac{1}{2}\sum_{i=1}^{n} a_iX_i \ .
\]
The following statement is a corollary of this formula and
Lemma~\ref{var:der:lem}.
\begin{prop}\label{var:der:prop}
If $\lambda=(H,\Ll,\mu)$ is a complex-valued measured geodesic
lamination and $x,y$ are in $\mathring H$ then
\begin{equation}\label{var:der:eq1}
    \frac{\d E_{z\lambda}(x,y)}{\,z}|_0=\frac{1}{2}\int_{[x,y]} X(t)\d\mu(t)
\end{equation}
where $X(t)$ is so defined:
\[
\left\{\begin{array}{ll}
        X(t)=X(l) & \textrm{if } t\in\Ll\textrm{ and }l\textrm{ is the leaf
        through } t\\
        X(t)=0    & \textrm{otherwise \ .}
       \end{array}\right.
\]
\end{prop}
\cvd 
\smallskip

Now we can compute the derivative of $h_t^\kappa$ at $0$.  Recall that
the canonical isomorphism between $\sG\lG(2,\mr)$ and $\mx_0$ sends
$X(l)$ to the unit spacelike vector orthogonal to $l$ and giving the
right orientation to $l$ (see Remark \ref{mod:ort:rm}).

\begin{cor}\label{var:der:DS:cor}
The derivative of $h^{1}_{t\lambda}$ at $0$ is an imaginary cocycle
in
\[
\coom1_{\Ad}(\Gamma,\sG\lG(2,\mc))=\coom1_{\Ad}(\Gamma,\sG\lG(2,\mr))\oplus
i\coom1_{Ad}(\Gamma,\sG\lG(2,\mr)) \ .
\]
Moreover, up to the identification of $\sG\lG(2,\mr)$ with $\mr^3$ we
have that
\[
    \dot h^1_{t\lambda}(0)=\frac{i}{2}\tau_\lambda
\]
where $\tau_\lambda\in\coom1(\Gamma,\mr^{3})$ is the translation part of
$h^0_\lambda$.
\end{cor}
\cvd 

In the same way we have the following statement
\begin{cor}\label{var:der:ADS:cor}
The derivative of $h^{(-1)}_{t\lambda}$ at $t=0$ is a pair of cocycles
$(\tau_-,\tau_+)\in\coom1(\Gamma,\sG\lG(2,\mr))
\oplus\coom1(\Gamma,\sG\lG(2,\mr))$. In particular, if $\tau_\lambda$ 
is the translation part of $h^0_\lambda$, then
\[
\begin{array}{l}
\tau_-=-\frac{1}{2}\tau_\lambda \ ,\\
\tau_+=\frac{1}{2}\tau_\lambda \ .
\end{array}
\]
\end{cor}
\cvd

\subsection{Derivatives of spectra}\label{derivateSPEC}
Let us denote by $\Cc$ the set of conjugacy classes
of hyperbolic elements of the group $\Gamma$.
For every $\kappa = 0,\pm 1$,
we associate to $[\gamma]\in \Cc$
two numerical ``characters'' 
$\ell_\lambda^\kappa([\gamma])$ and 
$\Mm_\lambda^\kappa([\gamma])$.

First consider $\kappa=0$.  
Define $\ell_\lambda^0([\gamma])$ to be
the translation length of $\gamma$. $\Mm_\lambda^0([\gamma])$ was
introduced by Margulis in~\cite{Marg}. Denote by
$\tau\in\Z^1(\Gamma,\mr^3)$ the translation part of
$h^0_\lambda$ (obtained by fixing a base point $x_0\in\mh^2$). 
Denote by $X\in\sG\lG(2,\mr)$ the
unit positive generator of the hyperbolic group containing $\gamma$.
Let $v\in\mr^3$ be, as above, the corresponding point in the Minkowski
space. Then we have
\[
      \Mm_\lambda^0([\gamma])=\E{v}{\tau(\gamma)} \,.
\] 
It is not hard to see that $\Mm^0_\lambda$ is well defined.\\
\smallskip

Consider now the case $\kappa=1$. Take $[\gamma]$ such that
 $h^1_\lambda(\gamma)$ is hyperbolic.  In this case
 $\ell^1_\lambda([\gamma])$ is the length of the simple closed
 geodesic $c$ in $\mh^3/h_\lambda^1(\gamma)$. On the other hand
 $\Mm_\lambda^1([\gamma])\in[-\pi,\pi]$ is the angle formed by a
 tangent vector $v$ orthogonal to $c$ at a point $x\in
 c\subset\mh^3/h_\lambda^1(\gamma)$ with the vector obtained by the
 parallel transport of $v$ along $c$.  A computation shows that
 \[
  \tr(h^1_\lambda(\gamma))=
2\ch(\frac{\ell^1_\lambda([\gamma])}{2}+i\frac{\Mm^1_\lambda([\gamma])}{2}) \ .
\]
In particular it follows that $h^1_\lambda(\gamma)$ is conjugated to
an element of $PSL(2,\R)$ if and only if
$\Mm^1_\lambda([\gamma])=0$.\\
   
Finally consider the case $\kappa=-1$. If $\gamma$ is hyperbolic, then
$h^{(-1)}_\lambda(\gamma)$ is a pair of hyperbolic transformations,
$(h_-(\gamma), h_+(\gamma))$ (in fact by choosing the base point on
the axis of $\gamma$, the axis of $\beta_-(x_0,\gamma x_0)$ intersects
the axis of $\gamma$).

There are exactly two spacelike lines $l_-,l_+$ invariant by 
$h^{(-1)}_\lambda(\gamma)$.
Namely, $l_+$ has endpoints
\[
\begin{array}{ll}
  p_-=(x^-_L,x^-_R)\ , & p_+=(x^+_L,x^+_R)
\end{array}
\]
and $l_-$ has endpoints
\[ 
\begin{array}{ll}
  q_-=(x^+_L,x^-_R)\ , & q_+=(x^-_L,x^+_R)
\end{array}
\]
where $x^\pm_L$ (resp. $x^\pm_R$) are the fixed points of
$h_-(\gamma)$ (resp. $h_+(\gamma)$).  Orient $l_+$ (resp. $l_-$)
from $p_-$ towards $p_+$ (resp. from $q_-$ towards $q_+$).  If $m,n$
are the translation lengths of $h_-(\gamma)$ and $h_+(\gamma)$
then  $h^{-1}_\lambda(\gamma)$ acts on $l_+$
by a positive translation of length equal to $\frac{m+n}{2}$ and on
$l_-$ by a translation of a length equal to $\frac{n-m}{2}$.  Thus let
us define
\[
\begin{array}{ll}
\ell^{-1}_\lambda([\gamma])=\frac{m+n}{2}\\
\Mm^{-1}_\lambda([\gamma])=\frac{n-m}{2} \ .
\end{array}
\]
\begin{prop}
If $\gamma$ is a hyperbolic element of $\Gamma$ then there exists
$t<1$ sufficiently small such that $h_{s\lambda}^1(\gamma)$ is
hyperbolic for $s<t$.  Moreover, for every choice of the curvature
$\kappa$ , the following formulae hold
\[
\begin{array}{l}
   \frac{\d \ell_{t\lambda}^\kappa([\gamma])}{\,\d t}|_0=0\\
   \frac{\d \Mm_{t\lambda}^\kappa([\gamma])}{\,\d t}|_0=
\Mm_\lambda^0([\gamma]) \ .
\end{array}
\]
\end{prop}
\Dim For $\kappa=0$ the statement is trivial.
\par Suppose $\kappa=1$. Denote by
$B_t$ the cocycle associated to the lamination $\lambda_t$
\[
    \tr(B_t(x_0,\gamma x_0)\gamma)=
2\ch(\frac{\ell^1_t([\gamma])+i\Mm^1_t([\gamma])}{2}) \ .
\]
By deriving at $0$ we obtain
\[
   \frac{1}{2}\tr(i X(\gamma)\gamma)=\sh(\ell([\gamma])/2)(\dot
   \ell^1([\gamma])|_0 + i \dot\Mm^1([\gamma])|_0)
\]
where $X(\gamma)$ is the element of $\sG\lG(2,\mr)$ corresponding to
$\tau(\gamma)\in\mr^3$ (where $\tau$ is the is the translation part of
$h^0_\lambda$ ).  Now if $Y\in\sG\lG(2,\mr)$ is the unit generator of
the hyperbolic group containing $\gamma$ we have
\[
     \gamma=\ch(\ell([\gamma])/2) I + \sh(\ell([\gamma])/2) Y \ .
\]
Thus we obtain
\[
  \dot \ell^1([\gamma])|_0 + i \dot\Mm^1([\gamma])|_0=i\Mm^0([\gamma]) \ . 
\]
An analogous computation shows the same result when $\kappa=-1$.
\cvd

\section {More compact Cauchy surfaces}\label{more_cocompact}
In this Section we focus on the case of compact Cauchy surfaces,
pointing out a few specific applications.  Throughout the section we
consider a cocompact group $\Gamma$, so that $F=\mh^2/\Gamma$ is
compact surface of genus $g\geq 2$. Moreover, a $\Gamma$-invariant
measured geodesic lamination $\lambda$ on $\mh^2$ is fixed.

\subsection {Derivative of $(\Uu^{-1}_{t\lambda})^*$}\label{derivate*}
Consider the family of AdS spacetimes
$$(\Uu^{-1}_{t\lambda})^* = \Uu^{-1}_{\lambda_t^*} \ .$$ obtained by
the $T$-symmetry.  We want to determine the \emph{derivative} at
$t=0$ of this family.

Recall that in such a case the set of $\Gamma$-invariant measured
geodesic laminations, say $\Mm\Ll(F)$, has a natural $\R$-linear
structure, induced by the identification of this space with $\coom
1(\Gamma,\mr^3)$. So it makes sense to consider $-\lambda$.  We have
(the meaning of the notations is as above)
\begin{prop}\label{der:inv:prop}
\[ \lim_{t \rightarrow 0}\  \frac{1}{t}\Uu^{-1}_{\lambda^*_t}=  
\Uu^0_{-\lambda}
\ .\]
\end{prop}
\Dim 

Let $(h'_t,h''_t)$ be the holonomy of $\ \Uu^{-1}_{t\lambda}$.  Denote
by $F^*_t$ the quotient of the past boundary of $\Kk_{t\lambda}$ by
$(h'_t,h''_t)$. Notice that $\lambda^*_t$ is a measured geodesic
lamination on $F^*_t$.  We claim that $(F^*_t,\lambda^*_t/t)$
converges to $(F,-\lambda)$ in $\Tt_g\times\Mm\Ll_g$ as $t\rightarrow
0$.  Before proving the claim we conclude the proof.

Choose a family of developing
maps
\[
   D_t:\tilde F\times\mr\rightarrow\Uu^0_{\lambda^*_t/t}\subset\mx_0
\]
such that $D_t$ converges to a developing map $D_0$ of
$\hat\Uu^0_{-\lambda}=\Uu^0_{-\lambda}/h^0_{-\lambda}$ as
$t\rightarrow 0$.  Denote by $k_t$ the flat Lorentzian metric on
$\tilde F\times\mr$ corresponding to the developing map $D_t$. We have
that $k_t$ converges to $k_0$ as $t\rightarrow 0$.  Moreover, if $T_t$
denotes the cosmological time on $F\times\mr$ induced by $D_t$, then
$T_t$ converges to $T_0$ in $\mathrm C^1(F\times\mr)$ as $t\rightarrow
0$.  Now, as in the proof of Proposition~\ref{lim:space},
$\Pp_{\lambda^*_t}$ is obtained by a Wick Rotation directed by the
gradient of $T_t$ with rescaling functions
\[
  \begin{array}{ll}
    \alpha=\frac{t^2}{1+(tT_t)^2}\,, \qquad
&\qquad\beta=\frac{t^2}{(1+(tT_t)^2)^2}\,.
  \end{array}
\]
By passing to the limit $t\rightarrow 0$ we get the statement.\\

Let us prove the claim.  First the set $\{(F^*_t,\lambda^*_t/t)|t\in
[0,1]\}$ is shown to be pre-compact in $\Tt_g\times\Mm\Ll_g$, and then
$(F,-\lambda)$ is proved to be the only possible limit of any sequence
$(F^*_{t_n},\lambda^*_{t_n}/t_n)$.\\

By Section~\ref{ge:quake}, $F'_t=\mh^2/h'_t(\Gamma)$
(resp. $F''_t=\mh^2/h''_t(\Gamma)$) is obtained by a right
(resp. left) earthquake on $F=\mh^2/\Gamma$ with shearing measured
lamination equal to $t\lambda$.  Thus, if $\lambda'_t$ is the measured
geodesic lamination of $F'_t$ corresponding to $t\lambda$ via the
canonical identification of $\Mm\Ll(F)$ with $\Mm\Ll(F'_t)$, we have
that $F''_t$ is obtained by a left earthquake on $F'_t$ along
$2\lambda'_t$ .

On the other hand let $(\lambda^*_t)'$ be the
measured geodesic lamination on $F'_t$ such that the right earthquake
along it sends $F'_t$ on $F''_t$.  Then, the quotient
$F^*_t$ of the past boundary of the convex core $\Kk_{t\lambda}$ is
obtained by a right earthquake along $F'_t$ with shearing lamination
$(\lambda^*_t)'$. Moreover, the bending locus $\lambda^*_t$ is the
lamination on $F^*_t$ corresponding to $2(\lambda^*_t)'$.

In order to prove that the family $\{(F^*_t,\lambda^*_t/t)|t\in [0,1]\}$
is pre-compact we will use some classical facts about $\Tt_g$. For the
sake of clarity we will recall them, referring to
~\cite{Pen,Ot} for details.

Denote by $\Cc$ the set of
conjugacy classes of $\Gamma$.
For $\lambda\in\Mm\Ll(S)$ we denote by $\iota_\gamma(\lambda)$ the
total mass of the closed geodesic curve corresponding to $[\gamma]$
with respect to the transverse measure given by $\lambda$.  The
following facts are well-known.
\begin{enumerate}
\item
Two geodesic laminations $\lambda$ on $S$ and $\lambda'$ on $S'$ are
identified by the canonical identification
$\Mm\Ll(S)\rightarrow\Mm\Ll(S')$ if and only if
$\iota_\gamma(\lambda)=\iota_\gamma(\lambda')$ for every
$[\gamma]\in\Cc$.
\item
A sequence $(F_n,\lambda_n)$ converges to $(F_\infty,\lambda_\infty)$
in $\Tt_g\times\Mm\Ll_g$ if and only if $F_n\rightarrow F_\infty$ and
$\iota_\gamma(\lambda_n)\rightarrow\iota_\gamma(\lambda_\infty)$ for
every $[\gamma]\in\Cc$.
\item
A subset $\{(F_i,\lambda_i)\}_{i\in I}$ of $\Tt_g\times\Mm\Ll_g$ is
pre-compact if and only if the base points $\{F_i\}$ runs in a compact
set of $\Tt_g$ and for every $[\gamma]\in\Cc$ there exists a constant
$C>0$ such that
\[
     \iota_\gamma(\lambda_i)<C\qquad\textrm{ for every }i\in I.
\]
\end{enumerate}
Clearly we have $F^*_t\rightarrow F$ as $t\rightarrow 0$.  Thus in
order to show that $(F_t^*,\lambda^*_t)$ is pre-compact it is sufficient
to find for every $[\gamma]\in\Cc$ a constant $C>0$ such that
\[
    \iota_\gamma(\lambda^*_t)<Ct
\]
for every $t\in [0,1]$.
  
The following lemma gives the estimate we need.
\begin{lem}\label{der:inv:lem}
For every compact set $K\subset\mh^2$ there exists a constant $M>0$
which satisfies the following statement.

If $\lambda$ is a
measured geodesic lamination on $\mh^2$ and $\beta$ is the right
cocycle associated to $\lambda$ then
\[
   ||\beta(x,y)-Id+\frac{1}{2}
\int_{[x,y]}X_\lambda(u)\d\lambda||\leq e^{M\lambda(x,y)}-1-M\lambda(x,y)
\]
where $X_\lambda(u)$ is defined as in~(\ref{var:der:eq1}), $x,y\in K$,
and $\lambda(x,y)$ is the total mass of the segment $[x,y]$.
\end{lem}
\Dim It is sufficient to prove the lemma when $\lambda$ is simplicial.
In this case denote by $l_1,\ldots,l_N$ the geodesics meeting the
segment $[x,y]$ with respective weights $a_1,\ldots,a_N$.  If
$X_i\in\sG\lG(2,\mr)$ is the unitary infinitesimal generator of the
positive translation along $l_i$ we have
\[
   \beta(x,y)= \exp(-a_1X_1/2)\circ\exp(-a_2X_2/2)\circ\cdots\circ\exp(-a_N
X_N/2)\ .
\]
Thus $\beta(x,y)$ is a real analytic function of $a_1,\ldots,a_n$.
If we write
\[
   \beta(x,y)=\sum_n A_n(a_1,\ldots,a_n)
\]
where $A_n$ is a matrix-valued homogenous polynomial in $x_1,\ldots,x_n$ of
degree $n$, then it is not difficult to see that
\[
     ||A_n||\leq (\sum_{i=1}^N a_i||X_i||)^n/n!\ .
\]
We have that
\[
   \beta(x,y)-Id+\frac{1}{2}
\int_{[x,y]}X_\lambda(u)\d\lambda=\sum_{i\geq 2} A_n(a_1,\ldots, a_N)
\]
Since the axes of transformations generated by $X_i$
cut $K$, there exists a constant $M>0$ (depending only on $K$) such
that $||X_i||<M$.  Thus
\[
   ||\beta(x,y)-Id+\frac{1}{2}
\int_{[x,y]}X_\lambda(u)\d\lambda||\leq e^{M\sum a_i}-1-M\sum a_i\ .
\]
\cvd
Let us go back to the proof of Proposition~\ref{der:inv:prop} .
Since
$\iota_\gamma(\lambda^*_t)=\iota_\gamma((\lambda^*_t)')$, we may
replace $\lambda^*_t$ with $(\lambda^*_t)'$.  Now let us put
$\gamma_t=h'_t(\gamma)$. We know that $\gamma_t$ is a
differentiable path in $ PSL(2,\mr)$ such that $\gamma_0=\gamma$ and
\[
   \dot\gamma(0)=-\frac{1}{2}\int_{[x,\gamma(x)]} X(u)\d\lambda(u)
\]
where $X(u)$ is defined as in~(\ref{var:der:eq1}).  On the other hand,
if $\beta_t$ is the right cocycle associated to the measured geodesic
lamination $2(\lambda^*_t)'$ we have
\[
  \beta_t(x,\gamma_tx)\gamma_t=\beta_{t\lambda}(x,\gamma x)\gamma
\]  
where $\beta_{t\lambda}$ is the left cocycle associated to $t\lambda$.
Thus $\beta_t(x,\gamma_tx)$ is a differentiable path and
\begin{equation}\label{der:rev:eq1}
  \lim_{t\rightarrow 0}
  \frac{\beta_t(x,\gamma_tx)-Id}{t}=\int_{[x,\gamma(x)]}X(u)\d\lambda\ .
\end{equation}
By Lemma~\ref{der:inv:lem}, there exists a
constant $C>0$ depending only on $\gamma$ such that
\[
    ||\beta_t(x,\gamma_tx)-Id||> ||\int_{[x,\gamma_tx]}
      X(u)\d\lambda^*_t|| - C \iota_\gamma((\lambda^*_t)')^2 .
\]
On the other hand, there exists a constant $L>0$ such that
\[
||\int_{[x,\gamma_tx]}
      X(u)\d\lambda^*_t||\geq L |\int_{[x,\gamma_tx]} X(u)\d\lambda^*_t|  
\]
where $|\cdot|$ denotes the Lorentzian norm of $\sG\lG(2,\mr)$.  Since
$X(u)$ are generators of hyperbolic transformations with disjoint axes
pointing in the same direction, the reverse of Schwarz inequality
inequality holds (see the proof of Lemma~\ref{ads:bend:lem})
\[
   \eta(X(u),X(v))^2\geq\eta(X(u),X(u))\eta(X(v),X(v))=1
\] 
and implies
\[
  ||\int_{[x,\gamma_tx]}
      X(u)\d(\lambda^*_t)'||\geq L\iota_\gamma((\lambda^*_t)') \ .
\]
From this inequality we obtain that
\[
     \left(L-C\iota_\gamma((\lambda^*_t)')\right)\iota_\gamma((\lambda^*_t)')
<||\beta_t(x,\gamma_tx)-Id|| \ .
\]
Dividing by $t$ the last inequality shows that
$\lambda^*_t(x,\gamma_tx)/t$ is bounded.  In particular we have proved
that $\{(F_t,\lambda^*_t)\}$ is pre-compact in
$\Tt_g\times\Mm\Ll_g$.\\

Now, let us set $\mu_t=\lambda^*_t/t$ and $\mu'_t=(\lambda^*)'_t/t$.
We have to show that if $\mu_{t_n}\rightarrow \mu_\infty$ then
$\mu_\infty=-\lambda$ in $\Mm\Ll(F)$.

Notice that $\mu'_{t_n}$ is
convergent and its limit is $\mu_\infty$.
Applying lemma~\ref{der:inv:lem} we get 
\[
   \lim_{t\rightarrow 0}\frac{\beta_t(x,\gamma_tx)-Id}{t}=
   -\int_{[x,\gamma(x)]}X_{\mu_\infty}(t)\d\mu_\infty
\]
By equation~(\ref{der:rev:eq1}) this limit is equal
to $\int_{[x,\gamma(x)]}X\d\lambda$ and this shows that
$\mu_\infty=-\lambda$.  \cvd

\subsection{Far away along a ray}\label{big-t}
Till now we have derived infinitesimal information at $t=0$. As
regards the behaviour along a ray for {\it big} $t$, let us just make
a qualitative remark.

We have noticed that, for every $t>0$, $\frac{1}{t}\Uu^0_{t\lambda}=
\Uu^0_\lambda$. Moreover, the flat spacetimes $\ \Uu^0_{t\lambda}\ $
are nice convex domains in $\mx_0$ which vary continuously and tamely
with $t$.  So, in the flat case, {\it apparently nothing substantially
new happens when $t>0$ varies.}  Similarly, this holds also for the
AdS past parts $\Pp_{t\lambda}\subset \Uu^{-1}_{t\lambda}$.  On the
other hand, radical qualitative changes do occur for $M_{t\lambda}$
(and $\Uu^{1}_{t\lambda}$) when $t$ grows. As $\lambda$ is
$\Gamma$-invariant for the cocompact group $\Gamma$, when $t$ is small
enough, we have a {\it quasi-Fuchsian} hyperbolic end. In particular,
the developing map is an embedding. When $t$ grows up, we find a first
value $t_0$ such that we are no longer in the quasi-Fuchsian region,
and for bigger $t$ the developing map becomes more and more ``wild''.
We believe that this different behaviour along a ray is conceptually
important: looking only at the flat Lorentzian sector, significant
critical phenomena should be lost; on the other hand, one could
consider the (flat or AdS Lorentzian towards hyperbolic geometry) Wick
rotations as a kind of ``normalization'' of the hyperbolic developing
map.
\smallskip

We give here a first simple application of these qualitative
considerations.

Assume that we are in the quasi-Fuchsian region.
So we have associated to $t\lambda$ three ordered pairs
of elements of the Teichm\"uller space $\Tt_g$. These are:
\smallskip

- the ``Bers parameter'' $(B_t^+,B_t^-)$
given by the conformal structure underlying the projective
asymptotic structures of the two ends of $Y_t$; 
\smallskip

- the hyperbolic structures $(C_t^+, C_t^-)$ of the boundary components
of the hyperbolic convex core of $Y_t$;
\smallskip

- the hyperbolic structures $(K_t^+,K_t^-)$ of the future and past
boundary components of the AdS convex core of $\Uu^{-1}_{t\lambda}$.
\smallskip

\noindent
It is natural to inquire about the relationship between
these pairs.
\smallskip

By construction, $K_t^+$ is isometric to $C_t^+$.  On the other hand
by Sullivan's Theorem (see \cite{Ep-M}), the Teichm\"uller distance of
$B_t^\pm$ from $C_t^\pm$ is uniformly bounded.  Now it is natural to
ask whether $C_t^-$ is isometric to $K_t^-$.  Actually it is not hard
to show that those spaces generally are not isometric.  In fact let us
fix a lamination $\lambda$ and let $t_0>0$ be the first time such that
the representation $h^1_{t_0\lambda}$ is not quasi-Fuchsian.  By Bers
Theorem~\cite{bers}   the family $\{(B^+_t,
B^-_t)\in\Tt_g\times\overline\Tt_g\}$ is not compact. Since $B^+_t$
converges to a conformal structure as $t$ goes to $t_0$ we have that
$\{B^-_t\}_{t\leq t_0}$ is a divergent family in $\overline\Tt_g$.  By
Sullivan's Theorem we have that $C_t^-$ is divergent too. On the other
hand $\{K_t^-\}_{t\leq t_0}$ is pre-compact.

\subsection{Volumes, areas and length of laminations}\label{volume}
Set $Y^0= \Uu^0_\lambda/h^0_\lambda(\Gamma)$,
$Y^1=\Uu^1_\lambda/h^1_\lambda(\Gamma)$ and
$Y^{-1}=\Pp_\lambda/h^{(-1)}_\lambda(\Gamma)$; \emph{i.e.} these
$Y^\kappa$, are the spacetime of constant curvature, with compact
Cauchy surface homeomorphic to $F$, related to each other via
equivariant canonical rescalings starting from the $\Gamma$-invariant
lamination $\lambda$.

In this Subsection we
compute the volume of any slab $Y^\kappa(\leq b)$ in terms of $F$ and
$\lambda$.

Let us outline the scheme of such a computation. We
first get a formula expressing the volume of $Y^\kappa(\leq b)$ in
terms of the areas of level surfaces $Y^\kappa(t)$ for $t\leq b$. Then
we compute these areas. Thanks to rescaling
formulae, it is sufficient to compute the area of $Y^0(t)$. When
the lamination is simplicial the computation is quite trivial. In the
general case, by using the continuity of the area of $Y^0(t)$ with
respect to the parameter $\lambda$, we express the area of $Y^0(t)$ in
terms of the well known notion of {\it length} of $\lambda$ 
(in fact of the induced lamination on $F$).

We use the following notation:
\medskip\par\noindent
- $V_{(F,\lambda)}(\kappa, b)$ is the volume of $Y^\kappa(\leq b)$;
\medskip\par\noindent
- $A_{(F,\lambda)}(\kappa, b)$ is the area of $Y^\kappa(b)$.
\begin{lem}
 The following formula holds
\begin{equation}\label{vol:eq}
    V(\kappa, b)=\int_0^b A(\kappa, t)\d t\,.
\end{equation}
\end{lem}
\Dim Let us fix some coordinates $x,y$ on the level surface
$Y^\kappa(1)$.  If $\varphi:Y\times\mr\rightarrow Y$ denotes the flow
of the gradient of the cosmological time, we get that
\[
  (x,y,t)\mapsto\varphi_t(x,y)
\]
furnishes a parameterization of $Y$. Moreover, since the gradient of
the cosmological time is a unitary vector, then the map
\[
  (x,y)\mapsto\varphi_t(x,y)
\]
takes values on the surface $Y^\kappa(1+t)$.

So the volume form of $Y^\kappa$, with respect to these coordinates,
takes the form
\[
   \Omega(x,y,t)=\omega_t(x,y)\d t
\]
(we are using again that the gradient of the cosmological time is unitary).
The formula~\ref{vol:eq} easily follows.
\cvd

When $\lambda$ is a weighted multicurve, the computation of the area
of $A_{(F,\lambda)}(0,b)$ is quite simple. Namely, if $l_1,\ldots,l_k$
are the leaves of $\lambda$ with weights $a_1,\ldots, a_k$, the
surface $\Uu(b)$ is obtained by rescaling the surface by $b$, and
replacing every $l_i$ by an Euclidean annulus of length $a_i$. So the
area of $Y^\kappa(b)$ is given by the formula
\begin{equation}\label{area1:eq}
   -2\pi\chi(F)b^2+b\sum_{i=1}^ka_i\ell_F(l_i)
\end{equation}
where $\ell_F(l_i)$ is the length of $l_i$.\par
In~\cite{Bon} it was shown that there exists a continuous function
\[
    \ell:\Tt_g\times\Mm\Ll\rightarrow\mr_+
\]
such that if $\lambda$ is a weighted multicurve then 
\[
  \ell(F,\lambda)=\sum_{i=1}^ka_i\ell_F(l_i)
\]
where $l_i$'s are the leaves of $\lambda$ and $a_i$'s are the
corresponding weights.  We call $\ell(F,\lambda)$ the length of
$\lambda$ with respect to $F$.

From (\ref{area1:eq}) we get
\begin{equation}\label{area2:eq}
   A_{(F,\lambda)}(0,b)=-2\pi b^2\chi(F)+b\ell(F,\lambda)
\end{equation}
whenever $\lambda$ is a weighted multicurve. 
The right hand of this expression continuously depends on $\lambda$ . 
In fact, by means of results of Sec. 6 of \cite{Bo} we have that
also $A_{(F,\lambda)}(0, b)$ varies continuously with $\lambda$.
Since weighted multicurve are dense in $\Mm\Ll_g$, formulae~(\ref{area2:eq})
holds for every measured geodesic lamination.\par
Since $Y^{-1}(b)$ (resp. $Y^1(b)$) is obtained by rescaling 
$Y^0(\tan b)$ (resp. $Y^0(\tgh b)$) by $\frac{1}{1+\tan^2 b}$
(resp. $\frac{1}{1-\tgh^2 b}$) we obtain the following formulae
\[
\begin{array}{l}
 A_{(F,\lambda)}(0,b)=-2\pi b^2\chi(F)\ +\ b\ell(F,\lambda)\,;\\
 A_{(F,\lambda)}(-1,b)= -2\pi\sin^2 b\chi(F)\ +\ \ell_F(\lambda)\sin
 b\cos b\,;\\ A_{(F,\lambda)}(1,b)=-2\pi\sh^2 b\chi(F) \ +\
 \ell_F(\lambda)\sh b\ch b\,.
\end{array}
\]
By~(\ref{vol:eq}), we have
\[
\begin{array}{l}
  V_{(F,\lambda)}(0,b)=-\frac{2\pi\chi(F)}{3}b^3\ +\
  \frac{\ell_F(\lambda)}{2}b^2\,;\\
  V_{(F,\lambda)}(-1,b)=-2\pi\chi(F)\frac{2b-\sin^2b}{4}\ +\
  \frac{1-\cos^2b}{4}\ell_F(\lambda)\,;\\
  V_{(F,\lambda)}(1,b)=-2\pi\chi(F)\frac{\sh^2b- 2b}{4}\ +\frac{\ch^2
  b-1}{4}\ell_F(\lambda).
\end{array}
\]
In particular we get that the volume of $Y^{-1}$ is given by
\[
   -\pi^2\chi(F)\, +\, \frac{1}{2}\ell_F(\lambda)\,.
\]
\begin{remark}{\rm
Given a maximal globally hyperbolic spacetime, the last formula allows
to compute the volume of its past part, in terms of parameters
obtained by looking at the future boundary of its convex
core. Clearly, by inverting the time-orientation, we obtain a similar
formula expressing the volume of the future part in terms of the past
boundary of the convex core.

Thus, the computation of the volume
of the whole  spacetime turns to be equivalent to the computation
of the volume of the convex core. Similar considerations hold
if $F/\Gamma$ is of finite area and we use laminations with compact
support.}
\end{remark}

\begin{remark}{\rm
Equation~(\ref{area1:eq}) and the remark that $A_{(F,\lambda)}(0,b)$
continuously varies with the pair $(F,\lambda)$ furnish another proof
of the existence and the continuity of the function $\ell$.}
\end{remark}

\section{Including particles}\label{particles}
We have seen in Section \ref{QD} that (globally hyperbolic)
$\Qq\Dd$-spacetimes contain in general ``particles'' that is conical
singularities along time-like lines (the ``world lines'' of the
particles). In that case the level surfaces of the canonical time were
{\it flat} surfaces with conical singularities (shortly: {\it flat
cone surfaces}), ``orthogonal'' to the particle world lines.  Starting
with {\it hyperbolic} cone surfaces and adapting $\Mm\Ll$ Wick
rotation-rescaling constructions, we could produce (correlated)
$\Mm\Ll$-{\it spacetimes with particles}. In the flat case, such a
kind of constructions were already considered for instance in
\cite{BG}(2); we refer to it also for some precision about the local
models at a particle. By the way, we recall that, up to a suitable normalization of the gravitational constant,
 the relation between the {\it cone
angle} $\beta=2\pi\alpha$ and the {\it ``mass''} $m$ of the particle
is given by
$$m=1-\alpha$$
hence the mass is positive only if $\beta < 2\pi$. On the other hand,
there is no natural geometric reason to impose a priori such a constraint
(for example the $\Qq\Dd$ spacetimes usually contain negative masses).
\smallskip

We do not intend here to fully develop such a generalization of our
theory. We limit ourselves to give some sketch in the case of compact
Cauchy surfaces.
\smallskip

Fix a base closed surface $S$ of genus $g$ with a set of $r>0$ marked
points $V=\{x_1,\dots,x_r\}$, such that $S\setminus V$ admits a
complete hyperbolic structure of finite area 
(that is $2-2g -r<0$). Denote by
$$\Tt_{g,r}$$ the corresponding Teichm\"uller space. Similarly to
$\Tt_g$ there is a canonically trivialized fiber bundle over
$\Tt_{g,r}$, such that the fiber over every $F\in \Tt_{g,r}$
consists of the measured geodesic laminations on $F$ with {\it compact
support}. Let us denote by
$$\Tt_{g,r}\times \Mm\Ll_{g,r}$$ such a (trivialized) bundle.

Fix $\Theta= (\beta_1,\dots,\beta_r)\in \mr^r_+$
and denote by 
$$\Cc_{g,r,\Theta}$$
the Theichm\"uller-like space of hyperbolic cone surface structures
on $(S,V)$ with assigned values $\Theta$ of the cone angles at $V$.
For every $F\in  \Cc_{g,r,\Theta}$, we denote by 
$$\Mm\Ll(F)$$
the set of ``measured geodesic laminations'' on $F$. We do no enter
the actual definition; note however that one should allow in general
{\it singular leaves} containing some cone points. 

\paragraph{Particles with ``big'' masses.}
Denote by $\Bb$ the subset of $\mr^r_+$ such that 
for every $j$, $0\leq \beta_j<\pi$. If $\Theta \in \Bb$ we say that it
corresponds  to ``big masses''. We have
\begin{prop} 
Assume $\Theta\in\Bb$, then:
\smallskip

(1) $\Cc_{g,r,\Theta}$ is naturally isomorphic to $\Tt_{g,r}$.
\smallskip

(2) For every $F\in  \Cc_{g,r,\Theta}$, every $\lambda \in \Mm\Ll(F)$
has compact support contained in $S\setminus V$; moreover $\Mm\Ll(F)$ 
is canonically isomorphic to  $\Mm\Ll_{g,r}$. 
\end{prop} \cvd
The first statement is due  Troyanov \cite{Tr}; the second follows 
by the same arguments used for $\Mm\Ll_g$ (see also \cite{BS}).

\noindent Hence  
$$\Tt_{g,r}\times \Mm\Ll_{g,r}$$ can be considered as a trivialized fiber
bundle over $\Cc_{g,r,\Theta}$, for every $\Theta$ corresponding to
big masses.

For every $(F,\lambda,\Theta)\in \Tt_{g,r}\times \Mm\Ll_{g,r}\times
\Bb$, denote by $(F_c,\lambda_c)$, $F_c\in \Cc_{g,r,\Theta}$,
$\lambda_c \in \Mm\Ll(F)$ the corresponding element. Then the
construction of Section \ref{ML_REGD} applies to $(F_c,\lambda_c)$,
far from the cone points of $F_c$, and produces a flat globally
hyperbolic spacetime $\hat\Uu^0(F,\lambda,\Theta)$ with particles. The
level surfaces of the canonical time with their conical points are homeomorphic to $(S,V)$, are 
orthogonal to the particle world lines and are (rescaled) hyperbolic
cone surfaces at the particles, with constant cone angles $\Theta$.
Similarly to the $\Qq\Dd$ case, the developing map on the complement
of the particles is no longer an embedding.

Now all the Wick rotation-rescaling formulas apply verbatim and either
produce cone hyperbolic manifolds or cone spacetimes, keeping the
cone angles $\Theta$.

Summing up, we have produced a structured family of $\Mm\Ll$-{\it spacetimes
with particles} of constant curvature $\kappa = 0, \pm 1$, that share
$$\Tt_{g,r}\times \Mm\Ll_{g,r}\times \Bb$$ as universal parameter spaces,
and are canonically correlated to each other. A natural question
(\cite{BS} is addressed to it) asks {\it to point out intrinsic
characterizations of this family}, in the spirit of the classification
Theorem \ref{FULL_CLASS}.

\paragraph{Allowing arbitrary masses.}
If $\Theta$ does not necessarily belong to $\Bb$ the situation is more
complicated and far to be understood. For example in \cite{BG}(2)
there are examples of (flat) spacetimes with particles obtained via so
called ``patchworking''. These are somewhat intriguing as they combine
features of both $\Qq\Dd$ and $\Mm\Ll$ ones. On the other hand, it is
natural to extend the $\Mm\Ll$ constructions to laminations with
singular leaves. In fact, one can realize that, at least in some case,
patchworking spacetimes are such generalized $\Mm\Ll$ ones.

\section{Open questions}\label{PROB}
In this final section we state a few open questions addressing 
further developments of the Wick rotation-rescaling theory.

\paragraph{(1) Characterization of AdS $\Mm\Ll(\mh^2)$-spacetimes.}
Give a characterization of AdS $\Mm\Ll(\mh^2)$-spacetimes and of
broken $T$-symmetry purely in terms of properties of the corresponding
curves at infinity (see Sections \ref{ge:quake} and \ref{3cusp}).
Related to it, find further sensible characterizations of projective
structures on surfaces, associated to $\Mm\Ll(\mh^2)$-spacetimes.

\paragraph{(2) On AdS canonical time.} 
Recall that the range of the canonical time of an AdS $\Mm\Ll$-spacetime
$\Uu^{-1}_\lambda$ is of the form $(0,a_0)$, $\pi/2<a_0 < \pi$.
Study $a_0=a_0(\lambda)$ as a function of $\lambda$.
We know that the canonical time is ${\rm C}^1$ on the past part $\Pp_\lambda$.
Study its lack of regularity on the slab $\Uu^{-1}_\lambda([\pi/2,a_0))$.

\paragraph{(3) Canonical versus CMC times.} We know since \cite{Mon}
that flat maximal globally spacetimes with compact Cauchy surface
admit a canonical constant mean curvature CMC time. In fact the
existence of canonical CMC time holds in general for maximal globally
spacetimes of constant curvature (see \cite{A-M-T, B-Z}). We could
study the behaviour of CMC time under the canonical Wick
rotation-rescaling. On the other hand, does there exist a Wick
rotation-rescaling theory {\it entirely} based on the CMC time,
instead of the canonical cosmological one? 
A partial (negative) answer has been given in~\cite{B-Z}:
in particular it has been shown that there is no Wick rotation 
with rescaling functions constant of the level surfaces of CMC time of a AdS spacetime, transforming it in a hyperbolic manifold.

\paragraph{(4) Ends of arbitrary tame hyperbolic $3$-manifolds.}
We have seen in Section \ref{END} that Wick rotation-rescaling applies
in a clean way on the ends of geometrically finite hyperbolic
$3$-manifolds. It would be interesting to use this machinery to treat
the ending geometry of more general topologically tame manifolds; in
particular to study limits of quasi-Fuchsian groups (like ones
occurring in the ``double limit theorem'' see \cite{Ot}) in terms of
the associated families of flat or AdS ending spacetimes (see also
Subsection \ref{big-t}).

\paragraph{(5) Wick cut locus.} We refer again to Section
\ref{END} of Chapter \ref{INTRO}.  Let $Y$ be a geometrically finite
hyperbolic $3$-manifold. Consider the ending Wick rotations of $Y$
towards (slabs of) AdS spacetimes. We would like to define and study a
somewhat canonical subset $\Ww(Y)$ of $Y$ (called its {\it Wick cut
locus}), such that $\Omega(Y)= Y\setminus \Ww(Y)$ is a maximal open
set which supports a C$^1$ and almost everywhere real analytic Wick
rotation towards a spacetime of constant curvature $\kappa = -1$,
extending the canonical ones on the ends. $\Omega(Y)$ should carry
almost everywhere (via Frobenius theorem and an analytic continuation
argument) a foliation by spacelike surfaces which extends the ending
one by the level surfaces of the canonical times. It might happen or
not that two ends of $Y$ are connected by (the closure of) timelike
curves orthogonal to the foliation. This induces a relation between
ends. We know from Bers parametrization that $Y$ is {\it
over-determined} by the family of its asymptotic projective
structures. One might wonder if the above relation reflects in some
way the implicit relationship existing among those projective
structures.

\paragraph{(6) Particles.} 
Fully develop a Wick rotation-rescaling theory on spacetimes with
particles (see Section \ref{particles}). 

\paragraph{(7) Higher dimensional Wick rotation-rescaling.}
The theory of flat regular domains of arbitrary dimension is developed
in \cite{Bo}. In particular, there is a notion of {\it measured
geodesic stratification} that generalizes the one of measured geodesic
laminations of the $(2+1)$ case. The subclass of so called {\it
simplicial} stratifications is particularly simple to dealing with,
and allows a very clean generalization of the results of Sections
\ref{laminations} - \ref{REGD_ML} of Chapter \ref{FGHST}.
Unfortunately, there is not a straightforward generalization of Wick
rotation-rescaling, not even for such a subclass.  Let us consider for
simplicity the $(3+1)$ case. Let $\Uu=\Uu_\lambda$ be a regular domain
in $\mm^4$ corresponding to a simplicial measured geodesic
stratification $\lambda$. The initial singularity $\Sigma_\Uu$ of
$\Uu$ is $2$-dimensional and, in general, it contains strata of
dimension $\geq 1$.  Let $Z$ be the union of $2$-dimensional
strata. Then the usual Wick rotation-rescaling formulas hold verbatim
on $\Uu\setminus r^{-1}(Z)$. However, they do not extend to the whole
of $\Uu$. In fact, consider for instance the Wick rotation.  We easily
realize that the portion of the level surface $\Uu(1)$ over a
$2$-stratum of $Z$ has a {\it flat} spacelike geometry. On the other
hand, the portion on a level surface of the distance function from the
hyperbolic boundary of the hyperbolic $4$-manifolds $M_\lambda$ (which
is still globally defined), that should naturally corresponds to it
has {\it spherical} spacelike geometry. So such a (conformal) global
Wick rotation, with universal rescaling functions, cannot exist.
Building a reasonably canonical Wick rotation-rescaling theory in
higher dimension is an interesting open question.



\end{document}